\documentclass[11pt,twoside,leqno,openany]{amsart}

\usepackage{amssymb,amsbsy,amsmath,amsfonts,amssymb,amscd,epsfig,
times,fancyhdr,graphics,color,xypic,wasysym,mathrsfs}

\let\mathcal\mathscr

\newcommand{\C}{\mathbb{C}}\newcommand{\D}{\mathbb{D}}

\newcommand{\N}{\mathbb{N}}\renewcommand{\P}{\mathbb{P}}
\newcommand{\Q}{\mathbb{Q}}\newcommand{\R}{\mathbb{R}}
\newcommand{\Z}{\mathbb{Z}}

\definecolor{blue}{cmyk}{1.,1.,0.,0.43}
\definecolor{red}{cmyk}{0.,1.,1.,0.43}
\definecolor{green}{cmyk}{1.,0.,1.,0.43}

\newcommand{\blue}{\textcolor{blue}}
\newcommand{\green}{\textcolor{green}}
\newcommand{\red}{\textcolor{red}}

\newtheorem{The}{Theorem}[section]

\newtheorem{MainTheorem}[The]{Main Theorem}
\newtheorem{Theorem}[The]{Theorem}
\newtheorem{Proposition}[The]{Proposition}
\newtheorem{Lemma}[The]{Lemma}
\newtheorem{Corollary}[The]{Corollary}

\theoremstyle{definition}

\newtheorem{Openproblem}[The]{Open problem}

\usepackage{fancyhdr}

\newcommand{\zero}[1]{\underline{#1}_{\circ}}
\newcommand{\zerozero}[1]{\underline{#1}_{\circ\circ}}

\begin{document}
\parindent 0.75cm

\pagenumbering{arabic}

\title[
Algebraic differential equations for entire holomorphic curves
]{
Algebraic differential equations
\\ 
for entire holomorphic curves
\\
in projective hypersurfaces of general type:
\\
optimal lower degree bound
}

\author{Jo\"el Merker}

\address{D\'epartement de math\'ematiques et applications, 
\'Ecole Normale Sup\'erieure, 45 rue d'Ulm, 75005 Paris, France}

\email{merker@dma.ens.fr}


\date{\number\year-\number\month-\number\day}

\begin{abstract}
Let $X = X^n \subset \P^{ n+1} ( \C)$ be a geometrically smooth
projective algebraic complex hypersurface. Using Green-Griffiths jets,
we establish the existence of nonzero global algebraic differential
equations that must be satisfied by every nonconstant entire
holomorphic curve $\C \to X$ if $X$ is of general type, namely if its
degree $d$ satisfies the optimal possible lower bound:
\[
d\geqslant 
n+3. 
\]
\end{abstract}

\maketitle

\vspace{-0.75cm}

\begin{center}
\begin{minipage}[t]{11.5cm}
\baselineskip =0.35cm
{\scriptsize

\centerline{\bf Table of contents}

{\bf 1.~Introduction\dotfill 
\pageref{Section-1}.}

{\bf 2.~Universal combinatorics of Green-Griffiths jets\dotfill 
\pageref{Section-2}.}

{\bf 3.~Euler-Poincar\'e characteristic of jet bundles and
multiple polylogarithms\dotfill 
\pageref{Section-3}.}

{\bf 4.~Exact Schur Bundle Decomposition of $\mathcal{ E}_{ 
\kappa, m}^{ GG} T_X^*$\dotfill 
\pageref{Section-4}.}

{\bf 5.~Asymptotic characteristic and asymptotic cohomology\dotfill 
\pageref{Section-5}.}

{\bf 6.~Emergence of basic numerical sums\dotfill 
\pageref{Section-6}.}

{\bf 7.~Asymptotic combinatorics of semi-standard Young 
tableaux\dotfill \pageref{Section-7}.}

{\bf 8.~Maximal length families\dotfill 
\pageref{Section-8}.}

{\bf 9.~Number of tight paths in semi-standard Young tableaux\dotfill 
\pageref{Section-9}.}

{\bf 10.~Bounded behavior of plurilogarithmic sums\dotfill 
\pageref{Section-10}.}

{\bf 11.~Algebraic sheaf theory and Schur bundles\dotfill 
\pageref{Section-11}.}

{\bf 12.~Asymptotic cohomology vanishing\dotfill 
\pageref{Section-12}.}

{\bf 13.~Observations about Demailly-Semple jet differentials\dotfill 
\pageref{Section-13}.}

}\end{minipage}
\end{center}

\medskip

\hfill
Shoshichi {\sc Kobayashi}${}^\dag$, Alain {\sc Lascoux}${}^\dag$, 
{\em In memoriam}\footnote{\,
The Main Theorem below has been presented at the {\sl
Memorial Symposium: "Geometry and Analysis on Manifolds"},
22-25 May 2013,
in honor of Professor Shoshichi Kobayashi${}^\dag$, 
organized by Takushiro
Ochiai, Keizo Hasegawa, Toshiki Mabuchi, Yoshiaki Maeda, Junjiro
Noguchi, Yoshihiko Suyama, Takashi Tsuboi, in the 
{\sl Graduate School of Mathematical Sciences} of the University of
Tokyo.}

\vspace{-0.15cm}

\markleft{Jo\"el Merker}
\markright{\sf \S1.~Introduction}
\section{\bf Introduction}
\label{Section-1}

Let $X$ be an $n$-dimensional ($n \geqslant 1$) compact complex
manifold and assume it to be {\sl of general type}, i.e., if as usual
$K_X = \Lambda^n T_X^*$ denotes its canonical line bundle, assume that
the dimension of the space of global pluricanonical sections:
\[
h^0\big(X,\,(K_X)^{\otimes m}\big) \geqslant {\sf Constant}\cdot
m^{\dim X} \ \ \ \ \ \ \ \ \ \ \ \ \ {\scriptstyle{({\sf
Constant}\,>\,0)}}
\] 
grows the fastest it can, as $m \to \infty$, namely the Kodaira
dimension of $X$ is maximal equal to $n$. According to a theorem due
to Kodaira, $X$ can then be embedded as a geometrically smooth
projective algebraic complex manifold in a certain complex projective
space $\P^N ( \C)$. Though it is somewhat delicate to select good
embeddings, it is algebraically convenient to view $X$ as being
projective {\em per se}.

In 1979, Green and Griffiths~\cite{gg1980} conjectured that there
should exist in $X$ a certain {\em proper} algebraic {\em sub}variety
$Y \subsetneqq X$ (possibly with singularities) inside which all
nonconstant entire holomorphic curves $f \colon \mathbb{ C} \to X$
must necessarily lie, without any such $f$ being allowed to wander
anywhere else in $X \backslash Y$. 

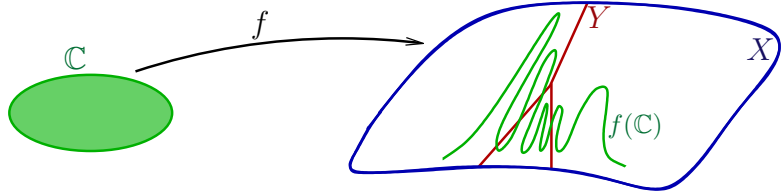
\begin{center}
\begin{picture}(0,0)%
\includegraphics{courbe-degenere.pstex}%
\end{picture}%
\setlength{\unitlength}{4144sp}%
\begingroup\makeatletter\ifx\SetFigFont\undefined%
\gdef\SetFigFont#1#2#3#4#5{%
  \reset@font\fontsize{#1}{#2pt}%
  \fontfamily{#3}\fontseries{#4}\fontshape{#5}%
  \selectfont}%
\fi\endgroup%
\begin{picture}(4645,1165)(871,-957)
\put(4341, 30){\makebox(0,0)[lb]{\smash{{\SetFigFont{12}{14.4}{\familydefault}{\mddefault}{\updefault}{\color[rgb]{.82,0,0}\red{$Y$}}%
}}}}
\put(5296,-158){\makebox(0,0)[lb]{\smash{{\SetFigFont{12}{14.4}{\familydefault}{\mddefault}{\updefault}{\color[rgb]{0,0,.69}\blue{$X$}}%
}}}}
\put(1212,-213){\makebox(0,0)[lb]{\smash{{\SetFigFont{12}{14.4}{\familydefault}{\mddefault}{\updefault}{\color[rgb]{0,.82,0}\green{$\C$}}%
}}}}
\put(4464,-598){\makebox(0,0)[lb]{\smash{{\SetFigFont{10}{12.0}{\familydefault}{\mddefault}{\updefault}{\color[rgb]{0,.82,0}\green{$f(\C)$}}%
}}}}
\put(2330,  5){\makebox(0,0)[lb]{\smash{{\SetFigFont{12}{14.4}{\familydefault}{\mddefault}{\updefault}{\color[rgb]{0,0,0}$f$}%
}}}}
\end{picture}%

\end{center}

According to a strategy of thought going back to Bloch, modernized by
Green-Griffiths and viewed in a new light by Siu, the `first half' of
this conjecture\,\,---\,\,so to say\,\,---\,\,consists in showing that
there exist some nonzero global algebraic jet differentials that must
be satisfied by every nonconstant entire holomorphic curve $f \colon
\C \to X$, {\em see}~\cite{ siu2004, mer2009, dmr2010} for aspects of
the `second half', not at all considered here. Furthermore, for a
systematic development of the Kobayashi hyperbolicity theory, the
reader is referred to Kobayashi's foundational book~\cite{ Kob1998},
and also beyond, to a recent monograph~\cite{NogWink2014} by Noguchi
and Winkelmann, both published in the {\em Grundlehren} Springer
series.

The principal theorem of this memoir is presented
specifically in the case where $X = X^n \subset \P^{ n+1} ( \C)$ is a
(geometrically smooth) $n$-dimensional {\em hypersurface}, because the
main mathematical difficulty is essentially to reach arbitrary
dimensions $n \geqslant 1$, as was shown recently by the complexity of
some of the formal computations sketched in~\cite{ mer2008b, dmr2010}
for the case of dimension $n = 4$. But because a substantial part of
our proof relies upon works of Br\"uckmann (\cite{ bru1972, bru1997,
brurac1990}) which hold in fact for complete intersections, it is very
likely that our results may be transferred to such a more general
context. Also, one could consider entire holomorphic maps $\C^p \to
X^n$ having maximal generic rank $p$ with for
any fixed $1 \leqslant p \leqslant n$, as did Pacienza and Rousseau
(\cite{ paro2008}) recently for $p = 2$ in the case of $X^3 \subset
\P^4 ( \C)$. Furthermore, we hope more generally that the techniques
developed here could in the future
enable us to handle any $X^n \subset \P^N (
\C)$ of general type having arbitrary codimension $N - n$, but
probably requiring more than just general type as a workable
assumption (\cite{ siu2009}). 
At least in codimension $1$, we are able to gain the
following optimal result toward the Green-Griffiths conjecture.

\begin{MainTheorem}
Let $X = X^n \subset \P^{ n+1} ( \C)$ be a geometrically smooth
$n$-dimensional projective algebraic complex hypersurface. If $X$ is
of general type, namely if its degree $d$ satisfies the optimal
lower bound:
\[
d\geqslant n+3, 
\]
then there exist global algebraic differential equations on $X$
that must be identically satisfied by every nonconstant entire
holomorphic curve $f \colon \C \to X$.

More precisely, if $\mathcal{ E}_{ \kappa, m}^{ GG} T_X^*$
denotes the bundle of Green-Griffiths jet
polynomials\footnote{\,
See Sections~2 and~3 for exact definitions. 
} 
of order $\kappa$ and of weight $m$ over $X$, then the following holds
true.

\smallskip
{\small\sf Firstly:} for the fixed ample line bundle $\mathcal{ A}
:= \mathcal{ O}_X (
1)$, one has:
\begin{equation}
\label{minoration-h-0}
\aligned
&
h^0
\big(X,\,\mathcal{E}_{\kappa,m}^{GG}T_X^*
\otimes 
\mathcal{A}^{-1}\big)
\geqslant
\\
&
\geqslant
\frac{m^{(\kappa+1)n-1}}{
(\kappa!)^n\,((\kappa+1)n-1)!}
\bigg\{
\frac{(\log\kappa)^n}{n!}
d(d-n-2)^n
-
{\sf Constant}_{n,d}\cdot(\log\kappa)^{n-1}
\bigg\}
\,-
\\
&
\ \ \ \ \ \ \ \ \ \ \ \ \ \ \ \ \ \ \ \ \ \ \ \ \ \ \ \ \ \ \ \ \ \ \
\ \ \ \ \ \ \ \ \ 
-\,
{\sf Constant}_{n,d,\kappa}\cdot
m^{(\kappa+1)n-2},
\endaligned
\end{equation}
and the right-hand side minorant visibly tends to $\infty$, as
soon as both $\kappa \geqslant \kappa_{ n, d}^0$ and $m \geqslant m_{
n, d, \kappa}^0$ are large enough.

\smallskip
{\small\sf Secondly:} 
If $P$ is any global section of $\mathcal{ E}_{ \kappa, m}^{ GG}
T_X^* \otimes \mathcal{ A}^{ -1}$, hence which vanishes on the ample
divisor associated to $\mathcal{ A}$, 
then every nonconstant entire holomorphic
curve $f \colon \C \to X$ must satisfy the corresponding algebraic
differential equation $P ( j^\kappa f) = 0$. 
\end{MainTheorem}

Since the late 1990's, after fundamental works of Bloch,
Green-Griffiths and Siu, the so-called Ahlfors-Schwarz lemma for entire
holomorphic curves was clarified in full generality, and the second
statement above is nowadays (well) known to be a consequence of the
first ({\em see} e.g. Section~7 in~\cite{ dem1997}).

The case $n = 2$ of this theorem dates back to Green-Griffiths 1979
(\cite{ gg1980}). In~\cite{ rou2006b}, Rousseau was the first
to study effective (Demailly-Semple) 
jet differentials in dimension 3, under
the conditions $d \geqslant 97$. In~\cite{ div2008}, Diverio
treated the next dimensions $n = 4$ and $n = 5$ (improving also $n =
3$ with $d \geqslant 74$), under the conditions $d \geqslant 298$ and
$d \geqslant 1222$. In~\cite{ mer2008b}, the author of the present
article improved for $n = 4$ the lower bound to $d \geqslant 259$.
In~\cite{ div2009}, 
Diverio showed the (noneffective) existence of a lower bound
degree $d_n$ such that $d \geqslant d_n$ insures existence of nonzero
global jet differentials. An effective $d_n$ was captured in~\cite{
dmr2010} ({\em see} $\widetilde{ d}_n^1$, p.~192):
\[
d\,
\geqslant\,
2^{n^4}\,
n^{4n^3}\,
3^{n^3}\,
n^{3n^2}\,
(n+1)^{n^2+1}\,
n^{2n}\,
12,
\]
far from the optimal $n + 3$. Afterwards, with an improved
approach based on equivariant cohomology, B\'erczi~\cite{ Berc2010,
Berc2012} was able to lower the bound to $d \geqslant n^{ 8 n} =
2^{ 8n \log_2 n}$, Demailly~\cite{dem2011} to $d \geqslant \frac{
n^4}{3}\, \big(n\,{\sf log}\, \big(n\, \big(\log(24\,n) \big)
\big)^n$ also using Green-Griffiths jets, and Darondeau~\cite{ Da2013,
Da2014} to $d \geqslant 5n^2\, n^n$
by exploring deeper the B\'erczi-Diverio-Merker-Rousseau approach.
Furthermore and notably, 
B\'erczi in the second part of~\cite{ Berc2010} ({\em see} also
Theorem~1.3 there),
by introducing a new compactification of the Demailly-Semple
invariant jet bundle inspired from a previous deep work by
B\'erczi-Szenes~\cite{BercSzen2012},
showed algebraic degeneracy
of entire holomorphic curves valued in a generic
projective hypersurface
$X^n \subset \P^{n+1}(\C)$ of degree $d \geqslant n^6$
under the assumption that a very plausible conjecture due to Rim\'anyi 
concerning the positivity of the coefficients
of the Thom polynomial of Morin singularities holds, and
under a certain assumption on the growth of its coefficients,
such an improvement $d \geqslant n^6$ 
on the degree bound being spectacular.
All the mentioned works focus on jet differentials of
order $\kappa = n$ equal to the dimension.

This Main Theorem above was presented in conference talks given in the CIRM
(June 2009), in the Hong-Kong University (August 2009) and also later
in some seminars (Paris, Marseille, Lyon), and appeared in May 2010 as
the preprint {\scriptsize\sf arxiv.org/abs/1005.0405/}.  In November
2010, Demailly ({\scriptsize\sf arxiv.org/abs/1011.3636/},
\cite{Demailly-2011}) was able to extend this result to any projective
manifold of general type, not necessarily being a hypersurface or a
complete intersection, also coming back to plain Green-Griffiths jets,
but developing completely different elaborate negative jet curvature
estimates which was inspired from an article of Cowen and Griffiths
(\cite{ Cowen-Griffiths-1976}).

It is certainly advisable to present the principal cornerstones of the
extended proof before entering its beautiful
core. Let therefore $X^n \subset \P^{ n+1} ( \C)$ be a complex
projective hypersurface of general type, that is to say having degree
$d \geqslant n + 3$, and let:
\[
\mathcal{E}_{\kappa,m}^{GG}T_X^*
\,\longrightarrow\,X
\]
denote the holomorphic vector bundle\,\,---\,\,over $X$\,\,---\,\,of
homogeneous polynomials of degree $m \geqslant 1$ in the jets of order
$\kappa \geqslant 1$ of (local) holomorphic maps $\D \to X$, where $\D
\subset \C$ is the unit disc.  Asymptotically, its Euler-Poincar\'e
characteristic is known, thanks to the seminal article~\cite{ gg1980}
by Green and Griffiths, to tend to infinity as:
\[
\aligned
&
\chi\big(
X,\,\mathcal{E}_{\kappa,m}^{GG}T_X^*
\big)
=
\\
&
=
\frac{m^{(\kappa+1)n-1}}{(\kappa!)^n\,
((\kappa+1)n-1)!}\,
\Big\{
{\textstyle{\frac{(\log\kappa)^n}{n!}}}
d(d-n-2)^n
-
{\sf Constant}_{n,d}\cdot(\log\kappa)^{n-1}
\Big\}
-
\\
&
\ \ \ \ \ \ \ \ \ \ \ \ \ \ \ \ \ \ \ \ \ \ \ \ \ \ \ \ \ \ \ \ \ \ \
\ \ \ \ \ \ \ \ \ 
-
{\sf Constant}_{n,d,\kappa}\cdot
m^{(\kappa+1)n-2},
\endaligned
\]
when $m \gg \kappa \gg 1$ both tend to $\infty$. 
Section~3 is devoted to reprove this formula in great details.
Moreover, it is known that when one 'twists', {\em i.e.} when one
{\em tensors} the Green-Griffiths jet bundle:
\[
\mathcal{ E}_{\kappa, m}^{GG} T_X^*
\otimes
\mathcal{A}^{-1}
\]
with a certain fixed ample line bundle $\mathcal{ A} \longrightarrow
X$ of the form $\mathcal{ A} = \mathcal{ O}_X (t)$ for a fixed integer
$t \geqslant 1$, then the asymptotic behavior of the Euler-Poincar\'e
characteristic remains exactly the same, for $\mathcal{ A}$ is in some
sense `submersed'. Hence the Main Theorem above states that the (most
interesting) vector space:
\[
H^0\big(X,\mathcal{E}_{\kappa,m}^{GG}T_X^*\big)
\]
of global (algebraic) holomorphic sections of this jet bundle has a
dimension which tends to $\infty$ in essentially the same asymptotic
way as its (much easier to estimate) characteristic.

In Section~2, we present in a self-contained elementary way all basic
properties of this jet bundle $\mathcal{ E}_{ \kappa, m}^{ GG}
T_X^*$. Specifically, we exhibit in great details the (known, \cite{
gg1980, dem1997}) natural filtration that $\mathcal{ E}_{ \kappa,
m}^{ GG} T_X^*$ possesses, {\em i.e.} a certain coordinate-invariant
nested sequence of holomorphic vector subbundles contained in it, and
we reprove there that the graded bundle ${\sf Gr}^\bullet \mathcal{ E}_{
\kappa, m}^{ GG} T_X^*$ associated to this filtration\,\,---\,\,which 
by definition consists of successive quotients
of sequence of nested subbundles in question\,\,---\,\,writes out:
\[
{\sf Gr}^\bullet
\mathcal{E}_{\kappa,m}^{GG}T_X^*
=
\bigoplus_{\ell_1+2\ell_2+\cdots+\kappa\ell_\kappa=m}\,
{\rm Sym}^{\ell_1}T_X^*
\otimes
{\rm Sym}^{\ell_2}T_X^*
\otimes
\cdots
\otimes
{\rm Sym}^{\ell_\kappa}T_X^*,
\]
where, according to a standard definition (\cite{ botu1982, ful1998,
fuha1991, macd1995}), for any integer $\ell \geqslant 1$, the {\sl
$\ell$-th symmetric tensor power} ${\rm Sym}^\ell T_X^*$ of the
cotangent bundle $T_X^*$ has, in local coordinates $(x_1, \dots, x_n)$
on $X$, a basis consisting of all:
\[
dx_{i_1}\odot
\cdots
\odot
dx_{i_\ell}
\ \ \ \ \ \ \ \ \ \ \ \ \ {\scriptstyle{(1\,\leqslant\,i_1\,
\leqslant\,\cdots\,\leqslant\,i_\ell\,\leqslant\,n)}},
\]
the symbol $\odot$ for the symmetric product being commutative,
contrary, of course, to the tensor product $\otimes$.

In the same Section~2, we also reprove the elementary fact that the Euler
Poincar\'e characteristic is unchanged:
\[
\chi\big(
X,\,\mathcal{E}_{\kappa,m}^{GG}T_X^*
\big)
=
\chi\big(
X,\,{\sf Gr}^\bullet\mathcal{E}_{\kappa,m}^{GG}T_X^*
\big),
\]
plus another quite crucial known fact\,\,---\,\,already used by Rousseau
in dimension 3 (\cite{ rou2006b, rou2007a})\,\,---\,\,according
to which positive
cohomology dimensions enjoy the agreeable majorations:
\[
\dim
H^q\big(
X,\,\mathcal{E}_{\kappa,m}^{GG}T_X^*
\big)
\,\leqslant\,
\dim\,H^q\big(
X,\,{\sf Gr}^\bullet\mathcal{E}_{\kappa,m}^{GG}T_X^*
\big)
\ \ \ \ \ \ \ \ \ \ \ \ \ {\scriptstyle{(q\,=\,1\,\cdots\,n)}},
\]
whence it instantly follows:
\[
\aligned
&
\dim
H^q\big(
X,\,\mathcal{E}_{\kappa,m}^{GG}T_X^*
\big)
\,\leqslant\,
\\
&
\,\leqslant\,
\sum_{\ell_1+2\ell_2+\cdots+\kappa\ell_\kappa=m}\,
\dim
H^q\big(X,\,
{\rm Sym}^{\ell_1}T_X^*
\otimes
{\rm Sym}^{\ell_2}T_X^*
\otimes
\cdots
\otimes
{\rm Sym}^{\ell_\kappa}T_X^*
\big),
\endaligned
\]
again for $q = 1, \dots, n$.

In this memoir, we develop (somewhat considerably) a strategy
successfully applied by Rousseau in dimension 3\,\,---\,\,also in a
logarithmic context, \cite{ rou2007b}, not covered here\,\,---, 
which consists in
majorating the right-hand side coholomogy dimensions $\dim H^q$, $q =
1, \dots, n$, by quantities that do not perturb too much the
explicitly known asymptotic positivity of the characteristic.  More
precisely, reminding that, for any holomorphic vector bundle $E
\longrightarrow X$, its Euler-Poincar\'e characteristic is the
alternating sum of its cohomology dimensions:
\[
\small
\aligned
\chi(X,E)
&
=
\dim H^0(X,E)
-
\dim H^1(X,E)
+
\dim H^2(X,E)
-
\dim H^3(X,E)
\,+
\\
&
\ \ \ \ \ 
+
\dim H^4(X,E)
-
-\cdots+
(-1)^n\,\dim H^n(X,E),
\endaligned
\]
a simple formula from which one trivially deduces the minoration:
\[
\footnotesize
\aligned
\dim H^0\big(X,\,
\mathcal{E}_{\kappa,m}^{GG}T_X^*\big)
\geqslant
\chi\big(X,\,
\mathcal{E}_{\kappa,m}^{GG}T_X^*\big)
\!-\!
\dim H^2\big(X,\,
\mathcal{E}_{\kappa,m}^{GG}T_X^*\big)
\!-\!
\dim H^4\big(X,\,
\mathcal{E}_{\kappa,m}^{GG}T_X^*\big)
\!-\cdots,
\endaligned
\]
the right-hand side subtracted terms being only (positive) {\em
even-dimensional} cohomology dimensions, if one is able, by means of
some argument, to show that the even sum of these dimensions : 
\[
\dim H^2 
+ 
\dim H^4 
+ 
\dim H^6
+
\cdots
\]
is somewhat smaller than $\chi$, when 
$m \gg \kappa \gg 1$ both tend to $\infty$, 
then one will deduce that the interesting space: 
\[
H^0
\big(X,\,\mathcal{E}_{\kappa,m}^{GG}T_X^*\big)
\]
of global holomorphic jet differentials on $X$ is of positive
dimension, a dimension which shall in fact also tends to $\infty$
when $m \gg \kappa \gg 1$.

As explained, by passing to the graded bundle ${\sf Gr}^\bullet
\mathcal{ E}_{\kappa, m}^{ GG} T_X^*$, one is led to estimate the
cohomology dimensions of the above sum of multi-tensored powers of the
symmetric powers of $T_X^*$.  Rousseau's approach in dimension
$3$\,\,---\,\,already suggested in Demailly's seminal memoir~\cite{
dem1997}\,\,---\,\,was in fact applied to the {\em sub}bundle:
\[
\mathcal{E}_{\kappa,m}^{DS}T_X^*
\,\subsetneqq\,
\mathcal{E}_{\kappa,m}^{GG}T_X^*
\]
of jet differentials that are invariant under reparametrization \footnote{\,
for results in dimension 4, {\em see}~\cite{ mer2008b}, for a survey
in dimensions 2, 3, 4 focused on invariant jets, see~\cite{ dmr2010,
Diverio-Rousseau-2011}, and for a
brief presentation, see Section~13, 
}, 
but the main feature, which concerns both
jet bundles and would apply
to any other holomorphic jet bundle as well, 
consists in {\em decomposing at first} this sum of
multi-tensored powers of symmetric powers of $T_X^*$ explicitly as a
certain direct sum of so-called {\sl Schur bundles}:
\[
\mathcal{S}^{(\ell_1,\dots,\ell_n)}T_X^*
\,\longrightarrow\,
X,
\]
about which we now make a brief reminder.

On a geometrically 
smooth projective algebraic hypersurface 
$X^n \subset \P^{ n+1} ( \C)$,  one classically studies a few
holomorphic vector bundles:

\medskip$\square$\ \ \
$T_X^*$;

\medskip$\square$\ \ \
$\Lambda^k T_X^*$ (Hodge theory); 

\medskip$\square$\ \ \
$K_X := \Lambda^n T_X^*$ canonical bundle;

\medskip$\square$\ \ \
$K_X^{\otimes m}$ its tensor powers (plurigenera); 

\medskip$\square$\ \ \
${\sf Sym}^k\, T_X^*$ (cotangential $k$-genus).

\medskip\noindent
All these are particular instances of the mentioned {\sl Schur bundles}:
\[
\mathcal{S}^{(\ell_1,\ell_2,\dots,\ell_n)}T_X^*
\,\longrightarrow\,
X,
\]
that are parametrized by decreasing sequences of nonnegative integers:
\[
\ell_1\geqslant\ell_2\geqslant\cdots\geqslant\ell_n\geqslant 0,
\]
and one recovers notably:
\[
\aligned
\Lambda^kT_X^*
&
=
\mathcal{S}^{(1,\dots,1,0,\dots,0)}T_X^*
\ \ \ \ \ \ \ 
\text{\rm with}\ \ k \ \ \text{\rm times}\ \ 1;
\\
{\rm Sym}^kT_X^*
&
=
\mathcal{S}^{(k,0,\dots,0)}T_X^*.
\endaligned
\]

The main definitional feature of these bundles is that they 
appear when one
decomposes in irreducible representations of ${\sf GL}_n (
\C)$ any $r$-th tensor power of the cotangent
bundle:
\[
\aligned
\underbrace{
T_X^*\otimes\cdots\otimes T_X^*}_{r\,\,\text{\rm times}}
=
\bigoplus_{(\ell)}\,
\Big[
\mathcal{S}^{(\ell_1,\dots,\ell_n)}T_X^*
\Big]^{\oplus N_{(\ell)}},
\endaligned
\]
with $\ell_1 \geqslant \cdots \geqslant \ell_n \geqslant 0$ where
$N(\ell) \in \N$ is a certain multiplicity. According to classical
representation theory (dating to the end of the
XIX\textsuperscript{th} Century), every representation (action) of
${\sf GL}_n ( \C)$ can be written as a certain direct sum of Schur
representations, which constitute the list of all possible {\em
irreducible} representations of ${\sf GL}_n(\C)$. Such splittings
pass in a natural way (one has to check) through ${\sf GL}_n(\C)$-valued
changes of
trivializations, whence a general algebraic fact is gently offered to
global complex geometry: Every holomorphic vector bundle $E$ over $X$,
on the fibers of which one can let ${\sf GL}_n ( \C)$ act, must in
principle decompose itself as a certain direct sum of Schur bundles,
which happen to be the elementary bricks with which one can
reconstitute any vector bundle in the so-called Grothendieck ring.

As mentioned above, in dimensions $2$ (Demailly) and $3$ (Rousseau),
for jet order $\kappa = 2$ and $\kappa = 3$ equal to the dimension, it
is not hard to obtain such a decomposition for the graded bundle
associated to the bundles of jets invariant under reparametrization:
\[
\aligned
{\sf Gr}^{\bullet}
E_{2,m}^{DS}T_X^*
&
=
\bigoplus_{a+3b=m}
\mathcal{S}^{(a+b,\,\,b)}T_X^*,
\\
{\sf Gr}^{\bullet}
E_{3,m}^{DS}T_X^{\ast}
&
=
\bigoplus_{a+3b+5c+6d=m}\mathcal{S}^{
(a+b+2c+d,\,\,b+c+d,\,\,d)}T_X^{\ast }.
\endaligned
\]
In dimension $4$, also for jet order $\kappa = 4$ equal to
the dimension, the author has obtained a more complex
decomposition (just below, 
the 41 subsets $\square_i$, $i = 1, 2, \dots, 41$, of
$\mathbb{N}^{ 14} \ni (a, b, \dots, l, m', n)$ are explicitly defined
in~\S12 of~\cite{mer2008b}):
\[
\aligned
\operatorname{Gr}^\bullet
&
E_{4,m}T_X^*
=
\bigoplus_{(a,b,c,d,e,f,g,h,i,j,k,l,m',n)\in\mathbb{N}^{14}\backslash
(\square_1\cup\cdots\cup\square_{41}) 
\atop
o+3a+5b+7c+6d+8e+10f+8g+10h+12i+14j+15k+17l+19m'+21n+10p\,=\,m}\,
\\
&
{\scriptsize
\mathcal{S}
\left(
\aligned
o+a+2b+3c+d+2e+3f+2g+2h+3i+4j+3k+3l+4m'+5n+p&
\\
a+b+c+d+e+f+2g+2h+2i+2j+2k+3l+3m'+3n+p&
\\
d+e+f+h+i+j+2k+2l+2m'+2n+p&
\\
p&
\endaligned
\right)}\,T_X^*,
\endaligned,
\]
the combinatorics of which is not easily devisable, and in higher
dimensions, no published effective result exists. There are some known
deep reasons why regularities in such objects may fail to be
discovered, or just, to exist, and our last Section~13 is devoted to
provide some evidence towards the necessity of abandoning the study of
the graded bundles of jets invariant under reparametrization:
\[
{\sf Gr}^\bullet\mathcal{E}_{\kappa,m}^{DS}T_X^*.
\]

In fact, the causal origin of the present memoir, dating back to
September 2008, lies in the unavoidable necessity of developing the
theory of jet bundles principally with the ancient, plain
Green-Griffiths jets.  In the published domain, concrete studies of
${\sf Gr}^\bullet \mathcal{ E}_{\kappa, m}^{ GG}T_X^*$ in dimensions
$3$, $4$, or even $5$, are completely missing, due to the fact that
experts (including the author) believed that 
the additional requirement of invariancy for jets would naturally
bring better positivity properties, as
was known in dimension $2$ without any
consciousness of the complexity of algebraic
invariant theory (again, {\em cf.}~Section~13
below). 

Undoubtedly, the algebraic complexity of the graded bundle 
${\sf Gr}^\bullet\mathcal{E}_{\kappa,m}^{DS}T_X^*$ of
invariant jets, already unwiedly for $n = \kappa = 4$
(\cite{mer2008b}),
forced us to stop in this direction, 

Then quite unexpectedly, and very serendipitously also, we realized in
the autumn 2008 that a decomposition of ${\sf Gr}^\bullet \mathcal{
E}_{\kappa, m}^{ GG}T_X^*$ is available not only in small dimensions
and for small jet order, but also in {\em arbitrary} dimension
and for jets of {\em any} order. Section~4 ends up with the
completion of the proof of a fundamental statement, the ingredients of
which were already discovered at the turn from the
XIX\textsuperscript{th} to the XX\textsuperscript{th} Century.

\begin{Theorem}
The graded vector bundle ${\sf Gr}^\bullet \mathcal{ E}_{ \kappa, m}^{
GG} T_X^*$ associated to the bundle $\mathcal{ E}_{ \kappa, m}^{ GG}
T_X^*$ of $\kappa$-th $m$-weighted Green-Griffiths jets identifies
to the following {\em exact direct sum} of Schur bundles:
\[
{\sf Gr}^\bullet\mathcal{E}_{\kappa,m}^{GG}T_X^*
=
\bigoplus_{\ell_1\geqslant\ell_2\geqslant\cdots\geqslant\ell_n
\geqslant 0}
\Big(
\mathcal{S}^{(\ell_1,\ell_2,\dots,\ell_n)}T_X^*
\Big)^{\oplus 
M_{\ell_1,\ell_2,\dots,\ell_n}^{\kappa,m}}, 
\]
with multiplicities $M_{\ell_1, \ell_2, \dots, \ell_n }^{ \kappa,
m} \in \N$ equal to the number of times a Young diagram ${\sf YD}_{
(\ell_1, \dots, \ell_n)}$ with row lengths equal to $\ell_1, \ell_2,
\dots, \ell_n$ can be filled in with positive integers $\lambda_i^j
\leqslant \kappa$ placed at its $i$-th row and $j$-th column so as to
constitute a semi-standard tableau:

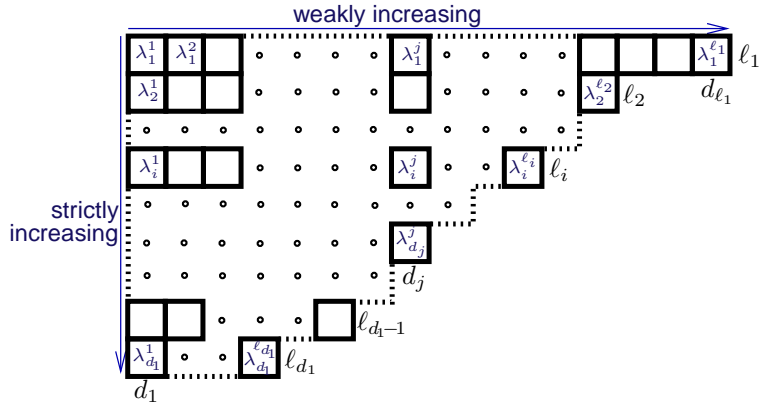
\begin{center}
\begin{picture}(0,0)%
\includegraphics{d1-fill-young-diagram.pstex}%
\end{picture}%
\setlength{\unitlength}{4144sp}%
\begingroup\makeatletter\ifx\SetFigFont\undefined%
\gdef\SetFigFont#1#2#3#4#5{%
  \reset@font\fontsize{#1}{#2pt}%
  \fontfamily{#3}\fontseries{#4}\fontshape{#5}%
  \selectfont}%
\fi\endgroup%
\begin{picture}(4364,2408)(-80,-1874)
\put(1619,426){\makebox(0,0)[lb]{\smash{{\SetFigFont{9}{10.8}{\familydefault}{\mddefault}{\updefault}{\color[rgb]{0,0,.69}\blue{\sf weakly increasing}}%
}}}}
\put(-80,-882){\makebox(0,0)[lb]{\smash{{\SetFigFont{9}{10.8}{\familydefault}{\mddefault}{\updefault}{\color[rgb]{0,0,.69}\blue{\sf increasing}}%
}}}}
\put(170,-743){\makebox(0,0)[lb]{\smash{{\SetFigFont{9}{10.8}{\familydefault}{\mddefault}{\updefault}{\color[rgb]{0,0,.69}\blue{\sf strictly}}%
}}}}
\put(685,-13){\makebox(0,0)[lb]{\smash{{\SetFigFont{7}{8.4}{\familydefault}{\mddefault}{\updefault}{\color[rgb]{0,0,.69}\blue{$\lambda_2^1$}}%
}}}}
\put(694,-473){\makebox(0,0)[lb]{\smash{{\SetFigFont{7}{8.4}{\familydefault}{\mddefault}{\updefault}{\color[rgb]{0,0,.69}\blue{$\lambda_i^1$}}%
}}}}
\put(913,209){\makebox(0,0)[lb]{\smash{{\SetFigFont{7}{8.4}{\familydefault}{\mddefault}{\updefault}{\color[rgb]{0,0,.69}\blue{$\lambda_1^2$}}%
}}}}
\put(673,209){\makebox(0,0)[lb]{\smash{{\SetFigFont{7}{8.4}{\familydefault}{\mddefault}{\updefault}{\color[rgb]{0,0,.69}\blue{$\lambda_1^1$}}%
}}}}
\put(2263,201){\makebox(0,0)[lb]{\smash{{\SetFigFont{7}{8.4}{\familydefault}{\mddefault}{\updefault}{\color[rgb]{0,0,.69}\blue{$\lambda_1^j$}}%
}}}}
\put(2242,-474){\makebox(0,0)[lb]{\smash{{\SetFigFont{7}{8.4}{\familydefault}{\mddefault}{\updefault}{\color[rgb]{0,0,.69}\blue{$\lambda_i^j$}}%
}}}}
\put(4033,199){\makebox(0,0)[lb]{\smash{{\SetFigFont{7}{8.4}{\familydefault}{\mddefault}{\updefault}{\color[rgb]{0,0,.69}\blue{$\lambda_1^{\ell_1}$}}%
}}}}
\put(3354,-32){\makebox(0,0)[lb]{\smash{{\SetFigFont{7}{8.4}{\familydefault}{\mddefault}{\updefault}{\color[rgb]{0,0,.69}\blue{$\lambda_2^{\ell_2}$}}%
}}}}
\put(2230,-904){\makebox(0,0)[lb]{\smash{{\SetFigFont{7}{8.4}{\familydefault}{\mddefault}{\updefault}{\color[rgb]{0,0,.69}\blue{$\lambda_{d_j}^j$}}%
}}}}
\put(2906,-477){\makebox(0,0)[lb]{\smash{{\SetFigFont{7}{8.4}{\familydefault}{\mddefault}{\updefault}{\color[rgb]{0,0,.69}\blue{$\lambda_i^{\ell_i}$}}%
}}}}
\put(659,-1600){\makebox(0,0)[lb]{\smash{{\SetFigFont{7}{8.4}{\familydefault}{\mddefault}{\updefault}{\color[rgb]{0,0,.69}\blue{$\lambda_{d_{\!1}}^1$}}%
}}}}
\put(1321,-1624){\makebox(0,0)[lb]{\smash{{\SetFigFont{7}{8.4}{\familydefault}{\mddefault}{\updefault}{\color[rgb]{0,0,.69}\blue{$\lambda_{d_{\!1}}^{\!\ell_{d_{\!1}}}$}}%
}}}}
\put(4057,-24){\makebox(0,0)[lb]{\smash{{\SetFigFont{10}{12.0}{\familydefault}{\mddefault}{\updefault}{\color[rgb]{0,0,0}$d_{\ell_1}$}%
}}}}
\put(3594,-56){\makebox(0,0)[lb]{\smash{{\SetFigFont{10}{12.0}{\familydefault}{\mddefault}{\updefault}{\color[rgb]{0,0,0}$\ell_2$}%
}}}}
\put(4284,183){\makebox(0,0)[lb]{\smash{{\SetFigFont{10}{12.0}{\familydefault}{\mddefault}{\updefault}{\color[rgb]{0,0,0}$\ell_1$}%
}}}}
\put(3151,-499){\makebox(0,0)[lb]{\smash{{\SetFigFont{10}{12.0}{\familydefault}{\mddefault}{\updefault}{\color[rgb]{0,0,0}$\ell_i$}%
}}}}
\put(666,-1825){\makebox(0,0)[lb]{\smash{{\SetFigFont{10}{12.0}{\familydefault}{\mddefault}{\updefault}{\color[rgb]{0,0,0}$d_1$}%
}}}}
\put(1568,-1645){\makebox(0,0)[lb]{\smash{{\SetFigFont{10}{12.0}{\familydefault}{\mddefault}{\updefault}{\color[rgb]{0,0,0}$\ell_{d_1}$}%
}}}}
\put(2001,-1418){\makebox(0,0)[lb]{\smash{{\SetFigFont{10}{12.0}{\familydefault}{\mddefault}{\updefault}{\color[rgb]{0,0,0}$\ell_{d_{\!1}\!-\!1}$}%
}}}}
\put(2269,-1140){\makebox(0,0)[lb]{\smash{{\SetFigFont{10}{12.0}{\familydefault}{\mddefault}{\updefault}{\color[rgb]{0,0,0}$d_j$}%
}}}}
\end{picture}%

\end{center}

\noindent
namely satisfying 
horizontal weak increase and vertical strict increase, 
with the further constraint that
the sum of all such integers:
\[
\aligned
m
&
=
\lambda_1^1+\cdots+\lambda_1^{\ell_n}
+\cdots+
\lambda_1^{\ell_2}+\cdots+\lambda_1^{\ell_1}
\\
&
+
\lambda_2^1+\cdots+\lambda_2^{\ell_n}
+\cdots+\lambda_2^{\ell_2}
\\
&
+\cdots\cdot\cdots\cdots\cdots\cdots+
\\
&
+
\lambda_n^1+\cdots+\lambda_n^{\ell_n}
\endaligned
\]
equals the prescribed weighted homogeneity degree $m$.
\end{Theorem}

Although these multiplicities $M_{\ell_1, \ell_2, \dots, \ell_n}^{
\kappa, m}$ of the appearing Schur bundles are described in an
apparently satisfactory combinatorial manner, the truth is
that, yet, their asymptotic behavior is not at all visible
in such a partially understood description, and
our Sections~6, 7, 8, 9, 10 will be devoted to develop
a presumably completely new study of asymptotic Young
diagrams which could in a near future have applications
to Probability Theory.

But before we present some of its aspects, setting temporarily aside
these multiplicities $M_{\ell_1, \ell_2, \dots, \ell_n}^{ \kappa, m}$,
we must as we know understand the cohomologies of the basic
irreducible bricks that are the Schur bundles $\mathcal{ S}^{ (\ell_1,
\dots, \ell_n)} T_X^*$. In Rousseau's paper~\cite{rou2006b} for
$X^3 \subset \P^4$ and in the joint paper~\cite{ dmr2010} for $X^4
\subset \P^5$, the following majorations of 
2\textsuperscript{nd} cohomology
dimensions were obtained:
\[
\footnotesize
\aligned
h^2\big(X,\mathcal{S}^{(\ell_{1},\ell_{2},\ell_{3})}T_X^*\big)
&
\,\leqslant\, 
d(d+13)\frac{3(\ell_{1}+\ell_{2}+\ell_{3})^3}{2}
(\ell_{1}-\ell_{2})(\ell_{1}-\ell_{3})
(\ell_{2}-\ell_{3})
+
{\rm O}(\vert\ell\vert^5),
\\
h^2\big(X,\,
\mathcal{S}^{(\ell_{1},\ell_{2},\ell_{3},\ell_{4})}T_X^*\big)
&
\,\leqslant\,
\frac{1}{80}\,d\,
(\ell_1-\ell_2)(\ell_1-\ell_3)(\ell_1-\ell_4)
(\ell_2-\ell_3)(\ell_2-\ell_4)(\ell_3-\ell_4)
\cdot
\\
&
\cdot
\big(\ell_1+\ell_2+\ell_3+\ell_4\big)^2
\Big[5\ell_2 \ell_1d^2+132\ell_2 \ell_1d
+132\ell_1 \ell_3d+5\ell_2 \ell_3d^2
\\
&
+132\ell_ 2 \ell_4 d+5\ell_2 d^2\ell_4
+132\ell_1 \ell_4 d+5\ell_3 \ell_4 d^2
+5\ell_1\ell_3d^2
\\
&
+132\ell_3 \ell_4d+132\ell_2 \ell_3 d
+1308\ell_2\ell_1+648\ell_2^2+648\ell_3^2
\\
&
+72\ell_3^2d+648\ell_1^2+72\ell_1^2d
+1308\ell_1\ell_4+5\ell_1d^2\ell_4
+1308\ell_2\ell_4
\\
&+1308\ell_2\ell_3+648\ell_4^2
+72\ell_2^2d+1308\ell_1\ell_3
+72\ell_4^2d+1308\ell_3\ell_4\Big]
\\
&
+
{\rm O}\big(\vert\ell\vert^9\big),
\endaligned
\]
respectively. In Section~5, we rather easily generalize such
majorations to arbitrary dimension.

\begin{Theorem}
Let $X = X^n \subset \P^{ n+1} ( \C)$ be a geometrically smooth
projective algebraic complex hypersurface of general type, i.e. of
degree $d \geqslant n+3$, and let $\ell = ( \ell_1, \dots, \ell_{ n-1},
\ell_n)$ with $\ell_1 \geqslant \cdots \geqslant \ell_{ n-1} \geqslant
\ell_n \geqslant 0$. If:
\[
\vert\ell\vert
=
\ell_1+\cdots+\ell_{n-1}+\ell_n
\geqslant
{\sf Constant}_{n,d},
\]
then for every $q = 1, 2, \dots, n$, the dimensions of the positive
cohomology groups of the Schur bundle $\mathcal{ S}^{ (\ell_1, 
\dots, \ell_{ n-1}, \ell_n)} T_X^*$ over $X$ satisfy a general
majoration of the form:
\[
\aligned
&
h^q\big(
X,\,\,
\mathcal{S}^{(\ell_1,\dots,\ell_{ n-1},\ell_n)}T_X^*
\big)
\leqslant
\\
&
\leqslant\,
{\sf Constant}_{n,d}\!\!\!
\prod_{1\leqslant i<j\leqslant n}
(\ell_i-\ell_j)\,
\bigg[
\sum_{\beta_1+\cdots+\beta_{n-1}+\beta_n=n}
\!\!\!
\ell_1^{\beta_1}\cdots\,\ell_{n-1}^{\beta_{n-1}}\ell_n^{\beta_n}
\bigg]
+
\\
&
\ \ \ \ \ \ \ \ \ \ \ \ \
\ \ \ \ \ \ \ \ \ \ \ \ \
\ \ \ \ \ \
+
{\sf Constant}_{n,d}\,
\bigg[
\sum_{\alpha_1+\dots+\alpha_n\leqslant\frac{n(n+1)}{2}-1}\,
\ell_1^{\alpha_1}\cdots\,\ell_n^{\alpha_n}
\bigg],
\endaligned
\]
with leading terms being homogeneous of degree $\frac{ n ( n+1)}{
2}$ with respect to the $\ell_i$ and divisible by all the differences
$(\ell_i - \ell_j)$, where $1 \leqslant i < j \leqslant n$.
\end{Theorem}

Of course, some appropriate values of the appearing constants can be
made explicit (such a task would be necessary only if one would desire
to make explicit a lower bound of $\kappa \geqslant \kappa_{n,d}^0$
insuring the conclusion of the Main Theorem).  However, the majoration
derived from a generalization of Rousseau's happens to be definitely
not sufficient for the purpose of reaching the Main Theorem, for one
realizes, with the tools developed in Section~10, that any monomial
which would, in the right-hand side majorant, appear under the form
(and there are some, in dimension $3$, in dimension $4$, in arbitrary
dimension $n$):
\[
\underbrace{{\sf Constant}}_{>\,0}
\cdot\,
\ell_n^n\,
\prod_{1\leqslant i<j\leqslant n}
(\ell_i-\ell_j)
\]
has the very annoying property that when one sums up: 
\[
\sum_{\ell_1\geqslant\cdots\geqslant\ell_n\geqslant 0}\,
M_{\ell_1,\dots,\ell_n}^{\kappa,m}\,
{\sf Constant}\,(\ell_n)^n\,
\prod_{1\leqslant i<j\leqslant n}
(\ell_i-\ell_j)
\]
taking
account of the multiplicities, then one obtains at the
end a term in $\log( \kappa )^n$:
\[
\frac{m^{(\kappa+1)n-1}}{(\kappa!)^n\,((\kappa+1)n-1)!}\,\,
\underbrace{\text{\footnotesize\sf New constant}}_{
\text{\sf again}\,>\,0}\,\,
\big(\log\kappa\big)^n,
\]
which unavoidably perturbs the leading term of the Euler-Poincar\'e
characteristic. Because such a majoration is not enough,
we devote Section~11 to a diagram chasing in long exact
cohomology sequences associated to short exact sequences
due to Br\"uckmann to obtain the following crucial cohomology
vanishing theorem.

\begin{Theorem}
Let $X = X^n \subset \P^{ n+1} ( \C)$ be a geometrically smooth
projective algebraic complex hypersurface of general type, i.e. of
degree $d \geqslant n+3$, and let $\ell = ( \ell_1, \dots, \ell_{ n-1},
\ell_n)$ with $\ell_1 \geqslant \cdots \geqslant \ell_{ n-1} \geqslant
\ell_n \geqslant 1$. If:
\[
\ell_n
\geqslant
{\textstyle{\frac{1}{d-n-2}}}
\big\{
n(d-1)
+
\ell_1-\ell_n
+
{\textstyle{\sum_{i=1}^{n-1}}}\,
(\ell_i-\ell_n)
\big\},
\]
then all the positive cohomologies vanish:
\[
0
=
H^q
\big(
X,\,
\mathcal{S}^{(\ell_1,\dots,\ell_{n-1},\ell_n)}T_X^*
\big)
\ \ \ \ \ \ \ \ \ \ \ \ \ 
{\scriptstyle{(q\,=\,1,\,\,2\,\cdots\,n)}}.
\]
\end{Theorem}

Performing the synthesis between these two theorems, we are able
to deduce the general majoration:
\begin{equation}
\label{prudent-majoration}
\aligned
&
\dim
H^q\big(X,\,\mathcal{S}^{(\ell_1,\dots,\ell_n)}T_X^*\big)
\leqslant\,
\\
&
\,\leqslant\,
{\sf Constant}_n\cdot
\big[1+d+\cdots+d^{n+1}\big]
\cdot
\prod_{1\leqslant i<j\leqslant n}\,
\big(\ell_i-\ell_j\big)\,\cdot
\\
&
\ \ \ \ \ \ \ \ \ \ \ \ \ \
\sum_{\beta_1+\cdots+\beta_{n-1}+\beta_n=n
\atop
\beta_n\leqslant n-1}
(\ell_1-\ell_2)^{\beta_1}\cdots
(\ell_{n-1}-\ell_n)^{\beta_{n-1}}
(\ell_n)^{\beta_n}
+
\\
&
\ \ \ \ \ \ \ \ \ \ \ \ \ \ \ \ \ \ \ \ \ \ \ \ \ 
\ \ \ \ \ \ \ \ \ \ \ \ \ \ \ \ \ \ \ \ \ \ \ \ \ 
+
{\sf Constant}_{n,d}\,
\big(
1+
\vert\ell\vert^{\frac{n(n+1)}{2}-1}
\big),
\endaligned
\end{equation}
where $\vert \ell \vert = \ell_1 + \cdots + \ell_n$, in which,
crucially, the exponent $\beta_n$ of $\ell_n$ is constrained to be
$\leqslant n-1$ so that after summation, no
perturbation in $( \log \kappa )^n$ occurs:
\[
\aligned
h^q\big(
X,\,\,
\mathcal{E}_{\kappa,m}^{GG}T_X^*
\big)
&
\leqslant
\sum_{\ell_1\geqslant\ell_2\geqslant\cdots\geqslant\ell_n\geqslant 0}\,
M_{\ell_1,\ell_2,\dots,\ell_n}^{\kappa,m}
\cdot
h^q\big(
X,\,\,
\mathcal{S}^{(\ell_1,\ell_2,\dots,\ell_n)}T_X^*
\big)
\\
&
\leqslant
{\textstyle{\frac{m^{(\kappa+1)n-1}}{(\kappa!)^n\,((\kappa+1)n-1)!}}}
\cdot
{\sf O}_{n,d}\big(\log(\kappa)^{n-1}\big)
+
{\sf O}_{n,d,\kappa}\big(
m^{(\kappa+1)n-2}
\big)
\\
&
\ \ \ \ \ \ \ \ \ \ \ \ \ \ \ \ \ \ \ \ \ \ \ \ \ \ \ \ \ \ \ 
{\scriptstyle{(q\,=\,1\,\cdots\,n)}}.
\endaligned
\]
Concerning remainders\footnote{\,
Throughout the paper, we shall sometimes write ${\sf O}_{n, d, \kappa}
\big(m^{ (\kappa+1) n - 2}\big)$ to denote a quantity which is
majorated by ${\sf Constant}_{ n, d, \kappa} \cdot m^{ (\kappa+1)n - 2
}$.
}, 
we indeed prove:

\begin{Theorem}
If $\alpha_1, \dots, \alpha_n$ are any nonnegative integers satisfying
$\alpha_1 + \cdots + \alpha_n \leqslant \frac{n(n+1)}{2} -1$, then:
\[
\small
\aligned
&
\sum_{\ell_1\geqslant\ell_2\geqslant\cdots\geqslant\ell_{n-1}
\geqslant\ell_n\geqslant 0}\,
M_{\ell_1,\ell_2,\dots,\ell_{n-1},\ell_n}^{\kappa,m}
\cdot
\big(\ell_1-\ell_2\big)^{\alpha_1}
\cdots\cdots
\big(\ell_{n-1}-\ell_n\big)^{\alpha_{n-1}}
\big(\ell_n\big)^{\alpha_n}
=
\\
&
=
\sum_{{\sf YT}\,{\sf semi-standard}
\atop
{\sf weight}({\sf YT})=m}
\big(\ell_1({\sf YT})-\ell_2({\sf YT})\big)^{\alpha_1}
\cdots
\big(\ell_{n-1}({\sf YT})-\ell_n({\sf YT})\big)^{\alpha_{n-1}}
\big(\ell_n({\sf YT})\big)^{\alpha_n}
\\
&
\leqslant
{\sf Constant}_{n,\kappa}
\cdot
m^{\alpha_1+\cdots+\alpha_n}
\cdot
m^{n\kappa-\frac{n(n-1)}{2}}
\\
&
\leqslant
{\sf Constant}_{n,\kappa}
\cdot
m^{(\kappa+1)n-2}.
\endaligned\,
\]
\end{Theorem}

{\small

\smallskip\noindent{\bf Acknowledgments.}
The author 
would like to thank Professor Ngaiming Mok, Director, and Professor
Yum-Tong Siu, C.V. Starr Visiting Professor, for inviting him at the
Hong Kong University's Department of Mathematics during August 4--14,
2009 when parts of the present memoir were finalized and presented
orally at a Summer conference entitled ``{\em Workshop on Complex
Geometry}''. Mainly, he must express his debt to the works~\cite{
rou2006a, rou2006b} by Erwan Rousseau for $n = \kappa = 3$ to whom
the present general Schur bundle strategy was borrowed, and also, he
would like to mention his debt to the works~\cite{ bru1972, bru1997}
by Peter Br\"uckmann. Besides, the present memoir has benefited of
enlightening discussions with Alain
Lascoux. Pascal Dingoyan indicated appropriate references. Jacky
Cresson and Michel Petitot provided information about multiple
polylogarithms. Lionel Darondeau suggested a few expositional
improvements. 

}

\markleft{Jo\"el Merker}
\markright{\sf \S2.~Universal combinatorics of Green-Griffiths jets}
\section{\bf Universal combinatorics of Green-Griffiths jets}
\label{Section-2}

\subsection{Jets of local holomorphic discs}
Let $X$ be a complex manifold
of dimension $n \geqslant 1$. 
In a local chart $(X, x_0) \simeq (\C^n,0)$ centered at a point $x_0
\in X$ equipped with $n$ complex coordinates $(x_1, \dots, x_n)$,
one looks at holomorphic discs passing through
$x_0$:
\[
f:(\D,0)\to(\C^n,0)\simeq(X,x_0),
\]
namely with $f ( 0) = x_0$, which possess of course $n$ components:
\[
\big(
f_1(\zeta),f_2(\zeta),\dots,f_n(\zeta)\big).
\]
For any integer $\kappa \geqslant 1$, the associated $\kappa$-jet map
of any such a holomorphic disc gathers all its $n\kappa$ derivatives
up to order $\kappa$ with respect to the (single) source variable
$\zeta \in \D$:
\[
j^\kappa f(\zeta)
=
\big(
f_1',\dots,f_n',f_1'',\dots,f_n'',\dots\dots,
f_1^{(\kappa)},\dots,f_n^{(\kappa)}
\big)
(\zeta).
\]
Accordingly, one is led to introduce $n \kappa$ new independent {\sl
jet coordinates} that will simply be denoted as:
\[
\big(
x_1',\dots,x_n',\,x_1'',\dots,x_n'',\,\dots\dots,\,
x_n^{(\kappa)},\dots,x_n^{(\kappa)}
\big),
\] 
so that $\big(x, x', x'', \dots, x^{ (\kappa)}\big)$ provide $n + n
\kappa$ coordinates on the space of {\em un}centered $\kappa$-jets of
maps $\D \to X$.

\subsection{Weighted homogeneous jet polynomials}
Above any fixed point $x_0 \in X$, Green-Griffiths (\cite{ gg1980})
introduced a certain ``fiber'' which consists of all polynomials in
these jet variables $x_i^{(\ell)}$, $1 \leqslant i \leqslant n$, $1
\leqslant \lambda \leqslant \kappa$, that are of the following type:
\begin{equation}
\label{gg-jets}
{\sf P}
\big(x',x'',\dots,x^{(\kappa)}
\big)
=
\!\!\!\!\!\!\!\!\!\!\!\!
\sum_{
\alpha_1,\alpha_2,\dots,\alpha_\kappa\in\N^n
\atop
\vert\alpha_1\vert+2\vert\alpha_2\vert
+\cdots+
\kappa\vert\alpha_\kappa\vert=m}\,
{\sf coeff}_{\alpha_1,\alpha_2,\dots,\alpha_\kappa}\cdot
\big(x'\big)^{\alpha_1}\big(x''\big)^{\alpha_2}
\cdots
\big(x^{(\kappa)}\big)^{\alpha_\kappa},
\end{equation}
where $m \geqslant 1$ is some integer, where $\alpha_\lambda = (
\alpha_{ \lambda, 1}, \dots, \alpha_{ \lambda, n} ) \in \N^n$ for $1
\leqslant \lambda \leqslant \kappa$ are multiindices of length $\vert
\alpha_\lambda \vert = \alpha_{ \lambda, 1} + \cdots + \alpha_{
\lambda, n}$, and where ${\sf coeff}_{ \alpha_1, \alpha_2, \dots,
\alpha_\kappa}$ are arbitrary complex coefficients, or equivalently
if written in greater length:
\[
\sum_{
\alpha_1,\alpha_2,\dots,\alpha_\kappa\in\N^n
\atop
\vert\alpha_1\vert+2\vert\alpha_2\vert
+\cdots+
\kappa\vert\alpha_\kappa\vert=m}
\!\!\!\!\!\!\!
{\sf coeff}_{\alpha_1,\alpha_2,\dots,\alpha_\kappa}\cdot
\prod_{1\leqslant i\leqslant n}
\big(x_i'\big)^{\alpha_{1,i}}
\prod_{1\leqslant i\leqslant n}
\big(x_i''\big)^{\alpha_{2,i}}
\cdots
\prod_{1\leqslant i\leqslant n}
\big(x_i^{(\kappa)}\big)^{\alpha_{\kappa,i}}.
\]
Visibly, such polynomials enjoy weighted homogeneity: 
\[
{\sf P}
\big(\delta x',\,\delta^2x'',\,\dots,\,\delta^\kappa x^{(\kappa)}\big)
=
\delta^m\,
{\sf P}\big(x',\,x'',\,\dots,\,x^{(\kappa)}\big)
\]
of the fixed weight $m$ with respect to the anisotropic complex jet
dilation defined by:
\[
\aligned
\delta\cdot\big(
x_{i_1}',x_{i_2}'',\dots,x_{i_\kappa}^{(\kappa)}\big)
:=
\big(
\delta\,x_{i_1}',\delta^2x_{i_2}'',\dots,
\delta^\kappa x_{i_\kappa}^{(\kappa)}\big),
\ \ \ \ \ \ \ \ \ \ \ \ \ \
\delta\in\C,
\endaligned
\]
whence for memory in all what follows one sees that:
\[
m
=
{\sf weight}
=
\text{\rm (fixed) total number of appearing ``primes''}. 
\]

\begin{Lemma}
As the point $x_0$ runs in $X$, these polynomial fibers organize
coherently as a {\em holomorphic vector} bundle $\mathcal{E}_{ \kappa,
m}^{ GG} T_X^*$ over $X$ of rank
equal to the number of arbitrary coefficients 
${\sf coeff}_{ \alpha_1, \dots, \alpha_\kappa}$, namely to:
\[
{\rm Card}\,
\big\{
(\alpha_1,\alpha_2,\dots,\alpha_\kappa)
\in
(\N^n)^\kappa
\colon
\vert\alpha_1\vert+2\vert\alpha_2\vert
+\cdots+
\kappa\vert\alpha_\kappa\vert
=
m
\big\}.
\]
\end{Lemma}

\proof
Consider an arbitrary change of (local) holomorphic chart on $X$:
\[
(x_1,\dots,x_n)
\longmapsto
(y_1,\dots,y_n)
=
\big(
\Psi_1(x_1,\dots,x_n),\,\dots\dots,\,
\Psi_n(x_1,\dots,x_n)
\big),
\]
understood as inducing a change of local trivialization for the
bundle. One must establish that the new coefficients of the
transformed jet polynomials express {\em linearly} in terms of the
coefficients of~\thetag{ \ref{gg-jets}}.

To begin with, the knowledge of how the $y_j^{ (\lambda)}$ express in
terms of the $x_i^{ (\tau )}$ is provided by an application of the
chain rule for the differentiation of a composed holomorphic disc
$\zeta \mapsto \Psi \big( f_1 ( \zeta), \dots, f_n ( \zeta)
\big)$. The closed combinatorial formula writes as follows, for any
$\lambda$ with $1 \leqslant \lambda \leqslant \kappa$ and for any $j$
with $1 \leqslant j \leqslant n$.

\begin{Theorem}
{\rm (see \cite{ cs1996, mer2008})}
The $\lambda$-jet of $\Psi_j ( f_1, \dots, f_n)$ is given by the
following multivariate Fa\`a di Bruno formula, written without the
argument $\zeta$:
\[
\scriptsize
\aligned
\big[
\Psi_j(f_1,\dots,f_n)
\big]^{(\lambda)}
&
=
\sum_{e=1}^\lambda\,
\sum_{1\leqslant\tau_1<\cdots<\tau_e\leqslant\lambda}\,
\sum_{\mu_1\geqslant 1,\dots,\mu_e\geqslant 1}\,
\sum_{\mu_1\tau_1+\cdots+\mu_e\tau_e=\lambda}\,
\frac{\lambda!}{
(\tau_1!)^{\mu_1}\mu_1!
\cdots
(\tau_e!)^{\mu_e}\mu_e!}\,
\\
&
\ \ \ \ \ \ \
\sum_{j_1^1,\dots,j_{\mu_1}^1=1}^n\,
\cdots\cdots\,
\sum_{j_1^e,\dots,j_{\mu_e}^e=1}^n\,
\frac{\partial^{\mu_1+\cdots+\mu_e}\Psi_j}{
\partial x_{j_1^1}\cdots\partial x_{j_{\mu_1}^1}
\cdots\cdots
\partial x_{j_1^e}\cdots\partial x_{j_{\mu_e}^e}}\cdot
\\
&
\ \ \ \ \ \ \ \ \ \ \ \ \ \ \ \ \ \ \ \ \ \ \ \ \ \ \ \ \ \ \ \ \ \ \ 
\ \ \ \ \ \ \ \ \ \ \ \ \ \ \ \ \ \ \ \ \ \ \ \ \ \ \ \ \ \ \ \ \ \ \ 
\cdot
f_{j_1^1}^{(\tau_1)}\cdots f_{j_{\mu_1}^1}^{(\tau_1)}
\cdots\cdots
f_{j_1^e}^{(\tau_e)}\cdots f_{j_{\mu_e}^e}^{(\tau_e)}.
\endaligned
\]
\end{Theorem}

To read this general formula, we comment it backward, understanding
it rather as a (polynomial, invertible) transformation between
independent jet variables:
\[
\footnotesize
\aligned
y_j^{(\lambda)}
&
=
\sum_{e=1}^\lambda\,
\sum_{1\leqslant\tau_1<\cdots<\tau_e\leqslant\lambda}\,
\sum_{\mu_1\geqslant 1,\dots,\mu_e\geqslant 1}\,
\sum_{\mu_1\tau_1+\cdots+\mu_e\tau_e=\lambda}\,
\frac{\lambda!}{
(\tau_1!)^{\mu_1}\mu_1!
\cdots
(\tau_e!)^{\mu_e}\mu_e!}\,
\\
&
\ \ \ \ \ \ \
\sum_{j_1^1,\dots,j_{\mu_1}^1=1}^n\,
\cdots\cdots\,
\sum_{j_1^e,\dots,j_{\mu_e}^e=1}^n\,
\frac{\partial^{\mu_1+\cdots+\mu_e}\Psi_j}{
\partial x_{j_1^1}\cdots\partial x_{j_{\mu_1}^1}
\cdots\cdots
\partial x_{j_1^e}\cdots\partial x_{j_{\mu_e}^e}}\cdot
\\
&
\ \ \ \ \ \ \ \ \ \ \ \ \ \ \ \ \ \ \ \ \ \ \ \ \ \ \ \ \ \ \ \ \ \ \ 
\ \ \ \ \ \ \ \ \ \ \ \ \ \ \ \ \ \ \ \ \ \ \ \ \ \ \ \ \ \ \ \ \ \ \ 
\cdot
x_{j_1^1}^{(\tau_1)}\cdots x_{j_{\mu_1}^1}^{(\tau_1)}
\cdots\cdots
x_{j_1^e}^{(\tau_e)}\cdots x_{j_{\mu_e}^e}^{(\tau_e)}.
\endaligned
\]
The general monomial $\prod x_\bullet^{ (\tau_1)} \, \prod x_\bullet^{
(\tau_2)} \cdots\, \prod x_\bullet^{ (\tau_2)}$ in the jet variables
gathers derivatives of increasing orders $\tau_1 < \tau_2 < \cdots <
\tau_e$ with $\mu_1, \mu_2, \dots, \mu_e$ counting their respective
numbers. Then $\Psi_j$ is subjected to a partial derivative of order
$\mu_1 + \mu_2 + \cdots + \mu_e$, the total number of
letters $x_\bullet^\bullet$ 
in the monomial in question. Because there are $n + 1$
variables $x_i$, the dots in the $x_\bullet^{ ( \tau_c)}$ should
receive indices, and in fact, there appear general sums $\sum_{ j_1^c,
\dots, j_{ \mu_c}^c = 1}^n$ over {\em all possible} such indices.

This precise closed combinatorial formula is not really needed for the
proof of our lemma, and instead, it is sufficient to know that each
$y_j^{ ( \lambda)}$ is a certain polynomial in the $x_i^{ (\tau_c)}$,
with coefficients depending linearly upon the $\lambda$-jet of $\Psi$,
the weight $\mu_1 \tau_1 + \cdots + \mu_e \tau_e$ of each appearing
monomial $x_{ j_1^1}^{ ( \tau_1)} \cdots x_{ j_{ \mu_1}^1}^{ (\tau_1)}
\cdots\cdots x_{j_1^e}^{ (\tau_e)} \cdots x_{ j_{ \mu_e}^e}^{
(\tau_e)}$ being constant equal to $\lambda$, and this fact is easily
proved by a rough induction argument. We can abbreviate this
as:
\[
y_j^{(\nu)}
=
\sum_{i_1=1}^n\,
\frac{\partial\Psi_j}{\partial x_{i_1}}
+\cdots+
\sum_{i_1,\dots,i_\lambda=1}^n\,
\frac{\partial^\lambda\Psi_j}{\partial x_{i_1}\cdots\,
\partial x_{i_\lambda}}
\]

Consequently, we must, as said, examine how an $y$-jet general
polynomial of weight $m$ like the $x$-jet polynomial ${\sf P}$
in~\thetag{ \ref{gg-jets}}:
\[
{\sf Q}
\big(y',\,y'',\,\dots,\,y^{(\kappa)}\big)
=
\sum_{\alpha_1,\,\alpha_2,\,\dots,\,\alpha_\kappa\in\N^n
\atop
\vert\alpha_1\vert+2\vert\alpha_2\vert
+\cdots+
\kappa\vert\alpha_\kappa\vert=m}
{\sf Q}_{\alpha_1,\alpha_2,\dots,\alpha_\kappa}
\cdot
\big(y'\big)^{\alpha_1}
\big(y''\big)^{\alpha_2}
\cdots
\big(y^{(\kappa)}\big)^{\alpha_\kappa}
\]
is transformed. From the theorem stated above (or from the rough
induction argument), we deduce that a general monomial of a weight
$m$:
\begin{equation}
\label{1-2-lambda-kappa}
\scriptsize
\aligned
\big(y'\big)^{\alpha_1}
\big(y''\big)^{\alpha_2}
&
\cdots
\big(y^{(\lambda)}\big)^{\alpha_\lambda}
\cdots
\big(y^{(\kappa)}\big)^{\alpha_\kappa}
=
\\
&
=
\bigg(
\sum_{i_1}\,\Psi_{x_{i_1}}\,x_{i_1}'
\bigg)^{\alpha_1}
\bigg(
\sum_{i_1}\,\Psi_{x_{i_1}}\,x_{i_1}''
+
\sum_{i_1,\,i_2}\,\Psi_{x_{i_1}x_{i_2}}\,x_{i_1}'x_{i_2}'
\bigg)^{\alpha_2}
\cdots
\\
&
\cdots
\bigg(
\sum_{i_1}\,\Psi_{x_{i_1}}\,x_{i_1}^{(\lambda)}
+\cdots+
\sum_{i_1,\dots,i_\lambda}\,\Psi_{x_{i_1}\dots
x_{i_\lambda}}\,x_{i_1}'\dots x_{i_\lambda}'
\bigg)^{\alpha_\lambda}
\cdots
\\
&
\cdots
\bigg(
\sum_{i_1}\,\Psi_{x_{i_1}}\,x_{i_1}^{(\kappa)}
+\cdots+
\sum_{i_1,\dots,i_\kappa}\,\Psi_{x_{i_1}\dots
x_{i_\kappa}}\,x_{i_1}'\dots x_{i_\kappa}'
\bigg)^{\alpha_\kappa}
\endaligned
\end{equation}
is clearly transformed to a jet polynomial of weight $m$:
\[
\aligned
\big(y'\big)^{\alpha_1}
\cdots
\big(y^{(\kappa)}\big)^{\alpha_\kappa}
=
\sum_{
\vert\beta_1\vert+\cdots+\kappa\vert\beta_\kappa\vert
=m}\,
{\sf H}_{\beta_1,\dots,\beta_\kappa}^{\alpha_1,\dots,\alpha_\kappa}
\big(j^\kappa\Psi\big)\cdot
(x')^{\beta_1}
\cdots
(x^{(\kappa)})^{\beta_\kappa}
\endaligned
\]
having coefficients that are certain universal polynomials in the
$\kappa$-jet of $\Psi$. It therefore follows that ${\sf Q} ( y',
\dots, y^{ (\kappa)})$ is transformed to:
\[
\aligned
&
\sum_{
\vert\alpha_1\vert+\cdots+\kappa\vert\alpha_\kappa\vert
=m}
\sum_{
\vert\beta_1\vert+\cdots+\kappa\vert\beta_\kappa\vert
=m}
{\sf Q}_{\alpha_1,\dots,\alpha_\kappa}\cdot
{\sf H}_{\beta_1,\dots,\beta_\kappa}^{\alpha_1,\dots,\alpha_\kappa}
\big(j^\kappa\Psi\big)\cdot
(x')^{\beta_1}
\cdots
(x^{(\kappa)})^{\beta_\kappa}
\\
&
=:
\sum_{
\vert\beta_1\vert+\cdots+\kappa\vert\beta_\kappa\vert
=m}
{\sf P}_{\beta_1,\dots,\beta_\kappa}\cdot
(x')^{\beta_1}
\cdots
(x^{(\kappa)})^{\beta_\kappa}
\endaligned
\]
with the following {\em linear} relationship between coefficients:
\[
{\sf P}_{\beta_1,\dots,\beta_\kappa}
=
\sum_{
\vert\alpha_1\vert+\cdots+\kappa\vert\alpha_\kappa\vert
=m}\,
{\sf H}_{\beta_1,\dots,\beta_\kappa}^{\alpha_1,\dots,\alpha_\kappa}
\big(j^\kappa\Psi\big)
\cdot
{\sf Q}_{\alpha_1,\dots,\alpha_\kappa}.
\]
This shows that $\mathcal{ E}_{ \kappa, m}^{ GG} T_X^*$ effectively is
a vector bundle, because the cocycle relations and the inverse
trivializations follow from the transitivity and from the
invertibility of change of local coordinates on $X$.
\endproof

\subsection{Symmetric pluri-tensor decomposition}
As is known in the domain (\cite{ gg1980, dem1997, dmr2010}),
a certain graded holomorphic vector bundle ${\sf Gr}^\bullet \mathcal{
E}_{ \kappa, m}^{ GG} T_X^*$ naturally associated to this
Green-Griffiths bundle $\mathcal{E}_{ \kappa, m}^{ GG} T_X^*$ happens
to decompose into the following direct sum of multi-tensored symmetric
powers of the cotangent bundle $T_X^*$ of $X$:
\[
\aligned
{\sf Gr}^\bullet
\mathcal{E}_{\kappa,m}^{GG}T_X^*
=
\bigoplus_{\ell_1+2\ell_2+\cdots+\kappa\ell_\kappa
=
m}\,
{\rm Sym}^{\ell_1}T_X^*
\otimes
{\rm Sym}^{\ell_2}T_X^*
\otimes\cdots\otimes
{\rm Sym}^{\ell_\kappa}T_X^*
\endaligned.
\]
Informally speaking, such a decomposition just relates to the fact
that the general $m$-weighted polynomial~\thetag{ \ref{gg-jets}}
looks like a linear combination of (tensor!) products of the
(individually symmetric!) monomials $(x')^{ \alpha_1}$, $(x'')^{
\alpha_2}$, \dots, $(x^{ (\kappa)})^{ \alpha_\kappa}$ with, say for a
good correspondence: 
\[
\ell_1\equiv\vert\alpha_1\vert,\ \ \ \ \
\ell_2\equiv\vert\alpha_2\vert,\ \ \ \ \ 
\dots\dots,\ \ \ \ \
\ell_\kappa\equiv\vert\alpha_\kappa\vert.
\]
But this view is not rigorous, so let us explain
with more details than in~\cite{ gg1980, dem1997, dmr2010} how one
builds the graded bundle ${\sf Gr}^\bullet \mathcal{E}_{ \kappa, m}^{
GG} T_X^*$.

Consider again the transformation~\thetag{ \ref{1-2-lambda-kappa}}. 
If $\alpha_{ \lambda+1} = \cdots = \alpha_\kappa = 0$
for some $\lambda$, one sees that
after expansion, the total number $\vert \alpha_1 \vert + 2 \vert
\alpha_2 \vert + \cdots + \lambda \vert \alpha_\lambda \vert$ of
primes remains unchanged. However, if there exists some $\mu$ with
$\lambda + 1 \leqslant \mu \leqslant \kappa$ such that $\alpha_\mu
\neq 0$, then in general the expansion of the factor $\big( \Psi (
x)^{ (\mu)} \big)^{ \alpha_\mu}$ adds a total of $\mu \vert \alpha_\mu
\vert$ further primes to the monomials in $x'$, $x''$, \dots, $x^{ (
\lambda)}$ that was already obtained by expanding the first $\lambda$
factors 
$\big( \Psi ( x)'\big)^{ \alpha_1}
\cdots\, \big( \Psi (x)^{ (\lambda)}\big)^{
\alpha_\lambda}$. 
Thus in all cases, after an arbitrary change of coordinates
$x \mapsto \Psi ( x)$, the $\lambda$-restricted weight $\vert \alpha_1
\vert + 2 \vert \alpha_2 \vert + \cdots + \lambda \vert \alpha_\lambda
\vert$ {\em can only increase}. Following~\cite{ gg1980}, one
may hence define for any $\lambda$ fixed in advance with $1 \leqslant
\lambda \leqslant \kappa$ a (decreasing) {\em filtered sequence}:
\[
\mathcal{E}_{\kappa,m}^{GG}T_X^*
=
\mathcal{F}_\lambda^0
\supset
\mathcal{F}_\lambda^1
\supset
\mathcal{F}_\lambda^2
\supset\cdots\supset
\mathcal{F}_\lambda^m
\supset
\{0\}
=
\mathcal{F}_\lambda^{m+1}
\]
of subbundles of $\mathcal{ E}_{ \kappa, m}^{ GG} T_X^*$ whose pieces
for any $q = 1, 2, \dots, m$ are naturally defined by:
\[
\footnotesize
\aligned
\mathcal{F}_\lambda^q
=
\mathcal{F}_\lambda^q
\big(
\mathcal{E}_{\kappa,m}^{GG}T_X^*
\big)
=
\left\{
\aligned
&
{\sf P}\big(x',\dots,x^{(\lambda)},\dots,x^{(\kappa)}\big)
\in\mathcal{E}_{\kappa,m}^{GG}T_X^*\
\text{\rm involving only monomials}
\\
&
(x')^{\alpha_1}\cdots (x^{(\lambda)})^{\alpha_\lambda}
\cdots (x^{(\kappa)})^{\alpha_\kappa}\
\text{\rm with}\
\vert\alpha_1\vert+\cdots+\lambda\vert\alpha_\lambda\vert
\geqslant
q
\endaligned
\right\}.
\endaligned
\]
Notice that $\mathcal{ F}_\lambda^m = \mathcal{ E}_{ \lambda, m}^{ GG}
T_X^*$. If we now set $\lambda = \kappa - 1$, the graded bundle
associated with this filtration:
\[
{\sf Gr}^\bullet
\mathcal{E}_{\kappa,m}^{GG}T_X^*
=
\big(
\mathcal{F}_{\kappa-1}^0
\big/
\mathcal{F}_{\kappa-1}^1
\big)
\oplus
\cdots
\oplus
\big(
\mathcal{F}_{\kappa-1}^{m-1}
\big/
\mathcal{F}_{\kappa-1}^m
\big)
\oplus
\mathcal{E}_{\kappa-1,m}^{GG}T_X^*
\]
is constituted of quotient factors:
\[
\mathcal{G}_{\kappa-1}^q
:=
\mathcal{F}_{\kappa-1}^q
\big/
\mathcal{F}_{\kappa-1}^{q+1}
\ \ \ \ \ \ \ \ \ \ \ \ \
{\scriptstyle{(q\,=\,0\,\cdots\,m\,-\,1)}}
\]
which consist of polynomials ${\sf P}$ as above for which:
\[
\vert\alpha_1\vert 
+\cdots+ 
(\kappa-1)\vert\alpha_{\kappa-1}\vert
=
q
\]
modulo polynomials for which $\vert \alpha_1 \vert +\cdots+ (\kappa-1)
\vert \alpha_{ \kappa-1} \vert \geqslant q + 1$. It follows at once
that $q + \kappa \vert \alpha_\kappa \vert = m$ in such polynomials,
that is to say:
\[
q
=
m-\kappa\ell_\kappa
\]
for some integer $\ell_\kappa \in \N$ with ${\sf Ent} 
\big[ \frac{ m}{ \kappa} \big]
\geqslant \ell_\kappa \geqslant 0$. In particular, this
quotient $\mathcal{ F}_{ \kappa - 1}^q \big / \mathcal{ F}_{ \kappa -
1}^{ q + 1}$ reduces to $\{ 0\}$ whenever $m - q$ is {\em not}
divisible by $\kappa$.

We now claim that:
\[
\mathcal{G}_{\kappa-1}^{m-\kappa\ell_\kappa}
\big(\mathcal{E}_{\kappa,m}^{GG}T_X^*\big)
\simeq
\mathcal{E}_{\kappa-1,\,m-\kappa\ell_\kappa}^{GG}T_X^*
\otimes
{\rm Sym}^{\ell_\kappa}T_X^*. 
\]
Indeed, under an arbitrary change of coordinates $x \mapsto \Psi ( x)
= y$, the constituents $x^{ (\kappa)}$ of the monomials of highest jet
$\big( x^{ (\kappa)} \big)^{ \alpha_\kappa}$ with $\vert \alpha_\kappa
\vert = \ell_\kappa$ are transformed to:
\[
y^{(\kappa)}
=
\sum_{i_1=1}^n\,\Psi_{x_{i_1}}\,x_{i_1}^{(\kappa)}
\ \
{\rm modulo}\ \
\big(x',\dots,x^{(\kappa-1)}\big),
\]
so that they visibly transform in exactly the same covariant way as
the covectors in $T_X^*$, namely: 
\[
d\Psi
=
\sum_{i_1=1}^n\, 
\Psi_{x_{i_1}}\, 
dx_{i_1}, 
\]
and so, the
$\big( x^{ (\kappa)} \big)^{ \alpha_\kappa}$ transform as ${\rm Sym}^{
\ell_\kappa} T_X^*$. The other constituents of $\mathcal{ F}_{
\kappa-1}^{ m - \kappa \ell_\kappa}$ depend only on the
$(\kappa-1)$-jet and are of the remaining weight $m - \kappa
\ell_\kappa$, whence the claimed isomorphism follows.

Putting together all these isomorphisms, we get: 
\[
\small
\aligned
{\sf Gr}^\bullet
\mathcal{E}_{\kappa,m}^{GG}T_X^*
&
=
\bigoplus_{{\sf Ent}[\frac{m}{\kappa}]\geqslant\ell_\kappa\geqslant 1}\,
\mathcal{G}_{\kappa-1}^{m-\kappa\ell_\kappa}
\bigoplus
\mathcal{E}_{\kappa-1,m}^{GG}T_X^*
\\
&
=
\bigoplus_{{\sf Ent}[\frac{m}{\kappa}]\geqslant\ell_\kappa\geqslant 1}\,
\Big(
\mathcal{E}_{\kappa-1,\,m-\kappa\ell_\kappa}^{GG}T_X^*
\otimes
{\rm Sym}^{\ell_\kappa}T_X^*
\Big)
\bigoplus
\big(
\mathcal{E}_{\kappa-1,m}^{GG}T_X^*
\otimes
{\rm Sym}^0T_X^*
\big)
\\
&
=
\bigoplus_{{\sf Ent}[\frac{m}{\kappa}]\geqslant\ell_\kappa\geqslant 0}\,
\mathcal{E}_{\kappa-1,\,m-\kappa\ell_\kappa}^{GG}T_X^*
\otimes
{\rm Sym}^{\ell_\kappa}T_X^*.
\endaligned
\]
Now an induction of this isomorphism applied to ${\sf Gr}^\bullet
\big( \mathcal{ E}_{ \kappa-1, \, m - \kappa \ell_\kappa}^{ GG} T_X^*
\big)$ yields the announced decomposition for ${\sf Gr}^\bullet
\mathcal{ E}_{\kappa, m}^{ GG} T_X^*$.

\begin{Theorem}
\label{graded-sym}
The holomorphic vector bundle $\mathcal{ E}_{ \kappa, m}^{ GG} T_X^*$
of Green-Griffiths polynomials of weight $m$ in the $\kappa$-jet of
local complex curves $\D \to X$ admits a natural filtration whose
associated graded bundle is isomorphic to the following direct sum
of multi-tensored symmetric powers of the cotangent bundle:
\[
\boxed{
\aligned
{\sf Gr}^\bullet
\mathcal{E}_{\kappa,m}^{GG}T_X^*
=
\bigoplus_{\ell_1+2\ell_2+\cdots+\kappa\ell_\kappa
=
m}\,
{\rm Sym}^{\ell_1}T_X^*
\otimes
{\rm Sym}^{\ell_2}T_X^*
\otimes\cdots\otimes
{\rm Sym}^{\ell_\kappa}T_X^*
\endaligned}\,.
\]
Furthermore, for every $q = 1, 2, \dots, n$, one has the inequalities:
\begin{equation}
\aligned
\label{first-inq-cohomology}
&
\dim H^q
\big(
X,\,\mathcal{E}_{\kappa,m}^{GG}T_X^*
\big)
\leqslant
\\
&
\leqslant
\sum_{\ell_1+2\ell_2+\cdots+\kappa\ell_\kappa=m}\,
\dim H^q
\Big(
X,\,\,
{\rm Sym}^{\ell_1}T_X^*
\otimes
{\rm Sym}^{\ell_2}T_X^*
\otimes\cdots\otimes
{\rm Sym}^{\ell_\kappa}T_X^*
\Big). 
\endaligned
\end{equation}
\end{Theorem}

To complete the proof, it only remains to check the inequality between
the cohomology dimensions. Let us consider instead in greater
generality the following situation, which clearly embraces the last
claim above.

\begin{Lemma}
Suppose a holomorphic
vector bundle $E \to X$ admits a filtration:
\[
\{0\}=E_{r+1}
\subset
E_r
\subset
E_{r-1}
\subset\cdots\subset
E_{k+1}
\subset
E_k
\subset
E_{k-1}
\subset\cdots\subset
E_1
\subset
E_0=E
\]
by nested holomorphic subbundles $E_k$, the associated graded
bundle being:
\[
{\rm Gr}^\bullet E
=
E_r
\oplus
\big(E_{r-1}/E_r\big)
\oplus\cdots\oplus
\big(E_k/E_{k+1}\big)
\oplus
\big(E_{k-1}/E_k\big)
\oplus\cdots\oplus
\big(E_0/E_1\big).
\]
where, for a good notational correspondence: ${\sf Gr}^k
E := E_k / E_{ k+1}$. Then for every $q = 0, 1, \dots, n$, the
following inequality between cohomological dimensions holds:
\[
\dim H^q\big(X,\,E\big)
\leqslant
\sum_{k=0}^r\,
\dim H^q\big(X,\,E_k/E_{k+1}\big).
\]
\end{Lemma}

\proof
To each obviously true short exact sequence:
\begin{equation}
\label{exact-graded-k}
0
\longrightarrow
{\sf Gr}^kE
\longrightarrow
E/E_{k+1}
\longrightarrow
E/E_k
\longrightarrow
0
\ \ \ \ \ \ \ \ \ \ \ \ \
{\scriptstyle{(k\,=\,0,\,1,\,\dots,\,r)}}
\end{equation}
is associated the long exact sequence between cohomology groups:
\begin{equation}
\label{long-exact}
\cdots
\longrightarrow
H^q
\big(
X,\,
{\sf Gr}^kE
\big)
\longrightarrow
H^q
\big(
X,\,
E/E_{k+1}
\big)
\longrightarrow
H^q
\big(
X,\,
E/E_k
\big)
\longrightarrow
\cdots,
\end{equation}
and the trivial majoration: $\dim B \leqslant \dim A + \dim C$ of the
dimension of any member $B$ of any long exact sequence of vector
spaces by the sum of the dimensions of its two immediate neighbors
gives us here:
\[
\dim H^q
\big(
X,\,
E/E_{k+1}
\big)
\leqslant
\dim H^q
\big(
X,\,
{\sf Gr}^kE
\big)
+
\dim H^q
\big(
X,\,
E/E_k
\big)
\ \ \ \ \ 
{\scriptstyle{(k\,=\,0,\,1,\,\dots,\,r)}}. 
\]
Starting from $k = 0$ for which ${\sf Gr}^0 E = E / E_1$ and $E / E_0
= \{ 0\}$, a plain summation up to $k = r$ of these inequalities
cancels out all factors involving $E / E_k$ except only one on the
left: $E / E_{ r+1} = E$, and we get the desired inequality:
\[
\dim H^q
\big(X,\,
E\big)
=
\dim H^q
\big(X,\,
E/E_{r+1}
\big)
\leqslant
\sum_{k=0}^r\,
\dim H^q
\big(X,\,
{\sf Gr}^kE
\big)
\]
which, when applied to the Green-Griffiths bundle, terminates 
our detailed restitution of the theorem.
\endproof

Furthermore, in specific situations where the dimensions of the first
cohomology groups of the graded pieces $E_k / E_{ k+1}$
do not vanish but happen to become asymptotically
(much) smaller than the dimensions of their 
zeroth cohomology groups, a
useful second lemma is as follows. 

\begin{Lemma}
Under the same assumptions and just for $i = 0$, in addition to the
above majoration $h^0 (X, E) \leqslant \sum_{ k=0}^r \, h^0 (X, \, E_k
/ E_{ k+1} )$, one has the {\em minoration:}
\[
h^0(X,\,E)
\geqslant
\sum_{k=0}^r\,
h^0\big(X,\,E_k/E_{k+1}\big)
-
\sum_{k=0}^r\,
h^1\big(X,\,E_k/E_{k+1}\big).
\]
\end{Lemma}

\proof
As a preliminary, we observe that any long exact sequence
\[
0
\longrightarrow
A
\overset{a}{\longrightarrow}
B
\overset{b}{\longrightarrow}
C
\overset{c}{\longrightarrow}
D
\overset{d}{\longrightarrow}
\cdots
\]
can be stopped at its fourth term by replacing it with 
the four-terms sequence: 
\[
0
\longrightarrow
A
\longrightarrow
B
\longrightarrow
C
\longrightarrow
{\rm Im}(c)
\longrightarrow
0
\]
which is easily checked to be also exact. Then by considering the
basic equality which comes out by alternately summing the dimensions
of its members:
\[
0
=
\dim(A)-\dim(B)+\dim(C)-\dim({\rm Im}(c)),
\]
we deduce from the trivial inequality $\dim ({\rm Im} ( c) ) \leqslant
\dim (D)$, the useful minoration:
\[
\dim(B)
\geqslant
\dim(A)+\dim(C)-\dim(D).
\]
Applying now such an inequality to the first four terms of the long
exact sequence associated to the $k$-th
quotient exact sequence~\thetag{ \ref{exact-graded-k}} above:
\[
\aligned
0
&
\longrightarrow
H^0\big(X,\,{\sf Gr}^kE\big)
\longrightarrow
H^0\big(X,\,E/E_{k+1}\big)
\longrightarrow
H^0\big(X,\,E/E_k\big)
\longrightarrow
\\
&
\longrightarrow
H^1\big(X,\,{\sf Gr}^kE\big)
\longrightarrow
\cdots,
\endaligned 
\]
we readily deduce: 
\[
h^0\big(X,\,E/E_{k+1}\big)
\geqslant
h^0\big(X,\,{\sf Gr}^kE\big)
+
h^0\big(X,\,E/E_k\big)
-
h^1\big(X,\,{\sf Gr}^kE\big).
\]
Starting then from $k = 0$ for which ${\sf Gr}^0 E = E / E_1$ and $E /
E_0 = \{ 0\}$, a plain summation of these inequalities up to $k = r$
cancels out all terms involving an $E / E_k$ except one: $E / E_{
r+1} = E$, and we get the announced minoration.
\endproof

As a corollary, by applying this second elementary lemma to the
Green-Griffiths bundle, we gain a possibly useful general minoration:
\begin{equation}
\footnotesize
\aligned
h^0\big(X,\,\mathcal{E}_{\kappa,m}^{GG}T_X^*\big)
&
\geqslant
\sum_{\ell_1+2\ell_2+\cdots+\kappa\ell_\kappa=m}\,
h^0
\Big(X,\,{\rm Sym}^{\ell_1}T_X^*\otimes{\rm Sym}^{\ell_2}T_X^*
\otimes\cdots\otimes{\rm Sym}^{\ell_\kappa}T_X^*
\Big)
-
\\
&
\ \ \ \ \
-
\sum_{\ell_1+2\ell_2+\cdots+\kappa\ell_\kappa=m}\,
h^1
\Big(X,\,{\rm Sym}^{\ell_1}T_X^*\otimes{\rm Sym}^{\ell_2}T_X^*
\otimes\cdots\otimes{\rm Sym}^{\ell_\kappa}T_X^*
\Big).
\endaligned
\end{equation}
Hence it is now clear that, in order to establish existence
nonzero global sections of $\mathcal{ E}_{ \kappa, m}^{ GG} T_X^*$ on
hypersurfaces of general type, it would suffice that the considered
sum of $h^1$'s grows less substantially than the sum of $h^0$'s, as
$\kappa$ tends to $\infty$ and as $m \gg \kappa$ tends to $\infty$
too. We conclude this section by recalling that the Euler-Poincar\'e
characteristic transfers better than cohomology dimensions
through exact sequences, namely without
inequalities.

\begin{Lemma}
\label{Gr-chi}
The Euler-Poincar\'e characteristic of the Green-Griffiths
bundle is equal to that of its associated graded bundle:
\[
\chi\big(X,\mathcal{E}_{\kappa,m}^{GG}T_X^*\big)
=
\chi\big(X,{\sf Gr}^\bullet
\mathcal{E}_{\kappa,m}^{GG}T_X^*\big).
\]
\end{Lemma}

\proof
In fact, in the general context of the previous two lemmas, because
for any exact sequence of finite-dimensional vector spaces:
\[
0
\longrightarrow
A_1
\longrightarrow
A_2
\longrightarrow
\cdots
\longrightarrow
A_N
\longrightarrow
0
\]
one may easily check that the alternating
sum of dimensions vanishes:
\[
0
=
\dim A_1
-
\dim A_2
+\cdots+
(-1)^N\,\dim A_N,
\]
one deduces from the long exact sequence of cohomology~\thetag{
\ref{long-exact}} that:
\[
\aligned
0
&
=
h^0(X,\,{\sf Gr}^kE)
-
h^0(X,\,E/E_{k+1})
+
h^0(X,\,E/E_k)
-
\\
&
\ \ \ \ \
-
h^1(X,\,{\sf Gr}^kE)
-
h^1(X,\,E/E_{k+1})
+
h^1(X,\,E/E_k)
+
\\
&
\ \ \ \ \ 
+\cdots\cdots\cdots\cdots\cdots\cdots\cdots\cdots\cdots
\cdots\cdots\cdots\cdots\cdots\cdots\cdots\cdots\cdot+
\\
&
\ \ \ \ \
+
(-1)^nh^n(X,\,{\sf Gr}^kE)
-
(-1)^nh^n(X,\,E/E_{k+1})
+
(-1)^n
h^n(X,\,E/E_k),
\endaligned
\]
or else after gathering terms column by column, that:
\[
0
=
\chi(X,\,{\sf Gr}^kE)
-
\chi(X,\,E/E_{k+1})
+
\chi(X,\,E/E_k).
\]
Finally, a plain summation $\sum_{ k=0}^n$ yields the formula claimed.
\endproof

\markleft{Jo\"el Merker}
\markright{\bf Euler-Poincar\'e characteristic computations}
\section{\bf Euler-Poincar\'e characteristic computations}
\label{Section-3}

\begin{Theorem}
\label{GG-characteristic}
{\rm (\cite{gg1980})}
On an arbitrary compact complex projective manifold of dimension 
$n \geqslant 1$, the Green-Griffiths jet bundle $\mathcal{ E}_{
\kappa, m}^{ GG} T_X^*$ has an Euler-Poincar\'e characteristic
asymptotically given by:
\[
\footnotesize
\aligned
\chi\big(X,\mathcal{E}_{\kappa,m}^{GG}T_X^*\big)
&
=
\chi\big(X,{\sf Gr}^\bullet
\mathcal{E}_{\kappa,m}^{GG}T_X^*\big)
\\
&
=
\sum_{\ell_1+2\ell_2+\cdots+\kappa\ell_\kappa=m}\,
\chi\Big(X,\,
{\rm Sym}^{\ell_1}T_X^*\otimes{\rm Sym}^{\ell_2}T_X^*
\otimes\cdots\otimes
{\rm Sym}^{\ell_\kappa}T_X^*
\Big)
\\
&
=
{\textstyle{\frac{m^{(\kappa+1)n-1}}{((\kappa+1)n-1)!\,(\kappa!)^{n}}}}
\Big\{
({\sf c}_1^*)^n
\frac{(\log\kappa)^n}{n!}
+
{\sf O}_n\big((\log\kappa)^{n-1}\big)
\Big\}
+
{\sf O}_{n,\kappa}\big(m^{(\kappa+1)n-2}\big),
\endaligned
\]
where ${\sf c}_1^* = {\sf c}_1 ( T_X^*) = - \, {\sf c}_1 ( T_X)$
is the first Chern class of $T_X^*$, a $(1, 1)$-cohomology class on
$X$, and where:

\smallskip{\bf (i)} the first remainder is a linear combination of
homogeneous\footnote{\,
As usual, we understand implicitly that each $(n, n)$-cohomology class
$({\sf c}_1^*)^{ \lambda_1} \cdots \, ( {\sf c}_n^* )^{ \lambda_n}$ is
{\em integrated} over $X$, hence represents the {\em numerical} value
$\int_X \, ({\sf c}_1^*)^{ \lambda_1} \cdots \, ( {\sf c}_n^* )^{
\lambda_n}$.}
terms $({\sf c}_1^*)^{ \lambda_1} ({\sf c}_2^*)^{\lambda_ 2} \cdots \,
( {\sf c}_n^* )^{ \lambda_n}$ with $\lambda_1 + 2\, \lambda_ 2+ \cdots
+ n\, \lambda_n = n$, with rational coefficients all bounded in
absolute value by ${\sf Constant}_n \, ( \log \kappa)^{ n-1}$;

\smallskip{\bf (ii)} 
the second remainder is a polynomial in $m$ of submaximal degree
$\leqslant ( \kappa + 1)n - 2$ whose coefficients are linear
combinations of the same $({\sf c}_1^*)^{ \lambda_1} ({\sf
c}_2^*)^{\lambda_ 2} \cdots \, ( {\sf c}_n^* )^{ \lambda_n}$ with
rational coefficients all also bounded in absolute value by ${\sf
Constant}_{ n, \kappa}$.
\end{Theorem}

In the case where $X \subset \P^{ n+1} ( \C)$ is a hypersurface of
degree $d \geqslant 1$, each homogeneous monomial $({\sf c}_1^*)^{
\lambda_1} ({\sf c}_2^*)^{\lambda_ 2} \cdots \, ( {\sf c}_n^* )^{
\lambda_n}$\,\,---\,\,implicitly integrated on $X$, as usually
understood in complex algebraic geometry\,\,---\,\,expresses in terms
of $n$ and $d$ by means of some universal formulas as follows. Let $h
:= {\sf c}_1\big( \mathcal{ O}_{ \mathbb{ P}^{ n+1}} (1)\big)$ denotes
the hyperplane $(1, 1)$-cohomology class, which satisfies $h^n =
\int_X h^n = d$. Then one may represent (\cite{ dmr2010}, p.~170):
\begin{equation}
\label{c-d}
\!\!\!\!
\left[
\aligned
{\sf c}_1
&
=
-h
\big(d-n-2\big)
\\
{\sf c}_2
&
=
h^2
\big(
d^2-
{\textstyle{\frac{(n+2)!}{(n+1)!\,\,1!}}}\,d
+
{\textstyle{\frac{(n+2)!}{n!\,\,2!}}}
\big)
\\
{\sf c}_3
&
=
-h^3
\big(
d^3
-
{\textstyle{\frac{(n+2)!}{(n+1)!\,\,1!}}}\,d^2
+
{\textstyle{\frac{(n+2)!}{n!\,\,2!}}}\,d
-
{\textstyle{\frac{(n+2)!}{(n-1)!\,\,3!}}}
\big)
\\
\,\cdot\cdot\,\,
&\,
\cdots\cdots\cdots\cdots\cdots\cdots\cdots\cdots\cdots\cdots\cdots
\cdots\cdots\cdots\cdots
\\
{\sf c}_n
&
=
(-1)^n\,h^n\big(
d^n
-
{\textstyle{\frac{(n+2)!}{(n+1)!\,\,1!}}}\,d^{n-1}
+\cdots+
(-1)^n\,{\textstyle{\frac{(n+2)!}{2!\,\,n!}}}\big).
\endaligned\right.
\end{equation}
Recall that the Chern classes of the tangent $T_X$ and of
the cotangent bundle $T_X^*$ of the hypersurface $X \subset
\P^{ n+1} ( \C)$ are linked together by the simple relation: 
\[
{\sf c}_k^*
:=
{\sf c}_k\big(T_X^*\big)
=
(-1)^k\,{\sf c}_k\big(T_X)
=
(-1)^k\,{\sf c}_k,
\] 
for $k = 1, 2, \dots, n$. It follows notably for instance that:
\[
({\sf c}_1^*)^n
=
(-1)^n
(-h)^n\,\big(d-n-2\big)^n
=
d\big(d-n-2\big)^n.
\]
Generally speaking, one easily convinces oneself that each homogeneous
degree $n$ monomial ${\sf c}_1^{ \tau_1} {\sf c}_2^{ \tau_2} \cdots
{\sf c}_n^{ \tau_n}$ identifies with a certain polynomial:
\[
{\sf c}_1^{\tau_1}{\sf c}_2^{\tau_2}\cdots{\sf c}_n^{\tau_n}
=
\sum_{k=1}^{n+1}\,
C_k^{\tau_1,\tau_2,\dots,\tau_n}\cdot
d^k
\] 
with respect to $d = \deg X$ having degree $\leqslant n+1$ with
integer coefficients $C_k^{ \tau_1, \tau_2, \dots, \tau_n} \in \Z$.
Furthermore, the constant coefficient
$C_0^{\tau_1,\tau_2,\dots,\tau_n} = 0$ is zero, because the factor
$h^{ \tau_1 + 2\tau_2 + \cdots + n \tau_n} = h^n = d$ is clearly
always present in every ${\sf c}_1^{ \tau_1} {\sf c}_2^{ \tau_2}
\cdots {\sf c}_n^{ \tau_n}$.

As a result, the Euler-Poincar\'e characteristic of $\mathcal{ E}_{
\kappa, m}^{ GG} T_X^*$:
\[
\small
\aligned
{\textstyle{\frac{m^{(\kappa+1)n-1}}{(\kappa!)^{n}((\kappa+1)n-1)!}}}
\Big\{
d\,(d-n-2)^n
\frac{(\log\kappa)^n}{n!}
+
{\sf O}_{n,d}\big((\log\kappa)^{n-1}\big)
\Big\}
+
{\sf O}_{n,d,\kappa}\big(m^{(\kappa+1)n-2}\big),
\endaligned
\]
visibly tends to $+\infty$ as $m$ and $\kappa$ tend to $+\infty$ as
soon as $X$ is of general type, that it is to say, as soon as:
\[
\deg X
\geqslant
\dim X+3.
\]

\proof
With more details, we redo Green-Griffiths' proof; fundamentals may be
found in~\cite{ ful1998} (pp.~50--59 plus Chap.~15), in~\cite{ botu1982}
and in~\cite{ hirz1966}.

To begin with, introduce the {\sl formal root decomposition}: 
\[
{\sf c}\big(T_X^*\big)
=
1+{\sf c}_1^*+{\sf c}_2^*+\cdots+{\sf c}_n^*
=
(1+{\sf a}_1^*)(1+{\sf a}_2^*)\cdots(1+{\sf a}_n^*)
\]
of the {\sl total Chern class}, namely of the sum ${\sf c} ( T_X^*)$
of the ${\sf c}_i^*$ so that ${\sf c}_i^*$ is the $i$-th elementary
symmetric function of the {\sl Chern roots} ${\sf a}_j^*$:
\[
{\sf c}_i^*
=
\sum_{1\leqslant j_1<j_2<\cdots<j_i\leqslant n}\,
{\sf a}_{j_1}^*{\sf a}_{j_2}^*
\cdots\,
{\sf a}_{j_i}^*
\ \ \ \ \ \ \ \ \ \ \ \ \ {\scriptstyle{(i\,=\,0,\,\,1\,\cdots\,n)}}.
\] 
Similarly, let the ${\sf a}_j$ denote the Chern roots of ${\sf c} (
T_X) = 1 + {\sf c}_1 + \cdots + {\sf c}_n$ and let the symbol $[ \, \,
\, ]_j$ denote projection to the $(j, j)$-cohomology class, so that
for instance $\big[ 1 + {\sf c}_1^* + \cdots + {\sf c}_n^* \big]_j =
{\sf c}_j^*$. To prove
the theorem, we must apply the Riemann-Roch-Hirzebruch
theorem~\cite{ hirz1966} which states that the Euler-Poincar\'e
characteristic:
\[
\chi
\big(X,\,\,
\mathcal{E}_{\kappa,m}^{GG}T_X^*
\big)
\overset{def}{=}
\sum_{0\leqslant q\leqslant n}\,
(-1)^q\,\dim
H^q
\big(X,\,\,
\mathcal{E}_{\kappa,m}^{GG}T_X^*
\big)
\]
is equal to the integral over $X$:
\[
\chi
=
\int_X\,
\big[
{\sf ch}
\big(\mathcal{E}_{\kappa,m}^{GG}T_X^*\big)\cdot
{\sf td}(T_X)
\big]_n
=
\int_X\,
\sum_{j=0}^n\,
\big[
{\sf ch}
\big(\mathcal{E}_{\kappa,m}^{GG}T_X^*\big)
\big]_{n-j}
\big[
{\sf td}(T_X)
\big]_j
\]
of the $(n, n)$-part of the product between the {\sl Chern character}
${\sf ch} \big( \mathcal{ E}_{ \kappa, m}^{ GG} T_X^* \big)$ of
$\mathcal{ E}_{ \kappa, m}^{ GG} T_X^*$ (to be computed in a while)
and the {\sl Todd class}\footnote{\,
Of course, all terms of degree $\geqslant n + 1$
in the ${\sf a}_j$ are $\equiv 0$, since
the associated cohomology classes vanish, as
does any form of bidegree $(p, q)$ with 
$p \geqslant n+1$ or $q \geqslant n+1$. 
} 
of $T_X$:
\[
\aligned
{\sf td}(T_X)
&
=
\frac{{\sf a}_1}{1-e^{-{\sf a}_1}}\,
\frac{{\sf a}_2}{1-e^{-{\sf a}_2}}\,
\cdots\,
\frac{{\sf a}_n}{1-e^{-{\sf a}_n}}
&
=
1
+
{\textstyle{\frac{1}{2}}}\,
{\sf c}_1
+
{\textstyle{\frac{1}{12}}}\,
\big[
{\sf c}_1^2+{\sf c}_2
\big]
+\cdots.
\endaligned
\]
As usual in asymptotic complex algebraic geometry (cf.~\cite{
dem2000, laza2004a, laza2004b}), for the
product: ${\sf chern} \cdot {\sf todd}$, picking cohomology classes of
positive degree $\geqslant 1$ in ${\sf td} ( T_X)$ forces to pick
classes in ${\sf ch} \big( \mathcal{ E}_{ \kappa, m}^{ GG} T_X^*
\big)$ of bidegree $\leqslant ( n - 1, \, n - 1)$, and then the
associated $m$-contributions are {\em smaller} than the maximal
possible: $m^{ ( \kappa + 1)n - 1}$. More precisely:

\begin{Lemma}
For every $j = 0, 1, \dots, n$, one has:
\[
\big[
{\sf ch}
\big(
\mathcal{E}_{\kappa,m}^{GG}T_X^*
\big)
\big]_{n-j}
=
{\sf O}_{n,\kappa}\big(m^{(\kappa+1)n-1-j}\big)
\]
and consequently all terms $\sum_{ j=1}^n$ are negligible
for our purposes, whence:
\[
\chi
\big(X,\,\,
\mathcal{E}_{\kappa,m}^{GG}T_X^*
\big)
=
\int_X\,
\big[
{\sf ch}
\big(\mathcal{E}_{\kappa,m}^{GG}T_X^*\big)
\big]_n
+
{\sf O}_{n,\kappa}\big(m^{(\kappa+1)n-2}\big).
\]
\end{Lemma}

\proof
As is known, the Chern character of the jet bundle equals that of its
graded decomposition, and the Chern character is both additive on
direct sums and multiplicative on tensor products, so that we
can write:
\[
{\sf ch}\big(
\mathcal{E}_{\kappa,m}^{GG}T_X^*\big)
=
\sum_{\ell_1+2\,\ell_2+\cdots+\kappa\ell_\kappa=m}\,
{\sf ch}
\big(
{\rm Sym}^{\ell_1}T_X^*
\big)
\,
{\sf ch}
\big(
{\rm Sym}^{\ell_2}T_X^*
\big)
\cdots\,
{\sf ch}
\big(
{\rm Sym}^{\ell_\kappa}T_X^*
\big). 
\] 
Furthermore, recall that the Chern character of an arbitrary symmetric
power of $T_X^*$ is given, in terms of the ${\sf a}_i^*$, by the known
formula:
\[
{\sf ch}\big({\rm Sym}^\ell T_X^*\big)
=
\sum_{
\substack{
x_1+\cdots+x_n=\ell
\\
x_1,\dots,x_n\in\N
}}\,
e^{x_1{\sf a}_1^*+\cdots+x_n{\sf a}_n^*}.
\]
We can therefore apply this to $\ell = \ell_\lambda$ for all $\lambda$
with $1 \leqslant \lambda \leqslant \kappa$:
\[
{\sf ch}\big({\rm Sym}^{\ell_\lambda}T_X^*\big)
=
\sum_{
x_{\lambda 1}+\cdots+x_{\lambda n}=\ell_\lambda}
\,
e^{x_{\lambda 1}\,{\sf a}_1^*
+\cdots+
x_{\lambda n}\,{\sf a}_n^*},
\]
where we have 
introduced nonnegative integers $x_{ \lambda i} \in \N$, $i = 1,
\dots, n$ parametrized by $\lambda$. When we expand the product of
all the $\kappa$ sums involved, the exponentiated terms add up and the
obtained sum together with the initial sum $\sum_{ \ell_1 + \cdots +
\kappa \ell_\kappa = m}$ unify as a single big sum:
\[
\footnotesize
\aligned
{\sf ch}\big(\mathcal{E}_{\kappa,m}^{GG}T_X^*\big)
=
\sum_{
\substack{
x_{11}+\cdots+x_{1n}
\\
+\cdots\cdots\cdots\cdots\cdots\cdots+
\\
+\kappa(x_{\kappa 1}+\cdots+x_{\kappa n})=m
}}
\exp
\big\{
(x_{11}+\cdots+x_{\kappa 1}){\sf a}_1^*
+\cdots+
(x_{1n}+\cdots+x_{\kappa n}){\sf a}_n^*
\big\}
\endaligned
\]
in which the $\ell_\lambda$ have been naturally removed, with the only
constraint that $\sum_{ \lambda = 1 }^\kappa\, \lambda\, (x_{ \lambda
1} + \cdots + x_{ \lambda n})$ be constant equal to $m$. Now we
observe the general summation rule:
\[
{\textstyle{\sum_{u_1+u_2+\cdots+u_\mu=m}}} 
\equiv 
{\textstyle{\sum_{u_2+\cdots+u_\mu\leqslant m}}}, 
\]
by simply taking $u_1 := m - u_2 - \cdots - u_\mu$, where the $u_j \in
\N$. Thus, we may eliminate $x_{ 11}$ in our argument of summation
and it follows at once for any $j = 0, 1, \dots, n$ that the quantity
we want to estimate is equal to:
\[
\footnotesize
\aligned
&
\big[
{\sf ch}\big(\mathcal{E}_{\kappa,m}^{GG}T_X^*\big)
\big]_{n-j}
=
\\
&
=
\!\!
\sum_{
\substack{
x_{12}+\cdots+x_{1n}
\\
+\cdots\cdots\cdots\cdots\cdots\cdots+
\\
+\kappa(x_{\kappa 1}+\cdots+x_{\kappa n})\leqslant m
}}
\!\!\!
\frac{1}{(n-j)!}\,
\big[
(\widehat{x_{11}}+x_{21}+\cdots+x_{\kappa 1}){\sf a}_1^*
+\cdots+
(x_{1n}+\cdots+x_{\kappa n}){\sf a}_n^*
\big]^{n-j},
\endaligned
\]
where the symbol $\widehat{ x_{ 11}}$ means that $x_{ 11}$ is replaced
by its value $m - x_{ 12} - \cdots - \kappa \, x_{ \kappa n}$.
Classically, by making the change of variables:
\[
y_{12}
:=
\frac{x_{12}}{m},\,
\dots,\,
y_{1n}
:=
\frac{x_{1n}}{m},\,\,
\dots\dots,\,\,
y_{\kappa 1}
:=
\frac{x_{\kappa 1}}{m},\,
\dots,\,
y_{\kappa n}
:=
\frac{x_{\kappa n}}{m},
\]
the discrete Riemann-like sum just obtained can be approximated by a
continuous integral performed on a $(\kappa n - 1)$-dimensional
simplex against the standard measure of $\R_+^{ \kappa n - 1}$:
\[
\footnotesize
\aligned
\big[
{\sf ch}\big(\mathcal{E}_{\kappa,m}^{GG}T_X^*\big)
\big]_{n-j}
&
=
m^{\kappa n-1+n-j}
\int_{
\substack{
y_{12}+\cdots+y_{1n}
\\
+\cdots\cdots\cdots\cdots\cdots\cdots+
\\
+\kappa(y_{\kappa 1}+\cdots+y_{\kappa n})\leqslant 1
}}\,
dy_{12}\cdots dy_{1n}
\cdots\cdots\,
dy_{\kappa 1}\cdots\,dy_{\kappa n}
\cdot
\\
&
\cdot
\frac{1}{(n-j)!}\,
\big[
(\widehat{y_{11}}+y_{21}+\cdots+y_{\kappa 1}){\sf a}_1^*
+\cdots+
(y_{1n}+\cdots+y_{\kappa n}){\sf a}_n^*
\big]^{n-j}
+
\\
&
\ \ \ \ \ \ \ \ \ \ \ \ \ \ \ \ \ \ \ \ \ \ \ \ \ \ \ \ \ \ \ \ \ \ \ 
\ \ \ \ \ \ \ \ \ \ \ \ \ \ \ \ \ \ \ \ \ \ \ \ \ \ \ \ \ \ \ \ \ \ \ 
\ \ \ \ \ \ \ \ \ \ 
+
{\sf O}_{n,\kappa}
\big(m^{(\kappa+1)n-j-2}\big),
\endaligned
\]
the remainder being automatically at most of the order of the
submaximal power of $m$. The integral remaining being visibly
independent of $m$, the conclusion is got.
\endproof

Consequently, in order to compute asymptotically our Euler-Poincar\'e
characteristic, we only have to estimate the integral above for $j =
0$, in which $\widehat{ y_{ 11}}$ is of course an abbreviation for $1
- y_{ 12} - \cdots - \kappa \, y_{ \kappa n}$. To this aim, we make
the multidilational change of variables: $y_{ \lambda i} \mapsto
\lambda\, y_{ \lambda i} =: z_{ \lambda i}$ and the asymptotic under
study becomes an integral over the {\em standard} $(n\kappa -
1)$-dimensional simplex:
\[
\footnotesize
\aligned
\chi
\big(
X,\,
\mathcal{E}_{\kappa,m}^{GG}T_X^*
\big)
&
=
\int_X\,
\big[{\sf ch}(\mathcal{E}_{\kappa,m}^{GG}T_X^*)\big]_n
+
{\sf O}_{n,\kappa}
\big(m^{(\kappa+1)n-2}\big)
\\
&
\equiv
\frac{m^{(\kappa+1)n\,-1}}{
n!\,\,(\kappa!)^n}\,
\int_{
\substack{
z_{21}+\cdots+z_{1n}
\\
+\cdots\cdots\cdots\cdots\cdots+
\\
+z_{\kappa 1}+\cdots+z_{\kappa n}
\leqslant 1
}}\,
dz_{12}\cdots dz_{1n}
\cdots\cdots\,
dz_{\kappa 1}\cdots\,dz_{\kappa n}
\cdot
\\
&\ \ \ \ \
\cdot
\bigg[
\bigg(\widehat{z_{11}}
+
\frac{z_{21}}{2}
+\cdots+
\frac{z_{\kappa 1}}{\kappa}
\bigg){\sf a}_1^*
+\cdots+
\bigg(
\frac{z_{1n}}{1}
+
\frac{z_{2n}}{2}
+\cdots+
\frac{z_{\kappa n}}{\kappa}
\bigg){\sf a}_n^*
\bigg]^n,
\endaligned
\]
where now the sign ``$\equiv$'' means modulo ${\sf O}_{n, \kappa}
\big( m^{ ( \kappa + 1 )n - 2} \big)$ and where $\widehat{ z_{ 11}} =
1 - z_{ 12} - \cdots - z_{ \kappa n}$. Applying now Newton's
multinomial formula:
\[
\footnotesize
\aligned
\big(Z_1+Z_2+\cdots+Z_n\big)^n
=
\sum_{q_1+q_2+\cdots+q_n=n}\,
\frac{n!}{q_1!\,q_2!\,\cdots\,q_n!}\,
(Z_1)^{q_1}(Z_2)^{q_2}\cdots(Z_n)^{q_n},
\endaligned
\]
we may expand the $n$-th power in the second line above, getting:
\[
\aligned
\chi
\big(
X,\,
\mathcal{E}_{\kappa,m}^{GG}T_X^*
\big)
&
=
\frac{m^{(\kappa+1)n-1}}{
\zero{n!}\,\,(\kappa!)^n}\,
\sum_{q_1+\cdots+q_n=n}
\frac{\zero{n!}}{q_1!\,\cdots\,q_n!}
\,({\sf a}_1^*)^{q_1}\cdots({\sf a}_n^*)^{q_n}
\cdot
\\
&
\ \ \ \ \ 
\cdot
\int_{
\substack{
z_{21}+\cdots+z_{1n}
\\
+\cdots\cdots\cdots\cdots\cdots+
\\
+z_{\kappa 1}+\cdots+z_{\kappa n}
\leqslant 1
}}\,
dz_{21}\cdots dz_{1n}
\cdots\cdots\,
dz_{\kappa 1}\cdots\,dz_{\kappa n}
\cdot
\\
&
\cdot
\bigg(
\widehat{z_{11}}
+
\frac{z_{21}}{2}
+\cdots+
\frac{z_{\kappa 1}}{\kappa}
\bigg)^{q_1}
\cdots\,
\bigg(
\frac{z_{1n}}{1}
+
\frac{z_{2n}}{2}
+\cdots+
\frac{z_{\kappa n}}{\kappa}
\bigg)^{q_n}.
\endaligned
\]
The $n!$ drops, a fact denoted with the symbol ``$\zero{ \,\, \,}$''.
Furthermore, in the integral\,\,---\,call it ${\sf I}_{ q_1, \dots,
q_n}$\,\,---\,\,which appears naturally in the last two lines, we yet
expand the $q_1$-th, \dots, the $q_n$-th powers:
\[
\footnotesize
\aligned
{\sf I}_{q_1,\dots,q_n}
&
=
\sum_{q_{11}+q_{21}+\cdots+q_{\kappa 1}=q_1}
\,\cdots\,
\sum_{q_{1n}+q_{2n}+\cdots+q_{\kappa n}}\,
\frac{q_1!}{q_{11}!\,q_{21}!\,\cdots\,q_{\kappa 1}!}
\,\cdots\,
\frac{q_n!}{q_{1n}!\,q_{2n}!\,\cdots\,q_{\kappa n}!}
\,\cdot
\\
&
\ \ \ \ \ \ \ \ \ \ \ \ \ \ \ \ \ \ \ \ \ \ \ \ \ \ \ \ \ \ \
\ \ \ \ \ \ \ \ \ \ \ \ \ \ \ \ \ \ \ \ \ \ \ \ \ \ \ \ \ \ \
\cdot
\frac{1}{(2)^{q_{21}}\,\cdots\,(\kappa)^{q_{\kappa 1}}}
\,\cdots\,
\frac{1}{(1)^{q_{1n}}\,(2)^{q_{2n}}\,\cdots\,(\kappa)^{q_{\kappa n}}}
\,\cdot
\\
&
\cdot
\int_{
\substack{
z_{21}+\cdots+z_{1n}
\\
+\cdots\cdots\cdots\cdots\cdots+
\\
+z_{\kappa 1}+\cdots+z_{\kappa n}
\leqslant 1
}}\,
dz_{21}\cdots dz_{1n}
\cdots\cdots\,
dz_{\kappa 1}\cdots\,dz_{\kappa n}
\cdot
\\
&
\ \ \ \ \ \ \ \ \ \ \ \ \ \ \ \ \ \ \ \ \ \ \ \ \ \ 
\cdot
(\widehat{z_{11}})^{q_{11}}\,(z_{21})^{q_{21}}
\,\cdots\,
(z_{\kappa 1})^{q_{\kappa 1}}
\,\,\cdots\cdots\,\,
(z_{1n})^{q_{1n}}\,(z_{2n})^{q_{2n}}
\,\cdots\,
z_{\kappa n}^{q_{\kappa n}}.
\endaligned
\]

\begin{Lemma}
\label{j1-j2-jp}
For any integer $p \geqslant 2$ and for any nonnegative integer
exponents $j_1, j_2, \dots, j_p \in \N$, one has:
\[
\footnotesize
\aligned
\int_{u_2+\cdots+u_p\leqslant 1
\atop
u_2\geqslant 0,\,\,\dots,\,\,u_p\geqslant 0}\,
[1-u_2-\cdots-u_p]^{j_1}u_2^{j_2}\cdots u_p^{j_p}\,
du_2\cdots du_p
=
\frac{j_1!\,\,j_2!\,\,\cdots\,\,j_p!}{
(j_1+j_2+\cdots+j_p+p-1)!}.
\endaligned
\]
\end{Lemma}

\proof
By decomposing the integrations, we may write this integral as:
\[
\int_0^1
\!\!\!
u_2^{j_2}\,du_2
\int_0^{1-u_2}
\!\!\!\!\!\!\!\!
u_3^{j_3}\,du_3
\,\cdots\cdots\,
\int_0^{1-u_2-\cdots-u_{p-1}}
\!\!\!\!\!\!\!\!\!\!\!\!\!\!\!\!\!\!\!\!\!\!\!\!\!\!\!
(1-u_2-\cdots-u_{p-1}-u_p)^{j_1}\,u_p^{j_p}\,du_p
=:
{\sf J}_{j_1,j_2,j_3,\dots,j_p}^p.
\]
Taking $j_p + 1$ times the primitive of the first factor in the last
integral and integrating successively by parts, this integral in
question receives the value:
\[
\aligned
&
\Big[
-\,
{\textstyle{\frac{(1-u_2-\cdots-u_{p-1}-u_p)^{j_1+j_p+1}}{
(j_1+1)\cdots(j_1+j_p)(j_1+j_p+1)}}}\,
j_p!
\Big]_0^{1-u_2-\cdots-u_{p-1}}
=
\\
&
\ \ \ \ \ \ \ \ \ \ \ \ \ \ \ \ \ \ \ \ \
=
{\textstyle{\frac{j_1!}{(j_1+j_p+1)!}}}\,
(1-u_2-\cdots-u_{p-1})^{j_1+j_p+1}\,j_p!.
\endaligned
\]
Thus, the case $p = 2$ is settled. If $p \geqslant 3$, inserting this
value just computed:
\[
\aligned
{\sf J}_{j_1,j_2,\dots,j_{p-1},j_p}^p
&
=
{\textstyle{\frac{j_1!\,\,j_p!}{(j_1+j_p+1)!}}}\,
{\sf J}_{j_1+j_p+1,j_2,\dots,j_{p-1}}^{p-1}
\\
&
=
{\textstyle{\frac{j_1!\,\,j_p!}{\zero{(j_1+j_p+1)!}}}}\,
{\textstyle{\frac{\zero{(j_1+j_p+1)!}\,j_2!\,\cdots\,j_{p-1}!}{
(j_1+j_p+1+j_2+\cdots+j_{p-1}+p-2)!}}},
\endaligned
\]
we get without effort the general conclusion by induction on $p$.
\endproof

So applying this elementary lemma, we may finish to compute our
integral:
\[
\footnotesize
\aligned
{\sf I}_{q_1,\dots,q_n}
&
=
\sum_{q_{11}+q_{21}+\cdots+q_{\kappa 1}=q_1}
\!\!\!\cdots\!\!\!
\sum_{q_{1n}+q_{2n}+\cdots+q_{\kappa n}}\,
\frac{q_1!}{\zero{q_{11}!\,q_{21}!\,\cdots\,q_{\kappa 1}!}}
\,\cdots\,
\frac{q_n!}{\zero{q_{1n}!\,q_{2n}!\,\cdots\,q_{\kappa n}!}}
\,\cdot
\\
&
\ \ \ \ \ \ \ \ \ \ \ \ \ \ \ \ \ \ \ \ \ \ \ \ \ \ \ \ \ \ \
\ \ \ \ \ \ \ \ \ \ \ \ \ \ \ \ \ \ \ \
\cdot
\frac{1}{(2)^{q_{21}}!\,\cdots\,(\kappa)^{q_{\kappa 1}}}
\,\cdots\,
\frac{1}{(1)^{q_{1n}}\,(2)^{q_{2n}}\,\cdots\,(\kappa)^{q_{\kappa n}}}
\,\cdot
\\
&
\ \ \ \ \ \ \ \ \ \ \ \ \ \ \ \ \
\cdot
\frac{\zero{q_{11}!\,q_{21}\,\cdots\,q_{\kappa 1}!}
\,\,\cdots\cdots\,\,\zero{q_{1n}!\,q_{2n}!\,\cdots\,q_{\kappa n}!}}{
(q_{11}+q_{21}+\cdots+q_{\kappa 1}
+\cdots\cdots+
q_{1n}+q_{2n}+\cdots+q_{\kappa n}
+\kappa n-1)!}.
\endaligned
\]
Remarkably, all the factorials $q_{ \lambda i}!$ drop. Furthermore,
the big factorial in the denominator visibly simplifies as
\[
(q_1+\cdots+q_n+\kappa n-1)!
=
(n+\kappa n-1)!, 
\]
and we get a formula for ${\sf I}_{ q_1, \dots, q_n}$ in which it will
appear soon to be convenient to reconstitute a product of $n$
independent big sums, and to this
aim, we add in advance the innocuous factor $\frac{ 1}{ (1)^{
q_{ 11}}}$:
\[
\scriptsize
\aligned
&
{\sf I}_{q_1,\dots,q_n}
=
\\
&
=
\frac{q_1!\,\cdots\,q_n!}{
(q_1+\cdots+q_n+\kappa n-1)!}
\!\!
\sum_{q_{11}+\cdots+q_{\kappa 1}=q_1}\!\!
\cdots\!\!
\sum_{q_{1n}+\cdots+q_{\kappa n}=q_n}\!\!
\frac{1}{(1)^{q_{11}}\cdots (\kappa)^{q_{\kappa 1}}}\,
\cdots\,
\frac{1}{(1)^{q_{1n}}\cdots (\kappa)^{q_{\kappa n}}}
\\
&
=
\frac{q_1!\,\cdots\,q_n!}{
((\kappa+1)n-1)!}
\bigg(
\sum_{q_{11}+\cdots+q_{\kappa 1}=q_1}\,
\frac{1}{(1)^{q_{11}}\cdots (\kappa)^{q_{\kappa 1}}}\,
\bigg)
\cdots
\bigg(
\sum_{q_{1n}+\cdots+q_{\kappa n}=q_n}\,
\frac{1}{(1)^{q_{1n}}\cdots (\kappa)^{q_{\kappa n}}}
\bigg). 
\endaligned
\]
Now, when we plug this formula in the computation of $\chi \big( X, \,
\mathcal{ E}_{ \kappa, m}^{ GG} T_X^* \big)$ that we interrupted
before stating the lemma, all the factorials $q_1!, \dots, q_n!$
appear once at a numerator place and once at a denominator place, so
they drop all and we finally get:
\[
\scriptsize
\aligned
\chi\big(X,\,\,\mathcal{E}_{\kappa,m}^{GG}T_X^*\big)
&
=
\frac{m^{(\kappa+1)n-1}}{
(\kappa!)^n\,((\kappa+1)n-1)!}\,
\bigg[
\sum_{q_1+\cdots+q_n=n}\,
({\sf a}_1^*)^{q_1}\cdots({\sf a}_n^*)^{q_n}
\cdot
\\
&
\ \ \ \ \
\cdot\,
\bigg(
\sum_{q_{11}+\cdots+q_{\kappa 1}=q_1}\,
\frac{1}{(1)^{q_{11}}\cdots (\kappa)^{q_{\kappa 1}}}\,
\bigg)
\cdots
\bigg(
\sum_{q_{1n}+\cdots+q_{\kappa n}=q_n}\,
\frac{1}{(1)^{q_{1n}}\cdots (\kappa)^{q_{\kappa n}}}
\bigg)
\bigg]
+
\\
&
\ \ \ \ \ \ \ \ \ \ \ \ \ \ \ \ \ \ \ \ \ \ \ \ \ \ \ \ \ \ \ 
\ \ \ \ \ \ \ \ \ \ \ \ \ \ \ \ \ \ \ \ \ \ \ \ \ \ \ \ \ \ \ 
\ \ \ \ \ \ \ \ \ \ \ \ \ \ \ \ \ \ \ \ \ \ \ \ \ \ \ \ \ \ \ 
\ \ \ \ \ \ \ \ \ \ \ \ 
+
{\sf O}
\big(m^{(\kappa+1)n-2}\big).
\endaligned
\]
We therefore have to deal with the asymptotic character, as $\kappa
\to \infty$, of the poly-logarithmic sums of the type:
\[
\Sigma_1^\kappa(q)
:=
\sum_{q_1+\cdots+q_\kappa=q
\atop
q_1\geqslant 0,\,\,\cdots,\,\,q_\kappa\geqslant 0}\,
\frac{1}{(1)^{q_1}\cdots(\kappa)^{q_{\kappa}}},
\]
where $q \in \N$ is arbitrary.

\begin{Lemma}
\label{Sigma-1-kappa-q}
As $\kappa \to \infty$, one has:
\[
\Sigma_1^\kappa(q)
=
\frac{(\log\kappa)^q}{q!}
+
{\sf O}_n\big((\log\kappa)^{q-1}\big).
\]
\end{Lemma}

\proof
Easily re-doable,
and in fact also known in the literature on polylogarithms
(\cite{ cemp1999}).
\endproof

From this last lemma, it follows at once that:
\[
\footnotesize
\aligned
&
\chi\big(X,\,\,\mathcal{E}_{\kappa,m}^{GG}T_X^*\big)
=
\frac{m^{(\kappa+1)n-1}}{
((\kappa+1)n-1)!\,(\kappa!)^n}\,
\bigg[
\sum_{q_1+\cdots+q_n=n}\,
({\sf a}_1^*)^{q_1}\cdots({\sf a}_n^*)^{q_n}
\cdot
\\
&
\ \ \ \ \ \ \ \ \ \ \ \ \ \ \ \ \ \ \ \ \ \ \ \ \ \ \ \ \ \ \ \ \ \ \ 
\cdot
\frac{(\log\kappa)^{q_1}}{q_1!}
\,\cdots\,
\frac{(\log\kappa)^{q_n}}{q_n!}
+
{\sf O}_n\big((\log\kappa)^{n-1}\big)
\bigg]
+
{\sf O}_{n,\kappa}
\big(m^{(\kappa+1)n-2}\big)
\\
&
=
{\textstyle{\frac{m^{(\kappa+1)n-1}}{
(\kappa!)^n\,((\kappa+1)n-1)!}}}\,
\Big[
({\sf a}_1^*+\cdots+{\sf a}_n^*)^n\,
{\textstyle{\frac{(\log\kappa)^n}{n!}}}
+
{\sf O}_n\big((\log\kappa)^{n-1}\big)
\Big]
+
{\sf O}_{n,\kappa}
\big(m^{(\kappa+1)n-2}\big),
\endaligned
\]
so the asymptotic formula exhibited in the theorem is established. To
conclude the proof, one easily convinces oneself by inspecting the
remainders that they indeed have the form claimed in {\bf (i)} and
{\bf (ii)}.
\endproof

\begin{Openproblem} 
{\sl Applying the concepts and the
combinatorics partly achieved in~\cite{ cemp1999, bbbl2001,
wald2000}, find {\em
closed explicit formulas} firstly for the remainder terms ${\sf O}_n
\big( (\log \kappa)^{ n-1} \big)$, secondly, for the remainder terms
${\sf O}_{ n, \kappa} \big( m^{ (\kappa + 1)n - 2} \big)$. As an
accessible preliminary, study the $\Sigma_q ( \kappa)$ completely.}
\end{Openproblem}

\markleft{Jo\"el Merker}
\markright{\sf \S4.~Exact Schur Bundle Decomposition}
\section{\bf Exact Schur Bundle Decomposition}
\label{Section-4}

\subsection{Schur bundles and Pieri rule}
Thanks to the filtration provided by Theorem~\ref{graded-sym} 
and to the basic cohomology inequalities
reproved in Section~2, the study of the Green-Griffiths jet bundle can
in principle be led back to the study of multitensored symmetric
powers:
\[
{\rm Sym}^{\ell_1}T_X^*
\otimes
{\rm Sym}^{\ell_2}T_X^*
\otimes\cdots\otimes
{\rm Sym}^{\ell_\kappa}T_X^*
\]
of the cotangent bundle. But it is known since the works of Isai
Schur at the turn to the 20\textsuperscript{th} century that these
multitensored bundles can even be decomposed in more atomic
independent bricks.

Since the complex linear group ${\sf GL}_n ( \C)$ acts naturally on
$T_X^*$ and on all of its tensor powers $\big( T_X^*\big)^{\otimes r}$
as well ($r = 1, 2, 3, \dots$), then by fundamental facts of
representation theory (Schur's theorems), it follows that the (in fact
complicated) direct sum ${\sf Gr}^\bullet \mathcal{E}_{\kappa, m}^{
GG}T_X^*$ provided by Theorem~\ref{graded-sym} can in
principle be represented as a certain direct sum of the so-called {\sl
Schur bundles}:
\[
\mathcal{S}^{(\ell_1,\ell_2,\dots,\ell_n)}T_X^*, 
\]
in which $\ell_1 \geqslant \ell_2 \geqslant \cdots \geqslant \ell_n
\geqslant 0$; we employ the notation of~\cite{ rou2006b} and the
reader is referred to the works of Br\"uckmann~\cite{ bru1972,
brurac1990, bru1997} and to the monographs~\cite{ fuha1991, stu1993,
macd1995,
krpr1996, pro2007} for background material, or alternatively to
Section~\ref{Section-11} below. In order to determine which
$\mathcal{S}^{ ( \ell_1, \dots, \ell_n)} T_X^*$ appear in ${\sf
Gr}^\bullet \mathcal{E}_{ \kappa, m}^{ GG}T_X^*$, possibly with some
multiplicity $\geqslant 1$, two options present themselves.

The first option would be to apply step by step the so-called {\sl
Pieri formula} (\cite{ fuha1991}, p.~455) to the 
direct sum representation:
\[
\footnotesize
\aligned
{\sf Gr}^\bullet
\mathcal{E}_{\kappa,m}^{GG}T_X^*
=
\bigoplus_{\ell_1+2\ell_2+\cdots+\kappa\ell_\kappa=m}\,
\mathcal{S}^{(\ell_1,0,\dots,0)}T_X^*
\otimes
\mathcal{S}^{(\ell_2,0,\dots,0)}T_X^*
\otimes\,\cdots\,\otimes
\mathcal{S}^{(\ell_\kappa,0,\dots,0)}T_X^*.
\endaligned
\]
Pieri indeed provides a neat combinatorial rule for representing any
tensor product of a Schur bundle with a symmetric power as a certain
direct sum of well controlled Schur bundles over $X$:
\begin{equation}
\label{pieri-rule}
\mathcal{S}^{(t_1,\dots,t_n)}T_X^* 
\otimes
\mathcal{S}^{(\ell,0,\dots,0)}T_X^*
=
\sum_{s_1+\cdots+s_n=\ell+t_1+\cdots+t_n
\atop
s_1\geqslant t_1\geqslant s_2\geqslant t_2
\geqslant\cdots\geqslant s_n\geqslant t_n\geqslant 0}\,
\mathcal{S}^{(s_1,\dots,s_n)}T_X^*. 
\end{equation}
However, when one tries to induct on such a formula, the complexity
increases dramatically as soon the number $\kappa$ of tensor factors
in ${\sf Gr}^\bullet \mathcal{ E}_{ \kappa, m}^{ GG} T_X^*$ passes
above $\kappa = 5$, even in dimension $n = 2$, and apparently, nothing
really effective or exploitable for us exists in the literature.

\subsection{Invariant theory approach}
The second option, more direct and more suited to asymptotic
approximations, consists in interpreting the problem directly in terms
of classical invariant theory, starting with the original
definition~\thetag{ \ref{gg-jets}}
of Green-Griffiths jets. Indeed,
the general $n\times n$ complex unipotent matrix:
\[
{\sf u}
:=
\left(
\begin{array}{ccccc}
1 & 0 & 0 & \cdots & 0
\\
{\sf u}_{21} & 1 & 0 & \cdots & 0 
\\
{\sf u}_{31} & {\sf u}_{32} & 1 & \cdots & 0
\\
\cdot\cdot & \cdot\cdot & \cdot\cdot & \cdots & \cdot\cdot
\\
{\sf u}_{n1} & {\sf u}_{n2} & {\sf u}_{n3} & \cdots & 1
\end{array}
\right),
\]
where the ${\sf u}_{ ij} \in \C$ are arbitrary complex numbers, acts
naturally and linearly on all the jet variables in such a way that
for any jet level $\lambda$ with $1 \leqslant \lambda \leqslant
\kappa$, one sets in matrix notation:
\[
g^{(\lambda)}
:=
{\sf u}\cdot f^{(\lambda)}\,,
\] 
that is to say in greater length: 
\[
\left\{
\aligned
g_1^{(\lambda)}
&
:=
f_1^{(\lambda)}
\\
g_2^{(\lambda)}
&
:=
f_2^{(\lambda)}
+
{\sf u}_{21}\,f_1^{(\lambda)}
\\
g_3^{(\lambda)}
&
:=
f_3^{(\lambda)}
+
{\sf u}_{32}\,f_2^{(\lambda)}
+
{\sf u}_{31}\,f_1^{(\lambda)}
\\
\cdots\cdot
&
\cdots\cdots\cdots\cdots\cdots\cdots\cdots\cdots\cdots\cdots
\\
g_n^{(\lambda)}
&
=
f_n^{(\lambda)}
+
{\sf u}_{n,n-1}\,f_{n-1}^{(\lambda)}
+
\cdots
+
{\sf u}_{n1}\,f_1^{(\lambda)}. 
\endaligned\right.
\]
A general fact from the classical representation theory of ${\sf GL}_n
( \C)$ states that the so-called {\sl vectors of highest weight}
identify precisely to those that remain invariant by this unipotent
action, namely to jet polynomials ${\sf P} \big( j^\kappa f\big)$
which satisfy the invariancy condition:
\[
{\sf P}\big(j^\kappa g\big)
=
{\sf P}\big({\sf u}\cdot j^\kappa f\big)
\equiv
{\sf P}\big(j^\kappa f\big),
\]
for every unipotent matrix ${\sf u} \in {\sf U}_n ( \C)$. Furthermore
and most importantly, there is a one-to-one correspondence between the
vectors of highest weight and the Schur bundles appearing in the
decomposition of ${\sf Gr}^\bullet \mathcal{E}_{ \kappa, m}^{ GG}T_X^*$,
the rule being as follows. Precisely speaking, the vector
space of unipotent-invariant polynomials (vectors of highest weight)
is shown to decompose as a direct sum of (linearly independent)
one-dimensional spaces generated by vectors ${\sf Q} = {\sf Q} \big(
j^\kappa f \big)$ that are eigenvalues for the action ${\sf e} \cdot
f_i^{ (\lambda)} := {\sf e}_i \, f_i^{ (\lambda)}$ of all diagonal
matrices of the form:
\[
{\sf e}
:=
\left(
\begin{array}{cccc}
{\sf e}_1 & 0 & \cdots & 0
\\
0 & {\sf e}_2 & \cdots & 0 
\\
\cdot\cdot & \cdot\cdot & \cdots & \cdot\cdot
\\
0 & 0 & \cdots & {\sf e}_n
\end{array}
\right),
\]
where ${\sf e}_1, {\sf e}_2, \dots, {\sf e}_n$ are arbitrary complex
numbers, so that there are certain characteristic exponents $\ell_i$
with the property that:
\[
{\sf Q}\big({\sf e}\cdot j^\kappa f)
=
({\sf e}_1)^{\ell_1}({\sf e}_2)^{\ell_2}\cdots ({\sf e}_n)^{\ell_n}
\,
{\sf Q}\big(j^\kappa f\big). 
\]
One shows (\cite{ ful1991})
that $\ell_1 \geqslant \ell_2 \geqslant \cdots \geqslant
\ell_n \geqslant 0$ and that such an eigenvector ${\sf Q}$ (of highest
weight) is precisely linked to the Schur bundle $\mathcal{S}^{ (
\ell_1, \ell_2, \dots, \ell_n)} T_X^*$ which corresponds to an
irreducible representation on a fiber over a point $x \in X$. Of
course, a specific Schur bundle $\mathcal{S}^{ ( \ell_1, \ell_2,
\dots, \ell_n)} T_X^*$ could well occur several times in the sought
decomposition of ${\sf Gr}^\bullet \mathcal{E}_{ \kappa, m}^{ GG}
T_X^*$, hence have a certain multiplicity $\geqslant 2$, because some
different linearly independent ${\sf Q}$'s could share the same
characteristic exponents $\ell_i$. In fact, this will indeed be the
case below, and determining such multiplicities, at least
asymptotically as $\kappa \to \infty$, will be crucial for us.

\subsection{Serendipity}
The knowledge of the algebra of invariants of the full unipotent group
${\sf U}_n(\C) \subset {\sf GL}_n ( \C)$ dates back to the nineteenth
century. As a matter of fact, the following four basic statements
Theorems~\ref{Theorem-A},
\ref{Theorem-B}, \ref{Theorem-C} and~\ref{Theorem-C} 
below, which will precede a main starting
theorem specially designed for our future purposes, are essentially
known and they are established in various sources.

\begin{Theorem}
\label{Theorem-A} 
{\rm (\cite{ stu1993, krpr1996, mist2005, pro2007})}
The algebra of jet
polynomials invariant under the above action of the full unipotent
group ${\sf U}_n ( \C) \subset {\sf GL}_n ( \C)$ is {\em generated},
{\em as an algebra}, by the collection of all the determinants
(minors):
\[
\small
\aligned
&
\left\vert
f_1^{(\lambda_1)}
\right\vert
\!=:\!
\Delta_1^{\lambda_1},
\ \ \ 
\left\vert\!
\begin{array}{cc}
f_1^{(\lambda_1)} & f_2^{(\lambda_1)}
\\
f_1^{(\lambda_2)} & f_2^{(\lambda_2)}
\end{array}
\!\right\vert
\!=:\!
\Delta_{1,2}^{\lambda_1,\lambda_2},
\ \ \ 
\left\vert\!
\begin{array}{ccc}
f_1^{(\lambda_1)} & f_2^{(\lambda_1)} & f_3^{(\lambda_1)}
\\
f_1^{(\lambda_2)} & f_2^{(\lambda_2)} & f_3^{(\lambda_2)}
\\
f_1^{(\lambda_3)} & f_2^{(\lambda_3)} & f_3^{(\lambda_3)}
\end{array}
\!\right\vert
\!=:\!
\Delta_{1,2,3}^{\lambda_1,\lambda_2,\lambda_3},
\\
&
\cdots\cdots,\,\,\,
\left\vert\!
\begin{array}{cccc}
f_1^{(\lambda_1)} & f_2^{(\lambda_1)} & \cdots & f_n^{(\lambda_1)}
\\
f_1^{(\lambda_2)} & f_2^{(\lambda_2)} & \cdots & f_n^{(\lambda_2)}
\\
\cdots & \cdots & \cdots & \cdots
\\
f_1^{(\lambda_n)} & f_2^{(\lambda_n)} & \cdots & f_n^{(\lambda_n)}
\end{array}
\!\right\vert
\!=:\!
\Delta_{1,2,\dots,n}^{\lambda_1,\lambda_2,\dots,\lambda_n},
\endaligned
\] 
in which the jet orders $1 \leqslant \lambda_1 < \lambda_2 <
\dots < \lambda_n \leqslant \kappa$ are all 
arbitrary
\footnote{\,
It is only necessary to consider strictly increasing integers
$\lambda_l$, since for every $i$ with $1 \leqslant i \leqslant n$, and
for every permutation $\sigma$ of $\{ 1, 2, \dots, i\}$ one clearly has:
\[
\Delta_{1,2,\dots,i}^{\lambda_{\sigma(1)},\lambda_{
\sigma(2)},\dots,\lambda_{\sigma(i)}}=(-1)^{{\rm sign}(
\sigma)}\,\Delta_{1,2,\dots,i}^{\lambda_1,\lambda_2,\dots,
\lambda_i}.
\] 
}, 
and these determinants are visibly invariant with
respect to the ${\sf U}_n ( \C)$-action.
\end{Theorem}

However, although all the determinants in question happen to be
linearly independent, one cannot just pretend that the whole
unipotent-invariant algebra identifies with the plain polynomial
algebra:
\[
\C\big[\Delta_1^{\lambda_1},\Delta_{1,2}^{\lambda_1,
\lambda_2},\dots,\Delta_{1,2,\dots,n}^{\lambda_1,\dots,
\lambda_n}\big],
\]
because several elementary {\em non}linear
relations exist between all these determinants; for instance, there
exist the basic quadratic Pl\"ucker relations\footnote{\,
\, ---\,\,the knowledge of which surely goes back to the seventeenth
century theory, at a time when elimination was the main tool in the
search for solving algebraic equations of degrees $2$, $3$, $4$ and
$5$.
} 
of first and of second type:
\[
\aligned
0
&
\equiv
\Delta_1^{\lambda_1}\,
\Delta_{1,2}^{\lambda_2,\lambda_3}
+
\Delta_1^{\lambda_3}\,
\Delta_{1,2}^{\lambda_1,\lambda_2}
+
\Delta_1^{\lambda_2}\,
\Delta_{1,2}^{\lambda_3,\lambda_1}
\\
0
&
\equiv
\Delta_{1,2}^{\lambda_1,\lambda_2}\,
\Delta_{1,2}^{\lambda_3,\lambda_4}
+
\Delta_{1,2}^{\lambda_1,\lambda_4}\,
\Delta_{1,2}^{\lambda_2,\lambda_3}
+
\Delta_{1,2}^{\lambda_1,\lambda_3}\,
\Delta_{1,2}^{\lambda_4,\lambda_2},
\endaligned
\] 
so that the binomial in the first line 
$\Delta_1^{ \lambda_1} \Delta_{ 1, 2}^{
\lambda_2, \lambda_3}$, viewed in the plain polynomial algebra $\C
\big[ \Delta_1^{ \lambda_1}, \, \Delta_{ 1, 2}^{ \lambda_1, \lambda_2}
\big]$, would possess {\em two} distinct representations: itself, and:
\[
-\,
\Delta_1^{\lambda_3}\,\Delta_{1,2}^{\lambda_1,\lambda_2}
-
\Delta_1^{\lambda_2}\,\Delta_{1,2}^{\lambda_3,\lambda_1}.
\]
Fortunately, the ideal of all relations between these
$\Delta$-determinants is also completely known and
understood. However, presenting {\em explicitly} this ideal of all
relations requires a bit of preparation and a few more indices.

\subsection{Ideal of relations between all $\Delta$
jet-determinants}
We therefore consider the collection of all determinants $\Delta_{ 1,
2,\dots, i}^{ \lambda_1, \lambda_2, \dots, \lambda_i}$ for every $i
\in \{ 1, \dots, n\}$ and for every choice of $i$ jet-line indices
$\lambda_1, \lambda_2, 
\dots, \lambda_i \in \{ 1, \dots, \kappa\}$. At first, we equip this
collection with a {\em partial} order by declaring that:
\[
\Delta_{1,2,\dots,i}^{\lambda_1,\lambda_2,\dots,\lambda_i}
<_{one}
\Delta_{1,2,\dots,j}^{\mu_1,\mu_2,\dots,\mu_j}
\]
if firstly:
\[
i\geqslant j
\]
and if secondly all the following inequalities
hold: 
\begin{equation}
\label{inq-lambda-mu}
\lambda_1\leqslant\mu_1,\ \ \ \ \
\lambda_2\leqslant\mu_2,\ \ \ \ \
\dots\dots,\ \ \ \ \
\lambda_j\leqslant\mu_j.
\end{equation}
Not all determinants are comparable for this order, {\em e.g.}
$\Delta_{ 1, 2}^{ 1,4}$ and $\Delta_{ 1, 2}^{ 2, 3}$ are {\em
in}comparable, and similarly, $\Delta_{ 1, 2}^{ 1, 4}$ and $\Delta_{
1, 2, 3}^{ 2, 3, 4}$ are incomparable too.
We will now see that there is a one-to-one
correspondence between incomparable
$\Delta$-determinants and (generalized) Pl\"ucker relations. 

Thus, let us pick any two general determinants $\Delta_{ 1, \dots, i}^{
\lambda_1, \dots, \lambda_i}$ and $\Delta_{ 1, \dots, j}^{ \mu_1,
\dots, \mu_j}$ that are incomparable and distinct. Permuting the pair
if necessary, we may assume that $i \geqslant j$. Furthermore, if $i
= j$, we may also assume without loss of generality that $(\lambda_1,
\dots, \lambda_i)$ is smaller than $(\mu_1, \dots, \mu_{i=j})$ in the
lexicographic ordering, namely there exists an index $s
\in \{ 1, \dots, i = j\}$ such that: 
\[
\lambda_1=\mu_1,\ \ \ \ \
\dots\dots,\ \ \ \ \
\lambda_{s-1}=\mu_{s-1},\ \ \ \ \
\lambda_s<\mu_s.
\]
Therefore in both cases $i > j$ and $i = j$,
we at least insure by these preliminary choices that:
\[
\Delta_{1,\dots,i}^{\lambda_1,\dots,\lambda_i}
\not>_{one}
\Delta_{1,\dots,j}^{\mu_1,\dots,\mu_j}. 
\]
Since by assumption, these two determinants
are incomparable, the reverse inequality must also fail:
\[
\Delta_{1,\dots,i}^{\lambda_1,\dots,\lambda_i}
\not<_{one}
\Delta_{1,\dots,j}^{\mu_1,\dots,\mu_j},
\]
and hence in the two cases $i > j$ and $i = j$, there must exist a
smallest index $t \in \{ 1, \dots, j\}$ such that:
\[
\lambda_1\leqslant \mu_1,\ \ \ \ \
\dots\dots,\ \ \ \
\lambda_{t-1}\leqslant\mu_{t-1},\ \ \ \ \
\lambda_t
> 
\mu_t, 
\]
because if otherwise all the inequalities~\thetag{
\ref{inq-lambda-mu}} would hold, one would have $\Delta_{1, \dots,
i}^{ \lambda_1, \dots, \lambda_i} <_{one} \Delta_{ 1, \dots, j}^{
\mu_1, \dots, \mu_j}$. In the case $i = j$, it is clear that $t$ can
only be $\geqslant s + 1$.

Remind that in any circumstance, the jet-line indices of the
determinants are strictly increasing:
\[
\lambda_1<\cdots<\lambda_t<\cdots<\lambda_i
\ \ \ \ \ \ \ \
\text{\rm and}
\ \ \ \ \ \ \ \
\mu_1<\cdots<\mu_t<\cdots<\mu_j.
\]
Diagrammatically, we may then represent a set of inequalities with a
pivotal solder, at the index $t$, between the $\mu_i$ and the
$\lambda_i$:
\[
\mu_1<\cdots<\mu_t
\underset{\sf solder}{<}
\lambda_t<\cdots<\lambda_i,
\]
by exhibiting, in two adjusted lines, the vertical spot where the join
takes place:
\[
\aligned
&
\mu_1<\mu_2<\cdots<\,\mu_t\,\,\,\,\,\,\,\,\,<
\underline{\mu_{t+1}\!<\mu_{t+2}<\cdots<\mu_j}_{\sf \,FIX}
\\
&
\ \ \ \ \ \ \ \ \
\underline{\lambda_1<\cdots<\!\lambda_{t-1}}_{\sf \,FIX}
\!<\,\lambda_t\,<\lambda_{t+1}
<\cdots<\lambda_{j-1}<\lambda_j
<\cdots<\lambda_i.
\endaligned
\]
Letting now $\pi \in \mathfrak{ S}_{ i+1}$ be any permutation of the
set $\{ 1, 2, \dots, i, i+1\}$ with $i + 1$ elements, we shall let it
act on the $i+1$ elements that are {\em not} underlined, so that $\pi$
transforms the $i+1$ integers:
\[
\aligned
\mu_1<\mu_2<\cdots<\mu_t
&
<
\\
&
<\lambda_t<\lambda_{t+1}
<\cdots<\lambda_{j-1}<\lambda_j<\cdots<\lambda_i
\endaligned
\]
to the $i+1$ permuted integers (not anymore necessarily
ordered increasingly): 
\[
\aligned
\big(\pi(\mu_1),\,\pi(\mu_2),\,\cdots,\,\pi(\mu_t),\,
&
\\
&
\pi(\lambda_t),\,\pi(\lambda_{t+1}),\,\dots,\,
\pi(\lambda_{j-1}),\,\pi(\lambda_j),\,\dots,\,\pi(\lambda_i)
\big). 
\endaligned
\]
Since our $\Delta$-determinants are skew-symmetric with respect to any
permutation of their lines, it is convenient to restrict attention
only to those permutations that respect strict ordering in the two
blocks:
\[
\aligned
&
\pi(\mu_1)<\pi(\mu_2)<\cdots<\pi(\mu_t)
\\
\text{\rm and:}\ \ \ \ \
&
\pi(\lambda_t)<\pi(\lambda_{t+1})<\dots<
\pi(\lambda_{j-1})<\pi(\lambda_j)<\dots<\pi(\lambda_i).
\endaligned
\]
At last, we are in a position to write down the most general quadratic
Pl\"ucker relations that are fundamental for the subject.

\begin{Theorem}
\label{Theorem-B} 
{\rm (\cite{ stu1993, krpr1996, mist2005, pro2007})}
For any two determinants $\Delta_{ 1, \dots, i}^{ \lambda_1,
\dots, \lambda_i}$ and $\Delta_{ 1, \dots, j}^{ \mu_1, \dots, \mu_j}$
with $i \geqslant j$ that are incomparable with respect to the partial
ordering ``$<_{one}$'', namely which have the concrete properties that:

\begin{itemize}

\smallskip\item[$\bullet$]
when $i > j$, there exists an index $t \in \{ 1, \dots, j\}$
such that:
\[
\lambda_1\leqslant\mu_1,\ \ \ \ \
\dots\dots,\ \ \ \ \
\lambda_{t-1}\leqslant\mu_{t-1},\ \ \ \ \
\text{\em but:}\ \ \ \ \
\lambda_t>\mu_t;
\]

\smallskip\item[$\bullet$]
when $i = j$, there exist two indices $s \in \{ 1, \dots, 
j\}$ and $t \in \{ 1, \dots, j\}$ with $t \geqslant s + 1$ such
that:
\[
\aligned
\lambda_1=\mu_1,\ \
\dots,\ \
&
\lambda_{s-1}=\mu_{s-1},\ \
\lambda_s<\mu_s,\ \
\\
&
\lambda_{s+1}\leqslant\mu_{s+1},\ \
\dots,\ \
\lambda_{t-1}\leqslant\mu_{t-1},\ \ 
\text{\em but again:}\ \
\lambda_t>\mu_t;
\endaligned
\]

\end{itemize}\smallskip

\noindent
the following general
quadratic (Pl\"ucker) relation holds identically in
the ground ring $\C \big[ f_{i_1}', 
f_{i_2}'', \dots, f_{i_\kappa}^{ (\kappa)} \big]$:
\[
\boxed{
\aligned
0
&
\equiv
\sum_{\pi\in\mathfrak{S}_{i+1}}\,
\sum_{
\pi(\lambda_t)<\cdots<\pi(\lambda_i)
\atop
\pi(\mu_1)<\cdots<\pi(\mu_t)}
{\sf sign}(\pi)\cdot
\Delta_{1,\dots,t-1,t,t+1,\dots,j-1,j,\dots,i}^{
\lambda_1,\dots,\lambda_{t-1},\pi(\lambda_t),\pi(\lambda_{t+1}),\dots,
\pi(\lambda_{j-1}),\pi(\lambda_j),\dots,\pi(\lambda_i)}\cdot
\\
&
\ \ \ \ \ \ \ \ \ \ \ \ \ \ \ \ \ \ \ \ \ \ \ \ \ \ \ \ \ \ \ \ \ \ 
\ \ \ \ \ \ \ \ \ \ \ \ \ \ \ \ 
\cdot
\Delta_{1,2,\dots,t,t+1,t+2,\dots,j}^{
\pi(\mu_1),\pi(\mu_2),\dots,\pi(\mu_t),
\mu_{t+1},\mu_{t+2},\dots,\mu_j}
\endaligned}\,.
\]
\end{Theorem}

We will not reproduce the proof here, but extract instead from the
cited references the further important information that the ideal of
relations between all our $\Delta$ jet-determinants:
\[
\aligned
\Delta_{1,\dots,i}^{\lambda_1,\dots,\lambda_i}
=
\left\vert
\begin{array}{ccc}
f_1^{(\lambda_1)} & \cdots & f_i^{(\lambda_1)}
\\
\cdot\cdot & \cdots & \cdot\cdot
\\
f_1^{(\lambda_i)} & \cdots & f_i^{(\lambda_i)}
\end{array}
\right\vert
\\
\ \ \ \ \ \ \ \ \ \ \ \ \
{\scriptstyle{(i\,=\,1\,\cdots\,n\,;\,\,
1\,\leqslant\,\lambda_1\,<\,\cdots\,<\,\lambda_i\,\leqslant\,\kappa)}}
\endaligned
\]
is generated (as an ideal) by all the above quadratic Pl\"ucker
relations. Moreover, these relations written explicitly above do
constitute a {\sl Gr\"obner basis} for a certain term order, 
presented as follows.

Introduce first as many independent variables $\nabla^{ \lambda_1^i,
\dots, \lambda_i^i}$ as there are $\Delta$ jet-determinants and consider
the ring $\C \big[ \nabla^{ \lambda_1^1}, \dots, \nabla^{ \lambda_1^n,
\dots, \lambda_n^n} \big]$. Totally order these variables by declaring
that: 
\[
\nabla^{\lambda_1^i,\dots,\lambda_i^i} 
<_{two} 
\nabla^{\mu_1^j,\dots,\mu_j^j} 
\]
if either $i > j$ or else if $i = j$ {\em and} $(\lambda_1^i, \dots,
\lambda_i^i)$ comes before $(\mu_1^j, \dots, \mu_{j=i}^j)$ in the
lexicographic ordering, which simply means that there exists an index
$s \in \{ 1, \dots, j\}$ such that:
\[
\lambda_1^i=\mu_1^j,\ \ \ \ \
\dots\dots,\ \ \ \ \
\lambda_{s-1}^i=\mu_{s-1}^j,\ \ \ \ \
\lambda_s^i<\mu_s^j. 
\]
This total order extend the partial
order ``$<_{one}$''.
Finally, let also ``$<_{two}$'' denote the {\sl reverse
lexicographic}\footnote{\,
Generally, if $x_1 <_{two} \cdots <_{two} < x_n$, the reverse
lexicographic (total) ordering induced on monomials says that $x_1^{
\alpha_1} \cdots x_n^{ \alpha_n}$ is {\sl smaller} than $x_1^{
\beta_1} \cdots x_n^{ \beta_n}$ if, when reading exponents from right
to left, equality holds: $\alpha_n = \beta_n$, \dots, $\alpha_{ u+1} =
\beta_{ u+1}$ until a first difference occurs: $\alpha_u \neq \beta_u$
for which $\alpha_u > \beta_u$ is {\em bigger} than $\beta_u$.
} 
term ordering on $\C \big[ \nabla^{ \lambda_1^1}, \dots, \nabla^{
\lambda_1^n, \dots, \lambda_n^n} \big]$ that is induced by this
variable ordering $<_{two}$. The set of polynomials ${\sf R} \big(
\nabla^{ \lambda_1^1}, \dots, \nabla^{ \lambda_1^n, \dots,
\lambda_n^n} \big)$ which annihilate identically after replacement by
the determinants:
\[
0
\equiv
{\sf R}
\big(
\Delta_1^{\lambda_1^1},\dots,
\Delta_{1,\dots,n}^{\lambda_1^n,\dots,\lambda_n^n}
\big)
\]
constitutes clearly an {\em ideal} of $\C \big[ \nabla^{ \lambda_1^1},
\dots, \nabla^{ \lambda_1^n, \dots, \lambda_n^n} \big]$.

\begin{Theorem}
\label{Theorem-C}
{\rm (\cite{ stu1993, krpr1996, mist2005, pro2007})}
The ideal of relations $\text{\small\sf Id-rel} (\Delta)$ between all
$\Delta$ jet-determinants is generated by all the Pl\"ucker relations
written above. Moreover, the collection of all these Pl\"ucker
relations constitutes already {\em per se} a {\em Gr\"obner basis} for
$\text{\small\sf Id-rel} (\Delta)$ under the term ordering
``$<_{two}$''. Finally, the products:
\[
\Delta_{1,\dots,i}^{\lambda_1,\dots,\lambda_i}
\cdot
\Delta_{1,\dots,j}^{\mu_1,\dots,\mu_j}
\]
of all possible incomparable pairs generate the (monomial) ideal of
leading monomials of elements of $\text{\small\sf Id-rel} (\Delta)$.
\end{Theorem}

\subsection{Polynomials modulo relations}
Thanks to this statement, we will be able to find a basis of the
$\C$-vector space:
\[
\text{\small\sf all}\,
\Delta\text{\small\sf -polynomials}\,
\big/\,
\text{\small\sf modulo their relations}.
\]
This will be very useful, for we saw that basis vectors are in
one-to-one correspondence with Schur bundles $\mathcal{ S}^{ ( \ell_1,
\dots, \ell_n)} T_X^*$ (see also below). The general Gr\"obner basis
theory then tells us that this vector space is generated by all
$\Delta$-monomials that are {\em not} multiple of any product of
incomparable pairs (leading monomials). In order to describe
explicitly this quotient vector space, 
we need a classical combinatorial object.

\subsection{Young diagrams}
Let $d_1 \geqslant 1$ be an integer and let $\ell_1, \ell_2, \dots,
\ell_{ d_1}$ be any collection of $d_1$ nonnegative integers collected
in decreasing order: 
\[
\ell_1\geqslant\ell_2
\geqslant\cdots
\geqslant\ell_{d_1}
\geqslant 1. 
\]
The {\sl Young diagram} ${\sf YD}_{ (\ell_1, \dots, \ell_{ d_1})}$
associated to such a $d_1$-tuple $(\ell_1, \dots, \ell_{ d_1})$ sits
in the right-bottom quadrant $\{ x \geqslant 0, \, y \leqslant 0\}$ of
the plane $\R^2 = \R^2 ( x, y)$ and it consists, in the $i$-th
horizontal strip $\{ -i \leqslant y \leqslant -i+1 \}$ from above, for
$i = 1, \dots, d_1$, of the $\ell_i$ empty unit squares:
\[
\square_i^j
:=
\big\{(x,y)\in\R^2\colon\,
-i\leqslant y\leqslant-i+1,\, 
j-1\leqslant x \leqslant j
\big\}
\]
placed, for $j = 1, \dots, \ell_i$, one after the other and starting
from the vertical $y$-axis (left-justi\-fication).

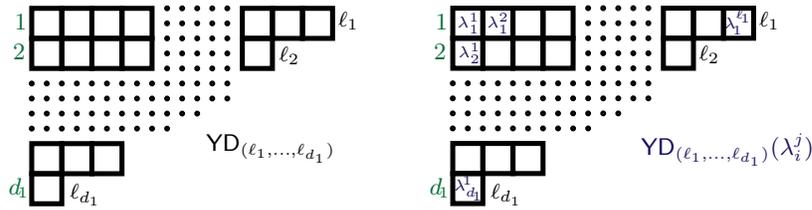
\begin{center}
\begin{picture}(0,0)%
\includegraphics{young-diagram.pstex}%
\end{picture}%
\setlength{\unitlength}{4144sp}%
\begingroup\makeatletter\ifx\SetFigFont\undefined%
\gdef\SetFigFont#1#2#3#4#5{%
  \reset@font\fontsize{#1}{#2pt}%
  \fontfamily{#3}\fontseries{#4}\fontshape{#5}%
  \selectfont}%
\fi\endgroup%
\begin{picture}(4492,1237)(808,-1174)
\put(1992,-826){\makebox(0,0)[lb]{\smash{{\SetFigFont{9}{10.8}{\familydefault}{\mddefault}{\updefault}{\color[rgb]{0,0,0}${\sf YD}_{(\ell_1,\dots,\ell_{d_1})}$}%
}}}}
\put(2427,-289){\makebox(0,0)[lb]{\smash{{\SetFigFont{8}{9.6}{\familydefault}{\mddefault}{\updefault}{\color[rgb]{0,0,0}$\ell_2$}%
}}}}
\put(2781,-81){\makebox(0,0)[lb]{\smash{{\SetFigFont{8}{9.6}{\familydefault}{\mddefault}{\updefault}{\color[rgb]{0,0,0}$\ell_1$}%
}}}}
\put(1175,-1109){\makebox(0,0)[lb]{\smash{{\SetFigFont{8}{9.6}{\familydefault}{\mddefault}{\updefault}{\color[rgb]{0,0,0}$\ell_{d_1}$}%
}}}}
\put(808,-1087){\makebox(0,0)[lb]{\smash{{\SetFigFont{9}{10.8}{\familydefault}{\mddefault}{\updefault}{\color[rgb]{0,.82,0}\green{$d_{\!1}$}}%
}}}}
\put(835,-276){\makebox(0,0)[lb]{\smash{{\SetFigFont{9}{10.8}{\familydefault}{\mddefault}{\updefault}{\color[rgb]{0,.82,0}\green{$2$}}%
}}}}
\put(844,-100){\makebox(0,0)[lb]{\smash{{\SetFigFont{9}{10.8}{\familydefault}{\mddefault}{\updefault}{\color[rgb]{0,.82,0}\green{$1$}}%
}}}}
\put(5300,-80){\makebox(0,0)[lb]{\smash{{\SetFigFont{8}{9.6}{\familydefault}{\mddefault}{\updefault}{\color[rgb]{0,0,0}$\ell_1$}%
}}}}
\put(3694,-1108){\makebox(0,0)[lb]{\smash{{\SetFigFont{8}{9.6}{\familydefault}{\mddefault}{\updefault}{\color[rgb]{0,0,0}$\ell_{d_1}$}%
}}}}
\put(3327,-1086){\makebox(0,0)[lb]{\smash{{\SetFigFont{9}{10.8}{\familydefault}{\mddefault}{\updefault}{\color[rgb]{0,.82,0}\green{$d_{\!1}$}}%
}}}}
\put(3354,-275){\makebox(0,0)[lb]{\smash{{\SetFigFont{9}{10.8}{\familydefault}{\mddefault}{\updefault}{\color[rgb]{0,.82,0}\green{$2$}}%
}}}}
\put(3363,-99){\makebox(0,0)[lb]{\smash{{\SetFigFont{9}{10.8}{\familydefault}{\mddefault}{\updefault}{\color[rgb]{0,.82,0}\green{$1$}}%
}}}}
\put(3474,-1064){\makebox(0,0)[lb]{\smash{{\SetFigFont{6}{7.2}{\familydefault}{\mddefault}{\updefault}{\color[rgb]{0,0,.69}\blue{$\lambda_{d_{\!1}}^{\!1}$}}%
}}}}
\put(3499,-268){\makebox(0,0)[lb]{\smash{{\SetFigFont{6}{7.2}{\familydefault}{\mddefault}{\updefault}{\color[rgb]{0,0,.69}\blue{$\lambda_2^1$}}%
}}}}
\put(3494,-88){\makebox(0,0)[lb]{\smash{{\SetFigFont{6}{7.2}{\familydefault}{\mddefault}{\updefault}{\color[rgb]{0,0,.69}\blue{$\lambda_1^1$}}%
}}}}
\put(3671,-87){\makebox(0,0)[lb]{\smash{{\SetFigFont{6}{7.2}{\familydefault}{\mddefault}{\updefault}{\color[rgb]{0,0,.69}\blue{$\lambda_1^2$}}%
}}}}
\put(5088,-92){\makebox(0,0)[lb]{\smash{{\SetFigFont{6}{7.2}{\familydefault}{\mddefault}{\updefault}{\color[rgb]{0,0,.69}\blue{$\lambda_{\!1}^{\ell_{\!1}}$}}%
}}}}
\put(4939,-285){\makebox(0,0)[lb]{\smash{{\SetFigFont{8}{9.6}{\familydefault}{\mddefault}{\updefault}{\color[rgb]{0,0,0}$\ell_2$}%
}}}}
\put(4594,-847){\makebox(0,0)[lb]{\smash{{\SetFigFont{9}{10.8}{\familydefault}{\mddefault}{\updefault}{\color[rgb]{0,0,.69}\blue{${\sf YD}_{(\ell_1,\dots,\ell_{d_1})}(\lambda_i^j)$}}%
}}}}
\end{picture}%

\end{center}

\noindent
It will be convenient to give names to the column lengths, say $d_j$
will denote that of the $j$-th, for $j = 1, \dots, \ell_1$. In summary
and for memory: 
\[
\ell_i
=
\text{\small\sf length of the $i$-th row;}
\ \ \ \ \ \ \ \ \ \ \ \
d_j
=
\text{\small\sf length of the $j$-th column}.
\]
Observe that the longest column lengths are equal to $d_1$ for all
indices $j$ between $1$ and $\ell_{ d_1}$, and more generally 
at any $i$-th row (see the 
zoom below), that the following
coincidence of column lengths holds:
\[
i
=
d_{1+\ell_{i+1}}
=\cdots=
d_{\ell_i}
\ \ \ \ \ \ \ \ \ \ \ \ \
{\scriptstyle{(1\,\leqslant\,i\,\leqslant\,d_1)}}.
\]
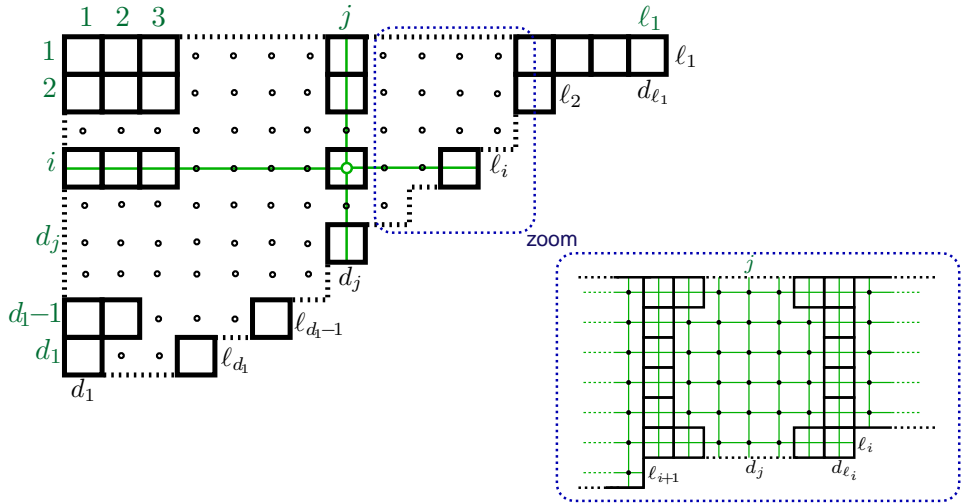
\begin{center}
\begin{picture}(0,0)%
\includegraphics{d1-reverse-young-diagram.pstex}%
\end{picture}%
\setlength{\unitlength}{4144sp}%
\begingroup\makeatletter\ifx\SetFigFont\undefined%
\gdef\SetFigFont#1#2#3#4#5{%
  \reset@font\fontsize{#1}{#2pt}%
  \fontfamily{#3}\fontseries{#4}\fontshape{#5}%
  \selectfont}%
\fi\endgroup%
\begin{picture}(5711,2997)(295,-2457)
\put(2277,420){\makebox(0,0)[lb]{\smash{{\SetFigFont{11}{13.2}{\familydefault}{\mddefault}{\updefault}{\color[rgb]{0,.82,0}\green{$j$}}%
}}}}
\put(495,204){\makebox(0,0)[lb]{\smash{{\SetFigFont{11}{13.2}{\familydefault}{\mddefault}{\updefault}{\color[rgb]{0,.82,0}\green{$1$}}%
}}}}
\put(500,-16){\makebox(0,0)[lb]{\smash{{\SetFigFont{11}{13.2}{\familydefault}{\mddefault}{\updefault}{\color[rgb]{0,.82,0}\green{$2$}}%
}}}}
\put(525,-460){\makebox(0,0)[lb]{\smash{{\SetFigFont{11}{13.2}{\familydefault}{\mddefault}{\updefault}{\color[rgb]{0,.82,0}\green{$i$}}%
}}}}
\put(3594,-56){\makebox(0,0)[lb]{\smash{{\SetFigFont{10}{12.0}{\familydefault}{\mddefault}{\updefault}{\color[rgb]{0,0,0}$\ell_2$}%
}}}}
\put(4284,183){\makebox(0,0)[lb]{\smash{{\SetFigFont{10}{12.0}{\familydefault}{\mddefault}{\updefault}{\color[rgb]{0,0,0}$\ell_1$}%
}}}}
\put(4057,-24){\makebox(0,0)[lb]{\smash{{\SetFigFont{10}{12.0}{\familydefault}{\mddefault}{\updefault}{\color[rgb]{0,0,0}$d_{\ell_1}$}%
}}}}
\put(720,409){\makebox(0,0)[lb]{\smash{{\SetFigFont{11}{13.2}{\familydefault}{\mddefault}{\updefault}{\color[rgb]{0,.82,0}\green{$1$}}%
}}}}
\put(1152,409){\makebox(0,0)[lb]{\smash{{\SetFigFont{11}{13.2}{\familydefault}{\mddefault}{\updefault}{\color[rgb]{0,.82,0}\green{$3$}}%
}}}}
\put(938,409){\makebox(0,0)[lb]{\smash{{\SetFigFont{11}{13.2}{\familydefault}{\mddefault}{\updefault}{\color[rgb]{0,.82,0}\green{$2$}}%
}}}}
\put(4064,408){\makebox(0,0)[lb]{\smash{{\SetFigFont{11}{13.2}{\familydefault}{\mddefault}{\updefault}{\color[rgb]{0,.82,0}\green{$\ell_1$}}%
}}}}
\put(1571,-1632){\makebox(0,0)[lb]{\smash{{\SetFigFont{10}{12.0}{\familydefault}{\mddefault}{\updefault}{\color[rgb]{0,0,0}$\ell_{d_{\!1}}$}%
}}}}
\put(669,-1812){\makebox(0,0)[lb]{\smash{{\SetFigFont{10}{12.0}{\familydefault}{\mddefault}{\updefault}{\color[rgb]{0,0,0}$d_1$}%
}}}}
\put(295,-1362){\makebox(0,0)[lb]{\smash{{\SetFigFont{11}{13.2}{\familydefault}{\mddefault}{\updefault}{\color[rgb]{0,.82,0}\green{$d_{\!1}\!\!-\!\!1$}}%
}}}}
\put(458,-1586){\makebox(0,0)[lb]{\smash{{\SetFigFont{11}{13.2}{\familydefault}{\mddefault}{\updefault}{\color[rgb]{0,.82,0}\green{$d_{1}$}}%
}}}}
\put(2009,-1422){\makebox(0,0)[lb]{\smash{{\SetFigFont{10}{12.0}{\familydefault}{\mddefault}{\updefault}{\color[rgb]{0,0,0}$\ell_{d_{\!1}\!-\!1}$}%
}}}}
\put(2272,-1141){\makebox(0,0)[lb]{\smash{{\SetFigFont{10}{12.0}{\familydefault}{\mddefault}{\updefault}{\color[rgb]{0,0,0}$d_j$}%
}}}}
\put(459,-901){\makebox(0,0)[lb]{\smash{{\SetFigFont{11}{13.2}{\familydefault}{\mddefault}{\updefault}{\color[rgb]{0,.82,0}\green{$d_{j}$}}%
}}}}
\put(3176,-471){\makebox(0,0)[lb]{\smash{{\SetFigFont{10}{12.0}{\familydefault}{\mddefault}{\updefault}{\color[rgb]{0,0,0}$\ell_i$}%
}}}}
\put(4704,-2271){\makebox(0,0)[lb]{\smash{{\SetFigFont{7}{8.4}{\familydefault}{\mddefault}{\updefault}{\color[rgb]{0,0,0}$d_j$}%
}}}}
\put(4702,-1051){\makebox(0,0)[lb]{\smash{{\SetFigFont{8}{9.6}{\familydefault}{\mddefault}{\updefault}{\color[rgb]{0,.82,0}\green{$j$}}%
}}}}
\put(5218,-2271){\makebox(0,0)[lb]{\smash{{\SetFigFont{7}{8.4}{\familydefault}{\mddefault}{\updefault}{\color[rgb]{0,0,0}$d_{\ell_i}$}%
}}}}
\put(5377,-2128){\makebox(0,0)[lb]{\smash{{\SetFigFont{7}{8.4}{\familydefault}{\mddefault}{\updefault}{\color[rgb]{0,0,0}$\ell_i$}%
}}}}
\put(4122,-2308){\makebox(0,0)[lb]{\smash{{\SetFigFont{7}{8.4}{\familydefault}{\mddefault}{\updefault}{\color[rgb]{0,0,0}$\ell_{i\!+\!1}$}%
}}}}
\put(3395,-912){\makebox(0,0)[lb]{\smash{{\SetFigFont{8}{9.6}{\familydefault}{\mddefault}{\updefault}{\color[rgb]{0,0,.69}\blue{\sf zoom}}%
}}}}
\end{picture}%

\end{center}

\subsection{Semi-standard Young tableaux} If $\lambda_i^j
\geqslant 1$ denote as many nonnegative integers as there are empty
squares $\square_i^j$, namely with $i = 1, \dots, d_1$ and $j = 1,
\dots, \ell_i$, a {\sl filling} ${\sf YD}_{ (\ell_1, \dots, \ell_{
d_1})} (\lambda_i^j)$ of the Young diagram ${\sf YD}_{ (\ell_1, \dots,
\ell_{ d_1})}$ by means of the $\lambda_i^j$ consists in putting each
$\lambda_i^j$ in each square $\square_i^j$. A {\sl semi-standard
(Young) tableau} is a filled Young diagram ${\sf YD}_{ (\ell_1, \dots,
\ell_{ d_1})} (\lambda_i^j)$ having the property that when reading its
full content:
\[
\aligned
\begin{array}{ccccccccccccc}
\lambda_1^1 & \cdots & \lambda_1^{\ell_{d_1}} & \cdots &
\lambda_1^{\ell_{d_1-1}} & \cdots\,\,\cdots & \lambda_1^{\ell_i} &
\cdots\,\,\cdots & \lambda_1^{\ell_2} & \cdots & \lambda_1^{\ell_1}
\\
\lambda_2^1 & \cdots & \lambda_2^{\ell_{d_1}} & \cdots &
\lambda_2^{\ell_{d_1-1}} & \cdots\,\,\cdots & \lambda_2^{\ell_i} &
\cdots\,\,\cdots & \lambda_2^{\ell_2} 
\\
\cdot & \cdots & \cdot & \cdots & \cdot &
\cdots\,\cdots & \cdot & \cdots\,\,\,\,\,\,\,\,\,\,\:
\\
\lambda_i^1 & \cdots & \lambda_i^{\ell_{d_1}} & \cdots &
\lambda_i^{\ell_{d_1-1}} & \cdots\,\,\cdots & \lambda_i^{\ell_i}
\\
\cdot & \cdots & \cdot & \cdots & \cdot & \cdots\,\,\,\,\,\,\,\,\,\,\
\\
\lambda_{d_1-1}^1 & \cdots & \lambda_{d_1-1}^{\ell_{d_1-1}} & \cdots &
\lambda_{d_1-1}^{\ell_{d_1-1}}
\\
\lambda_{d_1}^1 & \cdots & \lambda_{d_1}^{\ell_{d_1}},
\end{array}
\endaligned
\]
the integers $\lambda_i^j$ increase from top to bottom in each column,
and they are 
nondecreasing\footnote{\,
A so-called {\sl standard tableau} would require that the integers
$\lambda_i^j$ also increase along the rows.
} 
in each row from left to right, that is to say and more precisely:
\[
\aligned
&
\lambda_1^j<\lambda_2^j<
\cdots<\lambda_{d_j}^j
\ \ \ \ \ \ \ \ \ \ \ \ \
{\scriptstyle{(1\,\leqslant\,j\,\leqslant\,\ell_1)}}
\\
&
\lambda_i^1\leqslant\lambda_i^2\leqslant
\cdots\leqslant\lambda_i^{\ell_i}
\ \ \ \ \ \ \ \ \ \ \ \ \
{\scriptstyle{(1\,\leqslant\,i\,\leqslant\,d_1)}}.
\endaligned
\] 

\subsection{Vector space basis for the algebra 
of $\Delta$ jet-determinants} Coming back to our algebra of
determinants $\Delta_{ 1, 2, \dots, i}^{ \lambda_1, \lambda_2, \dots,
\lambda_i}$, the increasing sequence of their exponents
$\lambda_1 < \lambda_2 < \cdots < \lambda_i$ will sit in a column of
such a Young diagram. Since the row-size $i$ of any not
identically zero minor $\Delta_{ 1, 2, \dots, i}^{ \lambda_1,
\lambda_2, \dots, \lambda_i}$ must be $\leqslant n = {\rm rank}\,
(T_X^*)$, we will consider in fact only semi-standard tableaux whose
depth $d_1$ is always $\leqslant n$. Accordingly, when it happens that
$d_1 < n$ we shall adopt the natural convention that:
\[
\ell_1\geqslant\ell_2
\geqslant\cdots\geqslant
\ell_{d_1-1}\geqslant\ell_{d_1}
>
0=\ell_{d_1+1}=\cdots=\ell_n. 
\]
We are at last in a position to state the
starting point theorem. 

\begin{Theorem}
\label{Theorem-D}
{\rm (\cite{ stu1993, krpr1996, mist2005, pro2007})}
The infinite-dimensional quotient vector space:
\[
\text{\small\sf all}\,
\Delta\text{\small\sf -polynomials}\,
\big/\,
\text{\small\sf modulo their relations}
\]
possesses a basis over $\C$ consisting of all possible 
$\Delta$-monomials:
\[
\prod_{1\leqslant j\leqslant\ell_{d_1}}\!\!
\Delta_{1,\dots,d_1}^{\lambda_1^j,\dots,\lambda_{d_1}^j}\!\!
\prod_{\ell_{d_1}+1\leqslant j\leqslant\ell_{d_1-1}}\!\!
\Delta_{1,\dots,d_1-1}^{
\lambda_1^j,\dots,
\lambda_{d_1-1}^j}
\cdots\!\!
\prod_{\ell_2+1\leqslant j\leqslant\ell_1}\!\!
\Delta_1^{\lambda_1^j}
\]
such that the collection of appearing upper exponents $(\lambda_i^j)$
constitutes a semi-standard Young tableau:
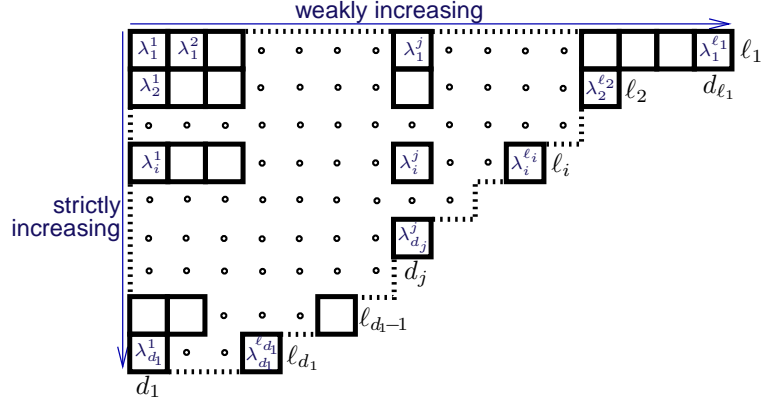
\begin{center}

\end{center}
\end{Theorem}

\subsection{Exact Schur bundle decomposition of ${\sf Gr}^\bullet 
\mathcal{E}_{ \kappa, m}^{ GG} T_X^*$} 
In order to apply this combinatorial information to our problem, we
may also represent the general $\Delta$-monomial written above more
concisely as:
\begin{equation}
\label{general-Delta-monomial}
\prod_{d_1\geqslant i\geqslant 1}\,\,
\prod_{1+\ell_{i+1}\leqslant j\leqslant\ell_i}\!\!
\Delta_{1,\dots,i}^{
\lambda_1^j,\dots,\lambda_i^j}. 
\end{equation}
First of all, every $\Delta$-determinant read off from such a product
happens to be an eigenvector for the action on jets of the diagonal
matrices ${\sf e} = {\rm diag} \big( {\sf e}_1, \dots, {\sf
e}_n\big)$:
\[
{\sf e}
\cdot
\Delta_{1,2,\dots,i}^{\lambda_1^j,\lambda_2^j,\dots,\lambda_i^j}
=
{\sf e}_1\,{\sf e}_2\cdots{\sf e}_i\,
\Delta_{1,2,\dots,i}^{\lambda_1^j,\lambda_2^j,\dots,\lambda_i^j},
\]
as is clear because the diagonal action just multiplies columns of
such a determinant by the quantities ${\sf e}_1, {\sf e}_2, \dots,
{\sf e}_j$:
\[
{\sf e}
\cdot
\Delta_{1,2,\dots,i}^{\lambda_1^j,\lambda_2^j,\dots,\lambda_i^j}
=
\left\vert
\begin{array}{ccccc}
{\sf e}_1f_1^{(\lambda_1^j)} & 
{\sf e}_2f_2^{(\lambda_1^j)} & \cdots & {\sf e}_if_i^{(\lambda_1^j)}
\\
{\sf e}_1f_1^{(\lambda_2^j)} & 
{\sf e}_2f_2^{(\lambda_2^j)} & \cdots & {\sf e}_if_i^{(\lambda_2^j)}
\\
\cdot\cdot & \cdot\cdot & \cdots & \cdot\cdot
\\
{\sf e}_1f_1^{(\lambda_i^j)} & 
{\sf e}_2f_2^{(\lambda_i^j)} & \cdots & 
{\sf e}_if_i^{(\lambda_i^j)}
\end{array}
\right\vert.
\]
Consequently, every general monomial in the $\Delta$-determinants
represented by the above arbitrary semi-standard tableau is also an
eigenvector:
\[
\aligned
{\sf e}\,\cdot
\bigg(
\prod_{d_1\geqslant i\geqslant 1}\,\,
\prod_{1+\ell_{i+1}\leqslant j\leqslant\ell_i}\!\!
&
\Delta_{1,\dots,i}^{
\lambda_1^j,\dots,\lambda_i^j}
\bigg)
=
{\sf e}\,\cdot
\big(
\text{\rm general $\Delta$-monomial}
\big)
\\
&
=
\prod_{d_1\geqslant i\geqslant 1}\,\,
\prod_{1+\ell_{i+1}\leqslant j\leqslant\ell_i}\!\!
{\sf e}_1\cdots{\sf e}_i
\cdot
\big(\text{\rm same $\Delta$-monomial}\big)
\\
&
=
\prod_{d_1\geqslant i\geqslant 1}\,\,
\big({\sf e}_1\cdots{\sf e}_i\big)^{\ell_i-\ell_{i+1}}
\cdot
\big(\text{\rm same $\Delta$-monomial}\big)
\\
&
=
({\sf e}_1)^{\ell_1}
({\sf e}_2)^{\ell_2}
\cdots
({\sf e}_n)^{\ell_n}\,
\cdot
\big(\text{\rm same $\Delta$-monomial}\big). 
\endaligned
\]
As a result, we deduce generally that: 
\[
\boxed{
\text{\small\sf Single semi-standard $\Delta$-monomial}
\longleftrightarrow
\text{\small\sf Unique Schur bundle}
}\,, 
\]
and more precisely, to the general monomial associated with a
semi-standard tableau ${\sf YD}_{( \ell_1, \dots, \lambda_n)}
(\lambda_i^j)$ corresponds bijectively the Schur bundle $\mathcal{
S}^{ ( \ell_1, \ell_2, \dots, \ell_n)} T_X^*$. Thus notably, the
related Schur bundle depends only on the diagram, and it {\em does
not depend on its filling by integers $\lambda_i^j$}.

Although essentially not new since it follows from 
Theorems~\ref{Theorem-A}, \ref{Theorem-B}, \ref{Theorem-C}
and~\ref{Theorem-D} above, 
the following basic statement appears nowhere as such in
the literature devoted the application of the jet bundle machinery to
the conjectures of Green-Griffiths and of Kobayashi, but it will
nonetheless constitute our basic starting point.

\begin{Theorem}
\label{exact-decomposition}
The graded vector bundle ${\sf Gr}^\bullet \mathcal{ E}_{ \kappa, m}^{
GG} T_X^*$ associated to the bundle $\mathcal{ E}_{ \kappa, m}^{ GG}
T_X^*$ of $\kappa$-th $m$-weighted Green-Griffiths jets identifies
to the following {\em exact direct sum} of Schur bundles:
\[
{\sf Gr}^\bullet\mathcal{E}_{\kappa,m}^{GG}T_X^*
=
\bigoplus_{\ell_1\geqslant\ell_2\geqslant\cdots\geqslant\ell_n
\geqslant 0}
\Big(
\mathcal{S}^{(\ell_1,\ell_2,\dots,\ell_n)}T_X^*
\Big)^{\oplus 
M_{\ell_1,\ell_2,\dots,\ell_n}^{\kappa,m}}, 
\]
with multiplicities $M_{\ell_1, \ell_2, \dots, \ell_n }^{ \kappa,
m} \in \N$ equal to the number of times a Young diagram ${\sf YD}_{
(\ell_1, \dots, \ell_n)}$ with row lengths equal to $\ell_1, \ell_2,
\dots, \ell_n$ can be filled in with positive integers $\lambda_i^j
\leqslant \kappa$ placed at its $i$-th line and $j$-th column so as to
constitute a semi-standard tableau, with the further constraint that
the sum of all such integers:
\[
\aligned
m
&
=
\lambda_1^1+\cdots+\lambda_1^{\ell_n}
+\cdots+
\lambda_1^{\ell_2}+\cdots+\lambda_1^{\ell_1}
\\
&
+
\lambda_2^1+\cdots+\lambda_2^{\ell_n}
+\cdots+\lambda_2^{\ell_2}
\\
&
+\cdots\cdot\cdots\cdots\cdots\cdots+
\\
&
+
\lambda_n^1+\cdots+\lambda_n^{\ell_n}
\endaligned
\]
equals the prescribed weighted homogeneity degree $m$.
\end{Theorem}

This apparently complete statement should not hide the fact that the
exact computation of the multiplicities $M_{ \ell_1, \ell_2, \dots,
\ell_n}^{ \kappa, m}$ is not provided in terms of $\kappa$, $m$ and
$\ell_1, \ell_2, \dots, \ell_n$. Manual attempts to find a usable,
closed and explicit formula for $M_{ \ell_1, \ell_2, \dots, \ell_n}^{
\kappa, m}$ showed us that the task could be hard, and we will proceed
differently, in an asymptotic manner, so as to avoid several
unnecessary computations which would anyway be inaccessible to us.

\begin{Corollary}
One has the following inequalities between the cohomology 
dimensions $h^q$ for all $q = 1, 2, \dots, n$:
\[
\footnotesize
\aligned
h^q
\big(
X,\,
\mathcal{E}_{\kappa,m}^{GG}T_X^*
\big)
&
\leqslant
\sum_{\ell_1+2\ell_2+\cdots+\kappa\ell_\kappa=m}\,
h^q
\Big(X,\,\,
{\rm Sym}^{\ell_1}T_X^*
\otimes
{\rm Sym}^{\ell_2}T_X^*
\otimes\cdots\otimes
{\rm Sym}^{\ell_\kappa}T_X^*
\Big)
\\
&
\leqslant
\sum_{\ell_1\geqslant\ell_2\geqslant\cdots\geqslant\ell_n\geqslant 0}\,
M_{\ell_1,\ell_2,\dots,\ell_n}^{\kappa,m}
h^q
\big(
X,\,\mathcal{S}^{(\ell_1,\ell_2,\dots,\ell_n)}T_X^*
\big).
\endaligned
\]
\end{Corollary}

\proof
The first one was already derived in~\thetag{ 
\ref{first-inq-cohomology}}. Then the decomposition into Schur
bundles of each tensored factor ${\rm Sym}^{ \ell_1 } T_X^* \otimes
\cdots \otimes {\rm Sym}^{ \ell_\kappa } T_X^*$ obtained {\em e.g.} by
an application of Pieri's rule~\thetag{ \ref{pieri-rule}}
enables one
to define a subfiltration to which the same reasoning as 
in~\thetag{ \ref{first-inq-cohomology}} applies.
\endproof

Thus, as in Rousseau's papers~\cite{ rou2006a, rou2006b} for $n = 3$
and $\kappa = 3$ and as in~\cite{ mer2008b} for $n = 4$ and $\kappa =
4$, the study of the cohomology of the Green-Griffiths bundle
$\mathcal{ E}_{ \kappa, m}^{ GG} T_X^*$ is led back to the study of
the cohomology of Schur bundles, which might in turn be complicated.

\markleft{Jo\"el Merker}
\markright{\S5.~Asymptotic characteristic and asymptotic cohomology}
\section{\bf Asymptotic characteristic 
\\
and asymptotic cohomology}
\label{Section-5}

\subsection{Giambelli determinants of Chern classes}
From Lemma~\ref{Gr-chi} and from Theorem~\ref{exact-decomposition}, 
we deduce at once from the
additivity of Euler-Poincar\'e characteristic that:
\[
\chi\big(
X,\,\,
\mathcal{E}_{\kappa,m}^{GG}T_X^*
\big)
=
\sum_{\ell_1\geqslant\ell_2\geqslant\cdots\geqslant\ell_n\geqslant 0}\,
M_{\ell_1,\ell_2,\dots,\ell_n}^{\kappa,m}
\cdot
\chi\big(
X,\,\,
\mathcal{S}^{(\ell_1,\ell_2,\dots,\ell_n)}T_X^*
\big).
\]
But there is a closed asymptotic general formula for: 
\[
\chi\big(
X,\,\mathcal{S}^{(\ell_1,\dots,\ell_n)}T_X^*\big)
=
(-1)^n\,
\chi\big(
X,\,\mathcal{S}^{(\ell_1,\dots,\ell_n)}T_X\big),
\]
where the $(-1)^n$ comes from ${\sf c}_k^* = ( -1)^k {\sf c}_k$.
Recall that a {\sl partition $(\nu_1, \nu_2, \dots, \nu_n)$} of $n$ is
just a collection of nonnegative integers $\nu_1 \geqslant \nu_2
\geqslant \cdots \geqslant \nu_n \geqslant 0$ whose sum $\nu_1 + \nu_2
+ \cdots + \nu_n$ equals $n$.

\begin{Theorem}
\label{E-P-characteristic}
{\rm (\cite{ mer2008b}\footnote{\,
After~\cite{ mer2008b} was posted on {\scriptsize{\sf arxiv.org}}, the
author was informed by E.~Rousseau that Br\"uckmann's Theorem~4
in~\cite{ bru1997} entails the above statement and
moreover, that it shows how to explicit the remainders. 
})} 
The terms of highest order with respect to $\vert \ell \vert =
\ell_1+ \cdots + \ell_n$ in the Euler-Poincar\'e characteristic of the
Schur bundle $\mathcal{ S}^{ ( \ell_1, \ell_2, \dots, \ell_n )} \,
T_X^*$ are homogeneous of order $\frac{ n ( n+1)}{ 2}$ and they are
given by a sum of determinants indexed by all the partitions $(\nu_1,
\dots, \nu_n)$ of $n$:
\[
\small
\aligned
&
(-1)^n\,
\chi\Big(X,\,\,
\mathcal{S}^{(\ell_1,\ell_2,\dots,\ell_n)}\,T_X^*\Big)
=
\\
&
=
\sum_{\nu\,\text{\rm partition of}\,\,n}\,
\frac{
{\sf C}_{\nu^c}}{(\nu_1+n-1)!\,\cdots\,\nu_n!}\,
\left\vert\!
\begin{array}{cccc}
{\ell_1'}^{\nu_1+n-1} & {\ell_2'}^{\nu_1+n-1} &
\cdots & {\ell_n'}^{\nu_1+n-1}
\\
{\ell_1'}^{\nu_2+n-2} & {\ell_2'}^{\nu_2+n-2} &
\cdots & {\ell_n'}^{\nu_2+n-2}
\\
\vdots & \vdots & \ddots & \vdots
\\
{\ell_1'}^{\nu_n} & {\ell_2'}^{\nu_n} & 
\cdots & {\ell_n'}^{\nu_n}
\end{array}
\!\right\vert
+
\\
&\ \ \ \ \ \ \ \ \ \
+
{\sf O}_n
\big(
\vert\ell\vert^{\frac{n(n+1)}{2}-1}
\big),
\endaligned
\]
where $\ell_i ' := \ell_i + n - i$ for notational brevity, with
coefficients ${\sf C}_{\nu^c}$ being expressed in terms of the
Chern classes ${\sf c}_k = {\sf c}_k\big(T_X\big)$ of $T_X$ by means
of {\em Giambelli's determinantal expression} depending upon the {\em
conjugate} partition $\nu^c$:
\[
{\sf C}_{\nu^c}
=
{\sf C}_{(\nu_1^c,\dots,\nu_n^c)}
=
\left\vert\!\!
\begin{array}{cccccc}
{\sf c}_{\nu_1^c} & {\sf c}_{\nu_1^c+1} &
{\sf c}_{\nu_1^c+2} & \cdots & {\sf c}_{\nu_1^c+n-1}
\\
{\sf c}_{\nu_2^c-1} & {\sf c}_{\nu_2^c} &
{\sf c}_{\nu_2^c+1} & \cdots & {\sf c}_{\nu_2^c+n-2}
\\
\vdots & \vdots & \vdots & \ddots & \vdots
\\
{\sf c}_{\nu_n^c-n+1} & {\sf c}_{\nu_n^c-n+2} &
{\sf c}_{\nu_n^c-n+3} & \cdots & 
{\sf c}_{\nu_n^c}
\end{array}
\!\!\right\vert,
\]
with the understanding, by convention, that ${\sf c}_k := 0$ for $k<
0$ or $k > n$, and that ${\sf c}_0 := 1$. Furthermore, the remainder
${\sf O}_n \big( \vert \ell \vert^{ \frac{ n ( n+1)}{ 2} } \big)$ is a
linear combination of homogeneous terms ${\sf c}_1^{ \tau_1} {\sf
c}_2^{ \tau_2} \cdots\, {\sf c}_n^{ \tau_n}$ with $\tau_1 + 2\tau_2 +
\cdots + n \tau_n = n$ each multiplied by some polynomial of degree
$\leqslant \frac{ n ( n+1)}{ 2} -1$ in the $\ell_i$ whose coefficients
are rational and bounded in absolute value by ${\sf Constant}_n$.
\end{Theorem}

Because it is elementarily checked that modulo ${\sf O}_n \big( \vert
\ell \vert^{ \frac{ n( n+1)}{2}-1} \big)$, one has:
\[
\small
\aligned
\left\vert\!
\begin{array}{cccc}
{\ell_1'}^{\nu_1+n-1} & {\ell_2'}^{\nu_1+n-1} &
\cdots & {\ell_n'}^{\nu_1+n-1}
\\
{\ell_1'}^{\nu_2+n-2} & {\ell_2'}^{\nu_2+n-2} &
\cdots & {\ell_n'}^{\nu_2+n-2}
\\
\vdots & \vdots & \ddots & \vdots
\\
{\ell_1'}^{\nu_n} & {\ell_2'}^{\nu_n} & 
\cdots & {\ell_n'}^{\nu_n}
\end{array}
\!\right\vert
\equiv
\left\vert\!
\begin{array}{cccc}
{\ell_1}^{\nu_1+n-1} & {\ell_2}^{\nu_1+n-1} &
\cdots & {\ell_n}^{\nu_1+n-1}
\\
{\ell_1}^{\nu_2+n-2} & {\ell_2}^{\nu_2+n-2} &
\cdots & {\ell_n}^{\nu_2+n-2}
\\
\vdots & \vdots & \ddots & \vdots
\\
{\ell_1}^{\nu_n} & {\ell_2}^{\nu_n} & 
\cdots & {\ell_n}^{\nu_n}
\end{array}
\!\right\vert,
\endaligned
\]
we may equivalently replace the $\ell_i'$-determinants by the
corresponding $\ell_i$-determinants in the formula of the theorem.
Then for coherence between the above theorem and the computation of
the Euler-Poincar\'e characteristic of $\mathcal{ E}_{ \kappa, m}^{ GG}
T_X^*$ conducted independently in Section~3, it should be true that
the sum of remainders attached to Schur bundles corresponds to the
last remainder of Theorem~\ref{GG-characteristic}:
\[
\sum_{\ell_1\geqslant\ell_2\geqslant\cdots\geqslant\ell_n\geqslant 0}\,
M_{\ell_1,\ell_2,\dots,\ell_n}^{\kappa,m}
\cdot
{\sf O}_n
\big(
\vert\ell\vert^{\frac{n(n+1)}{2}-1}
\big)
=
{\sf O}_{n,\kappa}
\big(
m^{(\kappa+1)n-2}
\big).
\]
This fact will be established later in
Proposition~\ref{majoration-general} below. 
Also, using~\thetag{ \ref{c-d}}, 
one should consider that all homogeneous products
of Chern classes ${\sf c}_1^{ \tau_1} \cdots {\sf c}_n^{ \tau_n}$ are
implicitly reexpressed in terms of $n$ and $d$, whence both remainders
are in fact of the form ${\sf O}_{n, d} \big( \vert \ell \vert^{
\frac{ n(n+1)}{2} - 1} \big)$ and ${\sf O}_{ n,d, \kappa} \big( m^{(
\kappa+1) n - 2} \big)$.

\subsection{Dimensions 2, 3 and 4} In greater length, 
let us for instance write down the expanded sums over partitions,
firstly in dimension $n = 2$, with two partitions $2 = 2 + 0 = 1 + 1$:
\[
-\,\chi\big(X,\,\mathcal{S}^{(\ell_1,\ell_2)}T_X^*\big)
=
\frac{{\sf c}_1^2-{\sf c}_2}{0!\,\,3!}\,
\left\vert\!\!
\begin{array}{cc}
\ell_1^3 & \ell_2^3
\\
1 & 1
\end{array}
\!\!\right\vert
+
\frac{{\sf c}_2}{1!\,2!}\,
\left\vert\!\!
\begin{array}{cc}
\ell_1^2 & \ell_2^2
\\
\ell_1 & \ell_2
\end{array}
\!\!\right\vert
+
{\sf O}\big(\vert\ell\vert^2\big);
\]
next in dimension $n = 3$, with three partitions $3 = 3 + 0 + 0 = 2 +
1 + 0 = 1 + 1 + 1$:
\[
\aligned
&
\chi\big(
X,\,
\mathcal{S}^{(\ell_1,\ell_2,\ell_3)}\,T_X^*
\big)
=
\\
&
=
\frac{{\sf c}_1^3-2\,{\sf c}_1{\sf c}_2+{\sf c}_3}{0!\,\,1!\,\,5!}\
\left\vert
\begin{array}{ccc}
\ell_1^5\, & \ell_2^5\, & \ell_3^5\,
\\
\ell_1\, & \ell_2\, & \ell_3\,
\\
1\, & 1\, & 1\,
\end{array}
\right\vert
+
\frac{{\sf c}_1{\sf c}_2-{\sf c}_3}{0!\,\,2!\,\,4!}\
\left\vert
\begin{array}{ccc}
\ell_1^4\, & \ell_2^4\, & \ell_3^4\,
\\
\ell_1^2\, & \ell_2^2\, & \ell_3^2\,
\\
1\, & 1\, & 1\,
\end{array}
\right\vert
+
\\
&
\ \ \ \ \ \ \ \
+
\frac{{\sf c}_3}{1!\,\,2!\,\,3!}\
\left\vert
\begin{array}{ccc}
\ell_1^3\, & \ell_2^3\, & \ell_3^3\,
\\
\ell_1^2\, & \ell_2^2\, & \ell_3^2\,
\\
\ell_1\, & \ell_2\, & \ell_3\,
\end{array}
\right\vert
+
{\sf O}\big(\vert\ell\vert^5\big).
\endaligned
\]
and finally in dimension $n = 4$, with
$5$ partitions 
$4 = 
4 + 0 + 0 + 0 = 
3 + 1 + 0 + 0 =
2 + 2 + 0 + 0 = 
2 + 1 + 1 + 0 = 
1 + 1 + 1 + 1$:
\[
\footnotesize
\aligned
&
\chi\big(
X,\,
\mathcal{S}^{(\ell_1,\ell_2,\ell_3,\ell_4)}\,T_X^*
\big)
=
\\
&
=
\frac{{\sf c}_1^4-3\,{\sf c}_1^2{\sf c}_2+{\sf c}_2^2
+2\,{\sf c}_1{\sf c}_3-{\sf c}_4}{0!\,\,1!\,\,2!\,\,7!}\
\left\vert
\begin{array}{cccc}
\ell_1^7\, & \ell_2^7\, & \ell_3^7\, & \ell_4^7\,
\\
\ell_1^2\, & \ell_2^2\, & \ell_3^2\, & \ell_4^2\,
\\
\ell_1^1\, & \ell_2^1\, & \ell_3^1\, & \ell_4^1\,
\\
1\, & 1\, & 1\, & 1\,
\end{array}
\right\vert
+
\\
&
+
\frac{{\sf c}_1^2{\sf c}_2-{\sf c}_2^2-{\sf c}_1{\sf c}_3+{\sf c}_4}
{0!\,\,1!\,\,3!\,\,6!}\
\left\vert
\begin{array}{cccc}
\ell_1^6\, & \ell_2^6\, & \ell_3^6\, & \ell_4^6\,
\\
\ell_1^3\, & \ell_2^3\, & \ell_3^3\, & \ell_4^3\,
\\
\ell_1\, & \ell_2\, & \ell_3\, & \ell_4\,
\\
1\, & 1\, & 1\, & 1\,
\end{array}
\right\vert
+
\frac{-{\sf c}_1{\sf c}_3+{\sf c}_2^2}
{0!\,\,1!\,\,4!\,\,5!}\
\left\vert
\begin{array}{cccc}
\ell_1^5\, & \ell_2^5\, & \ell_3^5\, & \ell_4^5\,
\\
\ell_1^4\, & \ell_2^4\, & \ell_3^4\, & \ell_4^4\,
\\
\ell_1\, & \ell_2\, & \ell_3\, & \ell_4\,
\\
1\, & 1\, & 1\, & 1\,
\end{array}
\right\vert
+
\\
&
+
\frac{{\sf c}_1{\sf c}_3-{\sf c}_4}{0!\,\,2!\,\,3!\,\,5!}\
\left\vert
\begin{array}{cccc}
\ell_1^5\, & \ell_2^5\, & \ell_3^5\, & \ell_4^5\,
\\
\ell_1^3\, & \ell_2^3\, & \ell_3^3\, & \ell_4^3\,
\\
\ell_1^2\, & \ell_2^2\, & \ell_3^2\, & \ell_4^2\,
\\
1\, & 1\, & 1\, & 1\,
\end{array}
\right\vert
+
\frac{{\sf c}_4}{1!\,\,2!\,\,3!\,\,4!}\
\left\vert
\begin{array}{cccc}
\ell_1^4\, & \ell_2^4\, & \ell_3^4\, & \ell_4^4\,
\\
\ell_1^3\, & \ell_2^3\, & \ell_3^3\, & \ell_4^3\,
\\
\ell_1^2\, & \ell_2^2\, & \ell_3^2\, & \ell_4^2\,
\\
\ell_1\, & \ell_2\, & \ell_3\, & \ell_4\,
\end{array}
\right\vert
+
{\sf O}\big(\vert\ell\vert^9\big).
\endaligned
\]

\subsection{Cohomology of Schur bundles}
One could be led to presume that the cohomology dimensions:
\[
h^q
=
\dim H^q
\big(
X,\,\,\mathcal{S}^{(\ell_1,\ell_2,\dots,\ell_n)}T_X^*
\big)
\ \ \ \ \ \ \ \ \ \ \ \ \ {\scriptstyle{(q\,=\,0,\,\,1\,\cdots\,n)}}
\]
of any Schur bundle over $X$ might be expressed similarly 
by means of a general formula of the kind:
\[
\footnotesize
\aligned
h^q
&
=
\sum_{\tau_1+2\tau_2+\cdots+n\tau_n=n}\,
{\sf c}_1^{\tau_1}{\sf c}_2^{\tau_2}\cdots\,
{\sf c}_n^{\tau_n}\,
\sum_{\alpha_1+\alpha_2+\cdots+\alpha_n\leqslant
\frac{n(n+1)}{2}}\,
\\
&
\ \ \ \ \
\sum_{\ell_1\geqslant\ell_2\geqslant\cdots\geqslant\ell_n\geqslant 0}\,
h_{\tau_1,\dots,\tau_n;\,\,\alpha_1,\dots,\alpha_n}^{
q;\,\,\ell_1,\dots,\ell_n}
\cdot
(\ell_1)^{\alpha_1}(\ell_2)^{\alpha_2}
\cdots\,
(\ell_n)^{\alpha_n}
\endaligned
\]
involving the Chern classes ${\sf c}_k$, the $\ell_i$ and certain
rational coefficients $h_{ \tau_1, \dots, \tau_n\,;\,\, \alpha_1,
\dots, \alpha_n}^{ q; \, \, \ell_1, \dots, \ell_n} \in \Q$, or
alternatively, after making the substitution~\thetag{ \ref{c-d}}, 
as follows:
\[
\footnotesize
\aligned
h^q
&
=
\sum_{k=1}^{n+1}\,d^k\,
\sum_{\alpha_1+\alpha_2+\cdots+\alpha_n\leqslant
\frac{n(n+1)}{2}}\,\,
\sum_{\ell_1\geqslant\ell_2\geqslant\cdots\geqslant\ell_n\geqslant 0}\,
h_{k;\,\,\alpha_1,\dots,\alpha_n}^{
q;\,\,\ell_1,\dots,\ell_n}
\cdot
(\ell_1)^{\alpha_1}(\ell_2)^{\alpha_2}
\cdots\,
(\ell_n)^{\alpha_n}.
\endaligned
\]
However, it turns out to be already known that purely algebraic
formulas are certainly impossible, only formulas with inequalities
discussing cases
can be hoped for. Indeed, Br\"uckmann computed in~\cite{ bru1972} the
exact cohomology dimensions:
\[
\dim
H^q
\big(X,\,\,\Lambda^rT_X^*\otimes\mathcal{O}_X(t)\big)
\]
for any $q = 0, 1, \dots, n$, any $r = 0, 1, \dots, n$ and any $t \in
\Z$, where $\Lambda^r T_X^*$ identifies with $\mathcal{ S}^{ (1,
\dots, 1, 0, \dots, 0)} T_X^*$ ($r$ times $1$), and it turns out that
the obtained formulas are only piecewise polynomial with respect
to the data $(n, d, q, r, t)$. In fact, making the convention
that $\Lambda^0 T_X^* \equiv \mathcal{ O}_X ( 0)$, it is at first
well known that:
\[
\aligned
\dim H^0\big(X,\,\mathcal{O}_X(t)\big)
&
=
{\textstyle{\binom{t+n+1}{n+1}}}
-
{\textstyle{\binom{t+n+1-d}{n+1}}},
\\
\dim H^q\big(X,\,\mathcal{O}_X(t)\big)
&
=
0
\ \ \ \ \ \ \ \ \ 
\text{\rm for all}\ q\ \text{\rm with}\
1\leqslant q\leqslant n-1,
\\
\dim H^n\big(X,\,\mathcal{O}_X(t)\big)
&
=
{\textstyle{\binom{d-n-2-t+n+1}{n+1}}}
-
{\textstyle{\binom{d-n-2-t+n+1-d}{n+1}}}.
\endaligned
\]
Using then $\Lambda^n T_X^* = K_X = \mathcal{ O}_X ( d - n - 2)$, one
deduces:
\[
\aligned
\dim H^0\big(X,\,\Lambda^nT_X^*\otimes\mathcal{O}_X(t)\big)
&
=
{\textstyle{\binom{d-n-2+t+n+1}{n+1}}}
-
{\textstyle{\binom{d-n-2+t+n+1-d}{n+1}}},
\\
\dim H^q\big(X,\,\Lambda^nT_X^*\otimes\mathcal{O}_X(t)\big)
&
=
0
\ \ \ \ \ \ \ \ \ 
\text{\rm for all}\ q\ \text{\rm with}\
1\leqslant q\leqslant n-1,
\\
\dim H^n\big(X,\,\Lambda^nT_X^*\otimes\mathcal{O}_X(t)\big)
&
=
{\textstyle{\binom{-t+n+1}{n+1}}}
-
{\textstyle{\binom{-t+n+1-d}{n+1}}}.
\endaligned
\]
On the other hand, for $1 \leqslant r \leqslant n - 1$, Br\"uckmann 
(\cite{ bru1972}) obtained complete dimension formulas:
\[
\aligned
\dim H^0\big(X,\,\Lambda^rT_X^*\otimes\mathcal{O}_X(t)\big)
&
=
{\textstyle{\binom{t-1}{r}}}
{\textstyle{\binom{t+n+1-r}{n+1-r}}},
\\
\dim H^q\big(X,\,\Lambda^rT_X^*\otimes\mathcal{O}_X(t)\big)
&
=
\delta_{q,r}\cdot\delta_{t,0}
\ \ \ \ \  
\text{\rm for}\,
1\leqslant q\leqslant n-1, q+r\neq n,
\\
\dim H^{n-r}\big(X,\,\Lambda^rT_X^*\otimes\mathcal{O}_X(t)\big)
&
=
\sum_{\mu=0}^{n+2}\,
(-1)^\mu\,
{\textstyle{\binom{n+2}{\mu}}}
{\textstyle{\binom{-t-rd-(\mu-1)(d-1)}{n+1}}}
+
\\
&
\ \ \ \ \
+
\delta_{n,2r}\cdot\delta_{t,0},
\\
\dim H^n\big(X,\,\Lambda^rT_X^*\otimes\mathcal{O}_X(t)\big)
&
=
{\textstyle{\binom{-t-1}{n-r}}}
{\textstyle{\binom{-t+n+1-2r}{n+1-2r}}}.
\endaligned
\]
Clearly, these formulas are only `semi-algebraic'. One does not find in
the literature complete formulas for cohomology dimensions of Schur
bundles having at least three distinct row lengths.

\subsection{Majorating the cohomology}
Rousseau's strategy developed in~\cite{ rou2006b} and in~\cite{
dmr2010} for dimensions 3 and 4 consists in avoiding exact, probably
unfeasible cohomology computations and in substituting for that
cohomology {\em inequalities}.

Let as before $X$ be a geometrically smooth projective algebraic
complex hypersurface in $\P^{ n+1} ( \C)$. Let $Fl ( T_X^*)$ denote
the (complete) flag manifold of $T_X^*$ which organizes as a
holomorphic vector bundle $\pi \colon Fl ( T_X^*) \to X$ of rank
$\frac{ n ( n +1)}{ 2}$ over $X$, the fiber of which above an
arbitrary point $x \in X$ consists of complete flags:
\[
0=E_{0,x}
\subset
E_{1,x}
\subset
\cdots
\subset
E_{n,x}=T_{X,x},
\]
where $\dim E_{ i, x} = i$. Let as before $\ell = (\ell_1, \ell_2,
\dots, \ell_n)$ with $\ell_1 \geqslant \ell_2 \geqslant \cdots
\geqslant \ell_n \geqslant 0$. According to Bott (\cite{ botu1982}),
there is a canonical {\em line} bundle $\mathcal{ B}^\ell ( T_X^*)$
over $Fl ( T_X^*)$ with the property that the Schur bundle $\mathcal{
S}^{ (\ell_1, \dots, \ell_n)} T_X^* \to X$ coincides with the direct
image $\pi_* ( \mathcal{ B}^\ell) \to X$ and whose fiber above an
arbitrary flag $E_x \in Fl ( T_X^*)$ is $\otimes_{i=1}^n\, \big(
\det(E_{ x, i}/E_{x, i-1}) \big)^{\otimes\ell_i}$. The fundamental
theorem of Bott
(\cite{ botu1982}) states that the two bundles $\mathcal{ S}^{ ( \ell_1,
\dots, \ell_n)} T_X^*$ and $\mathcal{ B}^\ell (T_X^*)$ have the same
cohomology, and it is therefore somewhat more convenient to deal with
$\mathcal{ B}^\ell (T_X^*)$, because {\em line} bundles are better
understood and more studied.

In fact, a certain control of the cohomology by means of inequalities
is available thanks to the so-called {\sl Holomorphic Morse
inequalities} due to Demailly which state as follows in a general
version (\cite{ dem2000}) suitable for applications devised by Trapani
(\cite{ trap1995}). Let $\mathcal{ E} \to X$ be a completely arbitrary
holomorphic vector bundle of rank $r \geqslant 1$ over a compact
Kähler manifold of dimension $n$, and let $\mathcal{ L} \to X$ be a
holomorphic line bundle subjected to the specific restriction that it
can be written as the difference: $\mathcal{ L} = \mathcal{ F} \otimes
\mathcal{ G}^{ -1}$ between two line bundles that are ample, or more
generally, numerically effective. Then about $\mathcal{ L}^k \otimes
\mathcal{ E}$ as $k \to \infty$, we have the following two collections
of asymptotic inequalities, firstly for plain cohomology dimensions,
secondly for their alternating sums:

\smallskip\noindent
{\small\sf $\bullet$\,\,Weak Morse inequalities:}
For any $q = 0, 1, \dots, n$, one has:
\[
\aligned
h^q\big(X,\,\,
\mathcal{L}^k\otimes\mathcal{E}\big)
&
\leqslant
r\,k^n\,\frac{1}{(n-q)!\,\,q!}\,
\int_X\,
{\sf c}_1(\mathcal{F})^{n-q}
\cdot
{\sf c}_1(\mathcal{G})^q
+
{\sf o}(k^n)
\\
&
\ \ \ \ \ \ \ \ \ \ \ \ \ \ \ \ \ \ \ \ 
{\scriptstyle{(q\,=\,0,\,\,1\,\cdots\,n)}}.
\endaligned
\]

\smallskip\noindent
{\small\sf $\bullet$\,\,Strong Morse inequalities:}
For any $q = 0, 1, \dots, n$, one has:
\[
\aligned
&
\sum_{0\leqslant q'\leqslant q}\,
(-1)^{q-q'}\,
h^{q'}\big(
X,\,\,
\mathcal{L}^k\otimes\mathcal{E}\big)
\leqslant
\\
&
\leqslant
r\,k^n\,
\sum_{0\leqslant q'\leqslant q}\,
\frac{(-1)^{q-q'}}{(n-q')!\,\,q'!}\,
\int_X\,
{\sf c}_1(\mathcal{F})^{n-q'}
\cdot
{\sf c}_1(\mathcal{G})^{q'}
+
{\sf o}(k^n).
\endaligned
\]

\smallskip\noindent
An algebraic proof of these inequalities (without $\mathcal{ E}$ and
for $X$ projective) by plain induction on dimension but not using any
tools from Analysis was given by Angelini in~\cite{ ange1996}. We then
borrow this scheme of proof, as it was applied by Rousseau (\cite{
rou2006b}) within the Schur bundle context. Weak type
inequalities will
suffice for us, and the goal is somehow to represent $\mathcal{
B}^\ell (T_X^*)$ as a difference between two line bundles that will be
positive, hence ample.

To begin with, since $T_X^* \otimes \mathcal{ O}_X ( 2)$ is generated
by its global sections, it is semi-positive. According to a general
property (\cite{ dem1987}), if a holomorphic vector bundle $\mathcal{
E} \to X$ is semi-positive, i.e. if $E \geqslant 0$, then the
corresponding line bundle $\mathcal{ B}^\ell ( \mathcal{ E})$ is also
semi-positive, i.e. $\mathcal{ B}^\ell ( \mathcal{ E}) \geqslant 0$.
Applying this to $\mathcal{ E} := T_X^* \otimes \mathcal{ O}_X ( 2)$,
we get, thanks to a natural isomorphism, that:
\begin{equation}
\label{semi-positive-B}
\mathcal{B}^\ell
\big(T_X^*\otimes\mathcal{O}_X(2)\big)
\simeq
\mathcal{B}^\ell(T_X^*)
\otimes
\pi^*\mathcal{O}_X(2\vert\ell\vert)
\geqslant 0
\end{equation}
is semi-positive, where $\vert \ell \vert = \ell_1 + \cdots + \ell_n$.
Tensoring then by $\pi^*
\mathcal{ O}_X ( \vert \ell \vert) > 0$, it thus
trivially follows that:
\[
\mathcal{B}^\ell(T_X^*)
\otimes
\pi^*\mathcal{O}_X(3\vert\ell\vert)
>
0
\]
is positive. Hence we can write (somehow artificially) $\mathcal{
B}^\ell ( T_X^*)$, which we will now write $\mathcal{ B}^\ell$ for
short, as the following difference:
\[
\mathcal{B}^\ell
=
\big[
\mathcal{B}^\ell\otimes\pi^*\mathcal{O}_X(3\vert\ell\vert)
\big]
\otimes
\big[
\pi^*\mathcal{O}_X(3\vert\ell\vert)
\big]^{-1}
\]
between two positive line bundles over $Fl ( T_X^*)$, with plainly:
\[
\mathcal{F}
:=
\mathcal{B}^\ell
\otimes
\pi^*\mathcal{O}_X(3\vert\ell\vert)
\ \ \ \ \ \ \ \ \ 
\text{\rm and}
\ \ \ \ \ \ \ \ \ 
\mathcal{G}
:=
\pi^*\mathcal{O}_X(3\vert\ell\vert),
\]
in the above notations for Morse inequalities.

Following Angelini and Rousseau, we need even more in order 
to force the
positive cohomologies $H^q \big( Fl( T_X^*), \, \mathcal{ F} \big)$,
$q = 1, \dots, n$, to be vanishing. We remind the Kodaira vanishing
theorem which stipulates that, on a projective algebraic complex
manifold $Z$, for every ample line bundle $\mathcal{ A} \to X$ one
has:
\[
0
=
H^q\big(Z,\,\mathcal{A}\otimes K_Z),
\]
for all $q = 1, \dots, n$. So on the flag manifold $Z := Fl ( T_X^*)$,
we not only need that $\mathcal{ F}$ be positive (hence ample), but we
need also, after decomposing in advance:
\[
\mathcal{F}
=
\big(\mathcal{F}\otimes (K_{Fl(T_X^*)})^{-1}\big)
\otimes
K_{Fl(T_X^*)},
\]
that $\mathcal{ A} := \mathcal{F}\otimes (K_{Fl(T_X^*)})^{-1}$
be positive (hence ample). 

For this, we recall at first the known 
isomorphisms (\cite{botu1982, dem1987}):
\[
\aligned
K_{Fl(T_X^*)}
&
\simeq
\big[
\mathcal{B}^{2n-1,\dots,3,1}
\big]^{-1}
\otimes
\pi^*(K_X)^{\otimes(n+1)}
\\
&
\simeq
\big[
\mathcal{B}^{2n-1,\dots,3,1}
\big]^{-1}
\otimes
\pi^*
\mathcal{O}_X
\big((n+1)(d-n-2)\big),
\endaligned
\]
from which we hence deduce:
\[
\aligned
&
\mathcal{F}
\otimes
(K_{Fl(T_X^*)})^{-1}
\simeq
\mathcal{B}^\ell
\otimes
\mathcal{B}^{2n-1,\dots,3,1}
\otimes
\pi^*
\mathcal{O}_X
\big(
3\vert\ell\vert
-
(n+1)(d-n-2)
\big)
\\
&
\simeq
\mathcal{B}^{\ell_1+2n-1,\dots,\ell_{n-1}+3,\ell_n+1}
\otimes
\pi^*
\mathcal{O}_X
\big(
3\vert\ell\vert
-
(n+1)(d-n-2)
\big).
\endaligned
\]
But similarly as in~\thetag{ \ref{semi-positive-B}} a short while ago,
we know that the bundle:
\[
\mathcal{B}^{\ell_1+2n-1,\dots,\ell_{n-1}+3,\ell_n+1}
\otimes
\pi^*
\mathcal{O}_X
\big(
2[\ell_1+2n-1+\cdots+\ell_{n-1}+3+\ell_n+1]
\big)
\]
is semi-positive, whence it is surely positive after it is tensored
only by $\pi^* \mathcal{ O}_X ( 1)$. Consequently, 
observing $2n-1 + \cdots + 3 + 1 = n^2$, our bundle
$\mathcal{ F} \otimes (K_{ Fl(T_X^*) })^{-1}$ will be positive when:
\[
3\vert\ell\vert
-
(n+1)(d-n-2)
\geqslant 
1
+
2(\vert\ell\vert+n^2), 
\]
that is to say when:
\[
\vert\ell\vert
\geqslant
1+2n^2
+
(n+1)(d-n-2),
\]
or with less effective information, when $\vert \ell \vert \geqslant
{\sf Constant}_{ n, d}$. Under this restriction concerning $\vert \ell
\vert$ which insures the applicability of Kodaira's vanishing theorem,
Rousseau's scheme of proof works in arbitrary dimension $n$
(cf.~\cite{ rou2006b} and also~\cite{ dmr2010} for the case $n = 4$), and
it yields the following majorations: 
\[
\footnotesize
\aligned
h^q
\big(X,\,
\mathcal{S}^{(\ell_1,\dots,\ell_n)}T_X^*
\big)
&
\leqslant
\chi\big(X,\,
\mathcal{S}^{(\ell_1,\dots,\ell_n)}T_X^*
\otimes
\mathcal{O}_X(3(q+1)\vert\ell\vert)
\big)
-
\\
&
\ \ \ \ \
-
{\textstyle{\binom{q}{1}}}\,
\chi\big(X,\,
\mathcal{S}^{(\ell_1,\dots,\ell_n)}T_X^*
\otimes
\mathcal{O}_X(3q\vert\ell\vert)
\big)
+
\\
&
\ \ \ \ \
+
\cdots\cdots\cdots\cdots\cdots\cdots
\cdots\cdots\cdots\cdots\cdots\cdots
+
\\
&
\ \ \ \ \
+(-1)^{q-1}\,
{\textstyle{\binom{q}{q-1}}}\,
\chi\big(X,\,
\mathcal{S}^{(\ell_1,\dots,\ell_n)}T_X^*
\otimes
\mathcal{O}_X(6\vert\ell\vert)
\big)
+
\\
&
\ \ \ \ \
+(-1)^q\,
{\textstyle{\binom{q}{q}}}\,
\chi\big(X,\,
\mathcal{S}^{(\ell_1,\dots,\ell_n)}T_X^*
\otimes
\mathcal{O}_X(3\vert\ell\vert)
\big),
\endaligned
\]
in terms of alternating sums of Euler-Poincar\'e characteristics.
Applying Br\"uckmann's formula for the explicit computation of the
appearing Euler-Poincar\'e characteristics (Theorem~4
in~\cite{bru1997}), we then get the following result.

\begin{Theorem}
\label{majoration-Rousseau}
Let $X = X^n \subset \P^{ n+1} ( \C)$ be a geometrically smooth
projective algebraic complex hypersurface of general type, i.e. of
degree $d \geqslant n+3$, and let $\ell = ( \ell_1, \dots, \ell_{ n-1},
\ell_n)$ with $\ell_1 \geqslant \cdots \geqslant \ell_{ n-1} \geqslant
\ell_n \geqslant 0$. If:
\[
\vert\ell\vert
=
\ell_1+\cdots+\ell_{n-1}+\ell_n
\geqslant
{\sf Constant}_{n,d},
\]
then for every $q = 1, 2, \dots, n$, the dimensions of the positive
cohomology groups of the Schur bundle $\mathcal{ S}^{ (\ell_1, 
\dots, \ell_{ n-1}, \ell_n)} T_X^*$ over $X$ satisfy a general
majoration of the form:
\[
\aligned
&
h^q\big(
X,\,\,
\mathcal{S}^{(\ell_1,\dots,\ell_{ n-1},\ell_n)}T_X^*
\big)
\leqslant
\\
&
\leqslant
{\sf Constant}_{n,d}
\prod_{1\leqslant i<j\leqslant n}
(\ell_i-\ell_j)\,
\bigg[
\sum_{\beta_1+\cdots+\beta_{n-1}+\beta_n=n}
\!\!\!
\ell_1^{\beta_1}\cdots\,\ell_{n-1}^{\beta_{n-1}}\ell_n^{\beta_n}
\bigg]
+
\\
&
\ \ \ \ \ \ \ \ \ \ \ \ \
\ \ \ \ \ \ \ \ \ \ \ \ \
\ \ \ \ \ \ \ \ \ \ 
+
{\sf Constant}_{n,d}\,
\bigg[
\sum_{\alpha_1+\dots+\alpha_n\leqslant\frac{n(n+1)}{2}-1}\,
\ell_1^{\alpha_1}\cdots\,\ell_n^{\alpha_n}
\bigg],
\endaligned
\]
with leading terms being homogeneous of degree $\frac{ n ( n+1)}{
2}$ with respect to the $\ell_i$ and divisible by all the differences
$(\ell_i - \ell_j)$, where $1 \leqslant i < j \leqslant n$.
\end{Theorem}

For the estimates that we will conduct in the next sections, we need
none of the three ${\sf Constant}_{ n, d}$ above to be
effective. Admitting this, raising if necessary the two ${\sf
Constants}_{ n, d}$ appearing in the right-hand side, it follows that
the majoration is in fact valid for every $\ell$, since the
restriction that $\vert \ell \vert$ be large enough can obviously be
absorbed by the ${\sf Constants}_{ n, d}$. Also, one must observe
that the Euler-Poincar\'e characteristic provided by 
Theorem~\ref{E-P-characteristic} satisfies the same kind of majoration,
hence {\em all} cohomology dimensions $h^0, h^1, \dots, h^n$ do the
same. However, we want to underline that, even with an effective
control on ${\sf Constant}_{ n, d}$ similar as in~\cite{ rou2006b,
dmr2010} for $n = 3$ and $n = 4$, the above kind of majoration cannot
at all conduct to the optimal degree bound $d \geqslant n + 3$ of the
Main Theorem, because we will see that the presence of the monomial
$\ell_n^n$ in $\sum_{ \beta} \, \ell^\beta$ forces to lose a nonzero
portion of the $(\log \kappa)^n$ entering in the Euler-Poincar\'e
characteristic when summing $\sum\, M_{\ell}^{\kappa, m}\, \mathcal{
S}^{ (\ell)}$ over Schur bundles.

\markleft{Jo\"el Merker}
\markright{\S6.~Emergence of basic numerical sums}
\section{\bf Emergence of basic numerical sums}
\label{Section-6}

\subsection{Expanding and rewriting}
At least, the explicit formula for the Euler characteristic and the
cohomology bounds for Schur bundles firmly motivates to consider basic
numerical sums of the form:
\[
\sum_{\ell_1\geqslant\ell_2\geqslant\cdots\geqslant\ell_n\geqslant 0}\,
M_{\ell_1,\ell_2,\dots,\ell_n}^{\kappa,m}\,
\ell_1^{\alpha_1}\ell_2^{\alpha_2}\cdots\ell_n^{\alpha_n},
\]
for any $n$-tuple of nonnegative integers $\alpha = (\alpha_1,
\alpha_2, \dots, \alpha_n) \in \N^n$ of length $\leqslant \frac{ n (
n+1)}{ 2} - 1$ if remainders are to be taken into consideration, or
else of length $\vert \alpha \vert = \alpha_1 + \cdots + \alpha_n$
constant equal to $\frac{ n ( n+1)}{ 2}$, if major terms are to be
studied.  After some reflections based on manuscript explorations and
on intense thought, it appears {\em a posteriori} convenient, if not
adequate, to express all quantities in terms of the successive
differences between horizontal lengths in the Young diagram:
\[
\ell_1-\ell_2,
\ \ \ \ \ \
\ell_2-\ell_3,
\ \ \ \ \ \
\cdots\cdots,
\ \ \ \ \ \
\ell_{n-1}-\ell_n,
\ \ \ \ \ \
\ell_n,
\]
that is to say, in terms of the horizontal lengths of the appearing
successive blocks of constant depths.
Thus accordingly, we may rewrite any appearing monomial
$\ell_1^{ \alpha_1} \cdots\, \ell_n^{ \alpha_n}$ by inserting
differences as follows:
\[
\aligned
\ell_1^{\alpha_1}\ell_2^{\alpha_2}\cdots
\ell_{n-1}^{\alpha_{n-1}}\ell_n^{\alpha_n}
=
\big(
\ell_1-\ell_2+\ell_2-\ell_3+\cdots+\ell_{n-1}-\ell_n+\ell_n
\big)^{\alpha_1}
\cdot
\\
\cdot
\big(
\ell_2-\ell_3+\cdots+\ell_{n-1}-\ell_n+\ell_n
\big)^{\alpha_2}
\cdot
\\
\cdot\cdots\cdots\cdots\cdots\cdots\cdots\cdots\cdot
\\
\cdot
\big(
\ell_{n-1}-\ell_n+\ell_n
\big)^{\alpha_{n-1}}
\cdot
\\
\cdot
\big(\ell_n\big)^{\alpha_n},
\endaligned
\]
and then we may simply expand all the appearing powers to obtain a
certain sum, with integer integer coefficients, of interesting
monomials of the specific form:
\[
\big(\ell_1-\ell_2\big)^{\alpha_1'}
\big(\ell_2-\ell_3\big)^{\alpha_2'}
\cdots\cdots
\big(\ell_{n-1}-\ell_n\big)^{\alpha_{n-1}'}
\big(\ell_n\big)^{\alpha_n'},
\]
the total degree in the $\ell_i$ being evidently preserved: 
\[
\alpha_1'+\alpha_2'+\cdots+\alpha_{n-1}'+\alpha_n'
=
\alpha_1+\alpha_2+\cdots+\alpha_n.
\]
Only multinomial coefficients being involved in the expansion,
we have a simple inequality of the form:
\begin{equation}
\label{l-l-l}
\aligned
&
\ell_1^{\alpha_1}
\cdots\,
\ell_{n-1}^{\alpha_{n-1}}\ell_n^{\alpha_n}
\leqslant
\\
&
\leqslant
{\sf Constant}_n\cdot
\!\!\!\!\!\!\!\!\!
\sum_{\alpha_1'+\cdots+\alpha_{n-1}'+\alpha_n'
=
\alpha_1+\cdots+\alpha_{n-1}+\alpha_n}
\!\!\!\!\!\!\!\!\!
\big(\ell_1-\ell_2\big)^{\alpha_1'}
\cdots
\big(\ell_{n-1}-\ell_n\big)^{\alpha_{n-1}'}
\big(\ell_n\big)^{\alpha_n'}.
\endaligned
\end{equation}

\subsection{Basic numerical sums}
As a consequence, in order to verify that the contribution in $\sum\,
M_\ell^{ \kappa, m} \cdot \ell^\alpha$ of any monomial
$\ell_1^{\alpha_1} \cdots \ell_{ n-1}^{ \alpha_{ n-1}} \ell_n^{
\alpha_n}$ of total degree:
\[
\alpha_1+\cdots+\alpha_{ n-1}+\alpha_n
\leqslant 
{\textstyle{\frac{n(n+1)}{2}}}-1
\]
which possibly appears in a general remainder of the form ${\sf O}_{
n, d} \big( \vert \ell \vert^{ \frac{ n ( n+1)}{ 2} - 1} \big)$ still
falls into the corresponding $m$-remainder ${\sf O}_{ n, \kappa} \big(
m^{ (\kappa + 1)n - 2} \big)$, we are led back to studying the
asymptotic behavior, as $m \to \infty$, of basic numerical sums of the
general form:
\[
\sum_{\ell_1\geqslant\ell_2\geqslant\cdots\geqslant\ell_{n-1}
\geqslant\ell_n\geqslant 0}\,
M_{\ell_1,\ell_2,\dots,\ell_{n-1},\ell_n}^{\kappa,m}
\cdot
\big(\ell_1-\ell_2\big)^{\alpha_1'}
\cdots\cdots
\big(\ell_{n-1}-\ell_n\big)^{\alpha_{n-1}'}
\big(\ell_n\big)^{\alpha_n'}, 
\]
with again $\alpha_1' + \cdots + \alpha_{ n-1} ' + \alpha_n' \leqslant
\frac{ n ( n+1)}{ 2} - 1$, and now, everything has become purely
combinatorial, that is to say, complex geometry concepts have entirely
disappeared.

On the other hand, after expanding any $\prod_{ i < j} ( \ell_i -
\ell_j) \, \ell^\beta$ with $\vert \beta \vert = n$ appearing both as
principal terms in the Euler-Poincar\'e characteristic of $\mathcal{
S}^{ (\ell)} T_X^*$ and in the cohomology majorations provided by 
Theorem~\ref{majoration-Rousseau}, and after performing the
rewriting in terms of $\ell_1 - \ell_2$, \dots, $\ell_{ n-1} -
\ell_n$, $\ell_n$ as above, we are led back to estimating the same
kind of basic numerical sums, but this time with $\alpha_1 ' + \cdots
+ \alpha_{ n-1}' + \alpha_n'= \frac{ n ( n+1)}{ 2}$.

Since we will not attempt to compute, even asymptotically, the
multiplicities $M_{\ell_1,\dots,\ell_n}^{\kappa,m}$ for which only
semi-algebraic formulas could exist as an examination for small values
of $n$ and $\kappa$ shows, we will rewrite such basic numerical sums
under the following more archetypal form\footnote{\,
The equality written follows immediately from the definitions: the
passage from the second line to the first line just consists in
counting the semi-standard Young Tableaux of weight $m$ which have the
same underlying Young diagram ${\sf YD}_{(\ell_1, \ell_2, \dots,
\ell_{ n-1}, \ell_n)}$, and their number is just what we denoted by
the multiplicity $M_{ \ell_1, \ell_2, \dots, \ell_{ n-1}, \ell_n}^{
\kappa, m}$. 
}: 
\begin{equation}
\label{M-YT}
\boxed{
\footnotesize
\aligned
&
\sum_{\ell_1\geqslant\ell_2\geqslant\cdots\geqslant\ell_{n-1}
\geqslant\ell_n\geqslant 0}\,
M_{\ell_1,\ell_2,\dots,\ell_{n-1},\ell_n}^{\kappa,m}
\cdot
\big(\ell_1-\ell_2\big)^{\alpha_1'}
\cdots\cdots
\big(\ell_{n-1}-\ell_n\big)^{\alpha_{n-1}'}
\big(\ell_n\big)^{\alpha_n'}
=
\\
&
=
\sum_{{\sf YT}\,{\sf semi-standard}
\atop
{\sf weight}({\sf YT})=m}
\big(\ell_1({\sf YT})-\ell_2({\sf YT})\big)^{\alpha_1'}
\cdots
\big(\ell_{n-1}({\sf YT})-\ell_n({\sf YT})\big)^{\alpha_{n-1}'}
\big(\ell_n({\sf YT})\big)^{\alpha_n'}
\endaligned}\,,
\end{equation}
where in the second line, ${\sf YT}$ runs over all the possible
semi-standard Young tableaux, where $\ell_i ( {\sf YT} )$ denote the
length of the $i$-th line of ${\sf YT}$, and where as before ${\sf
weight} ( {\sf YT} )$ denotes the total number of primes appearing in
the associated $\Delta$-monomial, that is to say, the sum of all the
$\lambda_i^j$ occupying the squares of ${\sf YT}$:
\[
\aligned
{\sf weight}
\big({\sf YT}\big)
&
=
{\sf weight}
\big(
{\sf YD}_{(\ell_1,\dots,\ell_n)}(\lambda_i^j)
\big)
\\
&
=
\sum_{1\leqslant j_1\leqslant\ell_1}\,
\lambda_1^{j_1}
+
\sum_{1\leqslant j_2\leqslant\ell_2}\,
\lambda_2^{j_2}
+\cdots+
\sum_{1\leqslant j_n\leqslant\ell_n}\,
\lambda_n^{j_n}.
\endaligned
\]
Then more tractable computations and partially explicit formulas will
come up as being somewhat available in the next sections.

Thus, assuming from now on 
that $\kappa \geqslant n$ is at least equal to the
dimension, our first main goal will be to establish 
(Corollary~\ref{corollary-remainder-1} below) that for every
$(\alpha_1', \dots, \alpha_{ n-1}', \alpha_n') \in \N^n$, the
following precise logarithmic-like majoration holds:
\[
\aligned
\sum_{{\sf YT}\,{\sf semi-standard}
\atop
{\sf weight}({\sf YT})=m}
&
\big(\ell_1({\sf YT})-\ell_2({\sf YT})\big)^{\alpha_1'}
\cdots
\big(\ell_{n-1}({\sf YT})-\ell_n({\sf YT})\big)^{\alpha_{n-1}'}
\big(\ell_n({\sf YT})\big)^{\alpha_n'}
\leqslant
\\
&
\leqslant
{\sf Constant}_{n,\kappa}\cdot
m^{\alpha_1'+\cdots+\alpha_{n-1}'+\alpha_n'}\cdot
m^{n\kappa-\frac{n(n-1)}{2}}.
\endaligned
\]
Applying these majorations when $\vert \alpha' \vert \leqslant \frac{
n ( n+1)}{ 2} - 1$, it will then follow in particular that the
right-hand side majorant is an ${\sf O}_{ n, \kappa} \big( m^{
(n+1)\kappa - 2} \big)$, whence remainders match through summation in
Euler-Poincar\'e characteristics, as was announced a bit
after Theorem~\ref{E-P-characteristic}.

Afterward, we will study what arises when $\alpha_1' + \cdots +
\alpha_{ n-1}' + \alpha_n' = \frac{ n ( n+1)}{ 2}$. In any case, we
need to analyze more deeply what comes out from the semi-standard
Young tableaux of weight $m$.

\markleft{Jo\"el Merker}
\markright{\S7.~Asymptotic combinatorics 
of semi-standard Young tableaux}
\section{\bf Asymptotic combinatorics 
\\
of semi-standard Young tableaux}
\label{Section-7}

\subsection{Repetitions in the $\Delta$-monomials}
In the general $\Delta$-monomial modulo the Pl\"ucker relations
given by~\thetag{ \ref{general-Delta-monomial}}:
\[
\prod_{d_1\geqslant i\geqslant 1}\,\,
\prod_{1+\ell_{i+1}\leqslant j\leqslant\ell_i}\!\!
\Delta_{1,\dots,i}^{
\lambda_1^j,\dots,\lambda_i^j},
\] 
there may exist (several) repetitions of a given determinant
$\Delta_{1,\dots,i}^{ \lambda_1^j,\dots,\lambda_i^j}$, since in the
semi-standard Young tableau, the increasing property enjoyed by the
$\lambda_i^j$ is only weak along its rows. So in ${\sf YD}_{ ( \ell_1,
\dots, \ell_n)} (\lambda_i^j)$, we should describe with more explicit
information the typical block of depth $i$:

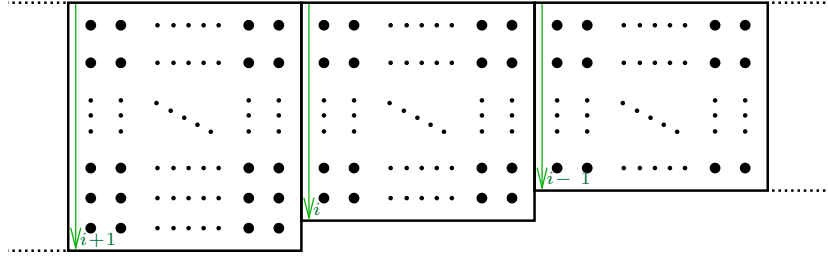
\smallskip
\begin{center}
\begin{picture}(0,0)%
\includegraphics{blocks.pstex}%
\end{picture}%
\setlength{\unitlength}{4144sp}%
\begingroup\makeatletter\ifx\SetFigFont\undefined%
\gdef\SetFigFont#1#2#3#4#5{%
  \reset@font\fontsize{#1}{#2pt}%
  \fontfamily{#3}\fontseries{#4}\fontshape{#5}%
  \selectfont}%
\fi\endgroup%
\begin{picture}(4949,1529)(744,-1298)
\put(1195,-1240){\makebox(0,0)[lb]{\smash{{\SetFigFont{7}{8.4}{\familydefault}{\mddefault}{\updefault}{\color[rgb]{0,.82,0}\green{$i\!+\!1$}}%
}}}}
\put(2591,-1064){\makebox(0,0)[lb]{\smash{{\SetFigFont{7}{8.4}{\familydefault}{\mddefault}{\updefault}{\color[rgb]{0,.82,0}\green{$i$}}%
}}}}
\put(3989,-874){\makebox(0,0)[lb]{\smash{{\SetFigFont{7}{8.4}{\familydefault}{\mddefault}{\updefault}{\color[rgb]{0,.82,0}\green{$i\!-\,1$}}%
}}}}
\end{picture}%

\end{center}

\noindent
which is naturally located between a block of depth $i + 1$ on its
left, and a block of depth $i - 1$ on its right. 
To this aim, let us rewrite
such a block as follows:
\[
\small
\aligned
\boxed{\,
\begin{array}{c}
\lambda_1^{1+\ell_{i+1}}
\\
\lambda_2^{1+\ell_{i+1}}
\\
\cdot\cdot
\\
\lambda_i^{1+\ell_{i+1}}
\end{array}
\begin{array}{c}
\cdots
\\
\cdots
\\
\cdot\cdot
\\
\cdots
\end{array}
\begin{array}{c}
\lambda_1^{\ell_i}
\\
\lambda_2^{\ell_i}
\\
\cdot\cdot
\\
\lambda_i^{\ell_i}
\end{array}
\!}
=
\boxed{
\left[\!\!
\begin{array}{c}
\mu_1^j
\\
\mu_2^j
\\
\cdot\cdot
\\
\mu_i^j
\end{array}
\!\!\right]^{\!\!a_{\mu_1^j,\mu_2^j,\dots,\mu_i^j}}
\!\!\!\cdots
\left[\!\!
\begin{array}{c}
\lambda_1^j
\\
\lambda_2^j
\\
\cdot\cdot
\\
\lambda_i^j
\end{array}
\!\!\right]^{\!\!a_{\lambda_1^j,\lambda_2^j,\dots,\lambda_i^j}}
\!\!\!\cdots
\left[\!\!
\begin{array}{c}
\nu_1^j
\\
\nu_2^j
\\
\cdot\cdot
\\
\nu_i^j
\end{array}
\!\!\right]^{\!\!a_{\nu_1^j,\nu_2^j,\dots,\nu_i^j}}
}\,.
\endaligned
\] 
Here firstly, looking at the two extreme (right and left) columns, we
changed the notation for later purposes, denoting $\mu_l^j :=
\lambda_l^{ 1 + \ell_{ i+1}}$ and $\nu_l^j := \lambda_l^{\ell_i}$ 
for any row
index $l = 1, 2, \dots, i$; secondly, the appearing exponents
$a_{*,*,\dots,*}$ are meant to denote repetitions of (bracketed) columns,
so that naturally their sum equals the horizontal length of the
initially considered $i$-th block:
\[
\ell_i
-
\ell_{i+1}
=
a_{\mu_1^j,\mu_2^j,\dots,\mu_i^j}
+
\cdots
+
a_{\lambda_1^j,\lambda_2^j,\dots,\lambda_i^j}
+
\cdots
+
a_{\nu_1^j,\nu_2^j,\dots,\nu_i^j}; 
\]
thirdly and lastly, the succession of columns now increases strictly
when one disregards the repetitions:
\[
\left[\!\!
\begin{array}{c}
\mu_1^j
\\
\mu_2^j
\\
\cdot\cdot
\\
\mu_i^j
\end{array}
\!\!\right]
<\cdots<
\left[\!\!
\begin{array}{c}
\lambda_1^j
\\
\lambda_2^j
\\
\cdot\cdot
\\
\lambda_i^j
\end{array}
\!\!\right]
<\cdots<
\left[\!\!
\begin{array}{c}
\nu_1^j
\\
\nu_2^j
\\
\cdot\cdot
\\
\nu_i^j
\end{array}
\!\!\right],
\]
where by definition we declare that a column $(\lambda_l' )_{ 1
\leqslant l \leqslant i}$ is smaller (strictly) than another column
$(\lambda_l'' )_{ 1 \leqslant l \leqslant i}$ if all its row elements
are smaller (weakly): $\lambda_l ' \leqslant \lambda_l''$ for $l = 1,
\dots, i$, {\em and if} 
there exists at least one row index $l_0$ for which
$\lambda_{ l_0} ' < \lambda_{ l_0} ''$. As a result, we have
represented our typical semi-standard Young tableau of depth $d_1$ as
follows by emphasizing precisely the column repetitions, all the
appearing columns being now pairwise distinct and ordered
increasingly:
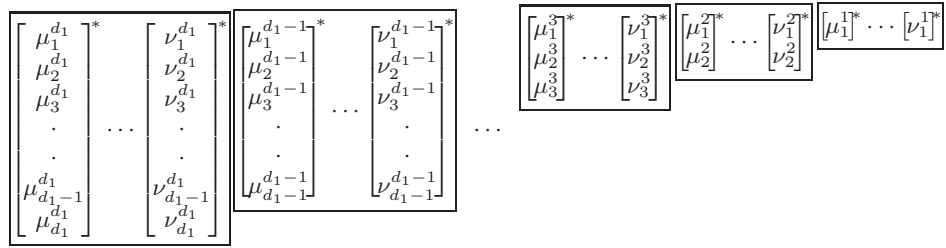
\[
\footnotesize
\aligned
\boxed{\!\!
\left[\!\!\!\!
\begin{array}{c}
\mu_1^{d_1}
\\
\mu_2^{d_1}
\\
\mu_3^{d_1}
\\
\cdot
\\
\cdot
\\
\mu_{d_1-1}^{d_1}
\\
\mu_{d_1}^{d_1}
\end{array}
\!\!\!\!\right]^{\!\!*}
\cdots
\left[\!\!\!\!
\begin{array}{c}
\nu_1^{d_1}
\\
\nu_2^{d_1}
\\
\nu_3^{d_1}
\\
\cdot
\\
\cdot
\\
\nu_{d_1-1}^{d_1}
\\
\nu_{d_1}^{d_1}
\end{array}
\!\!\!\!\right]^{\!\!*}\!\!}
\!\!\!\!\:\:\!
\begin{array}{c}
\boxed{\!\!
\left[\!\!\!\!
\begin{array}{c}
\mu_1^{d_1-1}
\\
\mu_2^{d_1-1}
\\
\mu_3^{d_1-1}
\\
\cdot
\\
\cdot
\\
\mu_{d_1-1}^{d_1-1}
\end{array}
\!\!\!\!\right]^{\!\!*}
\cdots
\left[\!\!\!\!
\begin{array}{c}
\nu_1^{d_1-1}
\\
\nu_2^{d_1-1}
\\
\nu_3^{d_1-1}
\\
\cdot
\\
\cdot
\\
\nu_{d_1-1}^{d_1-1}
\end{array}
\!\!\!\!\right]^{\!\!*}\!\!}
\\ 
\rule[-2.17pt]{0pt}{13pt}
\end{array}
\cdots
\begin{array}{c}
\boxed{\!\!
\left[\!\!\!\!
\begin{array}{c}
\mu_1^3
\\
\mu_2^3
\\
\mu_3^3
\end{array}
\!\!\!\!\right]^{\!\!*}
\cdots
\left[\!\!\!\!
\begin{array}{c}
\nu_1^3
\\
\nu_2^3
\\
\nu_3^3
\end{array}
\!\!\!\!\right]^{\!\!*}\!\!}
\\ 
\rule[-2.17pt]{0pt}{53pt}
\end{array}
\!\!\!\!\:\:\!\!\!\!
\begin{array}{c}
\boxed{\!\!
\left[\!\!\!\!
\begin{array}{c}
\mu_1^2
\\
\mu_2^2
\end{array}
\!\!\!\!\right]^{\!\!*}
\cdots
\left[\!\!\!\!
\begin{array}{c}
\nu_1^2
\\
\nu_2^2
\end{array}
\!\!\!\!\right]^{\!\!*}\!\!}
\\ 
\rule[-2.17pt]{0pt}{65pt}
\end{array}
\!\!\!\!\:\:\!\!\!\!
\begin{array}{c}
\boxed{\!\!
\left[\!\!\!\!
\begin{array}{c}
\mu_1^1
\end{array}
\!\!\!\!\right]^{\!\!*}
\cdots
\left[\!\!\!\!
\begin{array}{c}
\nu_1^1
\end{array}
\!\!\!\!\right]^{\!\!*}\!\!}
\\ 
\rule[-2.17pt]{0pt}{77pt}
\end{array}.
\endaligned
\]
Here, for reasons of space, the multiplicities $*$ are not written in
length, but as above, they should be read for a typical column as an
integer $a_{ \lambda_1^j, \lambda_2^j, \dots, \lambda_i^j }$ depending
on the column which is $\geqslant 1$; so we understand that the
multiplicities of appearing columns are always positive, but it may
well happen that some blocks of given depths are completely
missing\footnote{\,
However, in what we will be interested in later
for applications to the Green-Griffiths bundle of jets, we will have
to study only Young diagrams ${\sf YD}_{ ( \ell_1, \dots, \ell_n)}$
for which $\ell_1 - \ell_2$, \dots, $\ell_{ n-1} - \ell_n$ and
$\ell_n$ are all positive, and even in fact large, so letting
all blocks appear in diagrams is harmless.
}, 
so that at some
places, there are contacts between a block of depth, say $i+c$ on the
left for some $c \geqslant 2$, and a block of depth $i$ on the
right. Furthermore, inside any block, semi-standard inequalities must
hold, and between the two contacting columns of two neighboring
blocks, say of depth $i+1$ and of depth $i$:
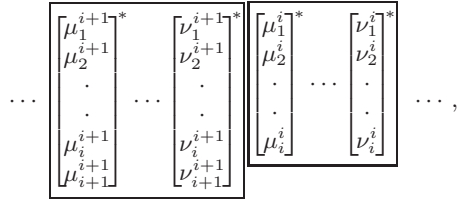
\[
\footnotesize
\aligned
\cdots\,\,
\boxed{\!\!
\left[\!\!\!\!
\begin{array}{c}
\mu_1^{i+1}
\\
\mu_2^{i+1}
\\
\cdot
\\
\cdot
\\
\mu_i^{i+1}
\\
\mu_{i+1}^{i+1}
\end{array}
\!\!\!\!\right]^{\!\!*}
\cdots
\left[\!\!\!\!
\begin{array}{c}
\nu_1^{i+1}
\\
\nu_2^{i+1}
\\
\cdot
\\
\cdot
\\
\nu_i^{i+1}
\\
\nu_{i+1}^{i+1}
\end{array}
\!\!\!\!\right]^{\!\!*}\!\!}
\!\!\!\!\:\:\!
\begin{array}{c}
\boxed{\!\!
\left[\!\!\!\!
\begin{array}{c}
\mu_1^i
\\
\mu_2^i
\\
\cdot
\\
\cdot
\\
\mu_i^i
\end{array}
\!\!\!\!\right]^{\!\!*}
\cdots
\left[\!\!\!\!
\begin{array}{c}
\nu_1^i
\\
\nu_2^i
\\
\cdot
\\
\cdot
\\
\nu_i^i
\end{array}
\!\!\!\!\right]^{\!\!*}\!\!}
\\ 
\rule[-2.17pt]{0pt}{12pt}
\end{array}
\cdots,
\endaligned
\]
there must of course in addition exist the semi-standard-like
truncated inequalities: 
\begin{equation}
\label{contact-inequalities}
\begin{array}{ccc}
\nu_1^{i+1} & \leqslant & \mu_1^i
\\
\nu_2^{i+1} & \leqslant & \mu_2^i
\\
\cdot & \leqslant & \cdot
\\
\nu_i^{i+1} & \leqslant & \mu_i^i
\\
\nu_{i+1}^{i+1},
\end{array}
\end{equation}
with nothing about the last element of the longest column; if more
generally, the contact holds between a nonvoid block of depth $i + c$
on the left and a nonvoid block of depth $i$ on the right, in the case
where blocks of the intermediate depths $i+c-1, \dots, i+1$ are
inextant, then the last $c$ elements $\nu_{ i+1}^{i+c}, \dots, \nu_{
i+c}^{ i+c}$ of the rightmost column of the longest block located on
the left are subjected to no constraint at all.

A notable example of such a semi-standard Young tableau representing a
$\Delta$-monomial written in such a way more informative than before
is the following:
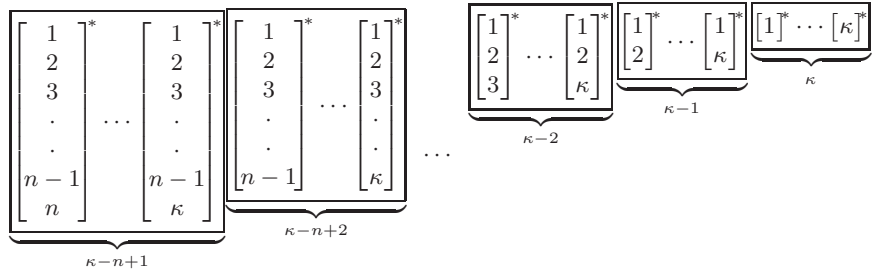
\[
\footnotesize
\aligned
\begin{array}{c}
\underbrace{\boxed{\!\!
\left[\!\!\!\!
\begin{array}{c}
1
\\
2
\\
3
\\
\cdot
\\
\cdot
\\
n-1
\\
n
\end{array}
\!\!\!\!\right]^{\!\!*}
\cdots
\left[\!\!\!\!
\begin{array}{c}
1
\\
2
\\
3
\\
\cdot
\\
\cdot
\\
n-1
\\
\kappa
\end{array}
\!\!\!\!\right]^{\!\!*}\!\!}}_{\kappa-n+1}
\\ 
\rule[-2.17pt]{0pt}{2pt}
\end{array}
\!\!\!\!\!\!
\begin{array}{c}
\underbrace{\boxed{\!\!
\left[\!\!\!\!
\begin{array}{c}
1
\\
2
\\
3
\\
\cdot
\\
\cdot
\\
n-1
\end{array}
\!\!\!\!\right]^{\!\!*}
\cdots
\left[\!\!\!
\begin{array}{c}
1
\\
2
\\
3
\\
\cdot
\\
\cdot
\\
\kappa
\end{array}
\!\!\!\right]^{\!\!*}\!\!}}_{\kappa-n+2}
\\ 
\rule[-2.17pt]{0pt}{22.5pt}
\end{array}
\cdots
\begin{array}{c}
\underbrace{\boxed{\!\!
\left[\!\!\!
\begin{array}{c}
1
\\
2
\\
3
\end{array}
\!\!\!\right]^{\!\!*}
\cdots
\left[\!\!\!
\begin{array}{c}
1
\\
2
\\
\kappa
\end{array}
\!\!\!\right]^{\!\!*}\!\!}}_{\kappa-2}
\\ 
\rule[-2.17pt]{0pt}{59pt}
\end{array}
\!\!\!\!\:\:\!\!\!\!
\begin{array}{c}
\underbrace{\boxed{\!\!
\left[\!\!\!
\begin{array}{c}
1
\\
2
\end{array}
\!\!\!\right]^{\!\!*}
\cdots
\left[\!\!\!
\begin{array}{c}
1
\\
\kappa
\end{array}
\!\!\!\right]^{\!\!*}\!\!}}_{\kappa-1}
\\ 
\rule[-2.17pt]{0pt}{71pt}
\end{array}
\!\!\!\!\:\:\!\!\!\!
\begin{array}{c}
\underbrace{\boxed{\!\!
\left[\!\!\!
\begin{array}{c}
1
\end{array}
\!\!\!\right]^{\!\!*}
\cdots
\left[\!\!\!
\begin{array}{c}
\kappa
\end{array}
\!\!\!\right]^{\!\!*}\!\!}}_{\kappa}
\\ 
\rule[-2.17pt]{0pt}{83pt}
\end{array}.
\endaligned
\vspace{-0.5cm}
\]
It has depth $d_1 = n$ and its first column on the left 
corresponds naturally to
the $n$-dimensional Wronskian:
\[
\Delta_{1,2,3,\dots,n-1,n}^{1,2,3,\dots,n-1,n}
=
\left\vert
\begin{array}{cccccc}
f_1' & f_2' & f_3 ' & \cdots & f_{n-1}' & f_n'
\\
f_1'' & f_2'' & f_3 '' & \cdots & f_{n-1}'' & f_n''
\\
f_1''' & f_2''' & f_3 ''' & \cdots & f_{n-1}''' & f_n'''
\\
\cdot & \cdot & \cdot & \cdots & \cdot & \cdot
\\
f_1^{(n-1)} & f_2^{(n-1)} & f_3^{(n-1)} & \cdots & f_{n-1}^{(n-1)} & 
f_n^{(n-1)}
\\
f_1^{(n)} & f_2^{(n)} & f_3^{(n)} & \cdots & f_{n-1}^{(n)} & 
f_n^{(n)}
\end{array}
\right\vert
\]
raised to a certain power $* = a_{ 1, 2, 3, \dots, n-1, n} \geqslant
1$; in its first block, the bottom indices of extant columns are $n,
n+1, \dots, \kappa -1, \kappa$ while all other indices above 
are constant horizontally;
in its second block, the bottom indices of extant columns are $n-1, n,
n+1, \dots, \kappa+1, \kappa$; and so on. Therefore, the number of
pairwise distinct columns is equal to\footnote{\,
Notice {\em passim}
that this number minus $1$ plus the (constant) degree of any
homogeneous monomial $\big( \ell_1 - \ell_2 \big)^{ \alpha_1} + \cdots
+ \big( \ell_{ n-1} - \ell_n \big)^{ \alpha_{ n-1}} \big( \ell_n
\big)^{ \alpha_n}$:
\[
n\kappa-{\textstyle{\frac{n(n-1)}{2}}}-1
+
{\textstyle{\frac{n(n+1)}{2}}}
=
(\kappa+1)n-1
\]
equals the exponent of $m$ in the formula of 
Theorem~\ref{GG-characteristic} about the
asymptotic behavior of the Euler-Poincar\'e characteristic. 
}: 
\[
(\kappa-n+1)+(\kappa-n+2)+\cdots+(\kappa-2)+(\kappa-1)+\kappa
=
n\kappa
-
{\textstyle{\frac{n(n-1)}{2}}}. 
\]
Observe that this example of semi-standard Young
tableau constitutes not
just one $\Delta$-monomial, but in fact it represents a whole specific
{\em family} of $\Delta$-monomials which are labelled by all possible
column exponents $* \geqslant 1$. Of course, the row lengths, or
equivalently and what is clearer, the differences of row lengths, may, 
as already seen, be expressed in terms of these exponents:
\[
\left\{
\aligned
\ell_n
&
=
a_{1,2,3,\dots,n-1,n}
+\cdots+
a_{1,2,3,\dots,n-1,\kappa}
\\
\ell_{n-1}-\ell_n
&
=
a_{1,2,3,\dots,n-1}
+\cdots+
a_{1,2,3,\dots,\kappa}
\\
\cdot\!\cdot
\ \ \ \ \ \ \,
\cdot\cdot\,
&
=
\cdots\cdots\cdots\cdots\cdots\cdots\cdots\cdot
\\
\ell_3-\ell_4
&
=
a_{1,2,3}
+\cdots+
a_{1,2,\kappa}
\\
\ell_2-\ell_3
&
=
a_{1,2}
+\cdots+
a_{1,\kappa}
\\
\ell_1-\ell_2
&
=
a_1
+\cdots+
a_\kappa.
\endaligned\right.
\]
One should also remember from Theorem~\ref{exact-decomposition} 
that the total number of primes in
any considered $\Delta$-monomial should be constant equal to $m$, a
condition that can now be read here as:
\[
\scriptsize
\aligned
m
&
=
\big[1+2+3+\cdots+n-1+n\big]a_{1,2,3,\dots,n-1,n}
+\cdots+
\big[1+2+3+\cdots+n-1+\kappa\big]a_{1,2,3,\dots,n-1,\kappa}
+
\\
&
+
\big[1+2+3+\cdots+n-1\big]a_{1,2,3,\dots,n-1}
+\cdots+
\big[1+2+3+\cdots+\kappa\big]a_{1,2,3,\dots,\kappa}
+
\\
&
+
\cdots\cdots\cdots\cdots\cdots\cdots\cdots\cdots\cdots\cdots\cdots
\cdots\cdots\cdots\cdots\cdots\cdots
+
\\
&
+
\big[1+2+3\big]a_{1,2,3}
+\cdots+
\big[1+2+\kappa\big]a_{1,2,\kappa}
+
\\
&
+
\big[1+2\big]a_{1,2}
+\cdots+
\big[1+\kappa\big]a_{1,\kappa}
+
\\
&
+
\big[1\big]a_1
+\cdots+
\big[\kappa\big]a_\kappa; 
\endaligned
\]
in this equation, each exponent $a_{1, 2, 3, \dots, i-1, \lambda}$, $i
\leqslant \lambda \leqslant \kappa$, is multiplied by the sum $1 + 2 +
3 + \cdots + i-1 + \lambda$ of its lower indices, an integer which
equals the total number of primes in the determinant $\Delta_{ 1, 2,
3, \dots, i-1, i}^{ 1, 2, 3, \dots, i-1, \lambda}$, so that the total
number of primes in the power $\big[ \Delta_{ 1, 2, 3, \dots, i-1,
i}^{ 1, 2, 3, \dots, i-1, \lambda} \big]^{ a_{ 1, 2, 3, \dots, i-1,
\lambda }}$ is indeed equal to the product:
\[
\big[1+2+3+\cdots+i-1+\lambda\big]
a_{1,2,3,\dots,i-1,\lambda}. 
\]

\subsection{Maximal number of pairwise distinct columns in 
${\sf YD}_{ (\ell_1, \dots, \ell_n)} ( \lambda_i^j)$} We now come back
to a general semi-standard Young tableau of depth $d_1$ with 
extreme columns in each block that can be arbitrary:
\begin{equation}
\label{arbitrary-semi-standard}
\footnotesize
\aligned
\begin{array}{c}
\underbrace{
\boxed{\!\!
\left[\!\!\!\!
\begin{array}{c}
\mu_1^{d_1}
\\
\mu_2^{d_1}
\\
\mu_3^{d_1}
\\
\cdot
\\
\cdot
\\
\mu_{d_1-1}^{d_1}
\\
\mu_{d_1}^{d_1}
\end{array}
\!\!\!\!\right]^{\!\!*}
\cdots
\left[\!\!\!\!
\begin{array}{c}
\nu_1^{d_1}
\\
\nu_2^{d_1}
\\
\nu_3^{d_1}
\\
\cdot
\\
\cdot
\\
\nu_{d_1-1}^{d_1}
\\
\nu_{d_1}^{d_1}
\end{array}
\!\!\!\!\right]^{\!\!*}\!\!}}_{
\#\,\leqslant\,D_{d_1}}
\\
\rule[-0.17pt]{0pt}{1pt}
\end{array}
\!\!\!\!\!\!
\begin{array}{c}
\underbrace{\boxed{\!\!
\left[\!\!\!\!
\begin{array}{c}
\mu_1^{d_1-1}
\\
\mu_2^{d_1-1}
\\
\mu_3^{d_1-1}
\\
\cdot
\\
\cdot
\\
\mu_{d_1-1}^{d_1-1}
\end{array}
\!\!\!\!\right]^{\!\!*}
\cdots
\left[\!\!\!\!
\begin{array}{c}
\nu_1^{d_1-1}
\\
\nu_2^{d_1-1}
\\
\nu_3^{d_1-1}
\\
\cdot
\\
\cdot
\\
\nu_{d_1-1}^{d_1-1}
\end{array}
\!\!\!\!\right]^{\!\!*}\!\!}}_{
\#\,\leqslant\,D_{d_1-1}}
\\ 
\rule[-2.17pt]{0pt}{24.5pt}
\end{array}
\cdots
\begin{array}{c}
\underbrace{\boxed{\!\!
\left[\!\!\!\!
\begin{array}{c}
\mu_1^3
\\
\mu_2^3
\\
\mu_3^3
\end{array}
\!\!\!\!\right]^{\!\!*}
\cdots
\left[\!\!\!\!
\begin{array}{c}
\nu_1^3
\\
\nu_2^3
\\
\nu_3^3
\end{array}
\!\!\!\!\right]^{\!\!*}\!\!}}_{
\#\,\leqslant\,D_3}
\\ 
\rule[-2.17pt]{0pt}{66pt}
\end{array}
\!\!\!\!\:\:\!\!\!\!\!\!\:
\begin{array}{c}
\underbrace{\boxed{\!\!
\left[\!\!\!\!
\begin{array}{c}
\mu_1^2
\\
\mu_2^2
\end{array}
\!\!\!\!\right]^{\!\!*}
\cdots
\left[\!\!\!\!
\begin{array}{c}
\nu_1^2
\\
\nu_2^2
\end{array}
\!\!\!\!\right]^{\!\!*}\!\!}}_{
\#\,\leqslant\,D_2}
\\ 
\rule[-2.17pt]{0pt}{78pt}
\end{array}
\!\!\!\!\:\:\!\!\!\!\!\!\:
\begin{array}{c}
\underbrace{\boxed{\!\!
\left[\!\!\!\!
\begin{array}{c}
\mu_1^1
\end{array}
\!\!\!\!\right]^{\!\!*}
\cdots
\left[\!\!\!\!
\begin{array}{c}
\nu_1^1
\end{array}
\!\!\!\!\right]^{\!\!*}\!\!}}_{
\#\,\leqslant\,D_1}
\\ 
\rule[-2.17pt]{0pt}{90pt}
\end{array}.
\endaligned
\vspace{-0.5cm}
\end{equation}
What is the maximal possible number of pairwise distinct $*$-ed 
columns? In the rightmost block, the number of entries in the single row
is clearly less than or equal to:
\[
D_1
:=
1+\nu_1^1-\mu_1^1.
\]
In full generality, there may well be several gaps from $\mu_1^1$ to
$\nu_1^1$ in the `cdots', so $D_1$ is an upper bound. If a general
$*$-ed column of depth $i$ with $1 \leqslant i \leqslant d_1$ has two
immediate neighbors, with no possible intermediate neighbors:
\[
\left[\!\!
\begin{array}{c}
\lambda_1'
\\
\lambda_2'
\\
\cdot\cdot
\\
\lambda_i'
\end{array}
\!\!\right]^{\!*}
<
\left[\!\!
\begin{array}{c}
\lambda_1
\\
\lambda_2
\\
\cdot\cdot
\\
\lambda_i
\end{array}
\!\!\right]^{\!*}
<
\left[\!\!
\begin{array}{c}
\lambda_1''
\\
\lambda_2''
\\
\cdot\cdot
\\
\lambda_i''
\end{array}
\!\!\right]^{\!*},
\]
then necessarily the two sums of horizontal differences: 
\[
\aligned
\lambda_1-\lambda_1'
+
\lambda_2-\lambda_2'
+\cdots+
\lambda_i-\lambda_i'
&
=
1
\\
\lambda_1''-\lambda_1
+
\lambda_2''-\lambda_2
+\cdots+
\lambda_i''-\lambda_i
&
=
1
\endaligned
\]
are both smallest possible, equal to $1$. It follows easily from this
observation and from the constraint of semi-standarcy that in the
penultimate block, in the antepenultimate block and in the general
$i$-th block, the numbers of pairwise distinct $*$-ed columns are at
most equal to, respectively:
\[
\aligned
D_2
&
:=
1+(\nu_1^2-\mu_1^2)+(\nu_2^2-\mu_2^2),
\\
D_3
&
:=
1+(\nu_1^3-\mu_1^3)+(\nu_2^3-\mu_2^3)+(\nu_3^3-\mu_3^3), 
\endaligned
\]
and to: 
\begin{equation}
\label{D-i}
D_i 
:= 
1+\sum_{l=1}^i\,(\nu_l^i-\mu_l^i). 
\end{equation}
Consequently, the total number of pairwise distinct columns in an
arbitrary semi-standard Young tableau is at most equal to $D_1 + D_2 +
D_3 + \cdots + D_{ d_1-1 } + D_{ d_1 }$, that is to say to:
\[
\scriptsize
\aligned
\begin{array}{ccccccccccccc}
& 1 & + & 1 & + & \cdots & + & 1 & + & 1 & + & 1 & +
\\
+ & \nu_1^{d_1}-\mu_1^{d_1} & + & \nu_1^{d_1-1}-\mu_1^{d_1-1} &
+ & \cdots & + & \nu_1^3-\mu_1^3 & + & \nu_1^2-\mu_1^2 & + &
\nu_1^1-\mu_1^1 & +
\\
+ & \nu_2^{d_1}-\mu_2^{d_1} & + & \nu_2^{d_1-1}-\mu_2^{d_1-1} &
+ & \cdots & + & \nu_2^3-\mu_2^3 & + & \nu_2^2-\mu_2^2 & +
\\
+ & \nu_3^{d_1}-\mu_3^{d_1} & + & \nu_3^{d_1-1}-\mu_3^{d_1-1} &
+ & \cdots & + & \nu_3^3-\mu_3^3 & + & 
\\
+ & \cdots\cdots\cdots & + & \cdots\cdots\cdots & + & \cdots & +
\\
+ & \!\!\!\nu_{d_1-1}^{d_1}-\mu_{d_1-1}^{d_1}\!\!\! & + &
\nu_{d_1-1}^{d_1-1}-\mu_{d_1-1}^{d_1-1} & +
\\
+ & \nu_{d_1}^{d_1}-\mu_{d_1}^{d_1}.
\end{array}
\endaligned
\]
But by permuting the order of appearance of $\nu$ and $\mu$ in each
subtraction of every line, this sum becomes:
\[
\footnotesize
\aligned
\begin{array}{cccccccccccccc}
 & & 1 & + & \cdots & + & 1 & + & 1 & + & 1 & -
\\
-\,\mu_1^{d_1} & + & \underline{\nu_1^{d_1}-\mu_1^{d_1-1}}_{\!<} & 
+ & \cdots & + & \underline{\nu_1^3-\mu_1^2}_{\!<} & + & 
\underline{\nu_1^2-\mu_1^1}_{\!<} & + & \nu_1^1
\\
-\,\mu_2^{d_1} & + & \underline{\nu_2^{d_1}-\mu_2^{d_1-1}}_{\!<} & 
+ & \cdots & + & \underline{\nu_2^3-\mu_2^2}_{\!<} & + & 
\nu_2^2 &
\\
\cdots & + & \cdots\cdots\cdots & + & \cdots
\\
-\,\mu_{d_1-1}^{d_1} & + & 
\underline{\nu_{d_1-1}^{d_1}-\mu_{d_1-1}^{d_1-1}}_{\!<} & + 
\\
-\,\mu_{d_1}^{d_1} & + & \nu_{d_1}^{d_1},
\end{array}
\endaligned
\]
and then by taking account just of the semi-standard-like
inequalities~\thetag{ \ref{contact-inequalities}}
(about the columns
of contact between two neighboring blocks), we see that all the pairs
that we have underlined above are $\leqslant 0$, whence:
\[
\footnotesize
\aligned
D_1+D_2+D_3+\cdots+D_{d_1-1}+D_{d_1}
&
\leqslant
d_1\cdot 1
+
(-\mu_1^{d_1}+\nu_1^1)
+
(-\mu_2^{d_1}+\nu_2^2)
+\cdots+
\\
&
\ \ \ \ \ \ \ \ \
+
(-\mu_{d_1-1}^{d_1}+\nu_{d_1-1}^{d_1-1})
+
(-\mu_{d_1}^{d_1}+\nu_{d_1}^{d_1}). 
\endaligned
\]
Finally, the strict inequalities $0 < \mu_1^{ d_1} < \mu_2^{ d_1} <
\cdots < \mu_{ d_1-1}^{ d_1} < \mu_{ d_1}^{ d_1}$ between the entries
of the first column yield trivial majorants:
\[
-\,\mu_1^{d_1}\leqslant -1,
\ \ \ \ \ 
-\,\mu_2^{d_1}\leqslant -2,
\dots\dots,
-\,\mu_{d_1-1}^{d-1}\leqslant -(d_1-1),
\ \ \ \ \ 
-\,\mu_{d_1}\leqslant -d_1,
\]
and since all the $\nu_i^j$ are $\leqslant \kappa$ anyway, we deduce
about any semi-standard Young tableau the following
majoration:
\[
\aligned
\text{\small\sf total number of distinct $*$-ed columns}
&
\leqslant
d_1+(-1+\kappa)+(-2+\kappa)+\cdots+
\\
&\ \ \ \ \ \ \ \ \
+(-(d_1-1)+\kappa)+(-d_1+\kappa)
\\
&
\leqslant
n+(-1+\kappa)+(-2+\kappa)+\cdots+
\\
&\ \ \ \ \ \ \ \ \
+(-(n-1)+\kappa)+(-n+\kappa)
\\
&
=
n\kappa
-
{\textstyle{\frac{n(n-1)}{2}}}. 
\endaligned
\]

\begin{Lemma}
\label{bound-distinct-columns}
The total number of pairwise distinct columns in a
semi-standard Young tableau of depth $\leqslant n$ filled in with
integers $\lambda_i^j \leqslant \kappa$ is in any case $\leqslant
n\kappa - \frac{ n ( n-1)}{ 2}$. 
\qed
\end{Lemma}

We now introduce several families of $\Delta$-monomials parametrized
by a fixed collection of pairs of columns (having various
multiplicities $* \geqslant 1$) that should occupy the left and right
extreme positions in blocks of decreasing depths.

\subsection{Main definition} 
Let $d_1$ be a positive integer $\leqslant n$ and let $\mu_l^i$ and
$\nu_l^i$ be integers defined for $i = 1, 2, 3, \dots, d_1-1, d_1$ and
$1 \leqslant l \leqslant i$ with $\mu_l^i \leqslant \kappa$ and
$\nu_l^i \leqslant \kappa$ which satisfy all the following
inequalities:
\label{main-definition}

\begin{itemize}

\smallskip\item[$\bullet$]
{\small\sf vertical downward increasing:} 
\[
0<\mu_1^i
<\cdots<
\mu_i^i 
\ \ \ \ \
\text{\rm and}
\ \ \ \ \ 
0<\nu_1^i<\cdots<\nu_i^i
\ \ \ \ \ \ \ \ \ \ \ \ \
{\scriptstyle{(i\,=\,1,\dots,\,d_1)}};
\]

\smallskip\item[$\bullet$]
{\small\sf weak increasing inside blocks:}
\[
\mu_l^i\leqslant\nu_l^i
\ \ \ \ \ \ \ \ \ \ \ \ \
{\scriptstyle{(i\,=\,1,\dots,\,d_1;\,\,
1\,\leqslant\,l\,\leqslant\,i)}};
\]

\smallskip\item[$\bullet$]
{\small\sf weak increasing for the contacts between neighboring
blocks\footnote{\,
The most general case where certain block lengths
$i$ are missing, so that block lengths sometimes jump for more than
one unit, is implicitly also embraced by such a definition, for it
suffices, about the indices $i$ of blocks that are thought to be
missing, to just prescribe somewhat arbitrarily some integers
$\mu_l^i$ and $\nu_l^i$ that violate the second condition ; the third
condition is then supposed to hold for the direct contacts between the
really extant neighboring blocks. Since in our later principal
considerations, we will not be dealing with semi-standard tableaux
having block gaps, it is not necessary to introduce further specific
notations here. 
}:} 
\[
\nu_l^{i+1}\leqslant\mu_l^i
\ \ \ \ \ \ \ \ \ \ \ \ \
{\scriptstyle{(i\,=\,1,\dots,\,d_1-1;\,\,
1\,\leqslant\,l\,\leqslant\,i)}}
\]

\end{itemize}\smallskip

\noindent
Then with such data, the family of semi-standard tableaux: 
\[
\boxed{{\sf YT}_{\!\kappa,m}\big(\mu_l^i, \nu_l^i\big)}
\]
is defined to consist of all possible concatenations:
\[
{\sf block}^{d_1}\big(\mu^{d_1},\nu^{d_1}\big)
\cdots\,
{\sf block}^i\big(\mu^i,\nu^i\big)
\cdots\,
{\sf block}^1\big(\mu^1,\nu^1\big)
\]
of {\em semi-standard} blocks\footnote{\,
By a {\sl semi-standard block} is of course meant a block in which
strict increase holds downward along columns, while weak increase
holds from left to right along rows.
} 
of the form:
\[
\footnotesize
\aligned
{\sf block}^i\big(\mu^i,\nu^i\big)
=
\boxed{\!
\left[\!\!\!
\begin{array}{c}
\mu_1^i
\\
\mu_2^i
\\
\cdot
\\
\cdot
\\
\mu_i^i
\end{array}
\!\!\!\right]^{\!a_{\mu_1^i,\dots,\mu_i^i}}
\cdots
\left[\!\!\!
\begin{array}{c}
\lambda_1^i
\\
\lambda_2^i
\\
\cdot
\\
\cdot
\\
\lambda_i^i
\end{array}
\!\!\!\right]^{\!a_{\lambda_1^i,\dots,\lambda_i^i}}
\cdots
\left[\!\!\!
\begin{array}{c}
\nu_1^i
\\
\nu_2^i
\\
\cdot
\\
\cdot
\\
\nu_i^i
\end{array}
\!\!\!\right]^{\!a_{\nu_1^i,\dots,\nu_i^i}}\!}\,,
\endaligned
\]
all $*$-ed columns being pairwise distinct and ordered increasingly from
left to right, with multiplicities:
\[
a_{\mu_1^i,\dots,\mu_i^i}
\geqslant 1,
\ \ \ \ \
\dots
\ \ \ \ \
a_{\lambda_1^i,\dots,\lambda_i^i}
\geqslant 1,
\ \ \ \ \
\dots
\ \ \ \ \
a_{\nu_1^i,\dots,\nu_i^i}
\geqslant 1
\]
which are {\em all} positive and with the further important
constraint that:
\[
\aligned
m
&
=
{\rm weight}\big({\sf block}^{d_1}(\mu^{d_1},\nu^{d_1})\big)
+\cdots+
\\
&
\ \ \ \ \
+
{\rm weight}\big({\sf block}^i(\mu^i,\nu^i)\big)
+\cdots+
{\rm weight}\big({\sf block}^1(\mu^1,\nu^1)\big),
\endaligned
\]
where according to an expectable, natural definition: 
\[
\footnotesize
\aligned
{\rm weight}\big({\sf block}^i(\mu^i,\nu^i)\big)
\overset{\text{\rm def}}{=}
&\,
\big(\mu_1^i+\cdots+\mu_i^i\big)
a_{\mu_1^i,\dots,\mu_i^i}
+\cdots+
\\
&
+
\big(\lambda_1^i+\cdots+\lambda_i^i\big)
a_{\lambda_1^i,\dots,\lambda_i^i}
+\cdots+
\big(\nu_1^i+\cdots+\nu_i^i\big)
a_{\nu_1^i,\dots,\nu_i^i}
\endaligned
\]
simply denotes the total number of primes (remember
Theorem~\ref{Theorem-A}) in the associated 
$\Delta$-monomial:
\[
\big(
\Delta_{1,\dots,i}^{\mu_1^i,\dots,\mu_i^i}
\big)^{a_{\mu_1^i,\dots,\mu_i^i}}
\cdots
\big(
\Delta_{1,\dots,i}^{\lambda_1^i,\dots,\lambda_i^i}
\big)^{a_{\lambda_1^i,\dots,\lambda_i^i}}
\cdots
\big(
\Delta_{1,\dots,i}^{\nu_1^i,\dots,\nu_i^i}
\big)^{a_{\nu_1^i,\dots,\nu_i^i}}.
\]
In a specific family ${\sf YT}_{\!\! \kappa, m} \big( \mu_l^i,
\nu_l^i\big)$, the freedom of variation lies: 1) in the choice of some
intermediate columns; 2) in the choice of the number of such
intermediate columns; 3) in the choice of the positive multiplicities
of all the columns.

\begin{Lemma}
The collection of all semi-standard Young tableaux ${\sf YT}$ of depth
$\leqslant n$ filled in with positive integers $\lambda_i^j \leqslant
\kappa$ whose weight equals $m$ is identical to the disjoint union:
\[
\bigcup_{\mu_l^i,\nu_l^i}
{\sf YT}_{\!\kappa,m}
\big(\mu_l^i,\nu_l^i\big)
\]
of the so-defined families of semi-standard tableaux.
\end{Lemma}

\proof
According to the preceding considerations, an arbitrary semi-standard
Young tableau looks like~\thetag{ \ref{arbitrary-semi-standard}},
hence belongs to ${\sf YT}_{\!\! \kappa, m} \big( \mu_l^i,
\nu_l^i \big)$ for some $\mu_l^i, \nu_l^i$. Disjointness follows from
the fact that the extreme column data $\big(\mu_l^i, \nu_l^i \big)$
are obviously pairwise distinct.
\endproof

By what has already been seen, the number of pairwise distinct columns
in any ${\sf block}^i (\mu^i, \nu^i)$ may well be equal to
zero\footnote{\,
This would correspond to the empty block case, {\em
cf.} a preceding footnote.
} 
and is always smaller than or equal to:
\[
D_i
:= 
1+
{\textstyle{\sum_{l=1}^i}}\,(\nu_l^i-\mu_l^i). 
\]
In order to fix ideas about the exact number of distinct columns, we
shall in addition split each (big) family ${\sf YT}_{ \! \kappa, m} (
\mu_l^i, \nu_l^i)$ in distinct, finer (sub)families as follows.

For every $i = 1, \dots, n$ and for every nonnegative integer:
\[
\tau_i\leqslant
D_i-1
=
{\textstyle{\sum_{l=1}^i}}\,
(\nu_l^i-\mu_l^i)
\]
less than the maximal possible number of distinct columns inside ${\sf
block}^i (\mu^i, \nu^i)$ {\em minus $1$}, let us choose a discrete
increasing path\footnote{\,
When a block of depth $i$ is inextant, we set $\tau^i := -1$ so that
the length $1 + \tau^i$ of any associated path $\gamma^i$ equals $0$:
possible paths $\gamma^i$ are thus necessarily empty in this case.
}: 
\[
\gamma^i\colon\,\,
\big\{0,1,2,\dots,\tau^i\big\}
\longrightarrow
\text{\small\sf downward increasing columns}
\in
\{1,\dots,\kappa\}^i
\]
from the $\mu^i$-column to the $\nu^i$-column, namely:
\[
\left[\!\!
\begin{array}{c}
\mu_1^i=\gamma_1^i(0)
\\
\mu_2^i=\gamma_2^i(0)
\\
\cdot\cdot
\\
\mu_i^i=\gamma_i^i(0)
\end{array}
\!\!\right]^{\!*}
<
\left[\!\!
\begin{array}{c}
\gamma_1^i(1)
\\
\gamma_2^i(1)
\\
\cdot
\\
\gamma_i^i(1)
\end{array}
\!\!\right]^{\!*}
<\cdots<
\left[\!\!
\begin{array}{c}
\gamma_1^i(s^i)
\\
\gamma_2^i(s^i)
\\
\cdot
\\
\gamma_i^i(s^i)
\end{array}
\!\!\right]^{\!*}
<\cdots<
\left[\!\!
\begin{array}{c}
\gamma_1^i(\tau^i)=\nu_1^i
\\
\gamma_2^i(\tau^i)=\nu_2^i
\\
\cdot\cdot
\\
\gamma_i^i(\tau^i)=\nu_i^i
\end{array}
\!\!\right]^{\!*},
\]
with $s^i = 0, 1, 2, \dots, \tau^i$ denoting the current (discrete)
time variable, such that the associated block:
\[
\footnotesize
\aligned
{\sf block}^i\big(\gamma^i\big)
:=
\boxed{\!
\left[\!\!\!
\begin{array}{c}
\gamma_1^i(0)
\\
\gamma_2^i(0)
\\
\cdot
\\
\gamma_i^i(0)
\end{array}
\!\!\!\right]^{\!*}
\left[\!\!\!
\begin{array}{c}
\gamma_1^i(1)
\\
\gamma_2^i(1)
\\
\cdot
\\
\gamma_i^i(1)
\end{array}
\!\!\!\right]^{\!*}
\cdots
\left[\!\!\!
\begin{array}{c}
\gamma_1^i(s^i)
\\
\gamma_2^i(s^i)
\\
\cdot
\\
\gamma_i^i(s^i)
\end{array}
\!\!\!\right]^{\!*}
\cdots
\left[\!\!\!
\begin{array}{c}
\gamma_1^i(\tau^i)
\\
\gamma_2^i(\tau^i)
\\
\cdot
\\
\gamma_i^i(\tau^i)
\end{array}
\!\!\!\right]^{\!*}\!}
\endaligned
\]
is semi-standard. Then with such data, the (sub)family of
semi-standard tableaux:
\[
\boxed{
{\sf YT}_{\!\kappa,m}
\big(\mu_l^i,\nu_l^i,\tau^i,\gamma^i(s^i)\big)}
\]
is defined to consist of all possible concatenations: 
\[
{\sf block}^n\big(\gamma^n\big)
\cdots
{\sf block}^i\big(\gamma^i\big)
\cdots
{\sf block}^1\big(\gamma^1\big)
\]
of the above specific blocks, with $*$-multiplicities:
\[
a_{\gamma_1^i(0),\dots,\gamma_i^i(0)}
\geqslant 1,
\ \ \ \ \
\dots
\ \ \ \ \
a_{\gamma_1^i(s^i),\dots,\gamma_i^i(s^i)}
\geqslant 1,
\ \ \ \ \
\dots
\ \ \ \ \
a_{\gamma_1^i(\tau^i),\dots,\gamma_i^i(\tau^i)}
\geqslant 1
\]
which are {\em all} positive, and with the further constraint, 
similar as before, that:
\[
m
=
{\rm weight}\big({\sf block}^n(\gamma^n)\big)
+\cdots+
{\rm weight}\big({\sf block}^i(\gamma^i)\big)
+\cdots+
{\rm weight}\big({\sf block}^1(\gamma^1)\big).
\]
Here of course, the weight of a general single column, namely having
with multiplicity $1$, equals:
\[
\gamma_1^i(s^i)+\gamma_2^i(s^i)+\cdots+\gamma_i^i(s^i),
\]
hence the $*$-ed column has weight:
\[
\big[
\gamma_1^i(s^i)+\gamma_2^i(s^i)+\cdots+\gamma_i^i(s^i)
\big]\,
a_{\gamma_1^i(s^i),\dots,\gamma_i^i(s^i)}.
\]
From now on, we shall denote the multiplicity of a general $*$-ed
column shortly by $a_{ s^i}^i$, instead of $a_{ \gamma_1^i ( s^i),
\dots, \gamma_i^i ( s^i)}$. The weight homogeneity condition therefore
reads:
\begin{equation}
\label{m-homogeneity-gamma}
\aligned
m
&
=
\sum_{i=1}^n\,
\sum_{0\leqslant s^i\leqslant\tau^i}
\big[\gamma_1^i(s^i)+\cdots+\gamma_i^i(s^i)\big]
a_{s^i}^i. 
\endaligned
\end{equation}

In a specific family ${\sf YT }_{\! \kappa, m} \big(
\mu_l^i, \nu_l^i, \tau^i, \gamma^i ( s^i) \big)$, the freedom of
variation now lies only in the multiplicities, since all the pairwise
distinct $*$-ed columns are fully prescribed in it. Notice that as $m$
is supposed to be quite large\footnote{\,
We will eventually let $m \to \infty$, similarly as in the
Euler-Poincar\'e characteristic of $\mathcal{ E}_{ \kappa, m}^{ GG}
T_X^*$.
} 
in comparison to $n$ and $\kappa$, then for any choice of pairwise
distinct column data $\big( \mu_l^i, \nu_l^i, \tau^i, \gamma^i ( s^i)
\big)$, the column weights $\big[\gamma_1^i(s^i) + \cdots +
\gamma_i^i(s^i) \big]$ being fixed and finite, there is still much
freedom for the multiplicities to fulfill the homogeneity condition in
question. We will in fact realize in a while that the number of
semi-standard Young tableaux of weight $m$ in any family ${\sf YT}_{\!
\kappa, m} \big( \mu_l^i, \nu_l^i, \tau^i, \gamma^i ( s^i) \big)$ is
an ${\sf O}_{ n, \kappa} \big( m^{ D - 1} \big)$, where $D = \sum_{
i=1}^n \, (1+ \tau^i)$ as before is the total number of pairwise
distinct columns.

By construction, it is clear that the union of the (sub)families
${\sf YT}_{\! \kappa,m} \big(\mu_l^i, \nu_l^i, \tau^i,
\gamma^i(s^i) \big)$ fills the previously introduced larger family:
\[
{\sf YT}_{\!\kappa,m}
\big(\mu_l^i,\nu_l^i\big)
=
\bigcup_{\tau^i,\gamma^i(s^i)}\,
{\sf YT}_{\!\kappa,m}
\big(\mu_l^i,\nu_l^i,\tau^i,\gamma^i(s^i)\big). 
\]

\begin{Lemma}
\label{lemma-disjoint-YT}
The collection of all semi-standard Young tableaux ${\sf YT}$ of depth
$\leqslant n$ filled in with positive integers $\lambda_i^j \leqslant
\kappa$ whose weight equals $m$ is identical to the {\em disjoint} 
union:
\[
\bigcup_{\mu_l^i,\nu_l^i,\tau^i,\gamma^i(s^i)}\,
{\sf YT}_{\!\kappa,m}
\big(\mu_l^i,\nu_l^i,\tau^i,\gamma^i(s^i)\big)
\]
of the so-defined families of semi-standard tableaux. Furthermore, the
number of possible such families ${\sf YT}_{\! \kappa,m}
\big( \mu_l^i, \nu_l^i, \tau^i, \gamma^i( s^i)\big)$ is smaller than
or equal to the (nonoptimal) constant:
\[
\prod_{i=1}^n\,
\big(
1+
{\textstyle{\frac{\kappa!}{(\kappa-i)!\,i!}}}
\big)^{1+i(\kappa-i)}
=
{\sf Constant}_{n,\kappa},
\]
independently of $m$.
\end{Lemma}

\proof
Disjointness (not yet argued)
of subfamilies inside a family ${\sf YT}_{\!
\kappa, m} \big( \mu_l^i, \nu_l^i\big)$ comes from the fact that any
collection of paths $(\gamma^1, \gamma^2, \dots, \gamma^n)$ prescribes
all the mutually distinct $*$-ed columns that are extant, their
multiplicities being all $\geqslant 1$.

In a block of depth $i$, a single $*$-ed column is either empty, or it
consists of $i$ numbers $\lambda_1, \dots, \lambda_i$ chosen without
repetition in $\{ 1, 2, \dots, \kappa\}$ and ordered afterward
increasingly. So the number of possible such columns (including the
empty one) is equal to $1 + \frac{ \kappa!}{ ( \kappa - i)!\,
i!}$. Since all $*$-ed columns are pairwise distinct, the maximal
number of $*$-ed columns that one may put in a semi-standard block of
depth $i$ will be attained for the blocks having the following two
extreme columns, which are the farthest ones for the ordering between
columns of depth $i$:
\[
\footnotesize
\aligned
\boxed{\!\!
\left[\!\!\!\!
\begin{array}{c}
1
\\
2
\\
\cdot
\\
i-1
\\
i
\end{array}
\!\!\!\!\right]^{\!\!*}
\cdots
\left[\!\!\!\!
\begin{array}{c}
\kappa-i+1
\\
\kappa-i+2
\\
\cdot
\\
\kappa-1
\\
\kappa
\end{array}
\!\!\!\!\right]^{\!\!*}\!\!}\,.
\endaligned
\]
It follows from~\thetag{ \ref{D-i}}
that one may put at most: 
\begin{equation}
\label{bound-tau-i}
1+(\kappa-i+1-1)+(\kappa-i+2-2)
+\cdots+
\kappa-i
=
1+i(\kappa-i)
\end{equation}
pairwise distinct $*$-ed columns in between so as to constitute a
semi-standard block. Without optimality, we then majorate the number
of possible semi-standard blocks of depth $i$ (including the empty
one) simply by the number $1 + \frac{ \kappa!}{ ( \kappa - i)!\, i!}$
of possible $*$-ed columns raised to a power equal to the maximal
number $1 + i ( \kappa - i)$ of pairwise distinct such $*$-ed
columns. What matters for the sequel is only that the obtained
majorant is independent of $m$.
\endproof

In summary, here is how we constitute our refined view of an arbitrary
semi-standard Young tableau: the data $(\mu_l^i, \nu_l^i)_{ 1
\leqslant l \leqslant i}$, subjected to the natural inequalities of
the Main definition in Subsection~\ref{main-definition}, prescribe the
left and right extreme $*$-ed column in all blocks of depth $i = 1, 2,
\dots, n$ (with possible block gaps); $\tau^i + 1$ is the number of
pairwise distinct $*$-ed columns in the $i$-th block, and these
columns are $\gamma^i ( 0 ), \dots, \gamma^i ( \tau^i )$; all
$*$-multiplicities of these columns are $\geqslant 1$.

\subsection{Asymptotically negligible families of
$\Delta$-monomials} By definition, for each semi-standard Young
tableau ${\sf YT}$ belonging to a fixed family: 
\[
{\sf YT}_{\!\kappa,m} 
\big(\mu_l^i,\nu_l^i,\tau^i,\gamma^i(s^i)\big), 
\]
the number $D$ of pairwise distinct $*$-ed columns is
equal to the sum of the lengths of the paths between two extreme
$*$-ed columns in every block\footnote{\,
By the preceding convention, inextant blocks contribute with $0$ to
this sum, {\em e.g.} all blocks of depths $d_1+1, \dots, n$ when the
depth $d_1$ of the tableau is $< n$.
}: 
\[
D
=
(1+\tau^1)+\cdots+(1+\tau^i)+\cdots+(1+\tau^n).
\]
However, the common horizontal
length of each of the $i$ rows in the block:
\[
\footnotesize
\aligned
{\sf block}^i\big(\gamma^i\big)
=
\boxed{\!
\left[\!\!\!
\begin{array}{c}
\gamma_1^i(0)
\\
\gamma_2^i(0)
\\
\cdot
\\
\gamma_i^i(0)
\end{array}
\!\!\!\right]^{\!a_0^i}
\left[\!\!\!
\begin{array}{c}
\gamma_1^i(1)
\\
\gamma_2^i(1)
\\
\cdot
\\
\gamma_i^i(1)
\end{array}
\!\!\!\right]^{\!a_1^i}
\cdots
\left[\!\!\!
\begin{array}{c}
\gamma_1^i(s^i)
\\
\gamma_2^i(s^i)
\\
\cdot
\\
\gamma_i^i(s^i)
\end{array}
\!\!\!\right]^{\!a_{s^i}^i}
\cdots
\left[\!\!\!
\begin{array}{c}
\gamma_1^i(\tau^i)
\\
\gamma_2^i(\tau^i)
\\
\cdot
\\
\gamma_i^i(\tau^i)
\end{array}
\!\!\!\right]^{\!a_{\tau^i}^i}\!}
\endaligned
\]
depends visibly on the multiplicities, and is equal to:
\[
a_0^i+a_1^i+\cdots+a_{s^i}^i+\cdots+a_{\tau^i}^i
=
\sum_{0\leqslant s^i\leqslant\tau^i}\,
a_{s^i}^i. 
\]
It follows that the lengths $\ell_1, \ell_2, \dots, \ell_{ n-1},
\ell_n$ of the first, second, \dots, $(n-1)$-th and $n$-th 
horizontal lines in
the semi-standard Young tableau ${\sf YT}$ are equal 
to\footnote{\,
Sums $\sum_{ 0 \leqslant s^i \leqslant \tau^i}\, a_{ s^i}^i$ for which
$\tau^i = -1$ (inextant blocks) are naturally thought to be inextant. 
}: 
\[
\small
\aligned
\ell_1
&
=
\sum_{0\leqslant s^n\leqslant\tau^n}\,
a_{s^n}^n
+
\sum_{0\leqslant s^{n-1}\leqslant\tau^{n-1}}\,
a_{s^{n-1}}^{n-1}
+
\cdots
+
\sum_{0\leqslant s^2\leqslant\tau^2}\,
a_{s^2}^2
+
\sum_{0\leqslant s^1\leqslant\tau^1}\,
a_{s^1}^1
\\
\ell_2
&
=
\sum_{0\leqslant s^n\leqslant\tau^n}\,
a_{s^n}^n
+
\sum_{0\leqslant s^{n-1}\leqslant\tau^{n-1}}\,
a_{s^{n-1}}^{n-1}
+
\cdots
+
\sum_{0\leqslant s^2\leqslant\tau^2}\,
a_{s^2}^2
\\
\cdot\cdot\cdot
&
=
\cdots\cdots\cdots\cdots\cdots\cdots\cdots\cdots\cdots\cdots\cdots\cdots
\cdots\cdots
\\
\ell_{n-1}
&
=
\sum_{0\leqslant s^n\leqslant\tau^n}\,
a_{s^n}^n
+
\sum_{0\leqslant s^{n-1}\leqslant\tau^{n-1}}\,
a_{s^{n-1}}^{n-1}
\\
\ell_n
&
=
\sum_{0\leqslant s^n\leqslant\tau^n}\,
a_{s^n}^n.
\endaligned
\]
As we already said in Section~6 above, it appears
{\em a posteriori} more adequate
to write everything in terms of the differences:
\[
\small
\aligned
\ell_1-\ell_2
=
\sum_{0\leqslant s^1\leqslant\tau^1}\,a_{s^1}^1,
\ \ 
\dots,
\ \
\ell_{n-1}-\ell_n
=
\sum_{0\leqslant s^{n-1}\leqslant\tau^{n-1}}\,a_{s^{n-1}}^{n-1},
\ \ \ \ \ \
\ell_n
=
\sum_{0\leqslant s^n\leqslant\tau^n}\,a_{s^n}^n
\endaligned
\]
of lengths of lines, which are nothing but row lengths of
blocks.

\begin{Proposition}
\label{majoration-general}
Fix $\mu_l^i$, $\nu_l^i$, $\tau^i$ and $\gamma^i ( s^i)$.
If $\alpha_1', \dots, \alpha_{ n-1}', \alpha_n'$ are arbitrary
nonnegative integers satisfying: 
\[
\alpha_1'+\cdots+\alpha_{n-1}'+\alpha_n'
\leqslant
{\textstyle{\frac{n(n+1)}{2}}},
\]
then there exists a positive quantity ${\sf
Constant}_{n, \kappa} > 0$ depending on $n$ and on $\kappa$ which is
independent of $m$ such that:
\[
\footnotesize
\aligned
\sum_{{\sf YT}\in{\sf YT}_{\!\kappa,m}
\left(\mu_l^i,\nu_l^i,\tau^i,\gamma^i(s^i)\right)}
&
\big(\ell_1({\sf YT})-\ell_2({\sf YT})\big)^{\alpha_1'}
\cdots
\big(\ell_{n-1}({\sf YT})-\ell_n({\sf YT})\big)^{\alpha_{n-1}'}
\big(\ell_n({\sf YT})\big)^{\alpha_n'}
\leqslant
\\
&
\leqslant
{\sf Constant}_{n,\kappa}
\cdot
m^{\alpha_1'+\cdots+\alpha_{n-1}'+\alpha_n'}
\cdot
m^{D-1},
\endaligned
\]
where $D = \sum_{ i=1}^n\, ( 1+ \tau^i)$ is the common number of
pairwise distinct $*$-ed columns shared by all Young tableaux ${\sf
YT} \in {\sf YT}_{\!\kappa, m} \big( \mu_l^i, \nu_l^i, \tau^i,
\gamma^i(s^i) \big)$.
\end{Proposition}

\proof
Substituting the values $\ell_i - \ell_{ i+1} = \sum_{ 0 \leqslant s^i
\leqslant \tau^i}\, a_{s^i}^i$ in the monomial:
\[
\big(\ell_1-\ell_2
\big)^{\alpha_1'}\cdots\big(\ell_{n-1}-\ell_n\big)^{
\alpha_{n-1}'}\big(\ell_n\big)^{\alpha_n'},
\]
and expanding the result, we may majorate:
\[
\aligned
&
\big(\ell_1-\ell_2\big)^{\alpha_1'}
\cdots
\big(\ell_{n-1}-\ell_n\big)^{\alpha_{n-1}'}
\big(\ell_n\big)^{\alpha_n'}
\leqslant
\\
&
\leqslant
{\sf Constant}_{\tau^1,\dots,\tau^n}\cdot
\sum_{\sum\,\beta_{s^1}^1+\cdots+\sum\,\beta_{s^n}^n
=
\alpha_1'+\cdots+\alpha_{n-1}'+\alpha_n'}
\bigg(
\prod_{i=1}^n\,
\prod_{0\leqslant s^i\leqslant\tau^i}\,
\big(a_{s^i}^i\big)^{\beta_{s^i}^i}
\bigg).
\endaligned
\]
Since according to~\thetag{ \ref{bound-tau-i}}
above, the $\tau^i \leqslant i ( \kappa - i)
\leqslant n \kappa$ are majorated in terms of $n$ and $\kappa$, we
have:
\[
{\sf Constant}_{\tau^1,\dots,\tau^n} 
\leqslant
{\sf Constant}_{ n,\kappa}.
\]
Moreover, since $\alpha_1 ' + \cdots + \alpha_{ n-1}' + \alpha_n'
\leqslant \frac{ n ( n+1)}{ 2}$, the number of terms in the sum:
\[
\sum_{\sum\,\beta_{s^1}^1+\cdots+\sum\,\beta_{s^n}^n
=
\alpha_1'+\cdots+\alpha_{n-1}'+\alpha_n'}
\big(\bullet\big)
\]
is also $\leqslant {\sf Constant}_{ n, \kappa}$.  Consequently, in
order to prove the proposition, it suffices to majorate by the same
claimed majorant:
\[
{\sf Constant}_{n,\kappa}
\cdot
m^{\alpha_1'+\cdots+\alpha_{n-1}'+\alpha_n'}
\cdot
m^{D-1}
\]
every single sum of the form:
\[
\sum_{{\sf YT}\in{\sf YT}_{\!\kappa,m}
\left(\mu_l^i,\nu_l^i,\tau^i,\gamma^i(s^i)\right)}
\bigg(
\prod_{0\leqslant s^1\leqslant\tau^1}\,
\big(a_{s^1}^1\big)^{\beta_{s^1}^1}
\cdots
\prod_{0\leqslant s^n\leqslant\tau^n}\,
\big(a_{s^n}^n\big)^{\beta_{s^n}^n}
\bigg),
\]
where the exponents $\beta_{ s^i}^i$, $i = 1, \dots, n$, $0 \leqslant
s^i \leqslant \tau^i$, are now fixed and subjected to the same
property that their sum equals:
\[
\sum_{i=1}^n\, 
\sum_{0\leqslant s^i\leqslant\tau^i}\,
\beta_{s^i}^i
=
\alpha_1'+\cdots+\alpha_{n-1}'+\alpha_n'.
\] 
Recall that Young tableaux in the family ${\sf YT}_{\!\kappa, m} \big(
\mu_l^i, \nu_l^i, \tau^i, \gamma^i(s^i) \big)$ have fixed set of
pairwise distinct columns, and that the freedom only lies in the
multiplicities $a_{ s^i}^i \geqslant 1$, $i = 1, \dots, n$, $0
\leqslant s^i \leqslant \tau^i$, of the columns. The considered sum:
\[
\sum_{{\sf YT}\in{\sf YT}_{\!\kappa,m}
\left(\mu_l^i,\nu_l^i,\tau^i,\gamma^i(s^i)\right)}
\big(\bullet\big)
\] 
coincides therefore with the sum: 
\[
\sum_{
\sum_{i=1}^n\,\sum_{0\leqslant s^i\leqslant\tau^i}\,
\left[\gamma_1^i(s^i)+\cdots+\gamma_i^i(s^i)\right]\,a_{s^i}^i=m}
\big(\bullet\big),
\]
which takes precisely account of the homogeneity constraint~\thetag{
\ref{m-homogeneity-gamma}}. 
Let us now set:
\[
b_{s^i}^i
:=
\left[\gamma_1^i(s^i)+\cdots+\gamma_i^i(s^i)\right]\,a_{s^i}^i,
\]
whence $a_{ s^i}^i \leqslant b_{ s^i}^i$ always\footnote{\,
Even in the case where the block of depth $i$ is inextant.
}, 
so that the sum in question now writes: 
\[
\small
\aligned
\sum_{
\sum\,b_{s^1}^1+\cdots+\sum\,b_{s^n}^n
=m}
&
\bigg(
\prod_{i=1}^n\,\prod_{0\leqslant s^i\leqslant\tau^i}\,
\big(a_{s^i}^i\big)^{\beta_{s^i}^i}
\bigg)
\\
&
\leqslant
\sum_{
\sum\,b_{s^1}^1+\cdots+\sum\,b_{s^n}^n
=m}
\bigg(
\prod_{i=1}^n\,\prod_{0\leqslant s^i\leqslant\tau^i}\,
\big(b_{s^i}^i\big)^{\beta_{s^i}^i}
\bigg).
\endaligned
\]
The number of nonzero variables $b_{ s^i}^i \in \N$ here is the same,
equal to $D$, as the number of nonzero exponents $a_{ s^i}^i$. The
conclusion now follows from the elementary general inequality:
\[
\sum_{b_1+\cdots+b_D=m\atop
b_1\geqslant 1,\dots,b_D\geqslant 1}\,
b_1^{\beta_1}\cdots b_D^{\beta_D}
\leqslant
{\sf Constant}_D
\cdot
m^{\beta_1+\cdots+\beta_D}
\cdot
m^{D-1},
\]
that can be established by approximating the sum by a Riemann
integral; of course, ${\sf Constant}_D \leqslant {\sf Constant}_{ n,
\kappa}$.
\endproof

From this proposition, we will deduce a few corollaries.  Firstly, as
announced earlier on at the end of Section~6, we have:

\begin{Corollary}
\label{corollary-remainder-1}
Let $\alpha_1', \dots, \alpha_{ n-1}', \alpha_n'$ be nonnegative
integers satisfying:
\[
\alpha_1'+\cdots+\alpha_{n-1}'+\alpha_n'
\leqslant
{\textstyle{\frac{n(n+1)}{2}}}
-
1.
\]
Then the following majoration holds for the summation over
{\em all} semi-standard Young tableaux of weight $m$:
\[
\aligned
\sum_{{\sf YT}\,{\sf semi-standard}
\atop
{\sf weight}({\sf YT})=m}
&
\big(\ell_1({\sf YT})-\ell_2({\sf YT})\big)^{\alpha_1'}
\cdots
\big(\ell_{n-1}({\sf YT})-\ell_n({\sf YT})\big)^{\alpha_{n-1}'}
\big(\ell_n({\sf YT})\big)^{\alpha_n'}
\leqslant
\\
&
\leqslant
{\sf Constant}_{n,\kappa}\cdot
m^{\alpha_1'+\cdots+\alpha_{n-1}'+\alpha_n'}\cdot
m^{n\kappa-\frac{n(n-1)}{2}}
\\
&
\leqslant
{\sf Constant}_{n,\kappa}\cdot
m^{(\kappa+1)n-2}.
\endaligned
\]
\end{Corollary}

\proof
According to Lemma~\ref{lemma-disjoint-YT}:
\[
\sum_{{\sf YT}\,{\sf semi-standard}
\atop
{\sf weight}({\sf YT})=m}\,
\big(\bullet\big)
=
\sum_{\mu_l^i,\,\,\nu_l^i}\,
\sum_{\tau_i}\,
\sum_{\gamma^i(s^i)}\,
\sum_{{\sf YT}\in
{\sf YT}(\mu_l^i,\nu_l^i,\tau^i,\gamma^i(s^i))}\,
\big(\bullet\big),
\]
and furthermore, the number of terms in the three first sums of the
right-hand side is $\leqslant {\sf Constant}_{ n, \kappa}$.  It
suffices then to apply the proposition which controls each fourth sum,
reminding from Lemma~\ref{bound-distinct-columns} that
each $D = \sum_{ i=1}^n \, 
(1 + \tau^i)$ 
is in any case $\leqslant n\kappa - \frac{ n ( n - 1)}{ 2}$.
\endproof

Secondly, from the simple inequality~\thetag{ \ref{l-l-l}}, 
we immediately deduce:

\begin{Corollary}
\label{corollary-remainder-2}
If $\alpha_1, \dots, \alpha_n$ are any nonnegative integers satisfying
$\alpha_1 + \cdots + \alpha_n \leqslant \frac{n(n+1)}{2} -1$, then:
\[
\aligned
\sum_{{\sf YT}\,{\sf semi-standard}
\atop
{\sf weight}({\sf YT})=m}
\big(\ell_1({\sf YT})\big)^{\alpha_1}
\cdots\,
\big(\ell_n({\sf YT})\big)^{\alpha_n}
&
\leqslant
{\sf Constant}_{n,\kappa}
\cdot
m^{\alpha_1+\cdots+\alpha_n}
\cdot
m^{n\kappa-\frac{n(n-1)}{2}}
\\
&
\leqslant
{\sf Constant}_{n,\kappa}
\cdot
m^{(\kappa+1)n-2}.
\endaligned
\]
\end{Corollary}

Lastly, we now introduce a certain collection of semi-standard Young
tableaux the contribution of which appears to also fall in the
remainder ${\sf O}_{ n, \kappa} \big( m^{ (\kappa+1)n - 2} \big)$:
we gather all the ones for which the number of
pairwise distinct columns is (strictly) less than the maximal possible
number $n\kappa - \frac{ n ( n-1)}{ 2}$:
\[
\boxed{
{\sf NGYT}_{\kappa,m}
:=
\bigcup_{
\mu_l^i,\nu_l^i,\tau^i,\gamma^i(s^i)
\atop
\sum_{i=1}^n(1+\tau^i)
\leqslant
n\kappa-\frac{n(n-1)}{2}-1}\,
{\sf YT}_{\!\kappa,m}
\big(\mu_l^i,\nu_l^i,\tau^i,\gamma^i(s^i)\big)}\,. 
\]
The number of appearing such families ${\sf YT}_{ \!\kappa, m} \big(
\mu_l^i, \nu_l^i, \tau^i, \gamma^i(s^i) \big)$ is of course less than
the majorant:
\[
\prod_{i=1}^n\,
\big(
1+
{\textstyle{\frac{\kappa!}{(\kappa-i)!\,i!}}}
\big)^{1+i(\kappa-i)}
\equiv
{\sf Constant}_{n,\kappa}
\]
provided by Lemma~\ref{lemma-disjoint-YT} 
for the total number of all families. 

\begin{Lemma}
For any integers $\alpha_1', \dots, \alpha_{ n-1}', \alpha_n'$ whose sum
equals $\frac{ n ( n+1)}{ 2}$, the contribution of:
\[
\aligned
\sum_{{\sf YT}\in{\sf NGYT}_{\kappa,m}}
&
\big(\ell_1({\sf YT})-\ell_2({\sf YT})\big)^{\alpha_1'}
\cdots
\big(\ell_{n-1}({\sf YT})-\ell_n({\sf YT})\big)^{\alpha_{n-1}'}
\big(\ell_n({\sf YT})\big)^{\alpha_n'}
\leqslant
\\
&
\leqslant
{\sf Constant}_{n,\kappa}\cdot
m^{(\kappa+1)n-2}
\endaligned
\]
is asymptotically negligible in comparison to the dominant power $m^{
(\kappa +1)n - 1}$.
\end{Lemma}

\proof
By what has been just seen, it suffices to verify that such a
majoration holds for a sum $\sum_{ {\sf YT} \in ( \bullet)}$ running
over a single family ${\sf YT}_{\!\kappa, m} \big( \mu_l^i, \nu_l^i,
\tau^i, \gamma^i ( s^i) \big)$ contained in ${\sf NGYT}_{ \kappa, m}$,
and this was already achieved by Proposition~\ref{majoration-general} 
above, since
$\frac{ n ( n+1)}{ 2} + D - 1 \leqslant
(\kappa+1)n-2$ always when $D \leqslant n\kappa - \frac{
n (n - 1)}{ 2} - 1$.
\endproof

\markleft{Jo\"el Merker}
\markright{\S8.~Maximal length families of semi-standard Young tableaux}
\section{\bf Maximal length families 
\\
of semi-standard Young tableaux}
\label{Section-8}

\subsection{Relevant families of $\Delta$-monomials}
In the remainder of the paper, we shall
now consider only exponents $\alpha_i'$ satisfying
$\alpha_1' + \cdots + \alpha_{ n-1}' + \alpha_n' = 
\frac{ n ( n+1)}{ 2}$. 
From the proposition just proved, we deduce that the complete sum:
\begin{equation}
\label{complete-sum}
\small
\aligned
\sum_{{\sf YT}\,{\sf semi-standard}
\atop
{\sf weight}({\sf YT})=m}
&
\big(\ell_1({\sf YT})-\ell_2({\sf YT})\big)^{\alpha_1'}
\cdots
\big(\ell_{n-1}({\sf YT})-\ell_n({\sf YT})\big)^{\alpha_{n-1}'}
\big(\ell_n({\sf YT})\big)^{\alpha_n'}
\endaligned
\end{equation}
splits up as a negligible sum plus a relevant sum that we should now
study:
\[
\sum_{{\sf YT}\in{\sf NGYT}_{\kappa,m}}
+
\sum_{{\sf YT}\not\in{\sf NGYT}_{\kappa,m}}.
\]
Hence, what remains to be studied is the collection of all families of
semi-standard Young tableaux:
\[
\boxed{
{\sf YT}_{\kappa,m}^{\rm max}
:=
\bigcup_{
\mu_l^i,\nu_l^i,\tau^i,\gamma^i(s^i)
\atop
\sum_{i=1}^n(1+\tau^i)
\,=\,
n\kappa-\frac{n(n-1)}{2}}\,
{\sf YT}_{\!\kappa,m}
\big(\mu_l^i,\nu_l^i,\tau^i,\gamma^i(s^i)\big)} 
\]
for which the number of pairwise distinct columns is maximal, equal to
$n\kappa - \frac{ n ( n-1)}{ 2}$. The following statement describes them
in great details.

\begin{Proposition}
\label{proposition-maximal}
The number $D$ of pairwise distinct columns in any semi-standard
Young tableau written as in~\thetag{ 
\ref{arbitrary-semi-standard}}
is in any case $\leqslant
n \kappa - \frac{ n ( n-1) }{ 2}$. Furthermore, a given semi-standard
Young tableau reaches the maximal number:
\[
D
=
n\kappa
-
{\textstyle{\frac{n(n-1)}{2}}}
\]
of pairwise distinct columns if and only if all the following
conditions are fulfilled:

\begin{itemize}

\smallskip\item[$\bullet$]
the depth of the tableau is maximal: $d_1 = n$; 

\smallskip\item[$\bullet$] 
nonvoid blocks of any depth $i = 1, 2, 3, \dots, n-1, n$ are all extant, 
so that the number of nonvoid blocks is maximal, equal to $n$; 

\smallskip\item[$\bullet$]
the leftmost $*$-ed column of the tableau corresponds to the
$n$-dimensional Wronskian $\Delta_{ 1, 2, 3, \dots, n-1, n}^{ 1, 2, 3,
\dots, n-1, n}$, reproduced a certain number $* \geqslant 1$ of times;

\smallskip\item[$\bullet$]
the bottom-right entry of every block is maximal:
\[
\nu_1^1=\nu_2^2=\nu_3^3=\cdots=\nu_{n-1}^{n-1}=\nu_n^n=\kappa; 
\]
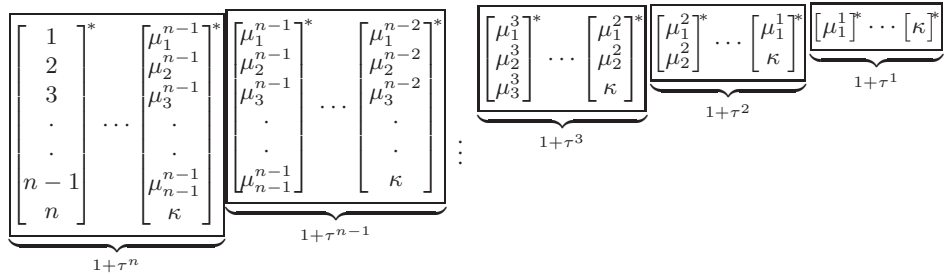
\[
\footnotesize
\aligned
\begin{array}{c}
\underbrace{\boxed{\!\!
\left[\!\!\!\!
\begin{array}{c}
1
\\
2
\\
3
\\
\cdot
\\
\cdot
\\
n-1
\\
n
\end{array}
\!\!\!\!\right]^{\!\!*}
\cdots
\left[\!\!\!\!
\begin{array}{c}
\mu_1^{n-1}
\\
\mu_2^{n-1}
\\
\mu_3^{n-1}
\\
\cdot
\\
\cdot
\\
\mu_{n-1}^{n-1}
\\
\kappa
\end{array}
\!\!\!\right]^{\!\!*}\!\!}}_{1+\tau^n}
\\ 
\rule[-2.17pt]{0pt}{2pt}
\end{array}
\!\!\!\!\!\!
\begin{array}{c}
\underbrace{\boxed{\!\!
\left[\!\!\!\!
\begin{array}{c}
\mu_1^{n-1}
\\
\mu_2^{n-1}
\\
\mu_3^{n-1}
\\
\cdot
\\
\cdot
\\
\mu_{n-1}^{n-1}
\end{array}
\!\!\!\right]^{\!\!*}
\cdots
\left[\!\!\!
\begin{array}{c}
\mu_1^{n-2}
\\
\mu_2^{n-2}
\\
\mu_3^{n-2}
\\
\cdot
\\
\cdot
\\
\kappa
\end{array}
\!\!\!\right]^{\!\!*}\!\!}}_{1+\tau^{n-1}}
\\ 
\rule[-2.17pt]{0pt}{22.5pt}
\end{array}
\vdots
\begin{array}{c}
\underbrace{\boxed{\!\!
\left[\!\!\!
\begin{array}{c}
\mu_1^3
\\
\mu_2^3
\\
\mu_3^3
\end{array}
\!\!\!\right]^{\!\!*}
\cdots
\left[\!\!\!
\begin{array}{c}
\mu_1^2
\\
\mu_2^2
\\
\kappa
\end{array}
\!\!\!\right]^{\!\!*}\!\!}}_{1+\tau^3}
\\ 
\rule[-2.17pt]{0pt}{59pt}
\end{array}
\!\!\!\!\:\:\!\!\!\!
\begin{array}{c}
\underbrace{\boxed{\!\!
\left[\!\!\!
\begin{array}{c}
\mu_1^2
\\
\mu_2^2
\end{array}
\!\!\!\right]^{\!\!*}
\cdots
\left[\!\!\!
\begin{array}{c}
\mu_1^1
\\
\kappa
\end{array}
\!\!\!\right]^{\!\!*}\!\!}}_{1+\tau^2}
\\ 
\rule[-2.17pt]{0pt}{71pt}
\end{array}
\!\!\!\!\:\:\!\!\!\!
\begin{array}{c}
\underbrace{\boxed{\!\!
\left[\!\!\!
\begin{array}{c}
\mu_1^1
\end{array}
\!\!\!\right]^{\!\!*}
\cdots
\left[\!\!\!
\begin{array}{c}
\kappa
\end{array}
\!\!\!\right]^{\!\!*}\!\!}}_{1+\tau^1}
\\ 
\rule[-2.17pt]{0pt}{83pt}
\end{array}.
\endaligned
\vspace{-0.5cm}
\]

\smallskip\item[$\bullet$]
the border column entries (excepting the last one, equal to $\kappa$,
of the longest column) of any pair of neighboring blocks are equal:
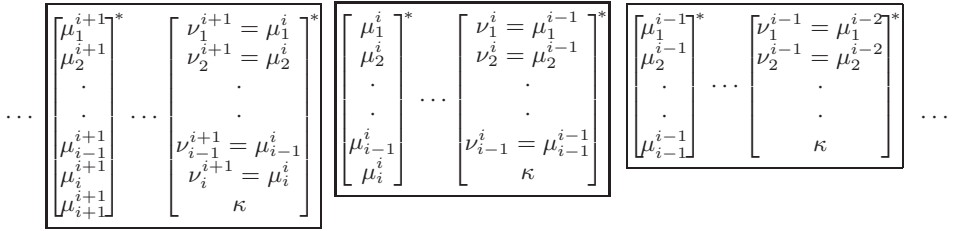
\[
\footnotesize
\aligned
\cdots\,\,
\boxed{\!\!
\left[\!\!\!\!
\begin{array}{c}
\mu_1^{i+1}
\\
\mu_2^{i+1}
\\
\cdot
\\
\cdot
\\
\mu_{i-1}^{i+1}
\\
\mu_i^{i+1}
\\
\mu_{i+1}^{i+1}
\end{array}
\!\!\!\!\right]^{\!\!*}
\cdots
\left[\!\!\!\!
\begin{array}{c}
\nu_1^{i+1}
=
\mu_1^i
\\
\nu_2^{i+1}
=
\mu_2^i
\\
\cdot
\\
\cdot
\\
\nu_{i-1}^{i+1}=\mu_{i-1}^i
\\
\nu_i^{i+1}=\mu_i^i
\\
\kappa
\end{array}
\!\!\!\!\right]^{\!\!*}\!\!}
\begin{array}{c}
\boxed{\!\!
\left[\!\!\!\!
\begin{array}{c}
\mu_1^i
\\
\mu_2^i
\\
\cdot
\\
\cdot
\\
\mu_{i-1}^i
\\
\mu_i^i
\end{array}
\!\!\!\!\right]^{\!\!*}
\cdots
\left[\!\!\!\!
\begin{array}{c}
\nu_1^i=\mu_1^{i-1}
\\
\nu_2^i=\mu_2^{i-1}
\\
\cdot
\\
\cdot
\\
\nu_{i-1}^i=\mu_{i-1}^{i-1}
\\
\kappa
\end{array}
\!\!\!\!\right]^{\!\!*}\!\!}
\\
\rule[-2.17pt]{0pt}{12pt}
\end{array}
\!\!\!\!\:\:\!
\begin{array}{c}
\boxed{\!\!
\left[\!\!\!\!
\begin{array}{c}
\mu_1^{i-1}
\\
\mu_2^{i-1}
\\
\cdot
\\
\cdot
\\
\mu_{i-1}^{i-1}
\end{array}
\!\!\!\!\right]^{\!\!*}
\cdots
\left[\!\!\!\!
\begin{array}{c}
\nu_1^{i-1}=\mu_1^{i-2}
\\
\nu_2^{i-1}=\mu_2^{i-2}
\\
\cdot
\\
\cdot
\\
\kappa
\end{array}
\!\!\!\!\right]^{\!\!*}\!\!}
\\ 
\rule[-2.17pt]{0pt}{22pt}
\end{array}
\cdots,
\endaligned
\]

\smallskip\item[$\bullet$]
the number of pairwise distinct columns in each block of depth $i$,
for $i=1, 2, 3, \dots, n-1, n$ is 
maximal\footnote{\,
By convention, we shall also call $\mu_1^n, \mu_2^n, \mu_3^n, \dots,
\mu_{ n-1}^n, \mu_n^n$ the entries $1, 2, 3, \dots, n-1, n$ of the
leftmost $*$-ed column.
}, 
equal to:
\[
\aligned
1
+
\tau^i
:=
&\,
1
+
(\mu_1^{i-1}-\mu_1^i)
+
(\mu_2^{i-1}-\mu_2^i)
+\cdots+
(\mu_{i-1}^{i-1}-\mu_{i-1}^i)
+
(\kappa-\mu_i^i)
\\
=
&\,
1+\kappa+
\sum_{l=1}^{i-1}\,\mu_l^{i-1}
-
\sum_{l=1}^i\,\mu_l^i,
\endaligned
\]
so that the total number of pairwise distinct columns is accordingly
indeed equal to:
\[
\small
\aligned
(1+\tau^1)+(1+\tau^2)+(1+\tau^3)
+\cdots+
(1+\tau^{n-1})+(1+\tau^n)
&
=
n
+
n\kappa
-
\sum_{l=1}^n\,\mu_l^n
\\
&
=
n\kappa
-
{\textstyle{\frac{n(n-1)}{2}}}. 
\endaligned
\]

\end{itemize}

\end{Proposition}

\proof
The majorant $n \kappa - \frac{ n ( n-1)}{ 2}$ has already been
obtained above, before the introduction of families of Young tableaux.
The remaining statements follow by thinking once again about what has
already been seen above.
\endproof

So the families of semi-standard Young tableaux having maximal number
$n\kappa - \frac{ n ( n-1)}{ 2}$ of pairwise distinct columns is
parameterized by all the collections of integers $\mu_i^j$ satisfying the
inequalities:
\begin{equation}
\label{strict-mu}
\aligned
&
1\leqslant\mu_1^1<\kappa
\\
&
1\leqslant\mu_1^2<\mu_2^2<\kappa
\\
&
1\leqslant\mu_1^3<\mu_2^3<\mu_3^3<\kappa
\\
&
\cdots\cdots\cdots\cdots\cdots\cdots\cdots
\\
&
1\leqslant\mu_1^{n-1}\!\!<\!\!\mu_2^{n-1}\!\!<\!\!\mu_3^{n-1}
<\cdots<\!
\mu_{n-1}^{n-1}\!\!<\kappa
\\
&
1\leqslant\mu_1^n<\,\mu_2^n<\,\mu_3^n
<\,\cdots\,<
\,\mu_{n-1}^n<\mu_n^n,
\endaligned
\end{equation}
together with the further semi-standard-like inequalities\footnote{\,
Diagramatically, this second set of inequalities reads as saying that
{\em vertically} in each column of the first array~\thetag{
\ref{strict-mu}} of inequalities, the integers $\mu_l^i$ are weakly
decreasing. In particular, $\kappa \geqslant \mu_n^n = n$, as
was assumed throughout earlier on.
}: 
\begin{equation}
\label{weak-mu}
\mu_l^l\geqslant\mu_l^{l+1}\geqslant\cdots\geqslant\mu_l^{n-1}
(\geqslant l).
\end{equation}
For brevity, let us write as:
\[
\boxed{
\mu_l^i\in\nabla_{\!n,\kappa}}
\]
the condition that the $\mu_l^i$ satisfy the two sets of
inequalities~\thetag{ \ref{strict-mu}} and~\thetag{
\ref{weak-mu}}. For any such choice of $\mu_l^i \in \nabla_{ \! n,
\kappa}$, we shall denote by:
\[
\boxed{{\sf YT}_{\!\kappa,m}^{\rm max}\big(\mu_l^i\big)}
\]
the family of semi-standard Young tableaux which consist of all
possible concatenations:
\[
{\sf block}^n\big(\gamma^n\big)
\cdots
{\sf block}^i\big(\gamma^i\big)
\cdots
{\sf block}^1\big(\gamma^1\big)
\]
of pathed blocks of the form: 
\[
\footnotesize
\aligned
{\sf block}^i\big(\gamma^i\big)
:=
\boxed{\!
\left[\!\!\!
\begin{array}{c}
\mu_1^i=\gamma_1^i(0)
\\
\mu_2^i=\gamma_2^i(0)
\\
\cdot\cdot
\\
\mu_{i-1}^i=\gamma_{i-1}^i(0)
\\
\mu_i^i=\gamma_i^i(0)
\end{array}
\!\!\!\right]^{\!a_0^i}
\left[\!\!\!
\begin{array}{c}
\gamma_1^i(1)
\\
\gamma_2^i(1)
\\
\cdot
\\
\gamma_{i-1}^i(1)
\\
\gamma_i^i(1)
\end{array}
\!\!\!\right]^{\!a_1^i}
\cdots
\left[\!\!\!
\begin{array}{c}
\gamma_1^i(s^i)
\\
\gamma_2^i(s^i)
\\
\cdot
\\
\gamma_{i-1}^i(s^i)
\\
\gamma_i^i(s^i)
\end{array}
\!\!\!\right]^{\!a_{s^i}^i}
\cdots
\left[\!\!\!
\begin{array}{c}
\gamma_1^i(\tau^i)=\mu_1^{i-1}
\\
\gamma_2^i(\tau^i)=\mu_2^{i-1}
\\
\cdot\cdot
\\
\gamma_{i-1}^i(\tau^i)=\mu_{i-1}^{i-1}
\\
\gamma_i^i(\tau^i)=\kappa
\end{array}
\!\!\!\right]^{\!a_{\tau^i}^i}\!}\,,
\endaligned
\]
where the lengths $1 + \tau^i$ of paths are maximal equal to:
\[
1+\tau^i
=
1+\kappa
+
\sum_{l=1}^{i-1}\,\mu_l^{i-1}
-
\sum_{l=1}^i\,\mu_l^i,
\]
so that between two successive $*$-ed columns:
\[
\left[\!\!
\begin{array}{c}
\gamma_1^i(s^i)
\\
\gamma_2^i(s^i)
\\
\cdot
\\
\gamma_i^i(s^i)
\end{array}
\!\!\right]^{\!*}
<
\left[\!\!
\begin{array}{c}
\gamma_1^i(s^i+1)
\\
\gamma_2^i(s^i+1)
\\
\cdot\cdot
\\
\gamma_i^i(s^i+1)
\end{array}
\!\!\right]^{\!*}
\]
one has the semi-standard inequalities $\gamma_l^i ( s^i) \leqslant
\gamma_l^i ( s^i + 1)$ for $l = 1, 2, \dots, i$ but the jump is
smallest possible, namely $\gamma_l^i ( s^i + 1) = \gamma_l^i ( s^i)$
for all $l$ except only one $l_0$ for which:
\[
\gamma_{l_0}^i(s^i+1)
=
1
+
\gamma_{l_0}^i(s^i). 
\] 
Such paths all of which jumps are unit
will be called {\sl tight paths}. With all these notations,
the initial complete sum~\thetag{
\ref{complete-sum}}
to be studied now writes:
\[
\small
\aligned
\sum_{{\sf YT}\,{\sf semi-standard}
\atop
{\sf weight}({\sf YT})=m}
(\bullet)
=
\sum_{{\sf YT}\in{\sf NGYT}_{\kappa,m}}
(\bullet)
+
\sum_{\mu_l^i\in\nabla_{\!n,\kappa}}\,
\sum_{{\sf YT}\in{\sf YT}_{\!\kappa,m}^{\rm max}
\left(\mu_l^i\right)}
(\bullet)
\endaligned,
\]
the first sum being negligible, in the sense that: 
\begin{equation}
\label{negl-max}
\boxed{
\footnotesize
\aligned
&
\sum_{{\sf YT}\,{\sf semi-standard}
\atop
{\sf weight}({\sf YT})=m}
\big(\ell_1({\sf YT})-\ell_2({\sf YT})\big)^{\alpha_1'}
\cdots
\big(\ell_{n-1}({\sf YT})-\ell_n({\sf YT})\big)^{\alpha_{n-1}'}
\big(\ell_n({\sf YT})\big)^{\alpha_n'}
=
\\
&
=
{\sf O}_{n,\kappa}\big(m^{(\kappa+1)n-2}\big)
+
\\
&
+
\sum_{\mu_l^i\in\nabla_{\!n,\kappa}}\,
\sum_{{\sf YT}\in{\sf YT}_{\!\kappa,m}^{\rm max}
\left(\mu_l^i\right)}
\big(\ell_1({\sf YT})-\ell_2({\sf YT})\big)^{\alpha_1'}
\cdots
\big(\ell_{n-1}({\sf YT})-\ell_n({\sf YT})\big)^{\alpha_{n-1}'}
\big(\ell_n({\sf YT})\big)^{\alpha_n'}
\endaligned}\,.
\end{equation}

\subsection{Grouping sums}
Let us therefore fix $\mu_l^i \in \nabla_{\! n, \kappa}$, let
us keep aside (and in mind) the first summation $\sum_{ \mu_l^i \in 
\nabla_{ \! n, \kappa}} ( \bullet)$, and let us
consider only the second summation:
\[
\small
\aligned
\sum_{{\sf YT}\in{\sf YT}_{\!\kappa,m}^{\rm max}
\left(\mu_l^i\right)}
\big(\ell_1({\sf YT})-\ell_2({\sf YT})\big)^{\alpha_1'}
\cdots
\big(\ell_{n-1}({\sf YT})-\ell_n({\sf YT})\big)^{\alpha_{n-1}'}
\big(\ell_n({\sf YT})\big)^{\alpha_n'}. 
\endaligned
\]
With $\tau^i := \kappa + \sum_{ l=1}^{ i-1}\, \mu_l^{ i-1} - \sum_{
l=1}^i \, \mu_l^i$, this apparently compact sum identifies in fact
with the multiple sums:
\[
\sum_{\gamma^1(s^1)}\cdots\sum_{\gamma^n(s^n)}\,\,
\sum_{a_{s^1}^1}\cdots\sum_{a_{s^n}^n}\,
(\bullet),
\]
the paths $\gamma^i (s^i)$ being all tight and the multiplicities $a_{
s^i}^i$ being constrained only by the weight condition:
\[
\sum_{i=1}^n\,\sum_{0\leqslant s^i\leqslant\tau^i}\,
\big[\gamma_1^i(s^i)+\cdots+\gamma_i^i(s^i)\big]\,a_{s^i}^i
=
m.
\]
Since the paths are tight, the sums (equal to the weights of columns): 
\[
\gamma_1^i(s^i)+\cdots+\gamma_i^i(s^i)
\]
jump of a single unit as $s^i = 0, 1, 2, \dots, \tau^i$, and to be
precise, they take the following exact values:
\[
\gamma_1^i(s^i)+\cdots+\gamma_i^i(s^i)
=
\left\{
\aligned
&
{\textstyle{\sum_{l=1}^i}}\,\mu_l^i\ \ \ \ \ \ \ \ \ \ \ \
\text{\rm for}\ \ s^i=0
\\
&
1+{\textstyle{\sum_{l=1}^i}}\,\mu_l^i\ \ \ \ \
\text{\rm for}\ \ s^i=1
\\
&
2+{\textstyle{\sum_{l=1}^i}}\,\mu_l^i\ \ \ \ \
\text{\rm for}\ \ s^i=2
\\
&
\cdots\cdots\cdots\cdots\cdots\cdots\cdots\cdots
\\
&
\kappa+{\textstyle{\sum_{l=1}^{i-1}}}\,\mu_l^i\ \ \ \ \
\text{\rm for}\ \ s^i=\tau^i,
\endaligned
\right.
\]
{\em independently of the paths}. In other words:
\[
\gamma_1^i(s^i)+\cdots+\gamma_i^i(s^i)
=
s^i+\mu_1^i+\cdots+\mu_i^i,
\]
with of course at the end (as already written): 
\[
\aligned
\gamma_1^i(\tau^i)+\cdots+\gamma_i^i(\tau^i)
&
=
\kappa+\sum_{l=1}^{i-1}\,\mu_l^{i-1}-\sum_{l=1}^i\,\mu_l^i
+\mu_1^i+\cdots+\mu_i^i
\\
&
=
\kappa+\mu_1^{i-1}+\cdots+\mu_{i-1}^{i-1}. 
\endaligned
\]
Recalling that\footnote{\,
By convention, $\ell_{ n+1} = \ell_{ n+1} \big( {\sf YT} \big) = 0$
always.
}: 
\[
\big(\ell_i({\sf YT})-\ell_{i+1}({\sf YT})\big)^{\alpha_i'}
=
\big(a_0^i+\cdots+a_{\tau^i}^i\big)^{\alpha_i'}
\ \ \ \ \ \ \ \ \ \ \ \ \
{\scriptstyle{(1\,\leqslant\,i\,\leqslant\,n)}},
\]
we may therefore write: 
\[
\small
\aligned
&
\sum_{{\sf YT}\in{\sf YT}_{\!\kappa,m}^{\rm max}
\left(\mu_l^i\right)}
\big(\ell_1({\sf YT})-\ell_2({\sf YT})\big)^{\alpha_1'}
\cdots
\big(\ell_{n-1}({\sf YT})-\ell_n({\sf YT})\big)^{\alpha_{n-1}'}
\big(\ell_n({\sf YT})\big)^{\alpha_n'}
=
\\
&
=
\sum_{\gamma^1(s^1)}\cdots\sum_{\gamma^n(s^n)}\,\,\,
\sum_{
\sum_{i=1}^n\sum_{0\leqslant s^i\leqslant\tau^i}\,
\left[s^i+\mu_1^i+\cdots+\mu_i^i\right]a_{s^i}^i=m}\,\,\,
\prod_{i=1}^n\,
\big(a_0^i+\cdots+a_{\tau^i}^i\big)^{\alpha_i'}. 
\endaligned
\]
Observing that the last written sum is independent of the paths, each
one of the first $n$ sums $\sum_{ \gamma^i ( s^i)}$ then collapses as
just multiplication by the number of considered paths $\gamma^i (
s^i)$. Thus, let $N_{ \mu_1^1}^\kappa$ denote the number\footnote{\,
In fact trivially, $N_{ \mu_1^1}^\kappa = 1$. }
of tight paths $\gamma^1 (s^1)$ from $\mu_1^1$ to $\kappa$; let $N_{
\mu_1^2, \mu_2^2}^{ \mu_1^1, \kappa}$ denote the number
of tight paths $\gamma^2 ( s^2)$ from the column ${}^t(\mu_1^2,
\mu_2^2)$, where ${}^t( \bullet )$ denotes transposition, to the
column ${}^t ( \mu_1^1, \kappa)$; and generally, let:
\[
\label{definition-N}
\boxed{
N_{\mu_1^i,\mu_2^i,\dots,\mu_{i-1}^i,\mu_i^i}^{
\mu_1^{i-1},\mu_2^{i-1},\dots,\mu_{i-1}^{i-1},\kappa}}
\]
denote the number of tight paths $\gamma^i ( s^i)$ from the column
${}^t \big( \mu_1^i, \mu_2^i, \dots, \mu_{ i-1}^i, \mu_i^i \big)$ to
the column ${}^t \big( \mu_1^{ i-1}, \mu_2^{ i-1}, \dots, \mu_{ i-1}^{
i-1}, \kappa \big)$, with the natural convention that, for $i = n$,
one has the notational equivalence:
\[
{}^t\big(\mu_1^n,\mu_2^n,\dots,\mu_{n-1}^n,\mu_n^n)
\equiv
{}^t\big(1,2,\dots,n-1,n\big). 
\]
Then with such notations, we may represent our sum as: 
\[
\small
\aligned
&
\sum_{\gamma^1(s^1)}\cdots\sum_{\gamma^n(s^n)}\,\,\,
\sum_{
\sum_{i=1}^n\sum_{0\leqslant s^i\leqslant\tau^i}\,
\left[s^i+\mu_1^i+\cdots+\mu_i^i\right]a_{s^i}^i=m}\,\,\,
\prod_{i=1}^n\,
\big(a_0^i+\cdots+a_{\tau^i}^i\big)^{\alpha_i'}
=
\\
&
=
N_{\mu_1^1}^\kappa\,
N_{\mu_1^2,\mu_2^2}^{\mu_1^1,\kappa}
N_{\mu_1^3,\mu_2^3,\mu_3^3}^{\mu_1^2,\mu_2^2,\kappa}
\cdots
N_{1,\dots,n-1,n}^{\mu_1^{n-1},\dots,\mu_{n-1}^{n-1},\kappa}\cdot
\\
&
\ \ \ \ \ \ \ \ \ \ 
\cdot
\sum_{
\sum_{i=1}^n\,\sum_{0\leqslant s^i\leqslant\tau^i}\,
\left[s^i+\mu_1^i+\cdots+\mu_i^i\right]a_{s^i}^i=m}\,\,\,
\prod_{i=1}^n\,
\big(a_0^i+\cdots+a_{\tau^i}^i\big)^{\alpha_i'}.
\endaligned
\]
In conclusion, remembering the dropped
$\sum_{ \mu_l^i \in \nabla_{ \! n, \kappa}}$, 
we have established that: 
\[
\small
\aligned
&
\sum_{\mu_l^i\in\nabla_{\!n,\kappa}}\,
\sum_{{\sf YT}\in{\sf YT}_{\!\kappa,m}^{\rm max}
\left(\mu_l^i\right)}
\big(\ell_1({\sf YT})-\ell_2({\sf YT})\big)^{\alpha_1'}
\cdots
\big(\ell_{n-1}({\sf YT})-\ell_n({\sf YT})\big)^{\alpha_2'}
\big(\ell_n({\sf YT})\big)^{\alpha_n'}
=
\endaligned
\]
\[
\small
\aligned
\\
&
=
\sum_{\mu_l^i\in\nabla_{\!n,\kappa}}\,
N_{\mu_1^1}^\kappa\,
N_{\mu_1^2,\mu_2^2}^{\mu_1^1,\kappa}
N_{\mu_1^3,\mu_2^3,\mu_3^3}^{\mu_1^2,\mu_2^2,\kappa}
\cdots
N_{1,\dots,n-1,n}^{\mu_1^{n-1},\dots,\mu_{n-1}^{n-1},\kappa}
\cdot
\\
&
\ \ \ \ \ \ \ \ \ \ \ \ \ \ \ \ \ \ \ \ 
\cdot
\sum_{
\sum_{i=1}^n\sum_{0\leqslant s^i\leqslant\tau^i}\,
\left[s^i+\mu_1^i+\cdots+\mu_i^i\right]a_{s^i}^i=m}\,
\prod_{i=1}^n\,
\big(a_0^i+\cdots+a_{\tau^i}^i\big)^{\alpha_i'}.
\endaligned
\]

\subsection{Approximating sums by integrals}
If we now set, similarly as in \S3: 
\[
b_{s^i}^i
:=
\frac{1}{m}\,
a_{s^i}^i, 
\]
the sum of the last line of the preceding equation: 
\[
{\sf S}_{n,\kappa,m}^{\alpha_1',\dots,\alpha_n'}\big(\mu_l^i\big)
:=
\sum_{
\sum_{i=1}^n\sum_{0\leqslant s^i\leqslant\tau^i}\,
\left[s^i+\mu_1^i+\cdots+\mu_i^i\right]a_{s^i}^i=m}\,
\prod_{i=1}^n\,
\big(a_0^i+\cdots+a_{\tau^i}^i\big)^{\alpha_i'}
\]
becomes: 
\[
\small
\aligned
{\sf S}_{n,\kappa,m}^{\alpha_1',\dots,\alpha_n'}\big(\mu_l^i\big)
&
=
\sum_{
\sum_{i=1}^n\sum_{0\leqslant s^i\leqslant\tau^i}\,
\left[s^i+\mu_1^i+\cdots+\mu_i^i\right]\frac{a_{s^i}^i}{m}=1}\,
m^{\alpha_1'+\cdots+\alpha_n'}\,
\prod_{i=1}^n\,
\bigg(
\sum_{0\leqslant s^i\leqslant\tau^i}\,
\frac{a_{s^i}^i}{m}
\bigg)^{\alpha_i'}
\\
&
=
m^{\frac{n(n+1)}{2}}\,
\sum_{
\sum_{i=1}^n\sum_{0\leqslant s^i\leqslant\tau^i}\,
\left[s^i+\mu_1^i+\cdots+\mu_i^i\right]\frac{a_{s^i}^i}{m}=1}\,
\prod_{i=1}^n\,
\bigg(
\sum_{0\leqslant s^i\leqslant\tau^i}\,
\frac{a_{s^i}^i}{m}
\bigg)^{\alpha_i'}.
\endaligned
\]
Approximating the so obtained sum by a 
Riemann integral\footnote{\,
A different view may be found in~\cite{ bero2007}.
}, 
we get
up to a negligible power of $m$: 
\[
\small
\aligned
{\sf S}_{n,\kappa,m}^{\alpha_1',\dots,\alpha_n'}\big(\mu_l^i\big)
&
=
{\sf O}_{n,\kappa}\big(m^{(\kappa+1)n-2}\big)
+
m^{\frac{n(n+1)}{2}}\cdot
m^{n\kappa-\frac{n(n-1)}{2}-1}\cdot
\\
&
\ \ \ \ \ \ \ \ \ \ \ \ \ \ \ \ \ \ \ \ \ \ \ \ \ \ \ 
\cdot
\int_{\sum_{i=1}^n\sum_{0\leqslant s^i\leqslant\tau^i}\,
\left[s^i+\mu_1^i+\cdots+\mu_i^i\right]b_{s^i}^i=1}\,\,
\prod_{i=1}^n\,
\bigg(
\sum_{0\leqslant s^i\leqslant\tau^i}\,
b_{s^i}^i
\bigg)^{\alpha_i'}
\endaligned
\]
Performing next the changes of variables: 
\[
c_{s^i}^i
:=
\big[s^i+\mu_1^i+\cdots+\mu_i^i\big]b_{s^i}^i,
\]
whence: 
\[
b_{s^i}^i
=
\frac{c_{s^i}^i}{s^i+\mu_1^i+\cdots+\mu_i^i},
\]
we transform the integral as: 
\[
\footnotesize
\aligned
{\sf S}_{n,\kappa,m}^{\alpha_1',\dots,\alpha_n'}\big(\mu_l^i\big)
&
=
m^{(\kappa+1)n-1}\cdot
\prod_{0\leqslant s^1\leqslant \tau^1}
\frac{1}{s^1+\mu_1^1}\,
\cdots
\prod_{0\leqslant s^n\leqslant\tau^n}
\frac{1}{s^n+\mu_1^n+\cdots+\mu_n^n}\,
\cdot
\\
&
\ \ \ \ \ \ \ \ 
\cdot
\int_{c_0^1+\cdots+c_{\tau^1}^1
+\cdots\cdots+
c_0^n+\cdots+c_{\tau^n}^n=1}\,\,
\prod_{i=1}^n\,
\bigg(
\sum_{0\leqslant s^i\leqslant\tau^i}\,
\frac{c_{s^i}^i}{s^i+\mu_1^i+\cdots+\mu_i^i}
\bigg)^{\alpha_i'}\,
dc'
+
\\
&
\ \ \ \ \ \ \ \ \ \ \ \ \ \ \ \ \ \ \ \ \ \ \ \ \ \ \ \ \ \ \ \ 
\ \ \ \ \ \ \ \ \ \ \ \ \ \ \ \ \ \ \ \ \ \ \ \ \ \ \ \ \ \ \ \ 
\ \ \ \ \ \ \ \ \ \ \ \ \ \ \ \ \ \ \ \ \ \ \ \ \ \ \ \ \ \ \ \ 
+
{\sf O}_{n,\kappa}\big(m^{(\kappa+1)n-2}\big),
\endaligned
\] 
where:
\[
dc'
:=
dc_0^1\cdots dc_{\tau^1}^1
\cdots\cdots
dc_0^n\cdots dc_{\tau^n}^n.
\]
Using the multinomial formula, we now expand the product of all the
$\alpha_i'$-th powers under the sum in the second line above:
\[
\footnotesize
\aligned
&
\prod_{i=1}^n\,
\bigg[
\sum_{0\leqslant s^i\leqslant\tau^i}\,
\frac{c_{s^i}^i}{s^i+\mu_1^i+\cdots+\mu_i^i}
\bigg]^{\alpha_i'}
=
\\
&
=
\prod_{i=1}^n\,
\bigg[
\sum_{q_0^i+\cdots+q_{\tau^i}^i=\alpha_i'}\,
\frac{\alpha_i'!}{q_0^i!\cdots q_{\tau^i}^i!}\,
\prod_{0\leqslant s^i\leqslant\tau^i}\,
\Big(
\frac{c_{s^i}^i}{s^i+\mu_1^i+\cdots+\mu_i^i}
\Big)^{q_{s^i}^i}
\bigg]
\endaligned
\]
\[
\aligned
&
=
\sum_{q_0^1+\cdots+q_{\tau^1}^1=\alpha_1'}\cdots
\sum_{q_0^n+\cdots+q_{\tau^n}^n=\alpha_n'}
\frac{\alpha_1'!}{q_0^1!\cdots q_{\tau^1}^1!}
\cdots
\frac{\alpha_n'!}{q_0^n!\cdots q_{\tau^n}^n!}\,\cdot
\\
&
\ \ \ \ \ \ \ 
\cdot
\prod_{0\leqslant s^1\leqslant\tau^1}\,
\Big(
\frac{c_{s^1}^1}{s^1+\mu_1^1}
\Big)^{q_{s^1}^1}
\cdots
\prod_{0\leqslant s^n\leqslant\tau^n}\,
\Big(
\frac{c_{s^n}^n}{s^n+\mu_1^n+\cdots+\mu_n^n}
\Big)^{q_{s^n}^n}
\endaligned
\]
\[
\aligned
&
=
\sum_{q_0^1+\cdots+q_{\tau^1}^1=\alpha_1'}\cdots
\sum_{q_0^n+\cdots+q_{\tau^n}^n=\alpha_n'}
\frac{\alpha_1'!}{q_0^1!\cdots q_{\tau^1}^1!}
\cdots
\frac{\alpha_n'!}{q_0^n!\cdots q_{\tau^n}^n!}\,
\cdot
\\
&
\ \ \ \ \ \ \ 
\cdot
\prod_{0\leqslant s^1\leqslant\tau^1}\,
\frac{1}{\big(s^1+\mu_1^1\big)^{q_{s^1}^1}}
\cdots
\prod_{0\leqslant s^n\leqslant\tau^n}\,
\frac{1}{\big(s^n+\mu_1^n+\cdots+\mu_n^n\big)^{q_{s^n}^n}}\,
\cdot
\\
&
\ \ \ \ \ \ \ 
\cdot
\prod_{0\leqslant s^1\leqslant\tau^1}\,
\big(
c_{s^1}^1
\big)^{q_{s^1}^1}
\cdots
\prod_{0\leqslant s^n\leqslant\tau^n}\,
\big(
c_{s^n}^n
\big)^{q_{s^n}^n}.
\endaligned
\]
After these expansions are done, in order to complete the computation
of ${\sf S}_{n, \kappa, m}^{ \alpha_1', \dots, \alpha_n'} \big( \mu_l^i
\big)$, we are left with the task of computing the integrals:
\[
\int_{c_0^1+\cdots+c_{\tau^1}^1
+\cdots\cdots+
c_0^n+\cdots+c_{\tau^n}^n=1}\,\,
\prod_{0\leqslant s^1\leqslant\tau^1}\,
\big(
c_{s^1}^1
\big)^{q_{s^1}^1}
\cdots
\prod_{0\leqslant s^n\leqslant\tau^n}\,
\big(
c_{s^n}^n
\big)^{q_{s^n}^n}\,
dc'.
\]
To this aim, we simply apply the elementary Lemma~\ref{j1-j2-jp}, 
and this then yields to us the desired value: 
\[
\small
\aligned
&
\int_{c_0^1+\cdots+c_{\tau^1}^1
+\cdots\cdots+
c_0^n+\cdots+c_{\tau^n}^n=1}\,
\prod_{0\leqslant s^1\leqslant\tau^1}\,
\big(
c_{s^1}^1
\big)^{q_{s^1}^1}
\cdots
\prod_{0\leqslant s^n\leqslant\tau^n}\,
\big(
c_{s^n}^n
\big)^{q_{s^n}^n}
=
\\
&
=
\frac{q_0^1!\,\cdots\,q_{\tau^1}^1!
\,\cdots\cdots\,
q_0^n!\,\cdots\,q_{\tau^n}^n!}{
(q_0^1+\cdots+q_{\tau^1}^1
+\cdots\cdots+
q_0^n+\cdots+q_{\tau^n}^n+
(1+\tau^1)+\cdots+(1+\tau^n)-1)!}
\endaligned
\]
\[
\aligned
&
=
\frac{q_0^1!\,\cdots\,q_{\tau^1}^1!
\,\cdots\cdots\,
q_0^n!\,\cdots\,q_{\tau^n}^n!}{
(\alpha_1'
+\cdots\cdots+
\alpha_n'
+n\kappa-\frac{n(n-1)}{2}-1)!}
\\
&
=
\frac{q_0^1!\,\cdots\,q_{\tau^1}^1!
\,\cdots\cdots\,
q_0^n!\,\cdots\,q_{\tau^n}^n!}{
\big((\kappa+1)n-1\big)!},
\endaligned
\]
since $q_0^i + \cdots + q_{ \tau^i}^i = \alpha_i'$ and since
$\alpha_1' + \cdots + \alpha_n' = \frac{ n ( n+1)}{ 2}$. 

Resuming what has been done, we therefore get:
\[
\small
\aligned
{\sf S}_{n,\kappa,m}^{\alpha_1',\dots,\alpha_n'}
\big(\mu_i^l\big)
&
=
m^{(\kappa+1)n-1}\cdot
\prod_{0\leqslant s^1\leqslant\tau^1}\,
\frac{1}{s^1+\mu_1^1}\,\cdots\,
\prod_{0\leqslant s^n\leqslant\tau^n}\,
\frac{1}{s^n+\mu_1^n+\cdots+\mu_n^n}\,\cdot
\\
&
\cdot
\sum_{q_0^1+\cdots+q_{\tau^1}^1=\alpha_1'}
\cdots
\sum_{q_0^n+\cdots+q_{\tau^n}^n=\alpha_n'}\,
\frac{\alpha_1'!}{\zero{q_0^1!\cdots q_{\tau^1}^1!}}
\cdots
\frac{\alpha_n'!}{\zero{q_0^n!\cdots q_{\tau^n}^n!}}
\cdot
\\
&
\cdot
\prod_{0\leqslant s^1\leqslant\tau^1}\,
\frac{1}{\big(s^1+\mu_1^1\big)^{q_{s^1}^1}}
\cdots
\prod_{0\leqslant s^n\leqslant\tau^n}\,
\frac{1}{\big(s^n+\mu_1^n+\cdots+\mu_n^n\big)^{q_{s^n}^n}}
\cdot
\\
&
\cdot
\frac{\zero{q_0^1!\,\cdots\,q_{\tau^1}^1!}
\,\cdots\cdots\,
\zero{q_0^n!\,\cdots\,q_{\tau^n}^n!}}{
\big((\kappa+1)n-1\big)!}
\,+
\\
&
\ \ \ \ \ \ \ \ \ \ \ \ \ \ \ \ \ \ \ \ \ \ \ \ \ \ \ \ \ \ \ \ 
\ \ \ \ \ \ \ \ \ \ \ \ \ \ \ \ \ \ \ \ \ \ \ \ \ \ \ \ \ \ \ \ 
\ \ \ \ \ \ \ \ \ \ \ \ \ \ \ \ \ \ \ 
+
{\sf O}_{n,\kappa}\big(m^{(\kappa+1)n-2}\big).
\endaligned
\]
The products of the factorials of the $q_{ s^i}^i$ disappear and a
reorganization gives:
\[
\footnotesize
\aligned
{\sf S}_{n,\kappa,m}^{\alpha_1',\dots,\alpha_n'}
\big(\mu_i^l\big)
&
=
\frac{m^{(\kappa+1)n-1}}{\big((\kappa+1)n-1\big)!}
\cdot
\\
&
\cdot
\prod_{0\leqslant s^1\leqslant\tau^1}\,
\frac{1}{s^1+\mu_1^1}\,\cdots\,
\prod_{0\leqslant s^n\leqslant\tau^n}\,
\frac{1}{s^n+\mu_1^n+\cdots+\mu_n^n}\,
\cdot
\\
&
\cdot
\alpha_1'!\,\cdots\,\alpha_n'!\,
\cdot
\sum_{q_0^1+\cdots+q_{\tau^1}^1=\alpha_1'}
\cdots
\sum_{q_0^n+\cdots+q_{\tau^n}^n=\alpha_n'}
\bigg(
\\
&
\bigg(
\prod_{0\leqslant s^1\leqslant\tau^1}\,
\frac{1}{\big(s^1+\mu_1^1\big)^{q_{s^1}^1}}
\cdots
\prod_{0\leqslant s^n\leqslant\tau^n}\,
\frac{1}{\big(s^n+\mu_1^n+\cdots+\mu_n^n\big)^{q_{s^n}^n}}
\bigg)
+
\\
&
\ \ \ \ \ \ \ \ \ \ \ \ \ \ \ \ \ \ \ \ \ \ \ \ \ \ \ \ \ \ \ \ 
\ \ \ \ \ \ \ \ \ \ \ \ \ \ \ \ \ \ \ \ \ \ \ \ \ \ \ \ \ \ \ \ 
\ \ \ \ \ \ \ \ \ \ \ \ \ \ \ \ \ \ \ 
+
{\sf O}_{n,\kappa}\big(m^{(\kappa+1)n-2}\big).
\endaligned
\]
Symbolically, instead of: 
\[
\prod_{0\leqslant s^1\leqslant\tau^1}\,
\frac{1}{s^1+\mu_1^1}\,
\prod_{0\leqslant s^2\leqslant\tau^2}\,
\frac{1}{s^2+\mu_1^2+\mu_2^2}
\,\cdots\,
\prod_{0\leqslant s^n\leqslant\tau^n}\,
\frac{1}{s^n+\mu_1^n+\cdots+\mu_n^n},
\]
we shall write without any risk of ambiguity: 
\[
\small
\aligned
\frac{1}{\kappa\cdots\mu_1^1}\,
\frac{1}{(\kappa+\mu_1^1)\cdots(\mu_2^2+\mu_1^2)}\,
\cdots\,
\frac{1}{(\kappa+\mu_{n-1}^{n-1}+\cdots+\mu_1^{n-1})
\cdots
(n+(n-1)+\cdots+1)},
\endaligned
\] 
the dots in the denominators meaning that one takes the product of all
integers, decreasingly, that are extant between the two written
extremal integers. In conclusion, we have established that:
\begin{equation}
\label{intermediate}
\boxed{
\footnotesize
\aligned
&
\sum_{{\sf YT}\,{\sf semi-standard}
\atop
{\sf weight}({\sf YT})=m}
\big(\ell_1({\sf YT})-\ell_2({\sf YT})\big)^{\alpha_1'}
\cdots
\big(\ell_{n-1}({\sf YT})-\ell_n({\sf YT})\big)^{\alpha_{n-1}'}
\big(\ell_n({\sf YT})\big)^{\alpha_n'}
=
\\
&
=
\frac{m^{(\kappa+1)n-1}}{\big((\kappa+1)n-1\big)!}
\sum_{\mu_l^i\in\nabla_{\!n,\kappa}}\,
\frac{N_{\mu_1^1}^\kappa}{
\kappa\cdots\mu_1^1}\,
\frac{N_{\mu_1^2,\mu_2^2}^{\mu_1^1,\kappa}}{
(\kappa+\mu_1^1)\cdots(\mu_2^2+\mu_1^2)}\,
\cdots\,
\\
&
\ \ \ \ \ \ \ \ \ \ \ \ \ \ \ \ \ \ \ \ \ \ \ \ \ \ \ \ \ \ \ \ \
\cdots
\frac{N_{\mu_1^n,\dots,\mu_{n-1}^n,\mu_n^n}^{
\mu_1^{n-1},\dots,\mu_{n-1}^{n-1},\kappa}}{
(\kappa+\mu_{n-1}^{n-1}+\cdots+\mu_1^{n-1})
\cdots
(\mu_n^n+\mu_{n-1}^n+\cdots+\mu_1^n)}\,\cdot
\\
&
\alpha_1'!\,\cdots\,\alpha_n'!
\cdot
\sum_{q_0^1+\cdots+q_{\tau^1}^1=\alpha_1'}\,
\cdots\,
\sum_{q_0^n+\cdots+q_{\tau^n}^n=\alpha_n'}
\\
&
\bigg(
\prod_{0\leqslant s^1\leqslant\tau^1}\,
\frac{1}{\big(s^1+\mu_1^1\big)^{q_{s^1}^1}}
\cdots
\prod_{0\leqslant s^n\leqslant\tau^n}\,
\frac{1}{\big(s^n+\mu_1^n+\cdots+\mu_{n-1}^n+\mu_n^n\big)^{q_{s^n}^n}}
\bigg)
+
\\
&
\ \ \ \ \ \ \ \ \ \ \ \ \ \ \ \ \ \ \ \ \ \ \ \ \ \ \ \ \ \ \ \ 
\ \ \ \ \ \ \ \ \ \ \ \ \ \ \ \ \ \ \ \ \ \ \ \ \ \ \ \ \ \ \ \ 
\ \ \ \ \ \ \ \ \ \ \ \ \ \ \ \ \ \ \ \ \ \ \ \ \ \ \ \ \ \ 
+
{\sf O}_{n,\kappa}\big(m^{(\kappa+1)n-2}\big),
\endaligned}
\end{equation}
where we recall for completeness that $\tau^i = \kappa +
\sum_{ l=1}^{ i-1}\, \mu_l^{ i-1} - \sum_{ l=1}^i\, \mu_l^i$
for convenient abbreviation.

\markleft{Jo\"el Merker}
\markright{\S9.~Number of tight paths in semi-standard Young tableaux}
\section{\bf Number of tight paths in semi-standard Young tableaux}
\label{Section-9}

\subsection{Summary} Thus, we are left with the task 
of computing or of 
majorating, for any $\alpha_1', \dots,
\alpha_n'$ with $\alpha_1' + \cdots + \alpha_n' = \frac{ n ( n+1)}{
2}$, sums:
\[
\footnotesize
\aligned
\square_{n,\kappa}^{\alpha_1',\dots,\alpha_n'}
:=
\sum_{\mu_l^i\in\nabla_{\!n,\kappa}}\,
&
(\kappa!)^n\cdot
\frac{N_{\mu_1^1}^\kappa}{
\kappa\cdots\mu_1^1}\,
\frac{N_{\mu_1^2,\mu_2^2}^{\mu_1^1,\kappa}}{
(\kappa+\mu_1^1)\cdots(\mu_2^2+\mu_1^2)}\,
\cdots\,
\\
&
\cdots
\frac{N_{\mu_1^n,\dots,\mu_{n-1}^n,\mu_n^n}^{
\mu_1^{n-1},\dots,\mu_{n-1}^{n-1},\kappa}}{
(\kappa+\mu_{n-1}^{n-1}+\cdots+\mu_1^{n-1})
\cdots
(\mu_n^n+\mu_{n-1}^n+\cdots+\mu_1^n)}\,\cdot
\\
&
\alpha_1'!\,\cdots\,\alpha_n'!
\cdot
\sum_{q_0^1+\cdots+q_{\tau^1}^1=\alpha_1'}\,
\cdots\,
\sum_{q_0^n+\cdots+q_{\tau^n}^n=\alpha_n'}
\\
&
\bigg(
\prod_{0\leqslant s^1\leqslant\tau^1}\,
\frac{1}{\big(s^1+\mu_1^1\big)^{q_{s^1}^1}}
\cdots
\prod_{0\leqslant s^n\leqslant\tau^n}\,
\frac{1}{\big(s^n+\mu_1^n+\cdots+\mu_{n-1}^n+\mu_n^n\big)^{q_{s^n}^n}}
\bigg)
\endaligned
\]
in which the weight $m$ has completely disappeared, while only the
dimension $n$ and the jet order $\kappa$ remain present.

At first, we would like to remind from Theorems~\ref{E-P-characteristic} 
and~\ref{majoration-Rousseau} that the basic numerical sums $\sum\,
M\cdot \ell^\alpha$ we must compute were in fact born when $\alpha_1
+ \cdots + \alpha_n$ is maximal equal to $\frac{ n ( n+1)}{ 2}$, after
rewriting in terms of $\ell_1 - \ell_2$, \dots, $\ell_{ n-1} - \ell_n$
and $\ell_n$ expressions of the form:
\[
\bigg[
\prod_{1\leqslant i<j\leqslant n}\,
(\ell_i-\ell_j)
\bigg]\cdot
\ell_1^{\beta_1}\cdots\ell_n^{\beta_n}
\]
that are multiple of the product of the $\ell_i - \ell_j$ with
$\beta_1 + \cdots + \beta_n = n$. Each $\ell_i - \ell_j$ then writes
as $\ell_i - \ell_{ i+1} + \cdots + \ell_{ j-1} - \ell_j$ with no
$\ell_n$ at all, and it follows that after the rewriting, the exponent
$\alpha_n'$ of $\ell_n$ is at most equal to $n$:
\begin{equation}
\label{beta-n-ell-ell}
\aligned
&
\bigg[
\prod_{1\leqslant i<j\leqslant n}\,
(\ell_i-\ell_j)
\bigg]\cdot
\ell_1^{\beta_1}\cdots\ell_{n-1}^{\beta_{n-1}}\ell_n^{\beta_n}
\\
&
\ \ \ \ \ \ \ \ \ \ 
\leqslant
{\sf Constant}_n\cdot
\sum_{\alpha_1'+\cdots+\alpha_{n-1}'+\alpha_n'=\frac{n(n+1)}{2}
\atop
\alpha_n'\leqslant n}\,
(\ell_1-\ell_2)^{\alpha_1'}\cdots\,(\ell_{n-1}-\ell_n)^{\alpha_{n-1}'}
\ell_n^{\alpha_n'}.
\endaligned
\end{equation}
Thus, the sums $\square_{ n, \kappa}^{ \alpha_1', \dots, \alpha_n'}$ 
we consider are such that $\alpha_1' + \cdots + \alpha_n ' = 
\frac{ n ( n+1)}{ 2}$ and $\alpha_n' \leqslant n$. 

\subsection{Logarithmic equivalents}
Next, we observe that for any integer $\alpha' \geqslant 1$,
as soon as $\tau \geqslant \alpha'$, one has:
\[
\aligned
&
\sum_{q_0+\cdots+q_\tau=\alpha'}\,
\frac{1}{(k)^{q_0}\cdots\,(k+\tau)^{q_\tau}}
=
\\
&
=
\frac{1}{\alpha'!}\,
\big[\log(k+\tau)-\log(k)\big]^{\alpha'}
+
{\sf O}_{\alpha'}
\big(
\big[\log(k+\tau)-\log(k)\big]^{\alpha'-1}
\big).
\endaligned
\]
If $\tau \leqslant \alpha' -1$, this sum is smaller than the written
power of a difference between two logarithms.  Since our goal now will
be to establish only an inequality of the form:
\begin{equation}
\label{square-n-kappa-inequality}
\boxed{
\square_{n,\kappa}^{\alpha_1',\dots,\alpha_n'}
\leqslant
{\sf Constant}_n\cdot
(\log\kappa)^{\alpha_n'}}\,,
\end{equation}
in which no particular knowledge about the ${\sf Constant}_n$ will
be required, again with $\alpha_1' + \cdots + \alpha_n' = \frac{ n (
n+1)}{ 2}$ and with $\alpha_n' \leqslant n$, it will even suffice to
observe that the last two lines in the definition of $\square_{ n,
\kappa}^{ \alpha_1', \dots, \alpha_n'}$ enjoy a majoration of the
sort:
\[
\footnotesize
\aligned
&
\alpha_1'!\,\alpha_2'!\,\cdots\,\alpha_n'!
\cdot
\sum_{q_0^1+\cdots+q_{\tau^1}^1=\alpha_1'}\,
\sum_{q_0^2+\cdots+q_{\tau^2}^2=\alpha_2'}\,
\cdots\,
\sum_{q_0^n+\cdots+q_{\tau^n}^n=\alpha_n'}\,\,\,
\prod_{0\leqslant s^1\leqslant\tau^1}\,
\frac{1}{\big(s^1+\mu_1^1\big)^{q_{s^1}^1}}\,
\\
&
\prod_{0\leqslant s^2\leqslant\tau^2}\,
\frac{1}{\big(s^2+\mu_1^2+\mu_2^2\big)^{q_{s^2}^2}}
\cdots
\prod_{0\leqslant s^n\leqslant\tau^n}\,
\frac{1}{\big(s^n+\mu_1^n+\cdots+\mu_{n-1}^n+\mu_n^n\big)^{q_{s^n}^n}}
\leqslant
\\
&
\leqslant
{\sf Constant}_n\cdot
\big[\log(\kappa)-\log(\mu_1^1)\big]^{\alpha_1'}
\big[\log(\kappa+\mu_1^1)-\log(\mu_2^2+\mu_1^2)\big]^{\alpha_2'}
\cdots
\\
&
\ \ \ \ \
\cdots
\big[
\log(\kappa+\mu_{n-1}^{n-1}+\cdots+\mu_1^{n-1})
-
\log(n+\cdots+2+1)
\big]^{\alpha_n'}.
\endaligned
\]
Consequently, we are left with establishing the following
proposition.

\begin{Proposition}
Let $\alpha_1', \dots, \alpha_n' \in \N$ with $\alpha_1' + \cdots +
\alpha_n' = \frac{ n( n+1)}{ 2}$ and with $\alpha_n' \leqslant
n$. Then for $\kappa \geqslant n$, one has the majoration:
\[
\footnotesize
\aligned
\widetilde{\square}_{n,\kappa}^{\alpha_1',\dots,\alpha_n'}
:=
\sum_{\mu_l^i\in\nabla_{\!n,\kappa}}\,
&
(\kappa!)^n\cdot
\frac{N_{\mu_1^1}^\kappa}{
\kappa\cdots\mu_1^1}\,
\frac{N_{\mu_1^2,\mu_2^2}^{\mu_1^1,\kappa}}{
(\kappa+\mu_1^1)\cdots(\mu_2^2+\mu_1^2)}\,
\cdots\,
\\
&
\cdots
\frac{N_{\mu_1^n,\dots,\mu_{n-1}^n,\mu_n^n}^{
\mu_1^{n-1},\dots,\mu_{n-1}^{n-1},\kappa}}{
(\kappa+\mu_{n-1}^{n-1}+\cdots+\mu_1^{n-1})
\cdots
(\mu_n^n+\mu_{n-1}^n+\cdots+\mu_1^n)}\,\cdot
\\
&
\big[\log(\kappa)-\log(\mu_1^1)\big]^{\alpha_1'}
\big[\log(\kappa+\mu_1^1)-\log(\mu_2^2+\mu_1^2)\big]^{\alpha_2'}
\cdots
\\
&
\ \ \ \ \
\cdots
\big[
\log(\kappa+\mu_{n-1}^{n-1}+\cdots+\mu_1^{n-1})
-
\log(n+\cdots+2+1)
\big]^{\alpha_n'}
\leqslant
\\
&
\!\!\!\!\!\!\!\!\!\!\!\!\!\!\!\!\!\!\!\!\!\!\!\!\!
\leqslant
{\sf Constant}_n\cdot(\log\kappa)^{\alpha_n'}.
\endaligned
\]
\end{Proposition}

\subsection{Tight paths} 
According to Proposition~\ref{proposition-maximal} and
to the definition made in~\thetag{ \ref{definition-N}}, the integer $N_{
\mu_1^1}^\kappa$ denotes the number of tight paths from the column
$\mu_1^1$ to the column $\kappa$, hence it is equal to $1$.  When the
dimension $n$ is equal to $2$, one may show that:
\[
N_{1,2}^{\mu_1^1,\kappa}
=
{\textstyle{\frac{(\kappa+\mu_1^1)!}{(\kappa-3)!\,(\mu_1^1-1)!}}}
-
{\textstyle{\frac{(\kappa+\mu_1^1-4)!}{(\kappa-1)!\,(\mu_1^1-3)!}}}.
\]
In higher dimensions, the exact computation of the numbers $N_{
\mu_1^2, \mu_2^2}^{ \mu_1^1, \kappa}$, $N_{ \mu_1^3, \mu_2^3,
\mu_3^3}^{ \mu_1^2, \mu_2^2, \kappa}$, \dots may certainly be done
and it involves only differences of multinomial coefficients, but very
many cases are to be considered according to certain inequalities
between the $\mu_i^j$. However, after some explorations, it appears
that in order to get the majoration claimed by the proposition, it
suffices to {\em majorate} these numbers uniformly
as follows.

\subsection{Majoration of the tight path numbers
$N_{ \mu_1^i, \dots, \mu_{i-1}^i, \mu_i^i}^{ \mu_1^{ i-1}, \dots,
\mu_{ i-1}^{ i-1}, \kappa}$} By definition, $N_{ \mu_1^i, \dots,
\mu_{i-1}^i, \mu_i^i}^{ \mu_1^{ i-1}, \dots, \mu_{ i-1}^{ i-1},
\kappa}$ counts the number of strictly increasing tight paths from the
column $\big[ \mu_1^i \,\, \cdots \, \mu_{ i-1}^i \,\, \mu_i^i
\big]^{\rm transposed}$ to the column $\big[ \mu_1^{ i-1} \,\, \cdots
\,\, \mu_{i-1}^{i-1} \,\, \kappa \big]^{\rm transposed}$ in the
$i$-dimensional lattice $\N^i$ with the supplementary constraint that
at each point $\big[ \gamma_1^i ( s^i) \,\, \cdots \,\, \gamma_{i-1}^i
( s^i) \,\, \gamma_i ( s^i) \big]^{ \rm transposed}$ of the path, the
inequalities $\gamma_1^i ( s^i) < \cdots < \gamma_{ i-1}^i ( s^i) <
\gamma_i ( s^i)$ (strict increase inside columns, downward) must be
satisfied.

If we relax this last constraint, there are clearly more paths.  But
the number of strictly increasing paths in a complete lattice is
elementarily computed.  Thus we deduce that:
\[
N_{\mu_1^i,\dots,\mu_{i-1}^i,\mu_i^i}^{
\mu_1^{i-1},\dots,\mu_{i-1}^{i-1},\kappa}
\leqslant
{\textstyle{
\frac{(\mu_1^{i-1}-\mu_1^i+\cdots+\mu_{i-1}^{i-1}-\mu_{i-1}^i
+\kappa-\mu_i^i)!}{
(\mu_1^{i-1}-\mu_1^i)!\,\cdots\,
(\mu_{i-1}^{i-1}-\mu_{i-1}^i)!\,
(\kappa-\mu_i^i)!}}}.
\]

\subsection{Removal of $\alpha_n'$} On the other hand, 
for any choice of $\mu_{ n-1}^{n-1}, \dots, \mu_1^{ n-1}$ as in 
Proposition~\ref{proposition-maximal}, the last difference
between logarithms:
\[
\big[
\log(\kappa+\mu_{n-1}^{n-1}+\cdots+\mu_1^{n-1})
-
\log(n+\cdots+2+1)
\big]^{\alpha_n'}
\]
enjoys, when $\kappa \gg n$, a doubly controlling inequality of the
form:
\[
\frac{1}{{\sf C}_n}\cdot
[\log(\kappa)]^{\alpha_n'}
\leqslant
\big[
\log(\kappa+\mu_{n-1}^{n-1}+\cdots+\mu_1^{n-1})
-
\log(n+\cdots+2+1)
\big]^{\alpha_n'}
\leqslant
{\sf C}_n\cdot
[\log(\kappa)]^{\alpha_n'}
\]
where the constant ${\sf C}_n > 1$ can be chosen arbitrarily close to
$1$ provided that $\kappa \geqslant \kappa_{{\sf C}_n} \gg n$ is large
enough.  Consequently, in order to establish the inequality~\thetag{
\ref{square-n-kappa-inequality}}, it suffices now to establish
the following concrete proposition, in which $\alpha_n' = 0$
has disappeared.

\begin{Proposition}
Let $\alpha_1', \dots, \alpha_{ n-1}' \in \N$ with $\alpha_1' + \cdots
+ \alpha_{ n-1}' \leqslant \frac{ n( n+1)}{ 2}$ and assume $\kappa \geqslant
n$. Then the following sum is bounded independently of $\kappa$:
\[
\aligned
&
\Delta_{n,\kappa}^{
\alpha_1',\dots,\alpha_{n-1}',0}
:=
\sum_{\mu_l^i\in\nabla_{n,\kappa}}\,
(\kappa!)^n\cdot
1\,
{\textstyle{\frac{(\mu_1^1-1)!}{\kappa!}}}
\cdot
{\textstyle{\frac{(\mu_1^1-\mu_1^2+\kappa-\mu_2^2)!}{
(\mu_1^1-\mu_1^2)!\,(\kappa-\mu_2^2)!}\,
\frac{(\mu_2^2+\mu_1^2-1)!}{(\kappa+\mu_1^1)!}}}
\cdot
\\
&
\ \ \ \ \ \ \ \ \ \ \ \ \ \ \ \ \ \ \ \ \ \ \ \ \ \
\cdot
{\textstyle{\frac{(\mu_1^2-\mu_1^3+\mu_2^2-\mu_2^3+\kappa-\mu_3^3)!}{
(\mu_1^2-\mu_1^3)!\,(\mu_2^2-\mu_2^3)!\,(\kappa-\mu_3^3)!}\,
\frac{(\mu_3^3+\mu_2^3+\mu_1^3-1)!}{(\kappa+\mu_2^2+\mu_1^2)!}}}
\cdots
\\
&
\ \ \ \ \ \ \ \ \ \ \ \ \ \ \ \ \ \ \ \ \ \ \ \ \ \
\cdots
{\textstyle{
\frac{(\mu_1^{n-2}-\mu_1^{n-1}+\cdots+\mu_{n-2}^{n-2}-\mu_{n-2}^{n-1}
+\kappa-\mu_{n-1}^{n-1})!}{
(\mu_1^{n-2}-\mu_1^{n-1})!\,\cdots\,(\mu_{n-2}^{n-2}-\mu_{n-2}^{n-1})!\,
(\kappa-\mu_{n-1}^{n-1})!}\,
\frac{(\mu_{n-1}^{n-1}+\mu_{n-2}^{n-1}+\cdots+\mu_1^{n-1}-1)!}{
(\kappa+\mu_{n-2}^{n-2}+\cdots+\mu_1^{n-2})!}}}
\cdot
\\
&
\ \ \ \ \ \ \ \ \ \ \ \ \ \ \ \ \ \ \ \ \ \ \ \ \ \
\cdot
{\textstyle{
\frac{(\mu_1^{n-1}-\mu_1^n+\cdots+\mu_{n-1}^{n-1}-\mu_{n-1}^n
+\kappa-\mu_n^n)!}{
(\mu_1^{n-1}-\mu_1^n)!\,\cdots\,(\mu_{n-1}^{n-1}-\mu_{n-1}^n)!\,
(\kappa-\mu_n^n)!}\,
\frac{(\mu_n^n+\mu_{n-1}^n+\cdots+\mu_1^n-1)!}{
(\kappa+\mu_{n-1}^{n-1}+\cdots+\mu_1^{n-1})!}}}
\cdot
\\
&
\ \ \ \ \ \ 
\cdot
\big[\log(\kappa)-\log(\mu_1^1)\big]^{\alpha_1'}
\cdot
\big[\log(\kappa+\mu_1^1)-\log(\mu_2^2+\mu_1^2)\big]^{\alpha_2'}
\cdot
\\
&
\ \ \ \ \ \ \
\cdot
\big[\log(\kappa+\mu_2^2+\mu_1^2)
-
\log(\mu_3^3+\mu_2^3+\mu_1^3)\big]^{\alpha_3'}
\cdots
\\
&
\ \ \ \ \ \ \
\cdots
\big[
\log(\kappa+\mu_{n-2}^{n-2}+\cdots+\mu_1^{n-2})
-
\log(\mu_{n-1}^{n-1}+\mu_{n-2}^{n-1}+\cdots+\mu_1^{n-1})
\big]^{\alpha_{n-1}'}
\leqslant
\\
&
\leqslant
{\sf Constant}_n.
\endaligned
\]
\end{Proposition}

\markleft{Jo\"el Merker}
\markright{\S10.~Bounded behavior of plurilogarithmic sums}
\section{\bf Bounded behavior of plurilogarithmic sums}
\label{Section-10}

\subsection{Simplifying the kernel}
To begin with, disregarding the logarithmic factors, or equivalently,
considering that $\alpha_1 ' = \alpha_2' = \alpha_3 ' = \cdots =
\alpha_{ n-1}' = 0$, we observe that the rational factor simplifies
a bit (a
factorial $\kappa!$ disappears) and can be majorated as follows: 
\[
\aligned
&
\zero{\kappa!}\,(\kappa!)^{n-2}\,
\underline{\kappa!}_{\bigcirc\!\!\!\!{\sf a}}
\cdot
1\,
{\textstyle{\frac{(\mu_1^1-1)!}{\zero{\kappa!}}}}
\cdot
{\textstyle{\frac{(\mu_1^1-\mu_1^2+\kappa-\mu_2^2)!}{
(\mu_1^1-\mu_1^2)!\,(\kappa-\mu_2^2)!}\,
\frac{(\mu_2^2+\mu_1^2-1)!}{(\kappa+\mu_1^1)!}}}
\cdots
\\
&
\ \ \ \ \ \ \ \ \ \ \ \ \ \ \ \ \ \ \ \ \ \ \ \ \ \
\cdot
{\textstyle{\frac{(\mu_1^2-\mu_1^3+\mu_2^2-\mu_2^3+\kappa-\mu_3^3)!}{
(\mu_1^2-\mu_1^3)!\,(\mu_2^2-\mu_2^3)!\,(\kappa-\mu_3^3)!}\,
\frac{(\mu_3^3+\mu_2^3+\mu_1^3-1)!}{(\kappa+\mu_2^2+\mu_1^2)!}}}
\cdots
\endaligned
\]
\[
\footnotesize
\aligned
&
\ \ \ \ \ \ \ \ \ \ \ \ \ \ \ \ \ \ \ \ \ \ \ \ \ \
\cdots
{\textstyle{
\frac{(\mu_1^{n-2}-\mu_1^{n-1}+\cdots+\mu_{n-2}^{n-2}-\mu_{n-2}^{n-1}
+\kappa-\mu_{n-1}^{n-1})!}{
(\mu_1^{n-2}-\mu_1^{n-1})!\,\cdots\,(\mu_{n-2}^{n-2}-\mu_{n-2}^{n-1})!\,
(\kappa-\mu_{n-1}^{n-1})!}\,
\frac{(\mu_{n-1}^{n-1}+\mu_{n-2}^{n-1}+\cdots+\mu_1^{n-1}-1)!}{
(\kappa+\mu_{n-2}^{n-2}+\cdots+\mu_1^{n-2})!}}}
\cdot
\\
&
\ \ \ \ \ \ \ \ \ \ \ \ \ \ \ \ \ \ \ \ \ \ \ \ \ \
\cdot
{\textstyle{
\frac{
\underline{(\mu_1^{n-1}-1+\cdots+\mu_{n-1}^{n-1}-(n-1)
+\kappa-n)!
}_{\bigcirc\!\!\!\!{\sf c}}
}{
(\mu_1^{n-1}-1)!\,\cdots\,(\mu_{n-1}^{n-1}-(n-1))!\,
\underline{(\kappa-n)!}_{\bigcirc\!\!\!\!{\sf a}}}\,
\frac{
\underline{(n+(n-1)+\cdots+1-1)!}_{
\bigcirc\!\!\!\!{\sf b}}
}{
\underline{(\kappa+\mu_{n-1}^{n-1}+\cdots+\mu_1^{n-1})!
}_{\bigcirc\!\!\!\!{\sf c}}
}}}
\leqslant
\endaligned
\]
\[
\footnotesize
\aligned
&
\leqslant
{\sf Constant}_n\cdot
(\kappa!)^{n-2}\cdot
(\mu_1^1-1)!\cdot
{\textstyle{\frac{(\mu_1^1-\mu_1^2+\kappa-\mu_2^2)!}{
(\mu_1^1-\mu_1^2)!\,(\kappa-\mu_2^2)!}\,
\frac{(\mu_2^2+\mu_1^2-1)!}{(\kappa+\mu_1^1)!}}}
\cdots
\\
&
\ \ \ \ \ \ \ \ \ \ \ \ \ \ \ \ \ \ \ \ \ \ \ \ \ \
\cdot
{\textstyle{\frac{(\mu_1^2-\mu_1^3+\mu_2^2-\mu_2^3+\kappa-\mu_3^3)!}{
(\mu_1^2-\mu_1^3)!\,(\mu_2^2-\mu_2^3)!\,(\kappa-\mu_3^3)!}\,
\frac{(\mu_3^3+\mu_2^3+\mu_1^3-1)!}{(\kappa+\mu_2^2+\mu_1^2)!}}}
\cdots
\endaligned
\]
\[
\footnotesize
\aligned
&
\ \ \ \ \ \ \ \ \ \ \ \ \ \ \ \ \ \ \ \ \ \ \ \ \ \
\cdots
{\textstyle{
\frac{(\mu_1^{n-2}-\mu_1^{n-1}+\cdots+\mu_{n-2}^{n-2}-\mu_{n-2}^{n-1}
+\kappa-\mu_{n-1}^{n-1})!}{
(\mu_1^{n-2}-\mu_1^{n-1})!\,\cdots\,(\mu_{n-2}^{n-2}-\mu_{n-2}^{n-1})!\,
(\kappa-\mu_{n-1}^{n-1})!}\,
\frac{(\mu_{n-1}^{n-1}+\mu_{n-2}^{n-1}+\cdots+\mu_1^{n-1}-1)!}{
(\kappa+\mu_{n-2}^{n-2}+\cdots+\mu_1^{n-2})!}}}
\cdot
\\
&
\ \ \ \ \ \ \ \ \ \ \ \ \ \ \ \ \ \ \ \ \ \ \ \ \ \
\cdot
{\textstyle{
\frac{1}{(\mu_1^{n-1}-1)!\cdots\,(\mu_{n-1}^{n-1}-(n-1))!}}}
\cdot
\\
&
\ \ \ \ \ \ \ \ \ \ \ \ \ \ \ \ \ \ \ \ \ \ \ \ \ \
\cdot
{\textstyle{
\frac{\kappa(\kappa-1)\cdots\,(\kappa-n+1)}{
(\kappa+\mu_{n-1}^{n-1}+\cdots+\mu_1^{n-1})\cdots\,
(\kappa+\mu_{n-1}^{n-1}+\cdots+\mu_1^{n-1}
-\frac{n(n+1)}{2}+1)}}},
\endaligned
\]
since the two pairs of terms underlined with ${\bigcirc\!\!\!\!{\sf
a}}$\, and ${\bigcirc\!\!\!\!{\sf c}}$\, appended can be put at the
end and simplified, while the pair of terms with
${\bigcirc\!\!\!\!{\sf b}}$\, appended, equal to the factorial
$(\frac{ n ( n+1)}{ 2}-1)!$, may be considered as just a ${\sf
Constant}_n$. But now, the last line is controlled as follows: 
\[
{\sf C}_n^{-1}\,\kappa^{-\frac{n(n-1)}{2}}
\leqslant
{\textstyle{
\frac{\kappa(\kappa-1)\cdots\,(\kappa-n+1)}{
(\kappa+\mu_{n-1}^{n-1}+\cdots+\mu_1^{n-1})\cdots\,
(\kappa+\mu_{n-1}^{n-1}+\cdots+\mu_1^{n-1}
-\frac{n(n+1)}{2}+1)}}}
\leqslant
{\sf C}_n\,\kappa^{-\frac{n(n-1)}{2}},
\]
for some constant ${\sf C}_n > 1$. Consequently, we are
reduced to the following proposition in which
we expressly mention that $\alpha_1' + \cdots + \alpha_{n-1}' \leqslant
\frac{ n(n+1)}{2}$, a restriction originating from 
our main concern, 
but its conclusion also holds true generally for
any integers $\alpha_1', \dots, \alpha_n' \in \N$, 
with a right-hand side majorating ${\sf Constant}_{\alpha_1', 
\dots, \alpha_n'}$ depending on them.

\begin{Proposition}
Let $\alpha_1', \dots, \alpha_{ n-1}' \in \N$ with $\alpha_1' + \cdots
+ \alpha_{ n-1}' \leqslant \frac{ n( n+1)}{ 2}$ and assume $\kappa \geqslant
n$. Then the following sum is bounded independently of $\kappa$:
\[
\aligned
&
{\sf K}_{\alpha_1',\dots,\alpha_{n-1}'}^n(\kappa)
:=
\sum_{\mu_l^i\in\nabla_{n,\kappa}}\,
\frac{1}{\kappa^{\frac{n(n-1)}{2}}}\cdot
(\kappa!)^{n-2}\cdot
(\mu_1^1-1)!
\cdot
\\
&
\ \ \ \ \ \ \ \ \ \ \ \ \ \ \ \
\cdot
{\textstyle{\frac{(\mu_1^1-\mu_1^2+\kappa-\mu_2^2)!}{
(\mu_1^1-\mu_1^2)!\,(\kappa-\mu_2^2)!}\,
\frac{(\mu_2^2+\mu_1^2-1)!}{(\kappa+\mu_1^1)!}}}
\cdot
{\textstyle{\frac{(\mu_1^2-\mu_1^3+\mu_2^2-\mu_2^3+\kappa-\mu_3^3)!}{
(\mu_1^2-\mu_1^3)!\,(\mu_2^2-\mu_2^3)!\,(\kappa-\mu_3^3)!}\,
\frac{(\mu_3^3+\mu_2^3+\mu_1^3-1)!}{(\kappa+\mu_2^2+\mu_1^2)!}}}
\cdots
\\
&
\ \ \ \ \ \ \ \ \ \ \ \ \ \ \ \ 
\cdots
{\textstyle{
\frac{(\mu_1^{n-2}-\mu_1^{n-1}+\cdots+\mu_{n-2}^{n-2}-\mu_{n-2}^{n-1}
+\kappa-\mu_{n-1}^{n-1})!}{
(\mu_1^{n-2}-\mu_1^{n-1})!\,\cdots\,(\mu_{n-2}^{n-2}-\mu_{n-2}^{n-1})!\,
(\kappa-\mu_{n-1}^{n-1})!}\,
\frac{(\mu_{n-1}^{n-1}+\mu_{n-2}^{n-1}+\cdots+\mu_1^{n-1}-1)!}{
(\kappa+\mu_{n-2}^{n-2}+\cdots+\mu_1^{n-2})!}}}
\cdot
\\
&
\ \ \ \ \ \ \ \ \ \ \ \ \ \ \ \ \ \ \ \ \ \ \ \ \ \
\ \ \ \ \ \ \ \ \ \ \ 
\cdot
{\textstyle{
\frac{1}{(\mu_1^{n-1}-1)!\cdots\,(\mu_{n-1}^{n-1}-(n-1))!}}}
\cdot
\\
&
\ \ \ \ \ 
\cdot
\big[\log(\kappa)-\log(\mu_1^1)\big]^{\alpha_1'}
\cdot
\big[\log(\kappa+\mu_1^1)-\log(\mu_2^2+\mu_1^2)\big]^{\alpha_2'}
\cdot
\\
&
\ \ \ \ \ \
\cdot
\big[\log(\kappa+\mu_2^2+\mu_1^2)
-
\log(\mu_3^3+\mu_2^3+\mu_1^3)\big]^{\alpha_3'}
\cdots
\\
&
\ \ \ \ \ \  
\cdots
\big[
\log(\kappa+\mu_{n-2}^{n-2}+\cdots+\mu_1^{n-2})
-
\log(\mu_{n-1}^{n-1}+\mu_{n-2}^{n-1}+\cdots+\mu_1^{n-1})
\big]^{\alpha_{n-1}'}
\leqslant
\\
&
\leqslant
{\sf Constant}_n.
\endaligned
\]
\end{Proposition}

Here, the summation $\sum_{ \mu_l^i \in \nabla_{ n, \kappa}}$
holds for $\mu_l^i$ satisfying the two collections of
inequalities~\thetag{ \ref{strict-mu}} and~\thetag{
\ref{weak-mu}}, and we may expand it symbolically 
using two symbols $\Sigma$ in order to notify well 
these two conditions:
\[
\sum_{\mu_l^i\in\nabla_{n,\kappa}}
\equiv
\sum_{
\substack{
1\leqslant\mu_1^1<\kappa
\\
1\leqslant\mu_1^2<\mu_2^2<\kappa
\\
1\leqslant\mu_1^3<\mu_2^3<\mu_3^3<\kappa
\\
\cdots\cdots\cdots\cdots\cdots\cdots\cdots\cdots\cdots\cdots
\cdots\cdots
\\
1\leqslant\mu_1^{n-2}<\mu_2^{n-2}<\mu_3^{n-2}
<\cdots<\mu_{n-2}^{n-2}<\kappa
\\
1\leqslant\mu_1^{n-1}<\mu_2^{n-1}<\mu_3^{n-1}<\cdots<
\mu_{n-2}^{n-1}<\mu_{n-1}^{n-1}<\kappa
}}
\sum_{
\substack{
\mu_1^1\geqslant\mu_1^2\geqslant\mu_1^3
\geqslant\cdots\geqslant
\mu_1^{n-2}\geqslant\mu_1^{n-1}
\\
\mu_2^2\geqslant\mu_2^3
\geqslant\cdots\geqslant
\mu_2^{n-2}\geqslant\mu_2^{n-1}
\\
\mu_3^3
\geqslant\cdots\geqslant
\mu_3^{n-2}\geqslant\mu_3^{n-1}
\\
\cdots\cdots\cdots\cdots\cdots\cdots\cdots
\\
\mu_{n-2}^{n-2}\geqslant\mu_{n-2}^{n-1}
}}
\]

\proof
In dimension $n = 2$, the sum writes:
\[
{\sf K}_{\alpha_1'}^2(\kappa)
=
\sum_{1\leqslant\mu_1^1<\kappa}\,
\frac{1}{\kappa}\cdot
\zero{(\mu_1^1-1)!}\cdot
\frac{1}{\zero{(\mu_1^1-1)!}}\cdot
\big[\log(\kappa)-\log(\mu_1^1)\big]^{\alpha_1'},
\]
and is seen to be an approximation of the Riemann integral:
\[
\int_0^1\,\big(-\log(x)\big)^{\alpha_1'}
=
\alpha_1'!,
\]
which is finite. 

In dimensions $n = 3$ and $n = 4$, the sum writes:
\[
\footnotesize
\aligned
{\sf K}_{\alpha_1',\alpha_2'}^3(\kappa)
&
=
\sum_{1\leqslant\mu_1^1<\kappa
\atop
1\leqslant\mu_1^2<\mu_2^2<\kappa}\,
\sum_{\mu_1^1\geqslant\mu_1^2}\,
\frac{1}{\kappa^3}\cdot
\kappa!\cdot
(\mu_1^1-1)!\cdot
\frac{(\mu_1^1-\mu_1^2+\kappa-\mu_2^2)!
}{
(\mu_1^1-\mu_1^2)!\,(\kappa-\mu_2^2)!}\,
\frac{(\mu_2^2+\mu_1^2-1)!}{(\kappa+\mu_1^1)!}
\cdot
\\
&
\cdot
\frac{1}{(\mu_1^2-1)!\,(\mu_2^2-2)!}\cdot
\big[\log\kappa-\log\mu_1^1\big]^{\alpha_1'}\,
\big[\log(\kappa+\mu_1^1)-\log(\mu_2^2+\mu_1^2)\big]^{\alpha_2'},
\endaligned
\]
and:
\[
\footnotesize
\aligned
{\sf K}_{\alpha_1',\alpha_2',\alpha_3'}^4(\kappa)
=
\sum_{
\substack{
1\leqslant\mu_1^1<\kappa
\\
1\leqslant\mu_1^2<\mu_2^2<\kappa
\\
1\leqslant\mu_1^3<\mu_2^3<\mu_3^3<\kappa
}}\,
&
\sum_{
\mu_1^1\geqslant\mu_1^2\geqslant\mu_1^3
\atop
\mu_2^2\geqslant\mu_2^3}\,
\frac{1}{\kappa^6}\cdot
(\kappa!)^2\cdot
(\mu_1^1-1)!\cdot
\\
&
\cdot
\frac{(\mu_1^1-\mu_1^2+\kappa-\mu_2^2)!
}{
(\mu_1^1-\mu_1^2)!\,(\kappa-\mu_2^2)!}\,
\frac{(\mu_2^2+\mu_1^2-1)!}{(\kappa+\mu_1^1)!}
\cdot
\\
&
\cdot
\frac{(\mu_1^2-\mu_1^3+\mu_2^2-\mu_2^3+\kappa-\mu_3^3)!}{
(\mu_1^2-\mu_1^3)!\,(\mu_2^2-\mu_2^3)!\,(\kappa-\mu_3^3)!}\,
\frac{(\mu_3^3+\mu_2^3+\mu_1^3-1)!}{(\kappa+\mu_2^2+\mu_1^2)!}
\cdot
\\
&
\cdot
\frac{1}{(\mu_1^3-1)!\,(\mu_2^3-2)!\,(\mu_3^3-3)!}
\cdot
\\
&
\cdot
\big[\log\kappa-\log\mu_1^1\big]^{\alpha_1'}\,
\big[\log(\kappa+\mu_1^1)-\log(\mu_2^2+\mu_1^2)\big]^{\alpha_2'}
\\
&
\cdot
\big[\log(\kappa+\mu_2^2+\mu_1^2)
-
\log(\mu_3^3+\mu_2^3+\mu_1^3)\big]^{\alpha_3'}.
\endaligned
\] 

The number of indiced quantities increasing much with the
dimension, let us treat in great details the dimension $n = 3$
of the general majoration:
\[
\small
\aligned
\frac{1}{\kappa^3}\,
\sum_{\lambda=1}^\kappa\,
\sum_{\mu=\lambda}^\kappa\,
\sum_{\nu=\lambda+1}^\kappa\,
&
\frac{\kappa!\,(\mu-1)!}{(\kappa+\mu)!}\,
\frac{(\kappa-\nu+\mu-\lambda)!}{(\kappa-\nu)!\,(\mu-\lambda)!}\,
\frac{(\nu+\lambda-1)!}{(\lambda-1)!\,(\nu-2)!}\,
\\
&
\bigg[
\log\,\frac{\kappa}{\mu}
\bigg]^\alpha\,
\bigg[
\log\,\frac{\kappa+\mu}{\lambda+\nu}
\bigg]^\beta
\,\leqslant\,
{\sf Constant}_{\alpha,\beta},
\endaligned
\]
using the somehow lighter unindiced notations:
\[
\lambda\,\equiv\,\mu_1^1,\ \ \ \ \ \ \ \ \ \ 
\mu\,\equiv\,\mu_1^2,\ \ \ \ \ \ \ \ \ \ 
\nu\,\equiv\,\mu_2^2.
\]
In fact, a similar majoration will hold true with just one small
modification:
\[
\small
\aligned
\frac{1}{\kappa^3}\,
\sum_{\lambda=1}^\kappa\,
\sum_{\mu=\lambda}^\kappa\,
\sum_{\nu=2}^\kappa\,
&
\frac{\kappa!\,(\mu-1)!}{(\kappa+\mu)!}\,
\frac{(\kappa-\nu+\mu-\lambda)!}{(\kappa-\nu)!\,(\mu-\lambda)!}\,
\frac{(\nu+\lambda-1)!}{(\lambda-1)!\,(\nu-2)!}\,
\\
&
\bigg[
\log\,\frac{\kappa}{\mu}
\bigg]^\alpha\,
\bigg[
\log\,\frac{\kappa+\mu}{\lambda+\nu}
\bigg]^\beta
\,\leqslant\,
{\sf Constant}_{\alpha,\beta},
\endaligned
\]
in which the third, last $\sum_{\nu = \lambda + 1}^\kappa$
is replaced by the larger sum $\sum_{ \nu = 2}^\kappa$, all
added terms in $\sum_{ \lambda = 1}^\kappa\,
\sum_{\mu = \lambda}^\kappa\, \sum_{ 2 \leqslant \nu \leqslant 
\lambda}$ being $\geqslant 0$, for the two logarithms act on rational
numbers all $\geqslant 1$.

Next, when $\beta = 1$, the second logarithm is classically
majorated by:
\[
\aligned
\log\,
\frac{\kappa+\mu}{\lambda+\nu}
&
=
\log\bigg(
1
+
\frac{\kappa-\nu+\mu-\lambda}{\lambda+\nu}
\bigg)
\\
&
\leqslant
\frac{\kappa-\nu+\mu-\lambda}{\lambda+\nu}
\\
&
\leqslant
\frac{\kappa-\nu+\mu-\lambda+1}{\lambda+\nu-1}
\endaligned
\]
(this latter trivial majoration being\,\,---\,\,a bit 
trickily\,\,---\,\,useful below), 
and since more generally for every $\beta \geqslant 1$, one has:
\[
\big[
\log(1+x)
\big]^\beta
\,\leqslant\,
{\sf Constant}_\beta\cdot 
x
\ \ \ \ \ \ \ \ \ \ \ \ \ {\scriptstyle{(x\,\geqslant\,0)}},
\]
dropping this last constant, it therefore suffices to show a majoration
of the kind:
\[
\aligned
\frac{1}{\kappa^3}\,
\sum_{\lambda=1}^\kappa\,
\sum_{\mu=\lambda}^\kappa\,
\sum_{\nu=2}^\kappa\,
&
\frac{\kappa!\,(\mu-1)!}{(\kappa+\mu)!}\,
\frac{(\kappa-\nu+\mu-\lambda)!}{(\kappa-\nu)!\,(\mu-\lambda)!}\,
\frac{(\nu+\lambda-1)!}{(\lambda-1)!\,(\nu-2)!}\,
\\
&
\bigg[
\log\,\frac{\kappa}{\mu}
\bigg]^\alpha\,
\frac{\kappa-\nu+\mu-\lambda+1}{\nu+\lambda-1}
\,\leqslant\,
{\sf Constant}_\alpha.
\endaligned
\]
Observing then that the numerator and the denominator
of the ultimate fraction can be absorbed in preceding factorials
(so was the trick arranged in advance), and
reorganizing the triple sum, we are studying the uniform
finiteness, as $\kappa \to \infty$, of:
\[
\small
\aligned
\frac{1}{\kappa^3}\,
\sum_{\lambda=1}\,\sum_{\mu=\lambda}^\kappa\,
\sum_{\mu=\lambda}^\kappa\,
\frac{\kappa!\,(\mu-1)!}{(\kappa+\mu)!}\,
\bigg[
\log\,\frac{\kappa}{\mu}
\bigg]^\alpha\,
\sum_{\nu=2}^\kappa\,
\frac{(\kappa-\nu+\mu-\lambda+1)!}{(\kappa-\nu)!\,(\mu-\lambda)!}\,
\frac{(\nu+\lambda-2)!}{(\lambda-1)!\,(\nu-2)!}.
\endaligned
\] 
The underlined simple sum then contracts thanks to a formula
which is elementarily proved and most probably well known.

\begin{Lemma}
For $1\leqslant \lambda \leqslant \mu \leqslant \kappa$ and
$\kappa \geqslant 2$, one has:
\[
\sum_{\nu=2}^\kappa\,
\frac{(\kappa-\nu+\mu-\lambda+1)!}{(\kappa-\nu)!\,(\mu-\lambda)!}\,
\frac{(\nu+\lambda-2)!}{(\lambda-1)!\,(\nu-2)!}
\,=\,
\lambda(\mu-\lambda+1)\,
\frac{(\kappa+\mu)!}{(\mu+2)!\,(\kappa-2)!}.
\qed
\]
\end{Lemma}

Hence after simplifications, we are left with showing the inequality:
\[
\small
\aligned
\frac{1}{\kappa^3}\,
\sum_{\lambda=1}^\kappa\,
\sum_{\mu=\lambda}^\kappa\,
\kappa(\kappa-1)\,
\lambda(\mu-\lambda+1)\,
\frac{1}{(\mu+2)(\mu+1)\mu}\,
\bigg[
\log\,\frac{\kappa}{\mu}
\bigg]^\alpha
\,\leqslant\,
{\sf Constant}_\alpha.
\endaligned
\]
But then, since it is very elementarily checked that:
\[
\bigg[
\log\,\frac{\kappa}{\mu}
\bigg]^\alpha
\,\leqslant\,
{\sf Constant}_\alpha\,
\sqrt{\frac{\kappa}{\mu}}
\ \ \ \ \ \ \ \ \ \ \ \ \ 
{\scriptstyle{(1\,\leqslant\,\mu\,\leqslant\,\kappa)}},
\]
it now only remains to convince oneself, keeping
only dominant terms, that:
\[
\frac{1}{\sqrt{\kappa}}\,
\sum_{\lambda=1}^\kappa\,
\sum_{\mu=\lambda}^\kappa\,
\frac{\lambda(\mu-\lambda)}{\mu^{3+\frac{1}{2}}}
\,\leqslant\,
{\sf Constant},
\]
independently of the bigness of $\kappa$, which is true.

The next cases of dimensions $n = 4$ and higher are treated 
similarly, with the same tools and arguments. Let us at least 
detail the case of dimension $4$ with $\alpha_1' = \alpha_2' = \alpha_3' = 0$,
the case of general $\alpha_1'$, $\alpha_2'$, $\alpha_3'$
being quite devisable from what has been done in dimension $3$.

Changing notation for summed quantities:
\[
\aligned
\lambda_3
&
\equiv
\mu_1^1,
\\
\lambda_2
&
\equiv
\mu_1^2,
\ \ \ \ \ \ \ \ \ \ \ \ \ \ 
\mu_2
\equiv
\mu_2^2,
\\
\lambda_1
&
\equiv
\mu_1^3,
\ \ \ \ \ \ \ \ \ \ \ \ \ \ 
\mu_1
\equiv
\mu_2^3,
\ \ \ \ \ \ \ \ \ \ \ \ \ \ 
\nu_1
\equiv
\mu_3^3,
\endaligned
\]
and also for the exponents of the logarithms:
\[
\alpha
\equiv
\alpha_1',
\ \ \ \ \ \ \ \ \ \ \ \ \ \ 
\beta
\equiv
\alpha_2',
\ \ \ \ \ \ \ \ \ \ \ \ \ \ 
\gamma
\equiv
\alpha_3',
\]
the sums we have to study write out:
\[
\footnotesize
\aligned
{\sf K}_{\alpha,\beta,\gamma}^4(\kappa)
=
\frac{1}{\kappa^6}
\sum_{
\substack{
1\leqslant\lambda_3<\kappa
\\
1\leqslant\lambda_2<\mu_2<\kappa
\\
1\leqslant\lambda_1<\mu_1<\nu_1<\kappa
}}\,
\sum_{
\lambda_3\geqslant\lambda_2\geqslant\lambda_1
\atop
\mu_2\geqslant\mu_1}\,
&
\kappa!\,\kappa!\,
(\lambda_3-1)!\,
\frac{(\lambda_3-\lambda_2+\kappa-\mu_2)!
}{
(\lambda_3-\lambda_2)!\,(\kappa-\mu_2)!}\,
\frac{(\mu_2+\lambda_2-1)!}{(\kappa+\lambda_3)!}
\,
\\
&
\,
\frac{(\lambda_2-\lambda_1+\mu_2-\mu_1+\kappa-\nu_1)!}{
(\lambda_2-\lambda_1)!\,(\mu_2-\mu_1)!\,(\kappa-\nu_1)!}\,
\frac{(\nu_1+\mu_1+\lambda_1-1)!}{(\kappa+\mu_2+\lambda_2)!}
\,
\\
&
\,
\frac{1}{(\lambda_1-1)!\,(\mu_1-2)!\,(\nu_1-3)!}
\,
\\
&
\,
\big[\log\kappa-\log\lambda_3\big]^{\alpha}\,
\big[\log(\kappa+\lambda_3)-\log(\mu_2+\lambda_2)\big]^{\beta}
\\
&
\,
\big[\log(\kappa+\mu_2+\lambda_2)
-
\log(\nu_1+\mu_1+\lambda_1)\big]^{\gamma}.
\endaligned
\] 
Now, we claim that this (already big) sum is majorated by the extended sum:
\[
\overline{\sf K}_{\alpha,\beta,\gamma}^4(\kappa)
:=
\frac{1}{\kappa^6}\,
\sum_{1\leqslant\lambda_1\leqslant\lambda_2
\leqslant\lambda_3\leqslant\kappa}\,
\sum_{2\leqslant\mu_1\leqslant\mu_2\leqslant\kappa}\,
\sum_{3\leqslant\nu_1\leqslant\kappa}
\big(
\text{\small\sf same quantities}
\big)
\]
a sum which, formally, just contains {\em more} terms that the
one we starting with, for the first collection of inequalities:
\[
\aligned
&
1\leqslant\lambda_3<\kappa,
\\
&
1\leqslant\lambda_2<\mu_2<\kappa,
\\
&
1\leqslant\lambda_1<\mu_1<\nu_1<\kappa,
\endaligned
\]
has plainly disappeared, and to be sure of such a claim, it would suffice
to observe that all the terms that are added in this larger
sum are nonnegative; but this fact is clearly okay, because by exactly the 
same argument as in dimension $3$ above,
all the three appearing logarithms bear on quantities that
are all $\geqslant 1$:

\medskip$\square$\,\,
$\kappa \geqslant \mu$ in the original and in the extended sum;

\medskip$\square$\,\,
$\kappa + \lambda_3 \geqslant \mu_2 + \lambda_2$ since 
$\kappa \geqslant \mu_2$ and since $\lambda_3 \geqslant \lambda_2$ in the
extended sum;

\medskip$\square$\,\,
$\kappa + \mu_2 + \lambda_2 \geqslant \nu_1 + \mu_1 + \lambda_1$ since
$\kappa \geqslant \nu_1$, since $\mu_2 \geqslant \mu_1$ and
since $\lambda_2 \geqslant \lambda_1$ in the extended sum.

\medskip\noindent
Hence thanks to the so justified inequality:
\[
{\sf K}_{\alpha,\beta,\gamma}^4(\kappa)
\,\leqslant\,
\overline{\sf K}_{\alpha,\beta,\gamma}^4(\kappa),
\]
it suffices to show that the larger sum in the right-hand side
is uniformly bounded as 
$\kappa \to \infty$.

Thus, let us restrict ourselves to the case $\alpha = \beta = \gamma = 0$,
which is the most significant one.
Writing everything in the sum, we must therefore compute:
\[
\footnotesize
\aligned
\overline{\sf K}_{0,0,0}^4(\kappa)
&
=
\frac{1}{\kappa^6}\,
\sum_{\lambda_1=1}^\kappa\,
\sum_{\lambda_2=\lambda_1}^\kappa\,
\sum_{\lambda_3=\lambda_2}^\kappa\,
\ \
\sum_{\mu_1=2}^\kappa\,
\sum_{\mu_2=\mu_1}^\kappa\,
\ \
\sum_{\nu_1=3}^\kappa
\\
&
\ \ \ \ \
\kappa!\,\kappa!\,(\lambda_3-1)!
\frac{(\kappa-\mu_2+\lambda_3-\lambda_2)!}{
(\kappa-\mu_2)!\,(\lambda_3-\lambda_2)!}\,\,
\frac{(\mu_2+\lambda_2-1)!}{(\kappa+\lambda_3)!}
\\
&
\ \ \ \ \
\frac{(\kappa-\nu_1+\mu_2-\mu_1+\lambda_2-\lambda_1)!}{
(\kappa-\nu_1)!\,(\mu_2-\mu_1)!\,(\lambda_2-\lambda_1)!}\,
\frac{(\nu_1+\mu_1+\lambda_1-1)!}{
(\kappa+\mu_2+\lambda_2)!}\,
\frac{1}{
(\lambda_1-1)!\,(\mu_1-2)!\,(\nu_1-3)!}\,
\\
&
\!\!\!\!\!\!\!\!\!\!\!\!\!\!\!\!\!\!\!\!\!\!
\big[\log(\kappa)-\log(\lambda_3)\big]^0\,
\big[\log(\kappa+\lambda_3)-\log(\mu_2+\lambda_2)\big]^0\,
\big[\log(\kappa+\mu_2+\lambda_2)-\log(\nu_1+\mu_1+\lambda_1)\big]^0.
\endaligned
\]

Naturally, we start out by extracting the summation with respect to $\nu_1$:
\[
\footnotesize
\aligned
&
\overline{\sf K}_{0,0,0}^4(\kappa)
=
\frac{1}{\kappa^6}\,
\sum_{\lambda_1=1}^\kappa\,
\sum_{\lambda_2=\lambda_1}^\kappa\,
\sum_{\lambda_3=\lambda_2}^\kappa\,
\ \
\sum_{\mu_1=2}^\kappa\,
\sum_{\mu_2=\mu_1}^\kappa\,
\\
&
\kappa!\,
\frac{\kappa!\,(\lambda_3-1)!}{(\kappa+\lambda_3)!}\,
\frac{(\kappa-\mu_2+\lambda_3-\lambda_2)!}{
(\kappa-\mu_2)!\,(\lambda_3-\lambda_2)!}\,
\frac{(\mu_2+\lambda_2-1)!}{
(\kappa+\mu_2+\lambda_2)!}\,
\frac{1}{
(\mu_2-\mu_1)!\,(\lambda_2-\lambda_1)!\,(\mu_1-2)!\,(\lambda_1-1)!}
\\
&
\ \ \ \ \
\sum_{\nu_1=3}^\kappa\,
\frac{(\kappa-\nu_1+\mu_2-\mu_1+\lambda_2-\lambda_1)!}{
(\kappa-\nu_1)!}\,
\frac{(\nu_1+\mu_1+\lambda_1-1)!}{(\nu_1-3)!},
\endaligned
\]
a summation which appears in the last line. The useful, elementary lemma, 
quite analogous to the previous one found in dimension $3$, is:

\begin{Lemma}
For $\kappa \geqslant 3$ and for $a, b\geqslant 0$, one has:
\[
\sum_{\nu=3}^\kappa\,
\frac{(\kappa-\nu+a)!}{(\kappa-\nu)!}\,
\frac{(\nu+b)!}{(\nu-3)!}
\,=\,
\frac{a!\,(b+3)!}{(a+b+4)!}\,
\frac{(\kappa+a+b+1)!}{(\kappa-3)!}.
\qed
\]
\end{Lemma}

A direct application then yields a replacement of the third line 
in question above:
\[
\footnotesize
\aligned
&
\overline{\sf K}_{0,0,0}^4(\kappa)
=
\frac{1}{\kappa^6}\,
\sum_{\lambda_1=1}^\kappa\,
\sum_{\lambda_2=\lambda_1}^\kappa\,
\sum_{\lambda_3=\lambda_2}^\kappa\,
\ \
\sum_{\mu_1=2}^\kappa\,
\sum_{\mu_2=\mu_1}^\kappa\,
\\
&
\kappa!\,
\frac{\kappa!\,(\lambda_3-1)!}{(\kappa+\lambda_3)!}\,
\frac{(\kappa-\mu_2+\lambda_3-\lambda_2)!}{
(\kappa-\mu_2)!\,(\lambda_3-\lambda_2)!}\,
\frac{(\mu_2+\lambda_2-1)!}{
\zero{(\kappa+\mu_2+\lambda_2)!}}\,
\frac{1}{
(\mu_2-\mu_1)!\,(\lambda_2-\lambda_1)!\,(\mu_1-2)!\,(\lambda_1-1)!}
\\
&
\ \ \ \ \
\frac{(\mu_2-\mu_1+\lambda_2-\lambda_1)!\,(\mu_1+\lambda_1+2)!}{
(\mu_2+\lambda_2+3)!}\,
\frac{\zero{(\kappa+\mu_2+\lambda_2)!}}{(\kappa-3)!},
\endaligned
\]
and the result, after a few elementary reorganizational
simplifications, becomes:
\[
\footnotesize
\aligned
\overline{\sf K}_{0,0,0}^4(\kappa)
&
=
\frac{\kappa(\kappa-1)(\kappa-2)}{\kappa^6}\,
\sum_{\lambda_1=1}^\kappa\,
\sum_{\lambda_2=\lambda_1}^\kappa\,
\sum_{\lambda_3=\lambda_2}^\kappa\,
\frac{\kappa!\,(\lambda_3-1)!}{(\kappa+\lambda_3)!}\,
\frac{1}{(\lambda_2-\lambda_1)!}\,\frac{1}{(\lambda_1-1)!}\,
\frac{1}{(\lambda_3-\lambda_2)!}\,
\\
&
\ \ \ \ \
\sum_{\mu_1=2}^\kappa\,
\sum_{\mu_2=\mu_1}^\kappa\,
\frac{(\kappa-\mu_2+\lambda_3-\lambda_2)!}{
(\kappa-\mu_2)!}\,
(\mu_2+\lambda_2-1)!\,
\\
&
\ \ \ \ \ \ \ \ \ \ \ \ \ \ \
\frac{1}{(\mu_2-\mu_1)!\,(\mu_1-2)!}\,
\frac{(\mu_2-\mu_1+\lambda_2-\lambda_1)!\,(\mu_1+\lambda_1+2)!}{
(\mu_2+\lambda_2+3)!}.
\endaligned
\]
Having reached this point, we cleverly exchange the two $\mu$-summations,
applying a straightforward discrete Fubini-like `theorem':
\[
\sum_{\mu_1=2}^\kappa\,\sum_{\mu_2=\mu_1}^\kappa
\,=\,
\sum_{\mu_2=2}^\kappa\,\sum_{\mu_1=2}^{\mu_2},
\]
and this immediately gives:
\[
\footnotesize
\aligned
&
\overline{\sf K}_{0,0,0}^4(\kappa)
=
\frac{\kappa(\kappa-1)(\kappa-2)}{\kappa^6}\,
\sum_{\lambda_1=1}^\kappa\,
\sum_{\lambda_2=\lambda_1}^\kappa\,
\sum_{\lambda_3=\lambda_2}^\kappa\,
\frac{\kappa!\,(\lambda_3-1)!}{(\kappa+\lambda_3)!}\,
\frac{1}{(\lambda_2-\lambda_1)!}\,\frac{1}{(\lambda_1-1)!}\,
\frac{1}{(\lambda_3-\lambda_2)!}\,
\\
&
\ \ \ \ \
\sum_{\mu_2=2}^\kappa\,
\frac{(\kappa-\mu_2+\lambda_3-\lambda_2)!}{
(\kappa-\mu_2)!}\,
\frac{(\mu_2+\lambda_2+1)!}{(\mu_2+\lambda_2+3)!}\,
\sum_{\mu_1=2}^{\mu_2}\,
\frac{(\mu_2-\mu_1+\lambda_2-\lambda_1)!}{
(\mu_2-\mu_1)!}\,
\frac{(\mu_1+\lambda_1+2)!}{(\mu_1-2)!}.
\endaligned
\]
But to close up the summation of the last line, there is
a general elementary formula that we already encountered in dimension
$3$:
\[
\sum_{\nu=2}^\kappa\,
\frac{(\kappa-\nu+a)!}{(\kappa-\nu)!}\,
\frac{(\nu+b)!}{(\nu-2)!}
\,=\,
\frac{a!\,(b+2)!}{(a+b+3)!}\,
\frac{(\kappa+a+b+1)!}{(\kappa-2)!},
\]
which, applied with the right $a$ and $b$ and with another $\kappa$, writes:
\[
\small
\aligned
\sum_{\mu_1=2}^{\mu_2}\,
\frac{(\mu_2-\mu_1+\lambda_2-\lambda_1)!}{
(\mu_2-\mu_1)!}\,
\frac{(\mu_1+\lambda_1+2)!}{(\mu_1-2)!}
\,=\,
\frac{(\lambda_2-\lambda_1)!\,(\lambda_1+4)!\,(\mu_2+\lambda_2+3)!}{
(\lambda_2+5)!\,(\mu_2-2)!},
\endaligned
\]
and inserting this in the last line above, we get:
\[
\footnotesize
\aligned
\overline{\sf K}_{0,0,0}^4(\kappa)
&
=
\frac{\kappa(\kappa-1)(\kappa-2)}{\kappa^6}\,
\sum_{\lambda_1=1}^\kappa\,
\sum_{\lambda_2=\lambda_1}^\kappa\,
\sum_{\lambda_3=\lambda_2}^\kappa\,
\frac{\kappa!\,(\lambda_3-1)!}{(\kappa+\lambda_3)!}\,
\frac{1}{\zero{(\lambda_2-\lambda_1)!}}\,\frac{1}{(\lambda_1-1)!}\,
\frac{1}{(\lambda_3-\lambda_2)!}\,
\\
&
\ \ \ \ \
\sum_{\mu_2=2}^\kappa\,
\frac{(\kappa-\mu_2+\lambda_3-\lambda_2)!}{
(\kappa-\mu_2)!}\,
\frac{(\mu_2+\lambda_2-1)!}{\zerozero{(\mu_2+\lambda_2+3)!}}\,
\frac{\zero{(\lambda_2-\lambda_1)!}\,(\lambda_1+4)!\,
\zerozero{(\mu_2+\lambda_2+3)!}}{
(\lambda_2+5)!\,(\mu_2-2)!},
\endaligned
\]
that is to say after simplifications:
\[
\footnotesize
\aligned
\overline{\sf K}_{0,0,0}^4(\kappa)
&
=
\frac{\kappa(\kappa-1)(\kappa-2)}{\kappa^6}\,
\sum_{\lambda_1=1}^\kappa\,
\sum_{\lambda_2=\lambda_1}^\kappa\,
\sum_{\lambda_3=\lambda_2}^\kappa\,
\frac{\kappa!\,(\lambda_3-1)!}{(\kappa+\lambda_3)!}\,
\frac{1}{(\lambda_1-1)!}\,
\frac{1}{(\lambda_3-\lambda_2)!}\,
\frac{(\lambda_1+4)!}{(\lambda_2+5)!}
\\
&
\ \ \ \ \
\sum_{\mu_2=2}^\kappa\,
\frac{(\kappa-\mu_2+\lambda_3-\lambda_2)!}{
(\kappa-\mu_2)!}\,
\frac{(\mu_2+\lambda_2-1)!}{(\mu_2-2)!}.
\endaligned
\]
Thanks to a second application of the same formula, the $\sum_{\mu_2}$ also
contracts:
\[
\footnotesize
\aligned
\overline{\sf K}_{0,0,0}^4(\kappa)
&
=
\frac{\kappa(\kappa-1)(\kappa-2)\kappa(\kappa-1)}{\kappa^6}\,
\sum_{\lambda_1=1}^\kappa\,
\sum_{\lambda_2=\lambda_1}^\kappa\,
\sum_{\lambda_3=\lambda_2}^\kappa\,
(\lambda_1+4)(\lambda_1+3)(\lambda_1+2)(\lambda_1+1)\lambda_1\,
\\
&
\ \ \ \ \
\frac{1}{(\lambda_2+5)(\lambda_2+4)(\lambda_2+3)(\lambda_2+2)}\,
\frac{1}{(\lambda_3+2)(\lambda_3+1)\lambda_3}
\\
&
=
\frac{\kappa(\kappa-1)(\kappa-2)\kappa(\kappa-1)}{\kappa^6}\,
\frac{1}{18}\,\kappa,
\endaligned
\]
the latter computation being elementary,
and as a result, we conclude that:
\[
\overline{\sf K}_{0,0,0}^4(\kappa)
\,\leqslant\,
\frac{1}{18},
\]
uniformly, especially when $\kappa \to \infty$.

Applying the same tricky majorations for the log-terms as in dimension
$3$, one realizes that the same kind of computation also goes
through. In arbitrary dimension $n \geqslant 2$, we will bypass
a complete indicial detailed writing, which would only use the same ideas.
\endproof

\subsection{Indirect majorations}
Before concluding this section, it is advisable to mention
that sums of the precise kind as that of the preceding
proposition already appeared implicitly before.  

Indeed, looking back at the Euler-Poincar\'e characteristic, in the
summation formula:
\[
\chi\big(X,\,
{\sf Gr}^\bullet\mathcal{E}_{\kappa,m}^{GG}T_X^*
\big)
=
\sum_{\ell_1\geqslant\ell_2\geqslant\cdots\geqslant\ell_n\geqslant 0}\,
M_{\ell_1,\ell_2,\cdots,\ell_n}\cdot
\chi\big(X,\,
\mathcal{S}^{(\ell_1,\ell_2,\dots,\ell_n)}T_X^*
\big),
\]
the coefficients of each Chern monomial ${\sf c}_1^{ \lambda_1} {\sc
c}_2^{ \lambda_2} \cdots\, {\sf c}_n^{ \lambda_n}$ must identify.  In
the Euler-Poincar\'e characteristic of the Schur bundle, the coefficient
of ${\sf c}_1^n$ is, up to a rational factor:
\[
\prod_{1\leqslant i<j\leqslant n}\,
(\ell_i-\ell_j)\,
\sum_{\beta_1+\cdots+\beta_{n-1}+\beta_n=n}\,
\ell_1^{\beta_1}\cdots\,\ell_{n-1}^{\beta_{n-1}}
\ell_n^{\beta_n}. 
\]
We then rewrite:
\[
\aligned
&
\sum_{\beta_1+\cdots+\beta_{n-1}+\beta_n=n}\,
\ell_1^{\beta_1}\cdots\,\ell_{n-1}^{\beta_{n-1}}
\ell_n^{\beta_n}
=
\\
&
=
\sum_{\beta_1'+\cdots+\beta_{n-1}'+\beta_n'}\,
C_{\beta_1',\dots,\beta_{n-1}',\beta_n'}\,
(\ell_1-\ell_2)^{\beta_1'}
\cdots\,
(\ell_{n-1}-\ell_n)^{\beta_{n-1}'}
\ell_n^{\beta_n'},
\endaligned
\]
with coefficients $C_{ \beta_1', \dots, \beta_{ n-1}', \beta_n'} \in
\N$. Notice that $C_{ 0, \dots, 0, n} = \binom{ n + n-1}{ n}$.
Identifying then the coefficients of ${\sf c}_1^n$, one
realizes then without any computation that:
\[
\boxed{
\footnotesize
\aligned
\sum_{{\sf YT}\,{\sf semi-standard}
\atop
{\sf weight}({\sf YT})=m}\,
\prod_{1\leqslant i<j\leqslant n}\,
(\ell_i-\ell_j)
\cdot
(\ell_n)^n
&
=
{\textstyle{\frac{m^{(\kappa+1)n-1}}{((\kappa+1)n-1)!}}}\,
1!\,2!\,\cdots\,(n-1)!\,
(\log\kappa)^n
+
\\
&
\ \ \ \ \
+
{\sf O}_n
\big(m^{(\kappa+1)n-1}
\cdot
(\log\kappa)^{n-1}\big)
+
{\sf O}_{n,\kappa}\big(m^{(\kappa+1)n-2}\big).
\endaligned
}
\]
Visibly, the power $(\ell_n)^n$ of $\ell_n$ corresponds to 
$(\log\kappa)^n$. 

Specifically, in dimension $n = 2$, looking at the coefficient
of ${\sf c}_2$ and making identification, one 
convinces oneself that one gets:
\[
\sum_{1\leqslant\lambda\leqslant\kappa}\,
\frac{1}{\lambda^2}
=
\sum_{1\leqslant\mu_1^1<\kappa}\,
(\kappa!)^2\,
{\textstyle{\frac{N_{\mu_1^1}^\kappa}{\kappa\cdots\,\mu_1^1}}}\,
{\textstyle{\frac{N_{1,2}^{\mu_1^1,\kappa}}{
(\kappa+\mu_1^1)\cdots\,(2+1)}}}\,
\bigg(
\sum_{q_0^1+\cdots+q_{\tau^1}^1=3}\,
\frac{1}{(\mu_1^1)^{q_0^1}\cdots\,(\kappa)^{q_{\tau^1}^1}}
\bigg),
\]
so one deduces without computation that the sum:
\[
\sum_{1\leqslant\mu_1^1<\kappa}\,
(\kappa!)^2\,
{\textstyle{\frac{N_{\mu_1^1}^\kappa}{\kappa\cdots\,\mu_1^1}}}\,
{\textstyle{\frac{N_{1,2}^{\mu_1^1,\kappa}}{
(\kappa+\mu_1^1)\cdots\,(2+1)}}}\,
\big[\log\kappa-\log\mu_1^1\big]^3
\]
is finite and bounded independently of $\kappa$. In dimensions $n =
3$ and higher, looking at the coefficient of ${\sf c}_n$, one sees
indirectly, without computations and without majorations that all the
sums $\square_{ n, \kappa}^{\alpha_1', \dots, \alpha_{n-1}', 0}$ which
appear after expressing:
\[
\prod_{1\leqslant i<j\leqslant n}\,
\sum_{\beta_1'+\cdots+\beta_{n-1}'=n}\,
\ell_1^{\beta_1'}\cdots\,\ell_{n-1}^{\beta_{n-1}'}
\]
in terms of $(\ell_1 - \ell_2)$, \dots, 
$(\ell_{ n-1} - \ell_n)$, $\ell_n$ are finite. 
These observations confirm what was delineated in the
previous paragraphs.  

\subsection{Summary}
In conclusion, either directly or indirectly by identification without
computations and without majorations, we have seen that for any
$\alpha_1', \dots, \alpha_{ n-1}', \alpha_n' \in \N$ with $\alpha_1' +
\cdots + \alpha_{ n-1}' + \alpha_n' \leqslant \frac{ n ( n+1)}{ 2}$,
the quantity $\Delta_{ n, \kappa}^{ \alpha_1', \dots, \alpha_{ n-1}',
0}$ is $\leqslant {\sf Constant}_n$, whence $\square_{ n, \kappa}^{
\alpha_1', \dots, \alpha_n'}$ is $\leqslant {\sf Constant}_n \cdot (
\log \kappa)^{ \alpha_n'}$, and from~\thetag{ \ref{intermediate}}, it
follows at the end that:
\begin{equation}
\label{finite-sums}
\footnotesize
\aligned
\sum_{{\sf YT}\,{\sf semi-standard}
\atop
{\sf weight}({\sf YT})=m}
&
\big(\ell_1({\sf YT})-\ell_2({\sf YT})\big)^{\alpha_1'}
\cdots
\big(\ell_{n-1}({\sf YT})-\ell_n({\sf YT})\big)^{\alpha_{n-1}'}
\leqslant
\\
&
\leqslant
\frac{m^{(\kappa+1)n-1}}{((\kappa+1)n-1)!\,(\kappa!)^n}
\cdot
{\sf Constant}_n
+
{\sf Constant}_{n,\kappa}\
\cdot
m^{(\kappa+1)n-2}.
\endaligned
\end{equation}

\markleft{Jo\"el Merker}
\markright{\S11.~Algebraic sheaf theory and Schur bundles}
\section{\bf Algebraic sheaf theory and Schur bundles}
\label{Section-11}

\subsection{Complex projective hypersurface and
line bundles $\mathcal{ O}_X ( k)$}
Let $X = X^n \subset \P^{n+1} ( \C)$ be a geometrically smooth complex
projective hypersurface of degree $d \geqslant 1$, defined in
homogeneous coordinates $z = [z_0 \colon z_1 \colon \cdots \colon z_n
\colon z_{ n+1} ]$ as the zero-set:
\[
X
=
\big\{
[z_0\colon z_1\colon\cdots\colon z_n\colon z_{n+1}]
\in
\P^{n+1}(\C)
\colon
P(z_0,z_1,\dots,z_n,z_{n+1})
=
0
\big\}
\]
of a certain holomorphic polynomial $P = P ( z) \in \C[ z_0, z_1,
\dots, z_n, z_{ n+1} ]$ which is homogeneous of a certain degree $d
\geqslant 1$ and whose differential $P_{z_0} dz_0 + \cdots + P_{
z_{n+1}} d z_{n+1}$ does not vanish at any point of $X$, 
so that $X$ has no singularities.
We will sometimes use the letter $N$ to denote $n+1$: 
\[
N
\overset{\sf notation}{\equiv}
n+1.
\]
The {\sl tautological line bundle} over $\P^N$ will be denoted by
$\mathcal{ O}_{ \P^N} ( -1)$ and its dual by $\mathcal{ O}_{ \P^N} (
1) := \mathcal{ O}_{ \P^N} ( -1)^*$. For various values of the
integer $k \in \Z$, the standard line bundles:
\[
\mathcal{O}_{\P^N}(k)
:= 
\mathcal{O}_{\P^N}(\pm 1)^{\otimes\vert k\vert},
\] 
where $\pm = {\rm sign} (k)$, will play a very decisive r\^ole in what
follows, as well as their restrictions to $X$, namely the bundles:
\[
\mathcal{O}_X(k)
:=
\mathcal{O}_{\P^N}(k)
\big\vert_X.
\]

\subsection{Canonical line bundles} 
For any $\P^N$, the (line) bundle of holomorphic differential forms of
maximal degree $N$ on $\P^N$:
\[
K_{\P^N}
=
\Lambda^NT_{\P^N}^*
\simeq
\mathcal{O}_{\P^N}(-N-1),
\]
is known, thanks to the {\sl adjunction formula}, to be isomorphic to
$\mathcal{ O}_{ \P^N} ( - N - 1)$. Similarly, the (line) bundle of
holomorphic differential forms of maximal degree $n$ on $X$:
\[
K_X
\overset{\sf notation}{\equiv}
\Lambda^nT_X^*
\simeq
\mathcal{O}_X(d-n-2)
\]
called the {\sl canonical bundle} of $X$ and central in complex algebraic
geometry, is known, again thanks to the adjunction formula, to be
isomorphic to $\mathcal{O}_X(d-n-2)$.

\subsection{Normal exact sequence}
To begin with, one has the so-called {\sl normal exact sequence}: 
\begin{equation}
\label{normal-exact}
\aligned
0
\longrightarrow
\mathcal{O}_{\P^{n+1}}(-d)
\overset{\text{\small\sf\em incl}}{\longrightarrow}
\mathcal{O}_{\P^{n+1}}(0)
\overset{\text{\small\sf\em rest}}{\longrightarrow}
\mathcal{O}_X(0)
\longrightarrow
0.
\endaligned 
\end{equation}
Here, the inclusion {\small\sf\em incl} is defined by multiplication
with the defining polynomial $P ( z_0, \dots, z_{ n+1})$ for $X$, and
the restriction {\small\sf\em rest}, of course from $\P^{ n+1}$ to
$X$, concerns functions, differential forms, bundles and sheaves.

\subsection{General sheaves of differential forms}
Let $r$ be an integer with $0 \leqslant r \leqslant n+1$ and consider
the bundle $\Lambda^r T_X^*$ of differential forms of degree $r$ on
$X$, with the convention that:
\begin{equation}
\label{convention-zero}
\Lambda^0T_X^*
\overset{\sf collapse}{\equiv}
\mathcal{O}_X(0).
\end{equation}
The functor $\mathcal{ F} \longmapsto \mathcal{ F} \otimes \mathcal{
G}$ is right exact, for any sheaf $\mathcal{ G}$, and is furthermore
also left exact when $\mathcal{ G}$ is locally free (in what follows,
only such sheaves will be considered). Here at any point $z \in X$,
the bundle $\Lambda^k T_X^*$ is, for any $k$ with $0 \leqslant k
\leqslant n$, a free $\mathcal{ O}_{X, z}$-module of rank $\binom{ n}{
k}$, hence by tensoring the above normal exact sequence, we obtain the
following exact sequence:
\[
\aligned
0
\longrightarrow
\Lambda^kT_{\P^{n+1}}^*
\otimes\mathcal{O}_{\P^{n+1}}(-d)
\longrightarrow
\Lambda^kT_{\P^{n+1}}^*
\otimes\mathcal{O}_{\P^{n+1}}(0)
\longrightarrow
\\
\longrightarrow
\Lambda^kT_{\P^{n+1}}^*
\otimes\mathcal{O}_X(0)
\longrightarrow
0.
\endaligned
\]

\subsection{Hook lengths of Young diagrams}
More generally, let $r \geqslant 1$ be any nonnegative integer and
let: 
\[
(\ell)
=
(\ell_1,\ell_2,\dots,\ell_n)
\]
be an arbitrary {\sl partition} of $r$ in at most $n$ parts, namely
the sum:
\[
\ell_1+\ell_2+\cdots+\ell_n 
= 
r
\] 
equals $r$, and the {\sl parts} $\ell_i$ are ordered decreasingly:
\[
\ell_1\geqslant\ell_2\geqslant\cdots\geqslant\ell_n
\geqslant 
0.
\]
Let $d_1, d_2, \dots, d_{\ell_1}$ denote the column lengths of the
diagram consisting of $\ell_1$ blank squares above $\ell_2$ blank
squares, \dots, above $\ell_n$ blank squares.

\noindent
With a bit more precisions, we hence can denote our arbitrary
partition as:
\[
\left[
\aligned
(\ell)
&
=
\big(\ell_1,\ell_2,\dots,\ell_{d_1},0,\dots,0\big)
\\
\ell_1
&
\geqslant\ell_2\geqslant\cdots
\geqslant\ell_{d_1}
\geqslant
1,
\endaligned\right.
\]
and as in Section~4, we will denote by: 
\[
{\sf YD}_{(\ell)} 
= 
{\sf YD}_{(\ell_1,\ell_2,\dots,\ell_{d_1}, 
0,\dots,0)}
\] 
the associated Young diagram. The {\sl hook-length} $h_{ i,j}$ of the
diagram at the square of coordinates $(i, j)$ is equal to:
\[
h_{i,j}
:=
\ell_i-j+d_j-i+1. 
\] 
A preliminary combinatorial fact, useful soon, is as follows.

\begin{Theorem}
{\rm (\cite{ fuha1991})}
The number of ways to fill in the $r$ blank cases of the diagram
${\sf YD}_{ (\ell_1, \dots, \ell_n)}$ 
just with the first $r$ nonnegative integers $1, 2, 3, \dots,
r$ in such a way that the appearing integers do increase (strictly)
along each row and do also increase (strictly) along each column is
equal to the {\em integer:}
\[
\nu_{(\ell)}
:=
\frac{r!}{\prod_{i,j}\,h_{i,j}}.
\]
\end{Theorem}

\subsection{Schur bundles}
On every fiber $\big( T_{X,x}^* \big)^{ \otimes r}$ of the $r$-th
tensor bundle $\big( T_X^* \big)^{ \otimes r}$ over a point $x \in X$,
the full linear group ${\sf GL}_n ( \C) \ni {\sf w}$ acts in a natural
way:
\[
{\sf w}\cdot
v_{i_1}^*
\otimes
v_{i_2}^*
\otimes\cdots\otimes
v_{i_r}^*
:=
{\sf w}(v_{i_1}^*)
\otimes
{\sf w}(v_{i_2}^*)
\otimes\cdots\otimes
{\sf w}(v_{i_r}^*),
\]
if by $(v_1^*, v_2^*, \dots, v_n^*)$ one denotes any fixed basis of
$T_{ X, x}^*$. Since the works of Schur?? at the end of the
19\textsuperscript{th}, it is known ({\em see}~\cite{ fuha1991}) how
one may decompose this action into irreducible (nondecomposable)
representations which generate the Schur bundles $\mathcal{ S}^{
(\ell_1, \ell_2, \dots, \ell_n)} T_X^*$ that were already considered
in Section~4. Let us provide more information here.

A {\em Young tableau} ${\sf YT}_{1,2, \dots, r}$ is a filling of a
given Young diagram ${\sf YD}_{( \ell_1, \dots, \ell_n)}$ having $r =
\ell_1 + \cdots + \ell_n$ blank boxes precisely by means of the first
$r$ positive integers $1, 2, \dots, r$. Notice {\em passim} that only
a special kind of Young tableaux was considered in the theorem above,
namely those which enjoy strict increase both along lines and columns,
and such combinatorial objects are usually called {\sl standard Young
tableaux}.

\subsection{Idempotents in the group algebra of 
permutations} Introduce also the {\sl group algebra} $\Q \cdot
\mathfrak{ S}_r$ over the permutation group:
\[
\mathfrak{ S}_r
= 
{\sf Perm}
\big(
\{1,2,\dots,r\}\big),
\]
whose general element is a typical sum $\sum_{ \sigma \in \mathfrak{
S}_r} \, c_\sigma \cdot \sigma$ having arbitrary rational coefficients
$c_\sigma \in \Q$, the addition: 
\[
{\textstyle{\sum_{\sigma\in\mathfrak{S}_r}}}\,
c_\sigma\cdot\sigma
+
{\textstyle{\sum_{\sigma\in\mathfrak{S}_r}}}\,
d_\sigma\cdot\sigma
=
\sum_{\sigma\in\mathfrak{S}_r}\,
\big(
c_\sigma
+
d_\sigma
\big)\cdot\sigma
\]
being obvious and the ``multiplication'':
\[
\Big(
{\textstyle{\sum_{\sigma'\in\mathfrak{S}_r}}}\,
c_{\sigma'}\cdot\sigma'
\Big)
\circ
\Big(
{\textstyle{\sum_{\sigma''\in\mathfrak{S}_r}}}\,
c_{\sigma''}\cdot\sigma''
\Big)
=
\sum_{\sigma'\in\mathfrak{S}_r}\,
\sum_{\sigma''\in\mathfrak{S}_r}\,
c_{\sigma'}\,c_{\sigma''}\cdot
\sigma'\circ\sigma''
\]
corresponding naturally to the composition $\sigma ' \circ \sigma''$
of permutations. For a given standard Young tableau ${\sf YT}_{ 1, \dots, r}$
which shall also be denoted shortly by ${\sf T}$, one introduces the
following element:
\begin{equation}
\label{e-idempotent}
e_{\sf T}
:=
\frac{\nu_{(\ell)}}{r!}
\cdot
\bigg(
\sum_{q\in Q_{\sf T}}\,{\sf sgn}(q)\cdot q
\bigg)
\circ
\bigg(
\sum_{p\in P_{\sf T}}\,p
\bigg)
\end{equation}
of the group algebra $\Q \cdot \mathfrak{ S}_r$, where $Q_{\sf T}$
denotes the set of permutations that preserve the numbers present in
each column of ${\sf T}$, and where similarly $P_{\sf T}$ denotes the
set of permutations that preserve the numbers present in each row of
${\sf T}$.

\begin{Theorem}
{\rm (\cite{ fuha1991})}
This element $e_{\sf T}$ is an idempotent:
\[
e_{\sf T}
\circ
e_{\sf T}
=
e_{\sf T},
\]
and the identity permutation ${\sf Id} \in \Q \cdot
\mathfrak{ S}_r$ decomposes as the sum of all
such idempotents:
\[
{\sf Id}
=
\sum_{
{\sf T}
=
\text{\sf standard Young tableau}
\atop
{\sf Card}({\sf T})
=
r
}
e_{\sf T}.
\]
\end{Theorem}

\subsection{Canonical decomposition of tensor powers of the
cotangent bundle}
The symmetric group $\mathfrak{ S}_r$ and therefore also the group
algebra $\Q \cdot \mathfrak{ S}_r$, act on $\big( T_X^* \big)^{
\otimes r}$ just by permuting the spots inside the tensor product:
\[
\sigma\cdot
v_1
\otimes
v_2
\otimes\cdots\otimes
v_r
:=
v_{\sigma^{-1}(1)}
\otimes
v_{\sigma^{-1}(2)}
\otimes\cdots\otimes
v_{\sigma^{-1}(r)}.
\]
The identity decomposition~\thetag{ \ref{e-idempotent}} then yields at
any point $x \in X$ the direct sum decomposition of 
the $r$-th tensor
power of the cotangent space:
\[
\big(T_{X,x}^*\big)^{\otimes r}
=
\bigoplus_{
{\sf T}=
\text{\sf Young tableau}
\atop
{\sf Card}({\sf T})=r}
=
\mathcal{S}^{{\sf T}}T_{X,x}^*
\ \ \ \ \ \
\text{\rm with}
\ \ \ \ \
\mathcal{S}^{{\sf T}}T_{X,x}^*
:=
e_{\sf T}
\cdot
\big(T_{X,x}^*\big)^{\otimes r}.
\]
This generalizes the simple well known case $r = 2$: 
\[
\big(T_{X,x}\big)^{\otimes 2}
=
\Lambda^2T_{X,x}^*
\oplus
{\rm Sym}^2T_{X,x}^*.
\]

\begin{Theorem}
{\rm (\cite{ fuha1991})}
For any Young Tableau ${\sf T}$, a basis of $\mathcal{ S}^{{\sf T}}
T_{ X, x}^*$ as a $\C$-vector space is constituted of all vectors of
the form:
\[
e_{\sf T}
\big(
v_{i_1}
\otimes
v_{i_2}
\otimes\cdots\otimes
v_{i_r}
\big),
\]
for any choice of integers $i_1, i_2, \dots, i_r \in \{ 1, \dots,
n\}$ having the property that the filling of the blank boxes of the
underlying Young diagram with the integers $i_1, \dots, i_{ \ell_1}$
in the first line, then with the integers $i_{ \ell_1 + 1}, \dots, i_{
\ell_1 + \ell_2}$ in the second line, and so on, provides at the end
a {\sl semi-standard Young tableau}, in the sense that integers are
always nondecreasing when read in each row from left to right, and are
always increasing (strictly) when read in each column from top to
bottom.
\end{Theorem}

It turns out (\cite{ fuha1991, bru1997, ful1998, dem1997, rou2006a})
that, if two arbitrary Young tableaux ${\sf T}$ and $\widetilde{ {\sf
T}}$ correspond to the same Young diagram, {\em i.e.} to the same
partition, then $\mathcal{ S}^{{\sf T}} T_{ X, x}^*$ and $\mathcal{
S}^{ \widetilde{ \sf T}} T_{ X, x}^*$ are {\em isomorphic}. Moreover,
for any ${\sf T}$, the linear action of ${\sf GL}_n ( \C)$ being
compatible with the changes of chart on $X$ and on $T_X^*$, one may
show that the various fibers $\mathcal{ S}^{ {\sf T}} T_{ X, x}^*$
organize coherently as a holomorphic bundle over $X$. In conclusion, a
fundamental Schur bundle decomposition theorem holds
which gives the complete generalization of, say: 
\[
\aligned
\big(T_X^*\big)^{\otimes 2}
&
=
\mathcal{S}^{(2,0,\dots,0)}T_X^*
{\scriptstyle{\,\bigoplus\,\,}}
\mathcal{S}^{(1,1,0,\dots,0)}T_X^*,
\\
\big(T_X^*\big)^{\otimes 3}
&
=
\mathcal{S}^{(3,0,\dots,0)}T_X^*
{\scriptstyle{\bigoplus}}\,
\big[\mathcal{S}^{(2,1,0,\dots,0)}T_X^*\big]^{\oplus 2}
{\scriptstyle{\bigoplus}}\,
\mathcal{S}^{(1,1,1,0,\dots,0)}T_X^*,
\endaligned
\]
provided $X$ is of dimension $\geqslant 3$; 
the last factor is dropped when $\dim X = 2$. 

\begin{Theorem}
{\rm (\cite{ fuha1991})}
For any integer $r \geqslant 1$, the $r$-th tensor power of the
cotangent bundle $T_X^*$ of an arbitrary $n$-dimensional complex
manifold $X$ splits up as a direct sum of so-called Schur bundles:
\[
\big(T_X^*\big)^{\otimes r}
=
\bigoplus_{
\ell_1\geqslant\ell_2\geqslant\cdots\geqslant\ell_n\geqslant 0
\atop
\ell_1+\ell_2+\cdots+\ell_n=r}
\big(
\mathcal{S}^{(\ell_1,\ell_2,\dots,\ell_n)}T_X^*
\big)^{\oplus\nu(\ell)}
\]
indexed by all the partitions $(\ell)$ of $r$. The rank of
$\mathcal{ S}^{ (\ell_1, \dots, \ell_n)} T_X^*$ as 
a complex vector bundle equals:
\[
{\rm rank}
\big(
\mathcal{S}^{(\ell_1,\ell_2,\dots,\ell_n)}T_X^*
\big)
=
\prod_{1\leqslant i<j\leqslant n}\,
\bigg(
\frac{\ell_i-\ell_j}{i-j}
+
1
\bigg),
\]
and the integer multiplicities:
\[
\nu_{(\ell)}
=
\frac{r!}{
\prod_{i,j}\,h_{i,j}}
\]
appearing in the decomposition are expressible in terms of the
hook lengths $h_{ i,j}$ of the concerned Young diagram ${\sf YD}_{
(\ell_1, \ell_2, \dots, \ell_n)}$.
\end{Theorem}

\subsection{Dividing by $K_X$}
Our main goal will now be to control the cohomology of the $\mathcal{
S}^{ (\ell_1, \dots, \ell_n)} T_X^*$ by a formula which will
complement the inequality of Theorem~\ref{majoration-Rousseau}, in the
case where $\ell_n$ is large (whence all the $\ell_i$ are so too). 
It is then natural to 
use the known formula:
\[
\mathcal{S}^{(\ell_1,\dots,\ell_n)}T_X^*
\otimes
K_X
=
\mathcal{S}^{(\ell_1,\dots,\ell_n)}T_X^*
\otimes
\mathcal{S}^{(1,\dots,1)}T_X^*
=
\mathcal{S}^{(\ell_1+1,\dots,\ell_n+1)}T_X^*
\]
under the subtraction form:
\[
\aligned
\mathcal{S}^{(\ell_1,\dots,\ell_{n-1},\ell_n)}T_X^*
&
=
\mathcal{S}^{(\ell_1-\ell_n,\dots,\ell_{n-1}-\ell_n,0)}T_X^*
\otimes
(K_X)^{\otimes\ell_n}
\\
&
=
\mathcal{S}^{(\ell_1-\ell_n,\dots,\ell_{n-1}-\ell_n,0)}T_X^*
\otimes
\mathcal{O}_X\big(\ell_n(d-n-2)\big),
\endaligned
\]
which underlines twisting by a certain $\mathcal{ O}_X ( t)$. 
On the occasion, it is
known thanks to analytical tools (cf. Section~6 in~\cite{ dem2000}) that
if $\mathcal{ E}$ is {\em any} holomorphic vector bundle on the
hypersurface $X \subset \P^{ n+1} ( \C)$ and if $\mathcal{ L}$ is an
ample (or even nef) {\em line} bundle on $X$, then:
\[
\dim H^q
\big(X,\,\,\mathcal{E}\otimes\mathcal{L}^{\otimes k}\big)
=
{\sf O}(k^{n-q}),
\]
for any $q = 0, 1, 2, \dots, n$. Using purely algebraic tools, what we
will do now is to make this estimate much more effective in the case
we are interested in, namely 
when $\mathcal{ E} = \mathcal{ S}^{ ( \ell_1 -
\ell_n, \dots, \ell_{ n-1} - \ell_n, 0)} T_X^*$ and when $\mathcal{
L} = \mathcal{ O}_X ( 1)$ on a general type hypersurface $X \subset
\P^{ n+1}$; in this case, $X$ is of degree $d \geqslant n+3$, whence:
\[
K_X 
= 
\mathcal{O}_X(d-n-2)
=
\mathcal{O}_X(1)^{\otimes(d-n-2)}
\]
is ample of course, so that the exponent $k := \ell_n ( d - n - 2)$
in:
\[
(K_X)^{\ell_n}
=
\big(\mathcal{O}_X(1)\big)^{\otimes(\ell_n(d-n-2))}
=
\mathcal{L}^{\otimes(\ell_n(d-n-2))}
\]
is positive and in fact large. However, 
the Landau-type estimate ``${\sf O}$'' above
provided by analytic techniques is not precise
enough and we need instead explicit {\em inequalities}.  To achieve
more effective estimates, three fundamental exact sequences of
holomorphic vector bundles due to Lascoux (\cite{ lasc1978}) and
to Br\"uckmann (\cite{ bru1997}) will be very helpful. Thus our goal 
is to study the cohomology of the twisted Schur bundles:
\[
\mathcal{S}^{(\ell_1',\dots,\ell_{n-1}',0)}T_X^*
\otimes
\mathcal{O}_X(t),
\]
when $t$ is large.

\subsection{First fundamental (long) exact sequence}
Dualizing the Euler exact sequence:
\[
0
\longrightarrow
\mathcal{O}_{\P^{n+1}}(0)
\longrightarrow
\mathcal{O}_{\P^{n+1}}(1)^{\oplus(n+2)}
\longrightarrow
T_{\P^{n+1}}
\longrightarrow
0,
\]
we get as a starter the exact sequence: 
\[
0
\longrightarrow
T_{\P^{n+1}}^*
\longrightarrow
\mathcal{O}_{\P^{n+1}}(-1)^{\oplus(n+2)}
\longrightarrow
\mathcal{O}_{\P^{n+1}}(0)
\longrightarrow
0.
\]
The procedure explained by Br\"uckmann in~\cite{ bru1997}
of taking the $r$-th tensor power of the extracted complex
composed of the last two bundles: 
\[
\cdots
\longrightarrow
0
\longrightarrow
0
\longrightarrow
\mathcal{O}_{\P^{n+1}}(-1)^{\oplus(n+2)}
\longrightarrow
\mathcal{O}_{\P^{n+1}}(0)
\longrightarrow
0
\]
and more generally, of taking any of its Schur powers, provides a
useful long exact sequence of holomorphic vector bundles which gives
a free resolution of $\mathcal{ S}^{ (\ell_1, \dots, \ell_n, \ell_{
n+1})} T_{\P^{ n+1}}^*$ 
on $\P^{ n+1}$. Instead of using the same letter $\mathcal{
S}$ for Schur bundles over $\P^{ n+1}$ and over $X$, we shall, in
order to underline a clearly visible distinction between $\P^{ n+1}$
and $X$, write:
\[
\mathcal{S}^{(\ell_1,\dots,\ell_n,\ell_{n+1})}
T_{\P^{n+1}}^*
\overset{\sf notation}{\equiv}
\Omega_{\P^{n+1}}^{(\ell_1,\dots,\ell_n,\ell_{n+1})},
\]
using the Greek letter\footnote{\,
Justification: in several articles, the letter $\Omega$ is
employed to denote the bundles $\Omega^k$ or $\Omega^k T_X^*$, $0
\leqslant k \leqslant n$, that we denoted by $\Lambda^k T_X^*$ above.
} 
$\Omega$ with `$\P^{ n+1}$' placed at the lower index place. 
 
Let now ${\sf T}$ be a Young tableau with $r$ boxes and with row
lengths $\ell_1\geqslant \ell_2 \geqslant \cdots\geqslant \ell_{ n+1}
\geqslant 0$, hence of depth $\leqslant n+1$. For convenient
abbreviation, we introduce the general notation:
\[
\Delta\big(\theta_1,\theta_2,\dots,\theta_K\big)
:=
\prod_{1\leqslant i<j\leqslant K}
\big(\theta_i-\theta_j\big)
\]
which is, up to sign, the value: 
\[
\left\vert\!
\begin{array}{ccccc}
1 & \theta_1 & \theta_1^2 & \cdots & \theta_1^{K-1}
\\
1 & \theta_2 & \theta_2^2 & \cdots & \theta_2^{K-1}
\\
\cdot\cdot & \cdot\cdot & \cdot\cdot & \cdots & \cdot\cdot
\\
1 & \theta_K & \theta_K^2 & \cdots & \theta_K^{K-1}
\end{array}\!
\right\vert
=
(-1)^{\frac{K(K-1)}{2}}
\Delta
\big(\theta_1,\theta_2,\dots,\theta_K\big)
\]
of a corresponding Van der Monde determinant. 

\begin{Theorem}
{\rm (\cite{bru1997})}
Let $d_1 = {\rm depth} ( {\sf T})$ be the depth of ${\sf T}$, which is
$\leqslant n+1$, let $r = \ell_1 + \cdots + \ell_{ d_1}$
be the number of boxes of ${\sf T}$, 
set:
\[
t_i
:=
r+\ell_i-i
\]
for all $i = 1, 2, \dots, n+1, n+2$ with of course:
\[
t_{d_1+1}=r-d_1-1,\,
\dots\dots,\,
t_{n+1}=r-n-1,\ \
t_{n+2}=r-n-2, 
\]
and define the rational number:
\[
b_0
:=
{\textstyle{\frac{1}{1!\,2!\,\cdots\,n!\,(n+1)!}}}
\cdot
\Delta
\big(
t_1,\dots,t_{n+1},t_{n+2}
\big),
\]
together with, for any $s = 1, 2, \dots, d_1$, the rational 
numbers:
\[
b_s
:=
{\textstyle{\frac{1}{1!\,2!\,\cdots\,n!\,(n+1)!}}}
\cdot
\sum_{1\leqslant i_1<\cdots<i_s\leqslant d_1}\,
\Delta
\big(
t_1,t_2,\dots,t_{i_1}-1,\dots,t_{i_s}-1,\dots,t_{n+1},t_{n+2}
\big). 
\]
Then there is a long exact sequence of holomorphic vector bundles
over $\P^{ n+1}$ of the form:
\[
\aligned
0
\longrightarrow
\Omega_{\P^{n+1}}^{(\ell_1,\dots,\ell_{d_1},0,\dots,0)}
\longrightarrow
\bigoplus_{b_0}
\mathcal{O}_{\P^{n+1}}(-r)
\longrightarrow
\bigoplus_{b_1}
\mathcal{O}_{\P^{n+1}}(-r+1)
\longrightarrow
\cdots
\\
\cdots
\longrightarrow
\bigoplus_{b_{d_1}}
\mathcal{O}_{\P^{n+1}}(-r+d_1)
\longrightarrow
0.
\endaligned
\]
\end{Theorem}

Then tensoring by $\mathcal{ O}_{ \P^{ n+1}} ( t)$ with an 
arbitrary $t \in \Z$, we get the useful: 

\begin{equation}
\label{first-exact-sequence}
\aligned
0
\longrightarrow
\Omega_{\P^{n+1}}^{(\ell_1,\dots,\ell_{d_1},0,\dots,0)}
\otimes
\mathcal{O}_{\P^{n+1}}(t)
\longrightarrow
\bigoplus_{b_0}
\mathcal{O}_{\P^{n+1}}(t-r)
\longrightarrow
\\
\longrightarrow
\bigoplus_{b_1}
\mathcal{O}_{\P^{n+1}}(t-r+1)
\longrightarrow
\cdots
\cdots
\longrightarrow
\bigoplus_{b_{d_1}}
\mathcal{O}_{\P^{n+1}}(t-r+d_1)
\longrightarrow
0.
\endaligned
\end{equation}

\subsection{Second fundamental (short) exact sequence}
Because any Schur bundle 
$\Omega_{ \P^{n+1}}^{ (\ell_1, \dots, \ell_n, \ell_{ n+1})}$ over
$\P^{ n+1}$ is, according to what precedes, a locally free sheaf of
$\mathcal{ O}_{\P^{ n+1}}$-modules of finite rank $\prod_{ 1 \leqslant
i < j \leqslant n+1}\, \big( \frac{ \ell_i - \ell_j}{ j - i} + 1
\big)$, a tensorisation of the normal exact sequence~\thetag{
\ref{normal-exact}} yields the general
short exact sequence
(\cite{ bru1972, rou2006a}):
\[
\aligned
0
\longrightarrow
\Omega_{\P^{n+1}}^{(\ell_1,\dots,\ell_n,\ell_{n+1})}
\otimes
\mathcal{O}_{\P^{n+1}}(-d)
&
\longrightarrow
\Omega_{\P^{n+1}}^{(\ell_1,\dots,\ell_n,\ell_{n+1})}
\otimes
\mathcal{O}_{\P^{n+1}}(0)
\longrightarrow
\\
&
\longrightarrow
\Omega_{\P^{n+1}}^{(\ell_1,\dots,\ell_n,\ell_{n+1})}
\otimes
\mathcal{O}_X(0)
\longrightarrow
0.
\endaligned
\]
Tensoring in addition again by $\mathcal{ O}_{ \P^{ n+1}} ( t)$ where
$t \in \Z$ is arbitrary, knowing $\mathcal{ O}_X ( 0 ) \otimes
\mathcal{ O}_{ \P^{ n+1}} ( t) = \mathcal{ O}_X ( t)$, we deduce the
general form of this (second, short) exact sequence that will be
useful below:
\begin{equation}
\label{second-exact-sequence}
\aligned
0
\longrightarrow
\Omega_{\P^{n+1}}^{(\ell_1,\dots,\ell_n,\ell_{n+1})}
\otimes
\mathcal{O}_{\P^{n+1}}(t-d)
\longrightarrow
\Omega_{\P^{n+1}}^{(\ell_1,\dots,\ell_n,\ell_{n+1})}
\otimes
\mathcal{O}_{\P^{n+1}}(t)
\longrightarrow
\\
\longrightarrow
\Omega_{\P^{n+1}}^{(\ell_1,\dots,\ell_n,\ell_{n+1})}
\otimes
\mathcal{O}_X(t)
\longrightarrow
0.
\endaligned
\end{equation}
Here, we make the convention similar to~\thetag{
\ref{convention-zero}} that when all
the $\ell_i$ are zero, $\Omega_{ \P^{ n+1}}^{ (0, \dots, 0, 0)}$
identifies to $\mathcal{ O}_{ \P^{ n+1}} ( 0)$, whence in this case
the written exact sequence reduces just to~\thetag{
\ref{normal-exact}}, tensored of course
by $\mathcal{ O}_{ \P^{ n+1}} ( t)$.

\subsection{Third fundamental exact sequence}
Lastly, starting from the cotangential normal exact sequence:
\begin{equation}
\label{O-X-d}
0
\longrightarrow
\mathcal{O}_X(-d)
\longrightarrow
T_{\P^{n+1}}^*\vert_X
\longrightarrow
T_X^*
\longrightarrow
0,
\end{equation}
(recall that $T_{ \P^{ n+1}}^* \vert_X = T_{\P^{ n+1}}^* \otimes
\mathcal{ O}_X ( 0)$), Br\"uckmann established that the Schur power of
the extracted complex:
\[
0
\longrightarrow
\mathcal{O}_X(-d)
\longrightarrow
T_{\P^{n+1}}^*\otimes\mathcal{O}_X(0)
\longrightarrow
0
\longrightarrow
0
\longrightarrow
\cdots
\]
provides a free resolution of $\mathcal{ S}^{ (\ell_1, \dots, 
\ell_n)} T_X^*$ (Theorem~ 3 in~\cite{ bru1997})
which may be written in great details as follows
when $\ell_n \geqslant 1$:
\[
\aligned
0
&
\longrightarrow
\bigoplus_{\delta_1+\cdots+\delta_n=n
\atop
\delta_i\,=\,0\,\,\text{\rm or}\,\,1}\,
\Omega_{\P^{n+1}}^{
(\ell_1,\dots,\ell_n,0)-(\delta_1,\dots,\delta_n,0)}
\otimes
\mathcal{O}_X(-nd)
\longrightarrow
\cdots
\\
\cdots
&
\longrightarrow
\bigoplus_{\delta_1+\cdots+\delta_n=k
\atop
\delta_i\,=\,0\,\,\text{\rm or}\,\,1}\,
\Omega_{\P^{n+1}}^{
(\ell_1,\dots,\ell_n,0)-(\delta_1,\dots,\delta_n,0)}
\otimes
\mathcal{O}_X(-kd)
\longrightarrow
\cdots
\\
\cdots
&
\longrightarrow
\Omega_{\P^{n+1}}^{(\ell_1,\dots,\ell_n,0)}
\otimes
\mathcal{O}_X(0)
\longrightarrow
\mathcal{S}^{(\ell_1,\dots,\ell_n)}T_X^*
\longrightarrow
0.
\endaligned
\]
Notice that the last upper entry $\ell_{ n+1}$ of each $\Omega$ is
zero. Of course, the direct sum for the first entry reduces just to
the single term:
\[
\Omega_{\P^{n+1}}^{(\ell_1-1,\dots,\ell_n-1,0)}
\otimes
\mathcal{O}_X(-nd).
\]
In full generality, if $d_1$ denotes the depth of the considered Young
diagram, hence if one has $\ell_1 \geqslant \cdots \geqslant \ell_{
d_1} \geqslant 1$ but $0 = \ell_{ d_1 + 1} = \cdots = \ell_n = \ell_{
n+1}$, the locally free resolution of $\mathcal{ S}^{ ( \ell_1, \dots,
\ell_{ d_1}, 0, \dots, 0)} T_X^*$ reads (\cite{ bru1997}):
\[
\aligned
0
&
\longrightarrow
\bigoplus_{\delta_1+\cdots+\delta_{d_1}=d_1
\atop
\delta_i\,=\,0\,\,\text{\rm or}\,\,1}\,
\Omega_{\P^{n+1}}^{
(\ell_1,\dots,\ell_{d_1},0,\dots,0,0)
-
(\delta_1,\dots,\delta_{d_1},0,\dots,0,0)}
\otimes
\mathcal{O}_X(-d_1d)
\longrightarrow
\cdots
\\
\cdots
&
\longrightarrow
\bigoplus_{\delta_1+\cdots+\delta_{d_1}=k
\atop
\delta_i\,=\,0\,\,\text{\rm or}\,\,1}\,
\Omega_{\P^{n+1}}^{
(\ell_1,\dots,\ell_{d_1},0,\dots,0,0)
-
(\delta_1,\dots,\delta_{d_1},0,\dots,0,0)}
\otimes
\mathcal{O}_X(-kd)
\longrightarrow
\cdots
\\
\cdots
&
\longrightarrow
\Omega_{\P^{n+1}}^{(\ell_1,\dots,\ell_{d_1},0,\dots,0,0)}
\otimes
\mathcal{O}_X(0)
\longrightarrow
\mathcal{S}^{(\ell_1,\dots,\ell_{d_1},0,\dots,0)}T_X^*
\longrightarrow
0,
\endaligned
\]
hence it just looks like a truncation of the preceding resolution
valid when $d_1 = n$.  Tensoring this by $\mathcal{ O}_X ( t)$, we
finally get what will be useful below:
\begin{equation}
\label{third-exact-sequence}
\aligned
0
&
\longrightarrow
\bigoplus_{\delta_1+\cdots+\delta_{d_1}=d_1
\atop
\delta_i\,=\,0\,\,\text{\rm or}\,\,1}\,
\Omega_{\P^{n+1}}^{
(\ell_1,\dots,\ell_{d_1},0,\dots,0,0)
-
(\delta_1,\dots,\delta_{d_1},0,\dots,0,0)}
\otimes
\mathcal{O}_X(t-d_1d)
\longrightarrow
\cdots
\\
\cdots
&
\longrightarrow
\bigoplus_{\delta_1+\cdots+\delta_{d_1}=k
\atop
\delta_i\,=\,0\,\,\text{\rm or}\,\,1}\,
\Omega_{\P^{n+1}}^{
(\ell_1,\dots,\ell_{d_1},0,\dots,0,0)
-
(\delta_1,\dots,\delta_{d_1},0,\dots,0,0)}
\otimes
\mathcal{O}_X(t-kd)
\longrightarrow
\cdots
\\
\cdots
&
\longrightarrow
\Omega_{\P^{n+1}}^{(\ell_1,\dots,\ell_{d_1},0,\dots,0,0)}
\otimes
\mathcal{O}_X(t)
\longrightarrow
\mathcal{S}^{(\ell_1,\dots,\ell_{d_1},0,\dots,0)}T_X^*
\otimes
\mathcal{O}_X(t)
\longrightarrow
0.
\endaligned
\end{equation}

\subsection{Cohomology of Schur bundles over 
$\P^{ n+1}$} In~\cite{ bru1997} too, using the first exact sequence
above plus further arguments, Br\"uckmann established the following
theorem which computes completely the dimensions of all the cohomology
groups of twisted Schur bundles over $\P^{ n+1}$.  As above, for fixed
$n+1 \geqslant 2$ and for fixed $\ell_1 \geqslant \cdots \geqslant
\ell_n \geqslant \ell_{ n+1} \geqslant 0$, we introduce the integers:
\[
t_i
:=
\ell_i
-
i
+
{\textstyle{\sum_{i=1}^{n+1}}}\,\ell_i
\ \ \ \ \ \ \ \ \ \ \ \ \ 
{\scriptstyle{(i\,=\,1\,\cdots\,n,\,\,n\,+\,1)}},
\]
which, visibly, are ordered decreasingly:
\[
t_1>t_2>\cdots>t_n>t_{n+1}.
\] 

\begin{Theorem}
\label{cohomology-Schur-P-n-1}
{\rm (\cite{ bru1997})} 
For any $t \in \Z$, the Euler-Poincar\'e characteristic of $\Omega_{
\P^{ n+1}}^{(\ell_1, \dots, \ell_n, \ell_{ n+1})}
\otimes \mathcal{ O}_{ \P^{ n+1}} ( t)$ is equal to:
\[
\aligned
\chi(t)
:=
&\,
\chi
\big(
\P^{n+1},\,\Omega_{\P^{n+1}}^{(\ell_1,\dots,\ell_n,\ell_{n+1})}
\!\otimes\!
\mathcal{O}_{\P^{n+1}}(t)
\big)
\\
=
&\,
{\textstyle{\frac{1}{1!\,2!\,\cdots\,n!\,(n+1)!}}}
\!\!
\prod_{1\leqslant i<j\leqslant n+1}
\!\!
(t_i-t_j)
\prod_{1\leqslant i\leqslant n}
\!
(t-t_i),
\endaligned
\]
whence it vanishes for $t$ equal to each one of the $t_i$.
Furthermore, as $t$ varies in $\Z$, at most one of the cohomology
dimensions:
\[
h^q(t)
:=
\dim
H^q\big( 
\P^{n+1},\,
\Omega_{\P^{n+1}}^{(\ell_1,\dots,\ell_n,\ell_{n+1})}
\otimes
\mathcal{O}_{\P^{n+1}}(t)\big)
\]
may be nonzero, and more precisely, $h^q(t)$ is nonzero and equal
to $(-1)^q \chi ( t)$ if and only if $t_{ q+1} + 1 \leqslant t
\leqslant t_{ q} -1$, while the other $h^{ q'} ( t)$ do vanish
for all $t$ in the same range. In
particular, for all:
\begin{equation}
\label{t-ell-sum-ell}
t\geqslant
\ell_1
+
{\textstyle{\sum_{i=1}^{n+1}}}\,\ell_i, 
\end{equation}
all the positive cohomology dimensions vanish:
\[
0
=
\dim H^q\big( 
\P^{n+1},\,
\Omega_{\P^{n+1}}^{(\ell_1,\dots,\ell_n,\ell_{n+1})}
\otimes
\mathcal{O}_{\P^{n+1}}(t)\big)
\ \ \ \ \ \ \ \ \ \ \ \ \ 
{\scriptstyle{(q\,=\,1,\,\,2\,\cdots\,n)}}.
\]
\end{Theorem}

\subsection{Applications}
For an application to the study of the cohomology of Schur bundles
over $X^n \subset \P^{ n+1}$, we shall apply the above theorems
specifically to the Young diagrams ${\sf YD}_{ (\ell_1, \dots, \ell_n,
0)}$ of depth $d_1 \leqslant n = \dim X$, with $\ell_{ n+1} = 0$
in order to gain the following complement
to Theorem~\ref{majoration-Rousseau}.

\begin{Theorem}
\label{Schur-vanishing}
Let $X = X^n \subset \P^{ n+1} ( \C)$ be a geometrically smooth
projective algebraic complex hypersurface of general type, i.e. of
degree $d \geqslant n+3$, and let $\ell = ( \ell_1, \dots, \ell_{ n-1},
\ell_n)$ with $\ell_1 \geqslant \cdots \geqslant \ell_{ n-1} \geqslant
\ell_n \geqslant 1$. If:
\[
\ell_n
\geqslant
{\textstyle{\frac{1}{d-n-2}}}
\big\{
n(d-1)
+
\ell_1-\ell_n
+
{\textstyle{\sum_{i=1}^{n-1}}}\,
(\ell_i-\ell_n)
\big\},
\]
then all the positive cohomologies vanish:
\[
0
=
H^q
\big(
X,\,
\mathcal{S}^{(\ell_1,\dots,\ell_{n-1},\ell_n)}T_X^*
\big)
\ \ \ \ \ \ \ \ \ \ \ \ \ 
{\scriptstyle{(q\,=\,1,\,\,2\,\cdots\,n)}}.
\]
\end{Theorem}

\proof
As anticipated above, after division by $(K_X)^{ \ell_n}$, we are lead
back to examining the cohomology of:
\[
\mathcal{S}^{(\ell_1-\ell_n,\dots,\ell_{n-1}-\ell_n,0)}T_X^*
\otimes
\mathcal{O}_X\big(\ell_n(d-n-2)\big).
\]
A bit more generally, using the second and the third exact 
sequences~\thetag{ \ref{second-exact-sequence}}
and~\thetag{ \ref{third-exact-sequence}}, 
we shall examine when the positive cohomologies of:
\[
\mathcal{S}^{(\ell_1',\dots,\ell_{n-1}',0)}T_X^*
\otimes
\mathcal{O}_X(t')
\]
do all vanish, and afterward, we shall set:
\[
\ell_1'
:=
\ell_1-\ell_n,
\,\,\dots,\,\,
\ell_{n-1}'
:=
\ell_{n-1}-\ell_n
\ \ \ \ \
\text{\rm and}
\ \ \ \ \ 
t'
:=
\ell_n(d-n-2).
\] 

We assume first that $\ell_{ n-1}' \geqslant 1$ and we shall discuss
the quite similar case $\ell_{ n-1}' = 0$ afterward.  The
consideration of the third exact sequence~\thetag{
\ref{third-exact-sequence}} with $d_1 = n-1$ and $(\ell_1', \dots, 
\ell_{ n-1}', 0)$ instead of $(\ell_1, \dots, \ell_{ n-1}, 0)$
then necessarily conducts us to 
the study of $\mathcal{ O}_X$-twisted Schur bundles over $\P^{ n+1}$:
\[
\Omega_{\P^{n+1}}^{(\ell_1'',\dots,\ell_{n-1}'',0)}
\otimes
\mathcal{O}_X(t'')
\]
whose Young diagram exponents $\ell_i''$ have values:
\[
(\ell_1'',\dots,\ell_{n-1}'',0)
=
(\ell_1',\dots,\ell_{n-1}',0)
-
(\delta_1',\dots,\delta_{n-1}',0)
\]
shifted a bit from the values of the $\ell_i'$, where $\delta_1 ' +
\cdots + \delta_{ n-1}' = k$ for $k = 0, 1, \dots, n - 1$, with of
course $\delta_i' = 0$ or $1$.  So to begin with, it is advisable to
study the cohomology of these $\mathcal{ O}_X$-twisted Schur bundles
over $\P^{ n+1}$.

To this aim, we look at the
second (short) exact sequence~\thetag{ \ref{second-exact-sequence}}
written with:
\[
(\ell_1,\dots,\ell_{n-1},\ell_n,\ell_{n+1})
:=
(\ell_1'',\dots,\ell_{n-1}'',0,0),
\]
for some arbitrary $\ell_1 '' \geqslant \cdots
\geqslant \ell_{ n-1}'' \geqslant 0$ and we abbreviate this
exact sequence as:
\[
0
\longrightarrow
\mathcal{P}
\longrightarrow
\mathcal{Q}
\longrightarrow
\mathcal{R}
\longrightarrow
0,
\]
where $\mathcal{ P} \to \P^{ n+1}$, $\mathcal{Q} \to \P^{ n+1}$ and
$\mathcal{R} \to X$ are the bundles:
\[
\aligned
\mathcal{P}
&
:=
\Omega_{\P^{n+1}}^{(\ell_1'',\dots,\ell_{n-1}'',0,0)}
\otimes
\mathcal{O}_{\P^{n+1}}(t''-d),
\\
\mathcal{Q}
&
:=
\Omega_{\P^{n+1}}^{(\ell_1'',\dots,\ell_{n-1}'',0,0)}
\otimes
\mathcal{O}_{\P^{n+1}}(t''),
\\
\mathcal{R}
&
:=
\Omega_{\P^{n+1}}^{(\ell_1'',\dots,\ell_{n-1}'',0,0)}
\otimes
\mathcal{O}_X(t''),
\endaligned
\]
so that all the cohomology dimensions of $\mathcal{ P}$ 
and of $\mathcal{ Q}$ are known thanks to 
Br\"uckmann's Theorem~\ref{cohomology-Schur-P-n-1}. 
Then in the long exact cohomology sequence associated to 
the short exact sequence:
\[
\aligned
0
&
\longrightarrow
H^0(\P^{n+1},\,\mathcal{P})
\longrightarrow
H^0(\P^{n+1},\,\mathcal{Q})
\longrightarrow
H^0(X,\,\mathcal{R})
\longrightarrow
\\
&
\longrightarrow
\zero{H^1(\P^{n+1},\,\mathcal{P})}
\longrightarrow
\zero{H^1(\P^{n+1},\,\mathcal{Q})}
\longrightarrow
H^1(X,\,\mathcal{R})
\longrightarrow
\\
&
\longrightarrow
\zero{H^2(\P^{n+1},\,\mathcal{P})}
\longrightarrow
\zero{H^2(\P^{n+1},\,\mathcal{Q})}
\longrightarrow
H^2(X,\,\mathcal{R})
\longrightarrow
\cdots
\\
\cdots
&
\longrightarrow
\zero{H^n(\P^{n+1},\,\mathcal{P})}
\longrightarrow
\zero{H^n(\P^{n+1},\,\mathcal{Q})}
\longrightarrow
H^n(X,\,\mathcal{R})
\longrightarrow
\\
\cdots
&
\longrightarrow
\zero{H^{n+1}(\P^{n+1},\,\mathcal{P})}
\longrightarrow
\zero{H^{n+1}(\P^{n+1},\,\mathcal{Q})}
\longrightarrow
0
\endaligned
\]
(the last $0$ because $\mathcal{R} \to X$ is a bundle over an
$n$-dimensional basis), all the underlined terms will vanish, 
namely one will have:
\[
\aligned
0
&
=
H^q
\big(
\P^{n+1},\,\,
\Omega_{\P^{n+1}}^{(\ell_1'',\dots,\ell_{n-1}'',0,0)}
\otimes
\mathcal{O}_{\P^{n+1}}(t''-d)
\big)
\ \ \ \ \ \ \ \ \ \ \ \ \ 
{\scriptstyle{(q\,=\,1,\,\,2\,\cdots\,n,\,\,n\,+\,1)}},
\\
0
&
=
H^q
\big(
\P^{n+1},\,\,
\Omega_{\P^{n+1}}^{(\ell_1'',\dots,\ell_{n-1}'',0,0)}
\otimes
\mathcal{O}_{\P^{n+1}}(t'')
\big)
\ \ \ \ \ \ \ \ \ \ \ \ \ 
{\scriptstyle{(q\,=\,1,\,\,2\,\cdots\,n,\,\,n\,+\,1)}},
\endaligned
\]
as soon as the following
two inequalities are satisfied by $t''$:
\[
\aligned
t''-d
&
\geqslant
\ell_1''
+
{\textstyle{\sum_{i=1}^{n-1}}}\,\ell_i'',
\\
t''
&
\geqslant
\ell_1''
+
{\textstyle{\sum_{i=1}^{n-1}}}\,\ell_i'',
\endaligned
\]
as is guaranteed by the
inequality~\thetag{ \ref{t-ell-sum-ell}} of 
Theorem~\ref{cohomology-Schur-P-n-1}.
But the first inequality obviously entails the second one, 
hence we deduce that all positive cohomologies:
\[
0
=
H^q
\big(
X,\,\,
\Omega_{\P^{n+1}}^{
(\ell_1'',\dots,\ell_{n-1}'',0,0)}
\otimes
\mathcal{O}_X(t'')
\big)
\ \ \ \ \ \ \ \ \ \ \ \ \ 
{\scriptstyle{(q\,=\,1,\,\,2\,\cdots\,n)}}
\]
of $\mathcal{ R}$ vanish as soon as:
\begin{equation}
\label{inequality-t-second}
t''
\geqslant
d
+
\ell_1''
+
{\textstyle{\sum_{i=1}^{n-1}}}\,\ell_i''.
\end{equation}
We observe that this fact is valid also when 
$\ell_{ n_1 + 1}'' = \cdots = \ell_{ n - 1} ''$ for 
some largest integer $n_1 \geqslant 0$ with 
$\ell_{ n_1}'' \geqslant 1$, because the second
exact sequence~\thetag{ \ref{second-exact-sequence}} we used
is subjected to no restriction.

\smallskip

We now come to dealing with the third exact sequence~\thetag{
\ref{third-exact-sequence}}. Cutting a long exact sequence in short
exact sequences, one may establish the following standard lemma, 
used e.g. in~\cite{ rou2006b}.

\begin{Lemma}
Consider a holomorphic vector bundle $\mathcal{ S} \to X$ equipped
with a free resolution of length $\leqslant n$ provided by a long
exact sequence of holomorphic vector bundles $\mathcal{ A}^0,
\mathcal{ A}^1, \dots, \mathcal{ A}^n$ over $X$:
\[
0
\longrightarrow
\mathcal{A}^n
\longrightarrow
\mathcal{A}^{n-1}
\longrightarrow
\cdots
\longrightarrow
\mathcal{A}^1
\longrightarrow
\mathcal{A}^0
\longrightarrow
\mathcal{S}
\longrightarrow
0.
\]
Then in order that all the positive cohomology groups vanish:
\[
0
=
H^1(X,\,\mathcal{S})
=
\cdots
=
H^n(X,\,\mathcal{S}),
\]
it suffices that:
\[
\aligned
0
=
H^1(X,\,\mathcal{A}^0)
=
H^2(X,\,\mathcal{A}^1)
=
H^3(X,\,\mathcal{A}^2)
=
\cdots
&
=
H^n(X,\,\mathcal{A}^{n-1})
\\
0
=
H^2(X,\,\mathcal{A}^0)
=
H^3(X,\,\mathcal{A}^1)
=
\cdots
&
=
H^n(X,\,\mathcal{A}^{n-2})
\\
0
=
H^3(X,\,\mathcal{A}^0)
=
\cdots
&
=
H^n(X,\,\mathcal{A}^{n-3})
\\
\cdots\cdots
&
\cdots\cdots\cdots\cdots\cdots
\\
0
&
=
H^n(X,\,\mathcal{A}^0).
\endaligned
\]
\end{Lemma}

So as said a short while ago, we aim to apply this lemma when looking
at the third exact sequence~\thetag{ \ref{third-exact-sequence}}
which, for the case we are interested in, writes precisely under the
form:
\[
\aligned
0
&
\longrightarrow
\bigoplus_{\delta_1'+\cdots+\delta_{n-1}'=n-1
\atop
\delta_i'\,=\,0\,\,\text{\rm or}\,\,1}\,
\Omega_{\P^{n+1}}^{
(\ell_1',\dots,\ell_{n-1}',0,0)
-
(\delta_1',\dots,\delta_{n-1}',0,0)}
\otimes
\mathcal{O}_X\big(t'-(n-1)d\big)
\longrightarrow
\cdots
\\
\cdots
&
\longrightarrow
\bigoplus_{\delta_1'+\cdots+\delta_{n-1}'=k
\atop
\delta_i'\,=\,0\,\,\text{\rm or}\,\,1}\,
\Omega_{\P^{n+1}}^{
(\ell_1',\dots,\ell_{n-1}',0,0)
-
(\delta_1',\dots,\delta_{n-1}',0,0)}
\otimes
\mathcal{O}_X\big(t'-kd\big)
\longrightarrow
\cdots
\\
\cdots
&
\longrightarrow
\Omega_{\P^{n+1}}^{(\ell_1',\dots,\ell_{n-1}',0,0)}
\otimes
\mathcal{O}_X(t')
\longrightarrow
\mathcal{S}^{(\ell_1',\dots,\ell_{n-1}',0)}
\mathcal{O}_X(t')
\longrightarrow
0.
\endaligned
\]
In the notations of the lemma, 
the resolution of:
\[
\mathcal{S}
:= 
\mathcal{S}^{(\ell_1',\dots,\ell_{n-1}',0)} 
\otimes 
\mathcal{O}_X(t') 
\]
is hence of length $n-1$ when we set:
\[
\aligned
\mathcal{A}^k
:=
\bigoplus_{\delta_1'+\cdots+\delta_{n-1}'=k
\atop
\delta_i'\,=\,0\,\,\text{\rm or}\,\,1}\,
&
\Omega_{\P^{n+1}}^{
(\ell_1',\dots,\ell_{n-1}',0,0)
-
(\delta_1',\dots,\delta_{n-1}',0,0)}
\otimes
\mathcal{O}_X\big(t'-kd\big)
\\
&
\ \ \ \ \ \ \ \ \ \ 
{\scriptstyle{(k\,=\,0,\,\,1\,\cdots\,n\,-\,1)}}.
\endaligned
\]
Then for the lemma to yield the vanishing of all the positive
cohomologies of $\mathcal{ S} = \mathcal{ S}^{ (\ell_1', \dots, \ell_{
n-1}', 0)} \otimes \mathcal{ O}_X ( t')$, it is evidently sufficient
that plainly all positive cohomologies of the $\mathcal{ A}^k$ vanish:
\[
0
=
H^q(X,\,\mathcal{A}^k)
\ \ \ \ \ \ \ \ \ \ \ \ \ 
{\scriptstyle{(q\,=\,1,\,\,2\,\cdots\,n\,;\,\,\,
k\,=\,0,\,\,1\,\cdots\,n\,-\,1)}}, 
\]
which is more than what is required in fact. But since each
$\mathcal{ A}^k$ is a direct sum, it even suffices that:
\[
\aligned
0
&
=
H^q
\big(
X,\,\,
\Omega_{\P^{n+1}}^{
(\ell_1',\dots,\ell_{n-1}',0,0)
-
(\delta_1',\dots,\delta_{n-1}',0,0)}
\otimes
\mathcal{O}_X(t'-kd)
\big)
\\
&
\ \ \ \ \ \ \ \ \ \ \ \ \ 
{\scriptstyle{(q\,=\,1,\,\,2\,\cdots\,n\,;\,\,\,
\delta_1'\,+\,\cdots\,+\,\delta_{n-1}'\,=\,k\,;\,\,\,
k\,=\,0,\,\,1\,\cdots\,n\,-\,1)}}.
\endaligned
\]
According to~\thetag{ \ref{inequality-t-second}}, this
holds true provided all the following inequalities are
satisfied:
\[
t'-kd
\geqslant
d
+
\ell_1'-\delta_1'
+
{\textstyle{\sum_{i=1}^{n-1}}}\,
(\ell_i'-\delta_i'),
\]
for every $k = 0, 1, \dots, n-1$ and every $\delta_1', \dots, \delta_{
n-1}' \in \{ 0, 1\}$ with $\delta_1' + \cdots + \delta_{ n-1}' =
k$. But since $- \delta_i' \leqslant 1$ and since $\sum ( - \delta_i')
= - k$, it suffices that, firstly for $k = 0, 1, \dots, n-2$:
\[
t'-kd
\geqslant
d
+
\ell_1'
+
{\textstyle{\sum_{i=1}^{n-1}}}\,
\ell_i'
-
k,
\]
and secondly for $k = n-1$, whence $- \delta_1' = - 1$ necessarily:
\begin{equation}
\label{t-prime-n-1}
t'-(n-1)d
\geqslant
d
+
\ell_1'-1
+
{\textstyle{\sum_{i=1}^{n-1}}}\,
\ell_i'
-
(n-1).
\end{equation}
But this last inequality, rewritten under the form:
\[
t'
\geqslant
n(d-1)
+
\ell_1'
+
{\textstyle{\sum_{i=1}^{n-1}}}\,
\ell_i'
\]
visibly entails all the inequalities for $k = 0, 1, \dots, n-2$.
Lastly, replacing $t' = \ell_n ( d - n - 2)$ and the $\ell_i' = \ell_i
- \ell_n$ by their values, we finally come to the numerical condition
claimed by the theorem for positive cohomologies of the $\mathcal{
S}^{ ( \ell_1, \dots, \ell_{ n-1}, \ell_n)} T_X^*$ to vanish.

To conclude the argument, it only remains to examine what happens with
the case, left aside, when $\ell_{ n-1}' = 0$.  In this case, there is
a nonnegative integer $n_1 \leqslant n-2$ with $\ell_1 ' \geqslant
\cdots \geqslant \ell_{ n_1}' \geqslant 1$ while $0 = \ell_{ n_1 +1}'
= \cdots = \ell_{ n-1}'$.  At first, if $n_1 = 0$, i.e. if all the
$\ell_i$ are equal to $\ell_n$, then $\mathcal{ S}^{ ( \ell_n, \dots,
\ell_n)} T_X^* = \mathcal{ O}_X \big( \ell_n ( d - n - 2)\big)$
reduces to a standard
line bundle $\mathcal{ O}_X ( t')$, and it is well known that:
\[
0
=
H^q
\big(
X,\,
\mathcal{O}_X(t')
\big)
\ \ \ \ \ \ \ \ \ \ \ \ \ 
{\scriptstyle{(q\,=\,1,\,\,2\,\cdots\,n)}}
\]
whenever $t' \geqslant 0$. 

Therefore, we may assume that $n_1$ satisfies $1
\leqslant n_1 \leqslant n-2$. As before, the subtraction of $(K_X)^{
\ell_n}$ yields:
\[
\ell_1'
=
\ell_1-\ell_n,
\,\,\dots,\,\,
\ell_{n_1'}
=
\ell_{n_1}-\ell_n
\ \ \ \ \
\text{\rm and}
\ \ \ \ \
0
=
\ell_{n_1+1}'
=\cdots=
\ell_{n-1}'
=
\ell_n',
\]
and again as always $t' = \ell_n ( d - n - 2)$. In the third
exact sequence, the factors then are:
\[
\aligned
\mathcal{A}^k
=
\bigoplus_{\delta_1'+\cdots+\delta_{n_1}'=k
\atop
\delta_i'\,=\,0\,\,\text{\rm or}\,\,1}\,
&
\Omega_{\P^{n+1}}^{
(\ell_1',\dots,\ell_{n_1}',0,\dots,0,0)
-
(\delta_1',\dots,\delta_{n_1}',0,\dots,0,0)}
\otimes
\mathcal{O}_X(t'-kd)
\\
&
\ \ \ \ \ \ \ \ \ \ \ \ \ 
{\scriptstyle{(k\,=\,0,\,\,1\,\cdots\,n_1)}},
\endaligned
\]
so the positive cohomologies vanish all provided that:
\[
t'-kd
\geqslant
d
+
\ell_1'-\delta_1'
+
{\textstyle{\sum_{i=1}^{n_1}}}\,
(\ell_i'-\delta_i')
\ \ \ \ \ \ \ \ \ \ \ \ \ 
{\scriptstyle{(k\,=\,0,\,\,1\,\cdots\,n_1\,;\,\,\,
\delta_1'\,+\,\cdots\,+\,\delta_{n_1}'\,=\,n_1)}},
\]
and because $k\leqslant n_1 \leqslant n-2$, these inequalities are all
less stringent than the one~\thetag{ \ref{t-prime-n-1}} we found
previously in the case when $\ell_{ n-1}' \geqslant 1$ (or
equivalently, when $n_1 = n-1$). This completes the proof of the
theorem.
\endproof

\markleft{Jo\"el Merker}
\markright{\S12.~Asymptotic cohomology vanishing}
\section{\bf Asymptotic cohomology vanishing}
\label{Section-12}

\subsection{Synthesis: uniform majoration for the cohomology
of Schur bundles} Two cohomology controls have been achieved. 
Firstly, according to Theorem~\ref{Schur-vanishing} stated above
and just proved, when:
\[
\ell_n 
\geqslant
{\textstyle{\frac{1}{d-n-2}}}\,
\big\{
n(d-1)+\ell_1-\ell_n
+
{\textstyle{\sum_{i=1}^{n-1}}}\,(\ell_i-\ell_n)
\big\},
\]
the positive cohomologies of Schur bundles vanish:
\[
0
=
h^q
\big(X,\,
\mathcal{S}^{(\ell_1,\ell_2,\dots,\ell_n)}T_X^*
\big)
\ \ \ \ \ \ \ \ \ \ \ \ \ 
{\scriptstyle{(q\,=\,1,\,2\,\cdots\,n)}}.
\]
Secondly, according to Theorem~\ref{majoration-Rousseau}, when:
\[
\vert\ell\vert
\geqslant
1+2n^2+(n+1)(d-n-2),
\]
the positive cohomologies enjoy a majoration of the shape:
\[
\footnotesize
\aligned
h^q
\big(X,\,
\mathcal{S}^{(\ell_1,\ell_2,\dots,\ell_n)}T_X^*
\big)
&
\leqslant
{\sf Constant}_{n,d}\cdot
\prod_{1\leqslant i<j\leqslant n}\,
(\ell_i-\ell_j)\,
\bigg\{
\\
&
\ \ \ \ \ \ \ \ \ \
\bigg\{
\sum_{\beta_1+\cdots+\beta_{n-1}+\beta_n=n}\,
(\ell_1-\ell_2)^{\beta_1}
\cdots\,
(\ell_{n-1}-\ell_n)^{\beta_{n-1}}
\ell_n^{\beta_n}
\bigg\}
+
\\
&
\ \ \ \ \
+
{\sf Constant}_{n,d}
\big(1+\vert\ell\vert^{\frac{n(n+1)}{2}-1}\big)
\ \ \ \ \ \ \ \ \ \ \ \ \ 
{\scriptstyle{(q\,=\,1,\,2\,\cdots\,n)}}.
\endaligned
\]
But then we may assume here that:
\[
\ell_n
<
{\textstyle{\frac{1}{d-n-2}}}\,
\big\{
n(d-1)+\ell_1-\ell_n
+
{\textstyle{\sum_{i=1}^{n-1}}}\,(\ell_i-\ell_n)
\big\},
\]
since otherwise the right-hand side majorant can be replaced by $0$,
and consequently, because it follows by exponentiation from such a
restriction on $\ell_n$ that:
\[
\ell_n^{\beta_n}
\leqslant
{\sf Constant}_{n,d}
\cdot
\sum_{\beta_1'+\cdots+\beta_{n-1}'\leqslant\beta_n}\,
(\ell_1-\ell_2)^{\beta_1'}
\cdots\,
(\ell_{n-1}-\ell_n)^{\beta_{n-1}'},
\]
we conclude that whenever $\vert \ell \vert \geqslant 1 + 2n^2 + (
n+1) ( d - n - 2)$, one has:
\[
\boxed{
\footnotesize
\aligned
&
h^q
\big(X,\,
\mathcal{S}^{(\ell_1,\dots,\ell_n)}T_X^*
\big)
\leqslant
\\
&
\leqslant
{\sf Constant}_{n,d}\cdot
\!\!\!
\prod_{1\leqslant i<j\leqslant n}\,
(\ell_i-\ell_j)
\bigg[
\sum_{\beta_1'+\cdots+\beta_{n-1}'=n}
\!\!\!
(\ell_1-\ell_2)^{\beta_1'}
\cdots\,
(\ell_{n-1}-\ell_n)^{\beta_{n-1}'}
\bigg]
+
\\
&
\ \ \ \ \
+
{\sf Constant}_{n,d}
\big(1+\vert\ell\vert^{\frac{n(n+1)}{2}-1}\big)
\ \ \ \ \ \ \ \ \ \ \ \ \ \ \ \ \ \ \ \ \ \ \ \ \ \ 
{\scriptstyle{(q\,=\,1,\,2\,\cdots\,n)}}.
\endaligned}
\]

\subsection{Application: cohomology control for
$\mathcal{ E}_{ \kappa, m}^{ GG} T_X^*$} Now, we make the following
observation: no Schur bundle $\mathcal{ S}^{ ( \ell_1, \dots, \ell_n)}
T_X^*$ for which $\vert \ell \vert < \frac{ m}{ \kappa}$ can appear in
the decomposition of ${\sf Gr}^\bullet \mathcal{ E}_{ \kappa, m}^{ GG}
T_X^*$ provided by Theorem~\ref{exact-decomposition},
just because all the integers $\lambda_i^j$ filling the Young diagram
${\sf YD}_{ (\ell_1, \dots, \ell_n)}$ satisfy all $1 \leqslant
\lambda_i^j \leqslant \kappa$, whence:
\[
\vert\ell\vert
\leqslant
m
\leqslant
\kappa\,\vert\ell\vert
\]
always. Thus, if we assume only that $\frac{ m}{ \kappa}$ is larger
than the above constant $1 + 2n^2 + (n+1) ( d - n -2)$, and we
certainly can assume this since both $m \gg \kappa$ and $\kappa \gg n$
are supposed to tend to infinity, then the cohomology majoration boxed
above can be applied to {\em all} Schur bundles entering the
decomposition of ${\sf Gr}^\bullet \mathcal{ E}_{ \kappa, m}^{ GG}
T_X^*$.

We are thus now in a position to accomplish the final series of
inequalities. For any $q = 1, 2, \dots, n$, reminding Sections~8, 9
and 10, we have:
\[
\footnotesize
\aligned
&
h^q
\big(X,\,
\mathcal{E}_{\kappa,m}^{GG}T_X^*
\big)
\leqslant
\sum_{\ell_1\geqslant\ell_2\geqslant\cdots\geqslant\ell_n\geqslant 0}\,
M_{\ell_1,\ell_2,\dots,\ell_n}^{\kappa,m}
\cdot
h^q\big(X,\,
\mathcal{S}^{(\ell_1,\ell_2,\dots,\ell_n)}T_X^*
\big)
\\
&
\leqslant
{\sf Constant}_{n,d}\,
\sum_{{\sf YT}\,{\sf semi-standard}
\atop
{\sf weight}({\sf YT})=m}\,
\prod_{1\leqslant i<j\leqslant n}
\big(\ell_i({\sf YT})-\ell_j({\sf YT})\big)\,
\bigg\{
\\
&
\ \ \ \ \ \ \ \ \ \ \ \ \ \ \ \ \ \ \ \ \ \ \ \ \ \
\ \ \ \ \ \ \ \ \ \ \ \ \
\bigg\{
\sum_{\beta_1'+\cdots+\beta_{n-1}'=n}\,
\big(\ell_1({\sf YT})-\ell_2({\sf YT}\big)^{\beta_1'}\cdots\,
\big(\ell_{n-1}({\sf YT})-\ell_n({\sf YT})\big)^{\beta_{n-1}'}
\bigg\}
+
\\
&
\ \ \ \ \ \ \ \ \ \ \ \ \
+
{\sf Constant}_{n,d}\,
\sum_{{\sf YT}\,{\sf semi-standard}
\atop
{\sf weight}({\sf YT})=m}\,
\sum_{\alpha_1+\cdots+\alpha_n\leqslant
\frac{n(n+1)}{2}-1}\,
\ell_1({\sf YT})^{\alpha_1}\cdots\,
\ell_n({\sf YT})^{\alpha_n}
\\
&
\leqslant
{\sf Constant}_{n,d}\,
\sum_{{\sf YT}\in{\sf YT}_{\kappa,m}^{\rm max}}\,
\sum_{\alpha_1'+\cdots+\alpha_{n-1}'=\frac{n(n+1)}{2}}\,
\\
&
\ \ \ \ \ \ \ \ \ \ \ \ \
\big(\ell_1({\sf YT})-\ell_2({\sf YT})\big)^{\alpha_1'}
\cdots\,
\big(\ell_{n-1}({\sf YT})-\ell_n({\sf YT})\big)^{\alpha_{n-1}'}
+
{\sf Constant}_{n,d,\kappa}\cdot m^{(\kappa+1)n-2}
\endaligned
\]
\[
\footnotesize
\aligned
&
\leqslant
{\sf Constant}_{n,d}\,
{\textstyle{\frac{m^{(\kappa+1)n-1}}{((\kappa+1)n-1)!\,\,(\kappa!)^n}}}
\sum_{\alpha_1'+\cdots+\alpha_{n-1}'=\frac{n(n+1)}{2}}
\bigg\{
\\
&
\bigg\{
\sum_{\mu_l^i\in\nabla_{n,\kappa}}\,
(\kappa!)^n
\cdot
{\textstyle{\frac{N_{\mu_1^1}^\kappa}{\kappa\cdots\,\mu_1^1}}}
\cdot
{\textstyle{\frac{N_{\mu_1^2,\mu_2^2}^{\mu_1^1,\kappa}}{
(\kappa+\mu_1^1)\cdots\,(\mu_2^2+\mu_1^2)}}}
\cdots\,
\\
&
\cdots\,
{\textstyle{\frac{
N_{\mu_1^{n-1},\dots,\mu_{n-2}^{n-1},\mu_{n-1}^{n-1}}^{
\mu_1^{n-2},\dots,\mu_{n-2}^{n-2},\kappa}}{
(\kappa+\mu_{n-2}^{n-2}+\cdots+\mu_1^{n-2})
\cdots\,
(\mu_{n-1}^{n-1}+\mu_{n-2}^{n-1}+\cdots+\mu_1^{n-1})}}}
\cdot 
{\textstyle{\frac{
N_{\mu_1^n,\dots,\mu_{n-1}^n,\mu_{n-1}^n}^{
\mu_1^{n-1},\dots,\mu_{n-1}^{n-1},\kappa}}{
(\kappa+\mu_{n-2}^{n-2}+\cdots+\mu_1^{n-2})
\cdots\,
(\mu_{n-1}^{n-1}+\mu_{n-2}^{n-1}+\cdots+\mu_1^{n-1})}}}
\cdot
\\
&
\cdot
\big[\log(\kappa)-\log(\mu_1^1)\big]^{\alpha_1'}\,
\big[\log(\kappa+\mu_1^1)-\log(\mu_2^2+\mu_1^2)\big]^{\alpha_2'}
\cdots\,
\\
&
\cdots
\big[
\log(\kappa+\mu_{n-2}^{n-2}+\cdots+\mu_1^{n-2})
-
\log(\mu_{n-1}^{n-1}+\mu_{n-2}^{n-1}+\cdots+\mu_1^{n-1})
\big]^{\alpha_{n-1}'}
\bigg\}
+
\\
&
+
{\sf Constant}_{n,d,\kappa}\cdot m^{(\kappa+1)n-2}
\endaligned
\]
\[
\footnotesize
\aligned
&
\leqslant
{\sf Constant}_{n,d}\,
{\textstyle{\frac{m^{(\kappa+1)n-1}}{((\kappa+1)n-1)!\,\,(\kappa!)^n}}}
\sum_{\alpha_1'+\cdots+\alpha_{n-1}'=\frac{n(n+1)}{2}}\,
\Delta_{n,\kappa}^{\alpha_1',\dots,\alpha_{n-1}',0}
+
\ \ \ \ \ \ \ \ \ \ \ \ \ \ \ \ \ \ \ \ \ \ \ \ \ \
\ \ \ \ \ \ \ \ \ \ \ \ \ \ \ \ \ \ \ \ \ \ \ \ \ \
\\
&
\ \ \ \ \ \ \ \ \ \ \ \ \
+
{\sf Constant}_{n,d,\kappa}\cdot m^{(\kappa+1)n-2}
\\
&
\leqslant
{\sf Constant}_{n,d}\,
{\textstyle{\frac{m^{(\kappa+1)n-1}}{((\kappa+1)n-1)!\,\,(\kappa!)^n}}}
+
{\sf Constant}_{n,d,\kappa}\cdot m^{(\kappa+1)n-2}.
\endaligned
\]
Lastly, in the trivial minoration:
\[
h^0
\geqslant
\chi-h^2-h^4-h^6-\cdots
\]
for $\mathcal{ E}_{ \kappa, m}^{ GG} T_X^*$, we may apply the
majorations just obtained with $q$ even and deduce that:
\[
\aligned
h^0\big(X,\,\mathcal{E}_{\kappa,m}^{GG}T_X^*\big)
\geqslant
\chi\big(X,\,\mathcal{E}_{\kappa,m}^{GG}T_X^*\big)
&
-
{\sf Constant}_{n,d}\,
{\textstyle{\frac{m^{(\kappa+1)n-1}}{((\kappa+1)n-1)!\,\,(\kappa!)^n}}}
(\log\kappa)^0
-
\\
&
-
{\sf Constant}_{n,d,\kappa}\cdot m^{(\kappa+1)n-2},
\endaligned
\]
so that we even get a better minoration of $h^0$ than
the one stated in the Main Theorem.

Lastly, with $\mathcal{ E}_{ \kappa, m}^{ GG} T_X^* \otimes
\mathcal{ O}_X ( -1)$ instead of $\mathcal{ E}_{ \kappa, m}^{ GG} T_X^*$,
the asymptotic Euler characteristic remains unchanged, and also, all
the previous estimations leading to a minoration of
$\dim H^0 \big(X, \mathcal{ E}_{\kappa, m}^{GG} T_X^* \big)$
remain unchanged as well, asymptotically
when $\kappa \to \infty$, $m \to \infty$,
as is usual and as follows from an inspection of,
{\em e.g.}, Br\"uckmann's three families of long 
exact cohomology sequences presented and applied in 
Section~\ref{Section-11}.
\qed

\markleft{Jo\"el Merker}
\markright{\S13.~Observations about invariant jet differentials}
\section{\bf Speculations about invariant jet differentials}
\label{Section-13}

\subsection{Demailly-Semple invariant jet differentials}
The group ${\sf G}_\kappa$ of $\kappa$-jets at the origin of local
reparametrizations $\phi ( \zeta) = \zeta + \phi'' ( 0) \, \frac{
\zeta^2}{ 2!} + \cdots + \phi^{ ( \kappa) } ( 0) \, \frac{
\zeta^\kappa }{ \kappa !} + \cdots$ of $(\C, 0)$ that are tangent to
the identity, namely which satisfy $\phi' ( 0) = 1$, may be seen to
act linearly on the $n \kappa$-tuples of jet variables $\big( f_{
j_1}', f_{ j_2}'', \dots, f_{ j_\kappa}^{ ( \kappa)} \big)$ by plain
matrix multiplication, {\em i.e.} when we set $g_i^{ ( \lambda)} :=
\big( f_i \circ \phi \big)^{ (\lambda)}$, a computation applying the
chain rule gives for each index $i$:
\[
\footnotesize
\aligned
\left(
\begin{array}{c}
g_i'
\\
g_i''
\\
g_i'''
\\
g_i''''
\\
\vdots
\\
g_i^{(\kappa)}
\end{array}
\right)
=
\left(
\begin{array}{ccccccc}
1 & 0 & 0 & 0 & \cdots & 0
\\
\phi'' & 1 & 0 & 0 & \cdots & 0
\\
\phi''' & 3\phi'' & 1 & 0 & \cdots & 0
\\
\phi'''' & 4\phi'''+3{\phi''}^2 & 6\phi'' & 1 & \cdots & 0
\\
\vdots & \vdots & \vdots & \vdots & \ddots & \vdots 
\\
\phi^{(\kappa)} & \cdots & \cdots & \cdots & \cdots & 1
\end{array}
\right)
\left(
\begin{array}{c}
f_i'\circ\phi
\\
f_i''\circ\phi
\\
f_i'''\circ\phi
\\
f_i''''\circ\phi
\\
\vdots
\\
f_i^{(\kappa)}\circ\phi
\end{array}
\right)
\ \ \ \ \ \ \ \ 
{\scriptstyle{(i\,=\,1\,\cdots\,n)}}.
\endaligned
\]
By definition
({\em see} \cite{ dem1997, rou2006a, mer2008b, dmr2010}), 
Demailly-Semple invariant jet polynomials
${\sf P} \big( j^\kappa f \big)$ 
satisfy, for some integer $m$: 
\[
{\sf P}\big(j^\kappa g)
=
{\sf P}\big(j^\kappa(f\circ\phi)\big) 
=
\phi'(0)^m\cdot
{\sf P}\big((j^\kappa f)\circ\phi\big)
=
{\sf P}\big((j^\kappa f)\circ\phi\big),
\] 
for any $\phi$. 

Then obviously when $\phi' (0) = 1$, the algebra ${\sf E}_\kappa^n$
just coincides with the algebra of invariants for the linear group
action represented by the group of matrices just written:
\[
{\sf P}\big(j^\kappa g\big)
=
{\sf P}\big({\sf M}_{\phi'',\phi''',\dots,\phi^{(\kappa)}}\cdot
j^\kappa f\big)
=
{\sf P}\big(j^\kappa f\big),
\]
with $\phi'', \phi''', \dots, \phi^{ ( \kappa)}$ interpreted as
arbitrary complex constants. Such a group clearly 
has dimension $\kappa - 1$.

This group of matrices is a subgroup of the full unipotent group,
hence it is {\em non-reductive}, and for this reason, it not immediate
to deduce finite generation, valid in the so developed invariant
theory of reductive actions, from Hilbert's averaging Reynold operator.
Moreover, though the invariants of the full unipotent group are
well understood (cf. Section~4), as soon as one looks at a {\em
proper} subgroup of it, formal harmonies seem to rapidly
disappear.

\subsection{Three challenging
questions about effectiveness that are, nevertheless, only preliminary} 
If one prefers to work with
Demailly-Semple jets (instead of working with plain Green-Griffiths
jets), then in order to reach the first stage which would correspond
to knowing the exact Schur bundle decomposition for $\mathcal{ E}_{
\kappa, m}^{ DS} T_X^*$ (instead of the one for $\mathcal{ E}_{
\kappa, m}^{ GG} T_X^*$ provided by the theorem stated at the end of
Section~4), one would have to answer {\em in an effective way} the
following three challenging questions, for which, step by step, we
explain why hidden difficulties would still remain.

\smallskip

{\sf Question~1:} {\sf\em Is the Demailly-Semple algebra 
finitely generated?}

\smallskip

At the end of 2010, B\'erczi and Kirwan~\cite{ BercKirw2012}
established the theorem that this algebra is indeed
finitely generated, for arbitray jet order $\kappa \geqslant 1$,
in any dimension $n \geqslant 1$.
Their approach is a remarkable prolongation on Mumford's geometric
invariant theory.

Concerning constructiveness, writing the 
$\kappa$-jet of an intrinsic holomorphic curve as:
\[
j^\kappa f
=
\big(
f_{i_1}',f_{i_2}'',\dots,f_{i_\kappa}^{(\kappa)}
\big)
\ \ \ \ \ \ \ \ \ \ \ \ \ 
{\scriptstyle{(1\,\leqslant\,i_1,\,\,i_2\,\cdots\,i_\kappa\,\leqslant\,n)}}
\]
they introduced a map of the form:
\[
\big(
f',f'',\dots,f^{(\kappa)}\big)
\longmapsto
\bigg(
{\sf r}
f',\,
{\sf r}\,
f''+{\sf r}\,(f')^2,\,\dots,\,
\sum_{\lambda_1+\cdots+\lambda_s=s}
{\sf r}
f^{(\lambda_1)}\cdots f^{(\lambda_s)},
\dots
\bigg),
\]
where ${\sf r}$ are certain various rational coefficients.
A synoptic matrix view is:
\[
\aligned
&
\ \ \ \ \ \ \
\overset{n}{\leftrightarrow}
\ \ \ \ \ \ \ \ \ \ \ \,
\overset{\binom{n+1}{2}}{\leftrightarrow}
\ \ \ \ \ \ \ \ \ \ \ \ \ \,
\overset{\binom{n+2}{3}}{\leftrightarrow}
\ \ \ \ \ \ \,
\overset{\binom{n+3}{4}}{\leftrightarrow}
\ \ \ \ \
\green{\cdots}
\\
&
\left(
\begin{array}{ccccc}
{\sf r}\,f' & 0 & 0 & 0 & \cdots
\\
{\sf r}\,f'' & {\sf r}\,(f')^2 & 0 & 0 & \cdots
\\
{\sf r}\,f'''& {\sf r}\,f'f'' & 
{\sf r}\,(f')^3 & 0 & \cdots
\\
{\sf r}\,f''''& \,\,{\sf r}\,f'f'''+{\sf r}\,(f'')^2 & 
\,\,{\sf r}\,f''(f')^2 &
\,\,\,{\sf r}\,(f')^3 & \cdots
\\
\vdots & \vdots & \vdots & \vdots & \ddots
\end{array}
\right),
\endaligned
\]
and the size of this matrix is:
\[
\kappa
\times
\Big[
{\textstyle{\binom{n+\kappa}{n}}}-1
\Big].
\]
For example, in the setting of Rousseau~\cite{rou2006a} where
$n = \kappa = 3$:
\[
\text{\rm size}
=
3\times 19.
\]
As a second example, in the setting of~\cite{mer2008b} where 
$n = \kappa = 4$:
\[
\text{\rm size}
=
4\times 69.
\]

The constructive aspect of B\'erczi-Kirwan's theorem is that one 
obtains generators of the Demailly-Semple algebra by 
taking all $\kappa \times \kappa$ minors of the above
matrix, getting a first algebra:
\[
\C\big[\Delta_{i_1,i_2,\dots,i_\kappa}\big],
\]
next by taking the {\em radical ideal} of this first minor algebra,
getting a second, {\em finitely generated algebra
which is full algebra of Demailly jets invariant
under reparametrization}.

But to really reach effectiveness
for applications to the determination
of global algebraic jet differentials,
if one would desire to recover,
with the B\'erczi-Kirwan approach,
the already quite
nontriviall case $n = \kappa = 4$ settled in~\cite{ mer2008b},
one should look at all $4 \times
4$ minors of a $4 \times 69$ matrix. This would make in sum a total
number of minors equal to:
\[
{\textstyle{\binom{69}{4}}}
=
{\textstyle{\frac{69\times 68\times 67\times 66}{1\times 2\times 3\times 4}}}
=
864\,501, 
\]
which would in any case require a considerable
amount of work. 

\smallskip

Suppose nevertheless that one can organize coherently some
understanding of the first algebra of minors. 
Then as a second step, one should also study and understand the
{\em radical ideal} of the algebra of minors.
And radical ideals are known, in effective algebra,
to hide complicated phenomena.

\smallskip
{\sf Question~2:} {\sf\em Is the ideal of relations between a set of
basic generating Demailly-Semple invariants finitely generated?}

\smallskip
Again, in order to be able to describe the {\em exact} Schur bundle
decomposition as was done in Section~4, it is absolutely necessary to
describe effectively and for arbitrary $n$, $\kappa$ the ideal of
relations. For jets of order $\kappa = 4$ in dimension $n = 4$, we
were unable to describe the full ideal of relations between the 2835
basic generating invariants listed
in~\cite{ mer2008b}, not to mention that we ignore what is the
minimal number of generators. We were saved in~\cite{ mer2008b} by the
fact that there are ``only'' 16 basic bi-invariants (minimal number)
and ``only'' 41 relations between them\footnote{\,
We believe that one could attack the seemingly accessible case
$n = 5$, $\kappa = 5$.
} 
(in a Gr\"obner basis for a
certain lexicographic order). 

All these speculative considerations
lead us {\em in fine} to the main metaphysical
question: {\em Are there observable, simple mathematical harmonies in
a certain set of generators and for all the relations between them?}
Without harmonies, there is absolutely 
no hope to treat the case where $n$ and
$\kappa$ are arbitrary.  For $n = 2$, $\kappa = 5$ and for $n = 4$,
$\kappa = 4$, we were unable, in~\cite{ mer2008b}, to discover any
combinatorially convincing global formal harmonies.  Nevertheless,
there could yet be some slight hope as follows.

\smallskip
{\sf Question~3:} {\sf\em Is the algebra of bi-invariants
Cohen-Macaulay?}

\smallskip
Cohen-Macaulayness would be nice. For reductive group actions, this is known
to be hold, but however, {\em almost never in a neat effective way}.
At least, one could dream that the Demailly-Semple algebra is
Cohen-Macaulay and that a basis of so-called {\em primary invariants}
presents some understandable harmonies. It is known, then, that the
effective calculations about Euler-Poincar\'e characteristic and
cohomologies with an adapted reduced Schur bundle decomposition become
much more tractable when one looks only at primary invariants. But for
$n = 4$, $\kappa = 4$ and for $n = 2$, $\kappa = 5$, going through the
mutually independent bi-invariant we exhibited in~\cite{ mer2008b} and
trying to change the generators, we were unable to see or to devise a
neat basis of primary invariants, though for $n = 2$, $\kappa = 4$,
one easily discovers such a basis at first glance. Hence, a
non-effective theorem claiming ``the algebra of Demailly-Semple is
Cohen-Macaulay'' would be useless toward the Green-Griffiths
conjecture because rather, one would really need to know the exact
description of a basis of primary invariants with all their weights in
order to start continuing working toward the Green-Griffiths
conjecture.  But if the algebra is not even Cohen-Macaulay, well, the
next tasks could be even more extremely challenging because, as we
already saw, the end of Section~4 opens several doors to other fields
of hard effective computations when one just deals with the much
simpler Green-Griffiths jets.

\smallskip

Last but not least, we would like to insist on the fact that in the
state of affairs which is current since the 19\textsuperscript{th}
Century, even for the most studied reductive action of ${\sf SL}_2 (
\C)$ on binary forms of degree $d$ in (only) two variables, the {\em
effective answers} to Questions~1, 2 and~3 is unknown for arbitrary
$d$, and is rather extremely challenging in fact. Cayley, Sylvester,
Gordan, Noether, Popov, Grosshans, Springer, Dixmier, Lazard,
Bedratyuk and others did not find any complete closed global tamed
combinatorial harmonies.

\smallskip
Thus in any case, the prohibitive complexity of any {\em effective,
applicable} description of the algebras of invariant jets clearly
prevented, already in 2008 once~\cite{ mer2008a, mer2008b} appeared,
to hope for reaching arbitrary dimension $n \geqslant 2$ and jet order
$\kappa \geqslant n$ with Demailly-Semple jets.

\smallskip
In 2008 also, following Demailly and using an
algebraic version of the {\sl Holomorphic Morse inequalities} 
due to Trapani, Diverio considered a certain {\em subbundle} 
of the bundle of
invariant jets already introduced in the fundamental paper~\cite{dem1997} 
of Demailly which,
in arbitrary dimension $n \geqslant 2$, pushes forward Demailly-El
Goul's~\cite{delg2000} and Siu-Yeung's Wronskians in dimension $2$. This,
for the first time after Siu, opened the door to arbitrary dimension
$n \geqslant 2$, but this was clearly not sufficient to reach the
first step towards the
Green-Griffiths-Lang conjecture, namely existence of differential
equations. In fact, when one examines~\cite{ dmr2010}, 
one realizes that on
hypersurfaces $X^n \subset \mathbb{ P}^{ n+1} ( \mathbb{ C})$, of
degree $d \geqslant 2^{ n^5}$,
differential equations exist with jet order $\kappa = n$ equal to the
dimension, but when one increases the jet order $\kappa \geqslant n$, 
a stabilization of the degree gain occurs, hence
it is impossible to reach the optimal $d \geqslant n+3$ for
the first step towards the
Green-Griffiths conjecture 
with Diverio's technique even for hypersurfaces
$X^n \subset \mathbb{ P}^{ n+1} ( \mathbb{ C})$.
This motivated to come back to the bundle of plain Green-Griffiths jets.

\smallskip

In May 2010, the present article appeared as {\footnotesize\sf
arxiv.org/abs/1005.0405}, and was never submitted to a
mathematics journal.

\smallskip

Six months later, also coming back to plain Green-Griffiths jets for
the reasons explained above, but developing different
negative jet curvature estimates inspired from 
Cowen and Griffiths~\cite{ Cowen-Griffiths-1976} but which had been
blocked for computational reasons by the
untractable algebraic complexity of invariant jets, Demailly~\cite{
dem2011} established the next significant advance towards the
conjecture by establishing, under the sole, optimal, assumption that
$X$ be of general type, that
nonconstant entire holomorphic curves $f \colon \mathbb{ C}
\longrightarrow X$ always satisfy nonzero differential
equations. In~\cite{Paun-2012}, Mihai Pa\u{u}n has provided
a useful guide to enter this computationally less complex field.

\vfill\end{document}